\documentclass[10pt]{article}
\usepackage[utf8]{inputenc}
\usepackage{amsmath}
\usepackage{amsfonts}
\usepackage{amssymb}
\usepackage{algorithm}
\usepackage{algpseudocode}
\usepackage{multirow}
\usepackage{float}
\usepackage{array}
\usepackage[left=0.75 in,right=0.75in,top=0.75in,bottom=0.75in]{geometry}
\usepackage{comment}
\usepackage{graphicx}
\usepackage{url}
\usepackage[authoryear]{natbib}
\usepackage{color}
\usepackage{pgfplots}
\usepackage{caption}
\usepackage{subcaption}
\usepackage{pdflscape}
\usepackage{booktabs} 
\usepackage{enumerate}
\usepackage{authblk}
\usepackage{blkarray}
\usepackage{bm}
\usepackage{marvosym} 
\usepackage{soul}
\usepackage{textcomp, xspace}
\usepackage{longtable}
\usepackage{lipsum} 
\usepackage{lineno}
\usepackage{fontawesome}
\usepackage{appendix}
\usepackage{titlesec}
\usepackage{amsthm}

\usepgfplotslibrary{dateplot}

\definecolor{greenxllite}{RGB}{196,215,155}
\definecolor{bluexllite}{RGB}{149,179,215}
\definecolor{redxllite}{RGB}{218,150,148}
\definecolor{greenxl}{RGB}{155, 187, 89}
\definecolor{bluexl}{RGB}{79,129,189}
\definecolor{redxl}{RGB}{192, 80, 77}
\definecolor{OliveGreen}{RGB}{70,136,52}

\definecolor{customblue}{RGB}{127,179,213}
\definecolor{customgray}{RGB}{153,163,164}

\newcommand{\zeroparkvar}[2]{x^{0}_{#1#2}}
\newcommand{\ninetyparkvar}[2]{x^{90}_{#1#2}}

\newcommand{\drivevar}[2]{y_{#1#2}}
\newcommand{\flowvar}[4]{f_{#1#2,#3#4}}
\newcommand{\flowvardash}[4]{g_{#1#2,#3#4}}
\newcommand{\countvar}[4]{z_{#1#2,#3#4}}

\newcommand{\lotset}{L}
\newcommand{\validzeroparkset}{P^{0}}
\newcommand{\validninetyparkset}{P^{90}}
\newcommand{\validdriveset}{D}

\newcommand{\existingdriveset}{Q}
\newcommand{\blockedset}{B}

\newcommand{\grid}{G}
\newcommand{\gridadjset}[2]{D_{#1#2}}
\newcommand{\gridinvadjset}[2]{D^{-1}_{#1#2}}
\newcommand{\gridarcset}{A}

\newcommand{\zeroneighborset}[2]{N_{#1#2}^{0}}

\newcommand{\ninetyneighborset}[2]{N_{#1#2}^{90}}

\newcommand{\invzeroparkset}[2]{\mathcal{P}_{#1#2}^{0}}
\newcommand{\invninetyparkset}[2]{\mathcal{P}_{#1#2}^{90}}
\newcommand{\invdriveset}[2]{\mathcal{D}_{#1#2}}

\newcommand{\numrows}{\mu}
\newcommand{\numcols}{\nu}
\newcommand{\parkwidth}{\omega}
\newcommand{\parklength}{\ell}
\newcommand{\drivewidth}{\delta}

\newcommand{\entranceset}{\{(\entrancex, \entrancey)\}}

\newcommand{\exitset}{\{(\exitx, \exity)\}}

\newcommand{\entrancex}{p}
\newcommand{\entrancey}{q}

\newcommand{\exitx}{r}
\newcommand{\exity}{s}

\newcommand{\vertexpartitions}[1]{\mathcal{V}_{#1}}

\newcommand{\vertexcutset}{V}
\newcommand{\twowaypartition}[1]{V_{#1}}
\newcommand{\twowaypartitioncomp}[1]{V_{#1}'}

\newcommand{\onewaypartitionsets}[1]{\mathcal{W}_{#1}}
\newcommand{\onewaypartition}{W}
\newcommand{\onewaypartitioncomp}{W'}

\newcommand{\invzeroparksetcut}[3]{\mathcal{P}^{0}_{#1#2}(#3)}
\newcommand{\invninetyparksetcut}[3]{\mathcal{P}_{#1#2}^{90}(#3)}
\newcommand{\invdrivesetcut}[3]{\mathcal{D}_{#1#2}(#3)}

\titlespacing*{\section}
{0pt}{0.1\baselineskip}{0.1\baselineskip}
\titlespacing*{\subsection}
{0pt}{0.1\baselineskip}{0.1\baselineskip}
\titlespacing*{\subsubsection}
{0pt}{0.1\baselineskip}{0.1\baselineskip}

\usetikzlibrary{patterns}
\newtheorem{theorem}{Theorem}
\newtheorem{proposition}{Proposition}
\theoremstyle{definition}
\newtheorem{remark}{Remark}

\allowdisplaybreaks

\newcommand{\Keywords}[1]{\par\noindent
{\small{\em \textbf{Keywords}\/}: #1}}

\usepackage[colorlinks]{hyperref}
\hypersetup{
    bookmarks=true,
    citecolor=blue,
    urlcolor=blue,
    linkcolor=blue
}

\title{A Branch-and-Cut Algorithm for the Optimal Design of Parking Lots with One-way and Two-way Lanes}

\author[1] {Helen Thomas*}
\author[2,3] {Tarun Rambha*\textsuperscript{\faEnvelopeO}\,}

\affil[1]{\small Department of Civil Engineering, Indian Institute of Science (IISc), Bengaluru, India}
\affil[2]{\small Centre for infrastructure, Sustainable Transportation, and Urban Planning (CiSTUP), Indian Institute of Science (IISc), Bengaluru, India}
\affil[3]{\small Robert Bosch Centre for Cyber Physical Systems (RBCCPS), Indian Institute of Science (IISc), Bengaluru, India}

\date{ }

\begin{document}
\maketitle
\let\thefootnote\relax\footnotetext{(*) Both authors contributed equally to this manuscript. (\faEnvelopeO~tarunrambha@iisc.ac.in) Corresponding author}
\vspace*{-5mm}
\begin{abstract}
We address the problem of maximizing the number of stalls in parking lots where vehicles park perpendicular to the driveways. Building on recent research on two-way driving lanes, we first formulate a mixed integer program to maximize the number of parking stalls using a flow-based approach. Parking lots are rasterized into a grid, and the proposed MIP model optimizes them in a generic manner, adapting to the grid resolution and stall size without requiring custom formulations. The constraints ensure the connectivity of parking stalls and driveways to the entrance/exit. This formulation is then extended to the case of one-way driving lanes. We then propose valid inequalities and a reformulation that can be solved using a branch-and-cut algorithm. This approach eliminates flow variables and big-M-type constraints, and improves solution times for medium-sized instances. The effectiveness of the suggested models is showcased on 325 parking lots from New York City. For instances where the flow version could be solved in 15 minutes, the branch-and-cut algorithm improved the median runtimes by 87.43\% for the one-way case and by 79.36\% for the two-way case, and achieved better optimality gaps than the baseline flow-based formulation for the other instances. Similar advantages were observed when run with a time budget of two hours. One-way configurations accommodated, on average, 18.63\% more vehicles on average than their two-way counterparts across all instances. Modifications to the proposed formulations that account for vehicle turning characteristics and the presence of multiple entrances and exits are also examined.

\vspace{3mm}
\end{abstract}

\section{Introduction}
In most cities, land for parking near central business districts and other high-demand points of interest is usually limited, making parking scarce and expensive. For example, owning a parking space in parts of New York City can cost a staggering six-digit figure \citep{nytimes2007, nytimes2019}. Such problems are not only acute for cars but also for other vehicles, such as buses, which need to park at depots located in dense urban areas (to avoid deadheading), and commercial heavy vehicles serving urban freight. Optimizing the capacity of commercial parking lots can also significantly improve the revenue of operators and reduce drivers' parking search times, especially during peak periods. For example, data from \cite{statista_parking_revenue_us} suggest that revenue from parking lots in New York has been rising, except during the COVID years, and was estimated to be approximately 1.37 billion USD in 2024. A similar trend can be observed across the United States, according to the Federal Reserve Economic Data database \citep{fred_parking_revenue}. From an occupancy standpoint, \cite{nelsonnygaard108th} studied three parking lots on West 108th Street in New York City and reported an average utilization of over 90\% during peak periods. Thus, it is essential to maximize the limited parking space available, particularly in small- to medium-sized off-street parking lots that tend to frequently reach capacity. Furthermore, the shapes of parking lots can be irregular and may contain obstacles, making it difficult to design these layouts manually. \cite{stephan2021layout}, using data from 177 real-world parking lots in Jena, Germany, found that manually designed lots typically have 25\% less capacity than those designed using Mixed Integer Programming (MIP). \cite{cudzik2024} developed heuristic methods that compared five real-world layouts from Europe and the United States and found that their method yielded lots with approximately 16\% more parking spaces than existing designs. These findings highlight that the design of parking lots involves several key decisions, which will be the focus of our paper. \textit{Where should parking stalls and driveways be located in lots with one-way or two-way driving lanes? How can we solve these optimization problems faster using exact methods? Can the layouts be adjusted to accommodate constraints arising from turn restrictions and multiple entrances and exits?} In the following subsections, we summarize related studies, highlight their gaps, and outline the contributions of our work.

\subsection{Background}
In many urban areas, drivers are often forced to park far from their destinations or spend considerable time searching for parking, which causes congestion. A study by \cite{hampshire2018share} estimated that in Stuttgart, the percentage of vehicles cruising for parking is about 15\%. This share varies widely across cities but is usually a significant fraction of the total traffic \citep{shoup2006cruising, lee2017cruising}. Several parking solutions have been explored in the literature, including measuring parking availability \citep{parking-ITS-DL, arora2019hard}, optimizing parking search times \citep{tang2014modeling, onstreet-parking-search}, and implementing policy measures \citep{parking-policies}. The reader can find in-depth surveys on intelligent parking systems in \cite{lin2017survey} and \cite{smart-parking-city}. At the parking lot level, recent research focuses on optimizing stall layouts for automated parking systems \citep{wu2019optimal, scheduling-unmanned-parking, automated-vehicles-parking, puzzle-parking-design, automated-valet-parking, wang2021, chen2023}. Problems of similar flavor appear in the planning and organization of automated warehouses and puzzle-based storage systems \citep{flake2002rush, gue2007puzzle, azadeh2019robotized, bukchin2022comprehensive, rosenfeld2022optimal, raviv2023optimal}. The two-dimensional bin packing problem \citep{lodi2002, pisinger2007} also involves fitting as many rectangular fields of different sizes as possible within a given region, but does not include reachability constraints. In this research, we focus on the problem of designing a car parking lot, but many of the methods and ideas can be trivially extended to other settings.

A crucial aspect of lot design involves determining the optimal configuration of parking stalls or fields. \cite{iranpour1989methodology} made one of the earliest attempts to solve this problem using a set of nonlinear equations to maximize capacity. They introduced the concept of an ease factor, which reflected the ease/difficulty associated with parking maneuvers and helped select the optimal layout among multiple solutions with equal parking capacity. \cite{abdelfatah2014parking} used MIP models and the LINDO solver to design optimal parking layouts at different angles. The total number of parking stalls was maximized by deciding how many and what types (based on the dimensions and angle of parking used) of rows of parking bays could be accommodated within the given parking lot. \cite{yildirim2019} employed a cutting-stock-like model to maximize the number of parking stalls within predefined bay patterns, while circulation driveways were pre-allocated to ensure accessibility. Other studies that optimize parking stalls are also restricted to rectangular- or triangular-shaped lots with predefined bay configurations \citep{syahrini2018, maina2025}. \cite{stiglmayr2020} presented one of the first approaches to the parking layout design problem without predefined bay–aisle structures or conventional configurations prescribed by design manuals. In their framework, the parking lot is discretized into a grid, with each cell assigned as either a parking space or a driving space. Distance-based and flow-based formulations are proposed to ensure accessibility, together with a heuristic solution method. Although accessibility is explicitly considered, the framework has limited applicability because it provides little control over the geometric dimensions of parking spaces and driveways. In research employing heuristic techniques, \cite{cudzik2024} assigned parking stalls along site boundaries and generated additional rows of stalls using directions derived from the boundary polylines. Multiple lot designs were generated, and the maximum number of stalls produced by the proposed methods was used for comparison with existing layouts. 

In summary, most of these methods lack generality. They typically address only selected components of the overall layout design problem or limit their attention to rectangular lots, and consequently do not provide practical value. Furthermore, many existing works do not model obstructions such as building elements, trees, sidewalks, and streetlights, and these practical constraints necessitate better optimization. \cite{stephan2021layout} addressed these drawbacks using a MIP and proposed heuristics to maximize the capacity of perpendicular parking lots. Their study rasterized a given ground plot into three different resolutions. \cite{TRISTAN2022website} suggested using cutting planes to reduce the time required to solve a similar MIP. 

\subsection{Gaps and contributions}
The results of field studies assessing vehicle capacity and maneuverability can serve as guidelines for parking design. However, manually designing parking lots that utilize space effectively, especially in irregularly shaped areas, is challenging and requires optimization. While a few automated solutions exist in the literature \citep{abdelfatah2014parking, syahrini2018, yildirim2019, cudzik2024, maina2025}, as discussed earlier, they are limited by the configurations, size, and shape of parking fields, among other factors. This work addresses the problem of designing parking stalls and driveways within a given plot, whether on land, at a single level of a multi-floor parking facility, or in a warehouse. We extend and improve recent MIP formulations \citep{stephan2021layout, TRISTAN2022website} that maximize the number of parking stalls oriented perpendicular to the driving lanes.

\begin{figure}[H]
  \centering
  \begin{subfigure}[b]{0.24\textwidth}
    \includegraphics[width=\textwidth]{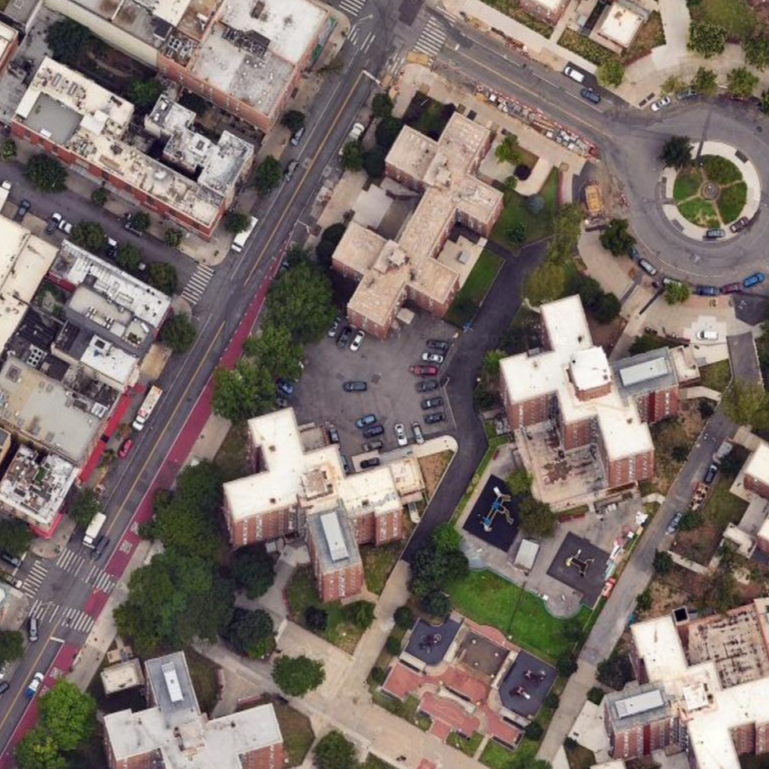}
    \caption{Google map view}
   \label{fig:gmaps}
  \end{subfigure}
  \begin{subfigure}[b]{0.24\textwidth}
    \includegraphics[width=\textwidth]{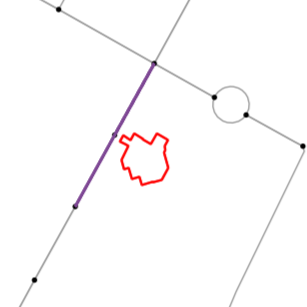}
    \caption{OSM network}
    \label{fig:osm}
  \end{subfigure}  
  \begin{subfigure}{0.24\textwidth}
    \centering
    \includegraphics[width=0.93\textwidth]{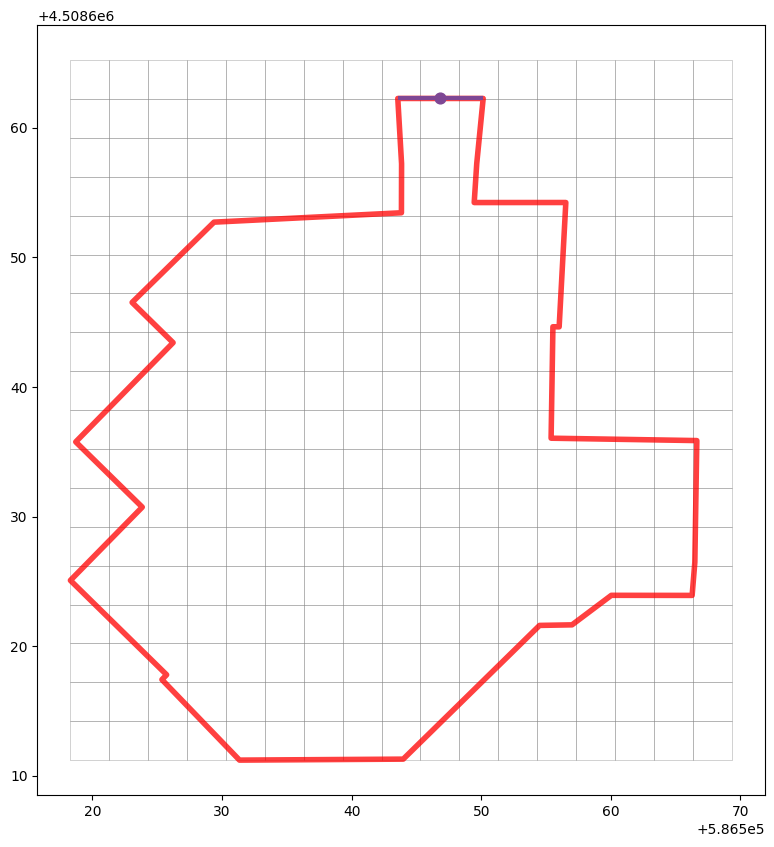}
    \caption{Grid with 3 m cells}
    \label{fig:3m_grid}
  \end{subfigure}
  \begin{subfigure}{0.24\textwidth}
    \centering
    \includegraphics[width=\textwidth]{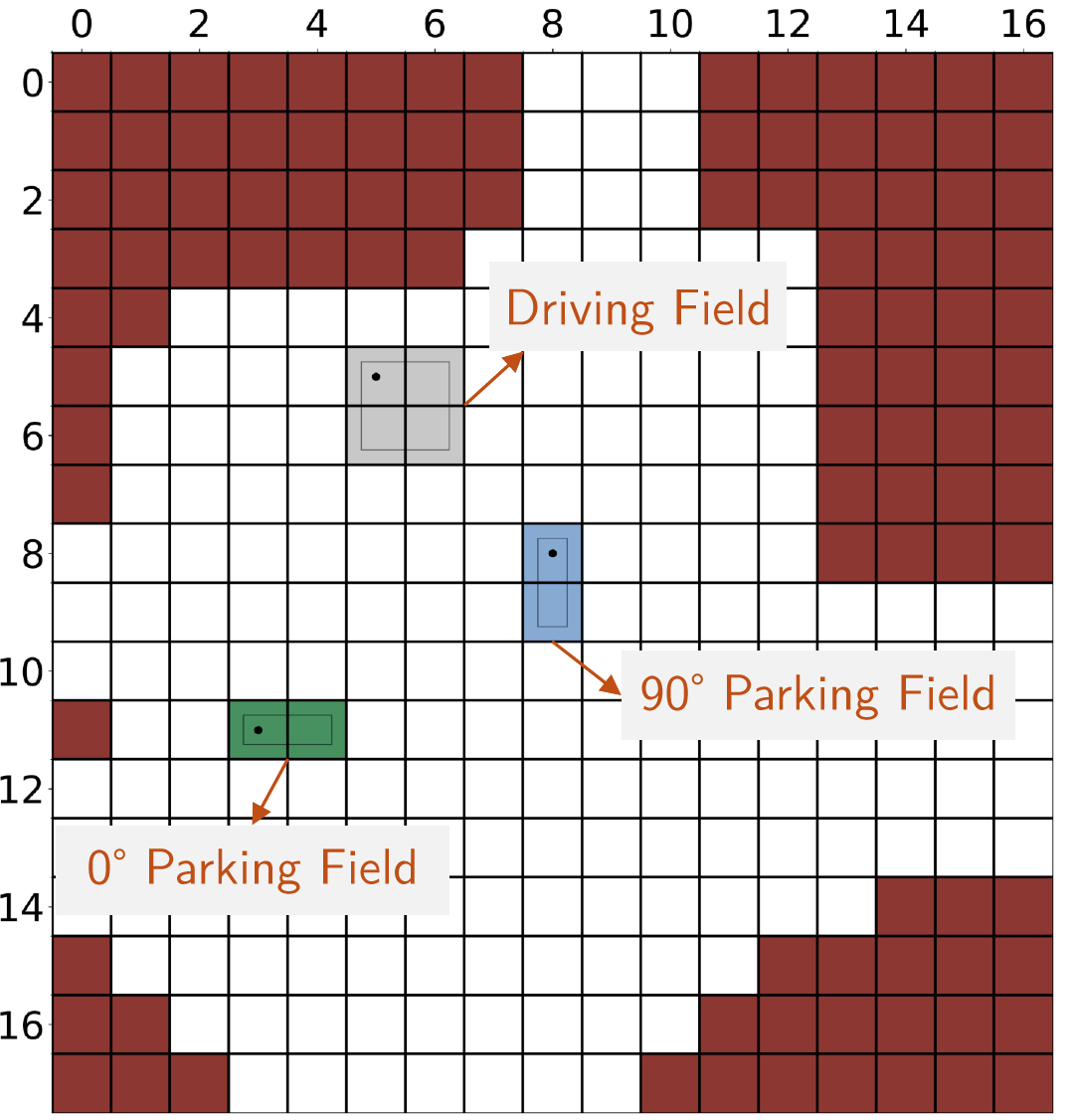}
    \caption{Sample fields}
    \label{fig:resolution_r2}
  \end{subfigure}
\caption{Processing a parking lot instance}
\end{figure}

A parking lot can be represented as a polygon with a grid of cells, some of which are blocked and some are available for parking. Access to the parking lot can be found by scanning a circular area (of a suitable radius) around the lot using, say, OpenStreetMaps (OSM) to identify the driving streets in its vicinity. The entrance to the parking lot may be known, or can be set as one of the edges of the polygon closest to the nearby streets. For example, Figure \ref{fig:gmaps} shows a satellite image of a parking lot, and its corresponding polygon is shown in Figure \ref{fig:osm}. The closest street is then used to set the entrance (see the purple segment in Figure \ref{fig:3m_grid}). The input to our formulations is a rasterized layout such as the one shown in Figure \ref{fig:resolution_r2}. These cells can accommodate parking stalls/fields shown in green and blue, depending on their orientation, and driving fields indicated in gray. Adjacent green and blue cells represent lateral ($0^{\circ}$) and longitudinal ($90^{\circ}$) parking fields, and each pair of cells accommodates a single vehicle. The red cells represent areas blocked either by obstructions or due to irregular boundaries.

\begin{figure}[H]
  \centering
  \begin{minipage}{0.75\textwidth}
  \begin{subfigure}{0.48\textwidth}
    \centering
    \includegraphics[width=0.75\textwidth]{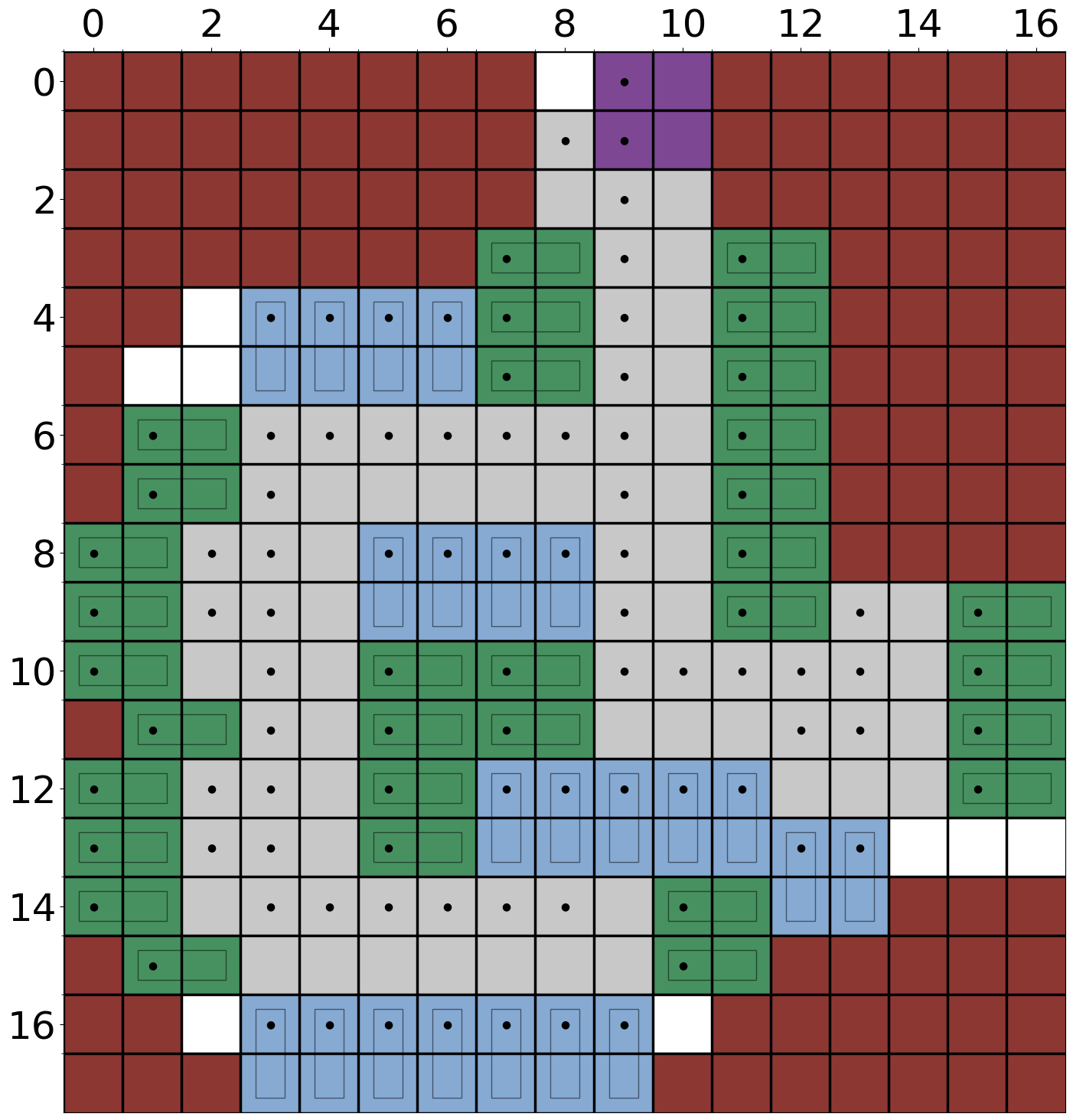}
    \caption{Two-way configuration with 54 stalls}
    \label{fig:example_two_way}
  \end{subfigure}
  \begin{subfigure}{0.48\textwidth}
    \centering
    \includegraphics[width=0.75\textwidth]{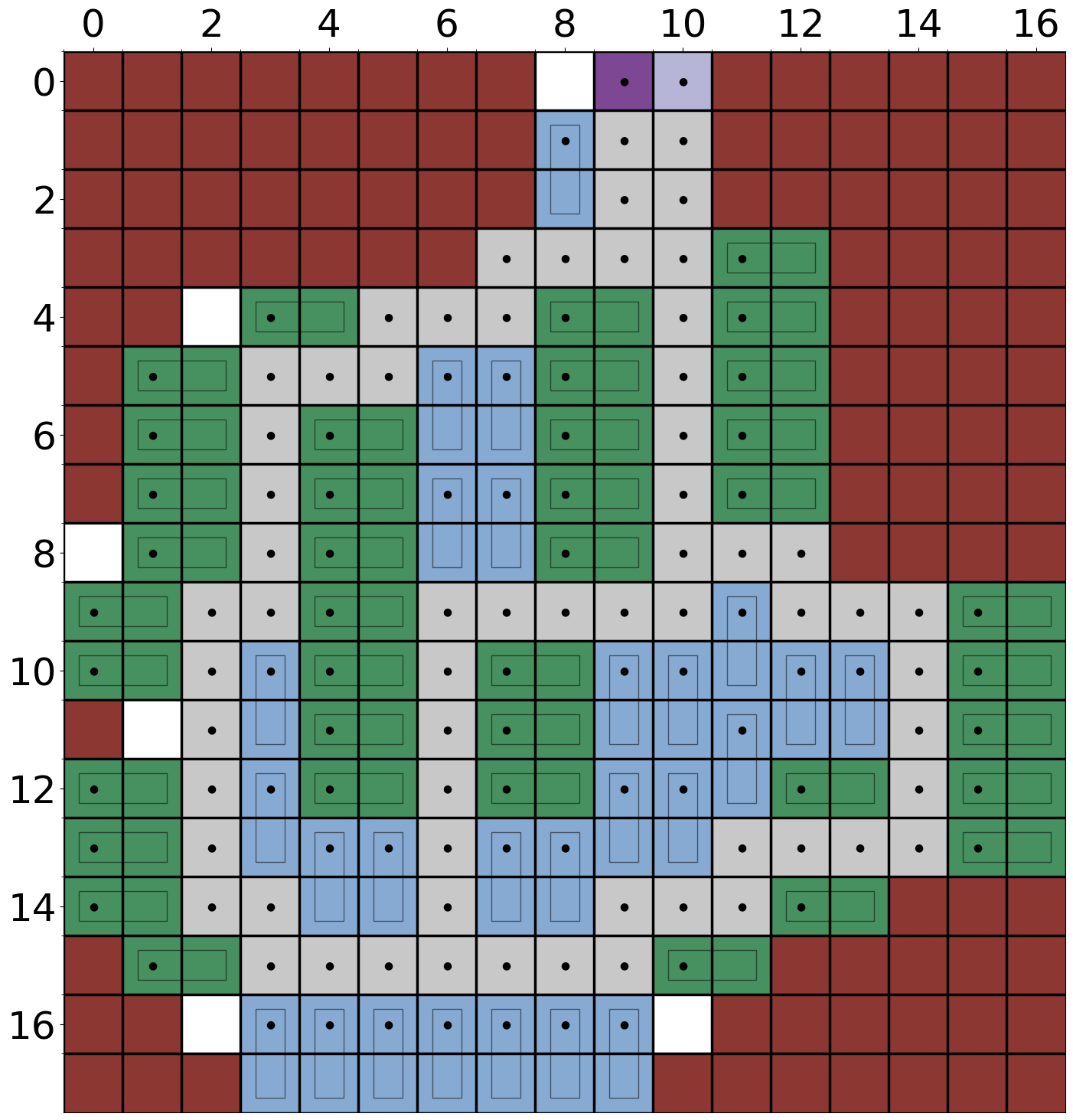}
    \caption{One-way configuration with 65 stalls}
    \label{fig:example_one_way}
  \end{subfigure}
\end{minipage}
\begin{minipage}{0.24\textwidth}
\flushleft
\vspace{-3mm}
\begin{tikzpicture}

    \definecolor{purple}{RGB}{126,71,148}
    \definecolor{blue}{RGB}{135,170,211}
    \definecolor{green}{RGB}{71,145,96}
    \definecolor{grey}{RGB}{200,200,200}
    \definecolor{lavender}{RGB}{182,181,216}
    \definecolor{red}{RGB}{141,55,51}
    \definecolor{white}{RGB}{255,255,255}


\filldraw[fill=white, draw=black] (-0.7,0.0) rectangle (-0.4,0.3);
\node[right=6pt] at (-0.4,0.15)
    {\footnotesize Anchor cell};

\fill (-0.55,0.15) circle (1pt);
    
\filldraw[fill=white, draw=black] (-0.7,0.55) rectangle (-0.4,0.85);
\node[right=6pt] at (-0.4,0.70)
    {\footnotesize Empty cell};

\filldraw[fill=lavender, draw=black] (-0.7,1.10) rectangle (-0.4,1.40);
\node[right=6pt] at (-0.4,1.25)
    {\footnotesize Exit cell};

\filldraw[fill=purple, draw=black] (-0.7,1.65) rectangle (-0.4,1.95);
\node[right=6pt] at (-0.4,1.80)
    {\footnotesize Entrance cell};

\filldraw[fill=grey, draw=black] (-0.7,2.20) rectangle (-0.4,2.50);
\node[right=6pt] at (-0.4,2.35)
    {\footnotesize Drive cell};

\filldraw[fill=blue, draw=black] (-0.7,2.75) rectangle (-0.4,3.05);
\node[right=6pt] at (-0.4,2.90)
    {\footnotesize Longitudinal park field cell};

\filldraw[fill=green, draw=black] (-0.7,3.30) rectangle (-0.4,3.60);
\node[right=6pt] at (-0.4,3.45)
    {\footnotesize Lateral park field cell};

\filldraw[fill=red, draw=black] (-0.7,3.85) rectangle (-0.4,4.15);
\node[right=6pt] at (-0.4,4.00)
    {\footnotesize Blocked cell};

\end{tikzpicture}
\end{minipage}
\caption{Optimal parking lot configurations for two-way and one-way driving lane configurations}
\label{fig:introduction_example}
\end{figure}

The MIP formulations presented in this paper are generic in nature. They solve the problem for homogeneous lots (e.g., passenger cars or unstacked cargo crates/containers in a yard or warehouse, as long as the stall size is the same throughout the lot) and model different resolutions using a single formulation. They also handle the case of designing lane configurations for one-way driving scenarios. Figure~\ref{fig:introduction_example} presents representative outputs obtained from the proposed models. The purple field functions as both the entrance and exit in two-way configurations. The role of anchor cells will be described in the next section.

The main contributions of this research are summarized below:

\begin{itemize}
\itemsep 0pt
\item We introduce a generalized optimization framework for real-world parking layouts that can model irregular lot geometries, various rasterized resolutions, and homogeneous vehicle sizes. Compared to manually designed lot configurations, which have been documented to yield lower capacities, these approaches can help planners optimize space utilization.
\item We extend existing MIP models to one-way lane configurations, modeling the direction of vehicle movement and adding constraints that connect the driveways to both the entrance and exit while avoiding dead ends.
\item We propose a family of valid inequalities that can be added to the model upfront to generate solutions with a small number of connected driving field components. We also develop a reformulation of the problem and a branch-and-cut algorithm in which feasibility cuts are dynamically added to ensure driveway connectivity, eliminating flow variables and big-M constraints.
\item We demonstrate computational savings using 325 medium-sized instances (up to half the size of a football field) from New York City. For two-way (one-way) configurations, flow-based formulations optimally solved 187 (149) of the 325 instances within a time limit of 15 min for each instance. On the other hand, our branch-and-cut algorithm, applied to a reformulation of the problem, solved 245 (288) instances optimally with median runtimes that were 79.36\% (87.43\%) faster when both methods found optimal solutions. One-way configurations resulted in an overall increase in the number of stalls by 18.63\% compared to lots with two-way driveways. 
\end{itemize}

\begin{remark}
Our MIP formulations generalize the models of \cite{stephan2021layout}, as illustrated by the following instances from their work. For example, from the right-most panel of Figure \ref{fig:two_way_resolutions}, one can generate a coarser problem instance by combining 2 × 2 units of cells to generate the data for the middle panel. Similarly, 2 × 2 red-cell units in the middle panel can be grouped to create the blocked cells in the left-most panel of the Figure. 

\begin{figure}[H]
  \centering
  \begin{subfigure}[t]{0.99\textwidth}
    \centering
    \includegraphics[width=0.27\textwidth]{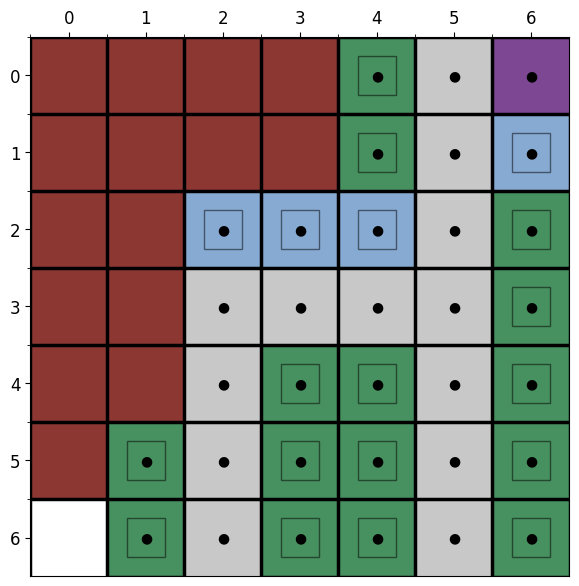}
    \hspace{2mm}
    \includegraphics[width=0.27\textwidth]{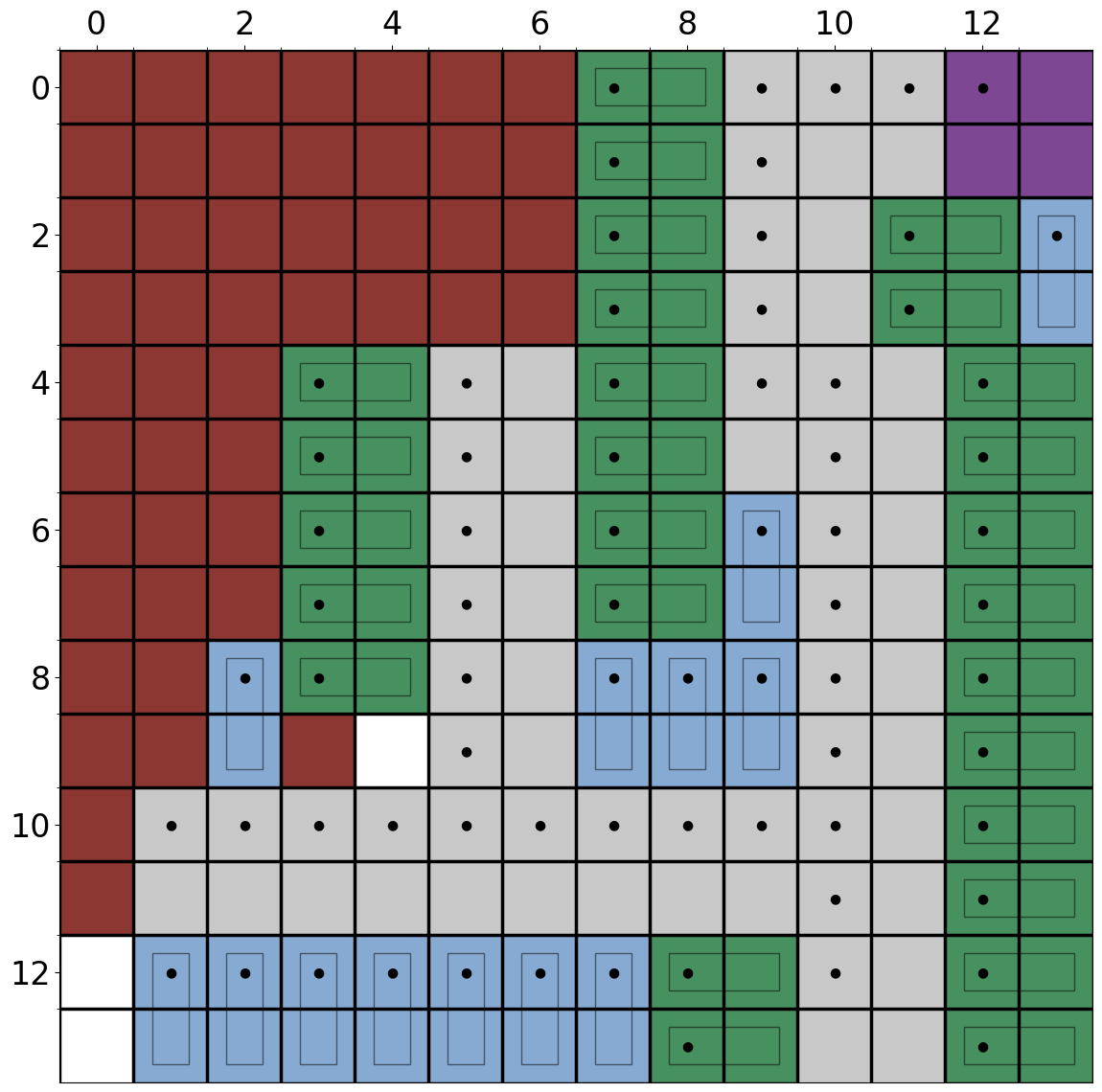}
    \hspace{2mm}
    \includegraphics[width=0.27\textwidth]{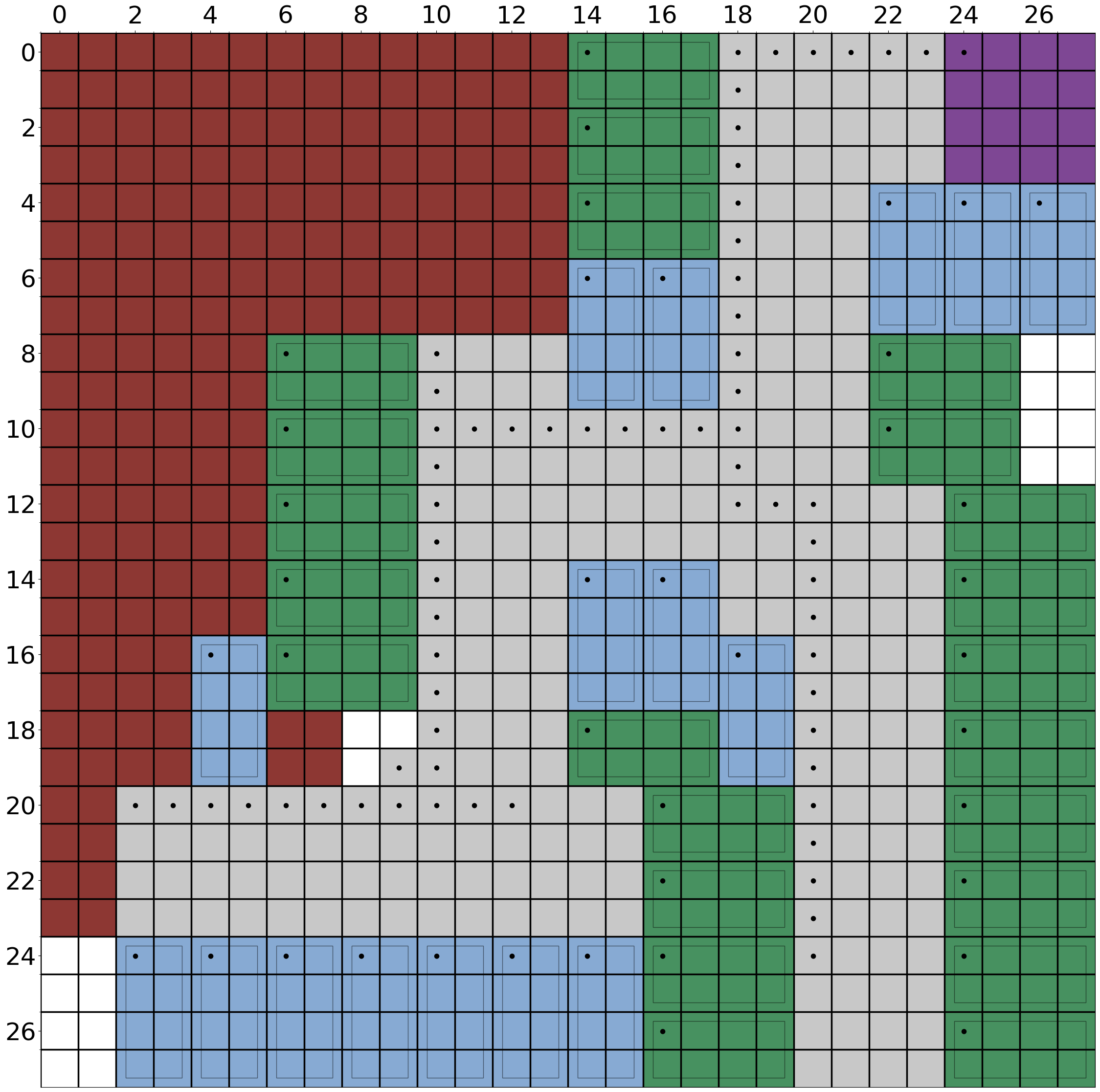}
    \caption{Two-way layouts with different resolutions}
    \label{fig:two_way_resolutions}
  \end{subfigure}

  \vspace{1em} 

  \begin{subfigure}[t]{0.99\textwidth}
    \centering
    \includegraphics[width=0.27\textwidth]{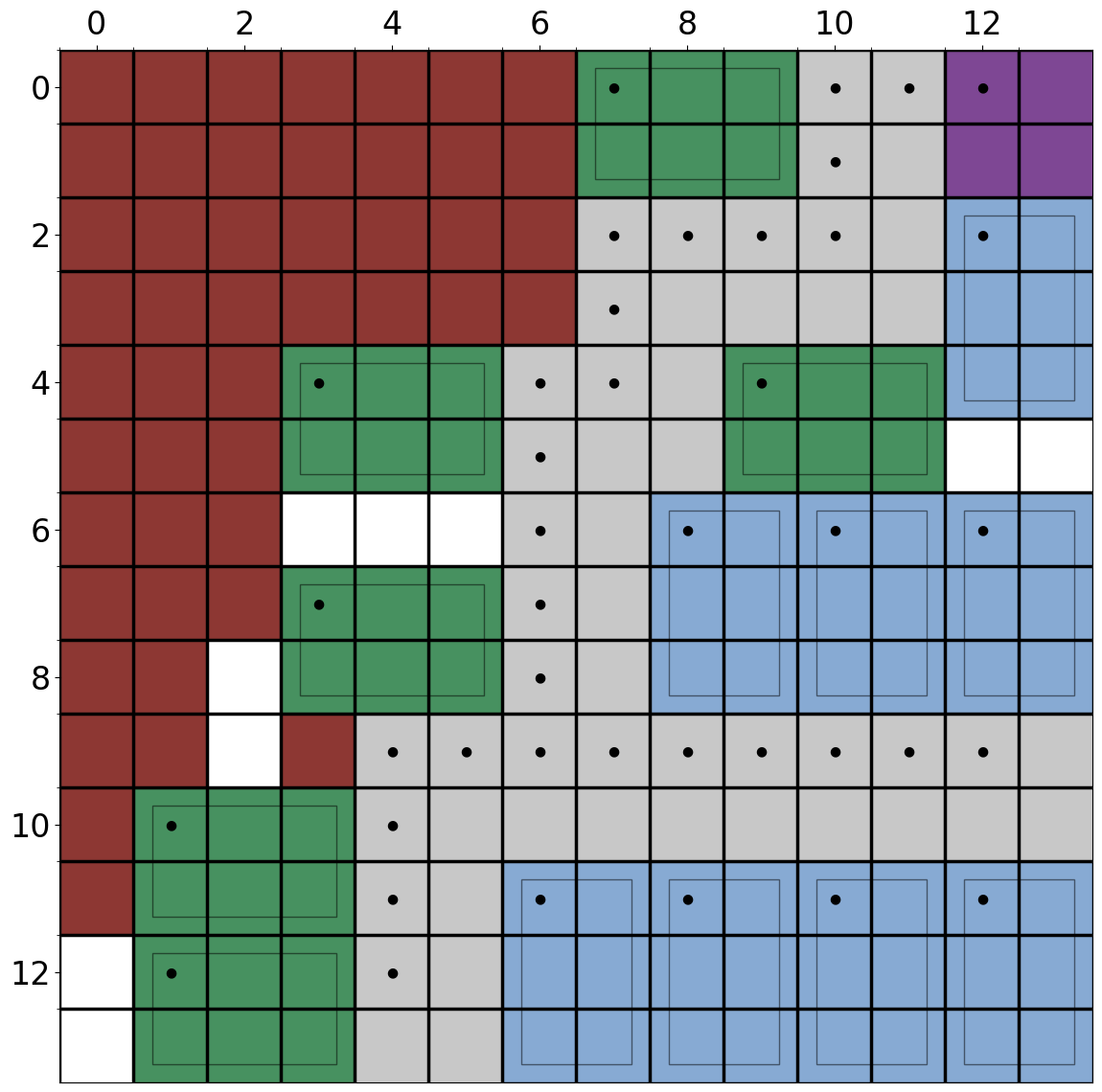}
    \hspace{4mm}
    \includegraphics[width=0.27\textwidth]{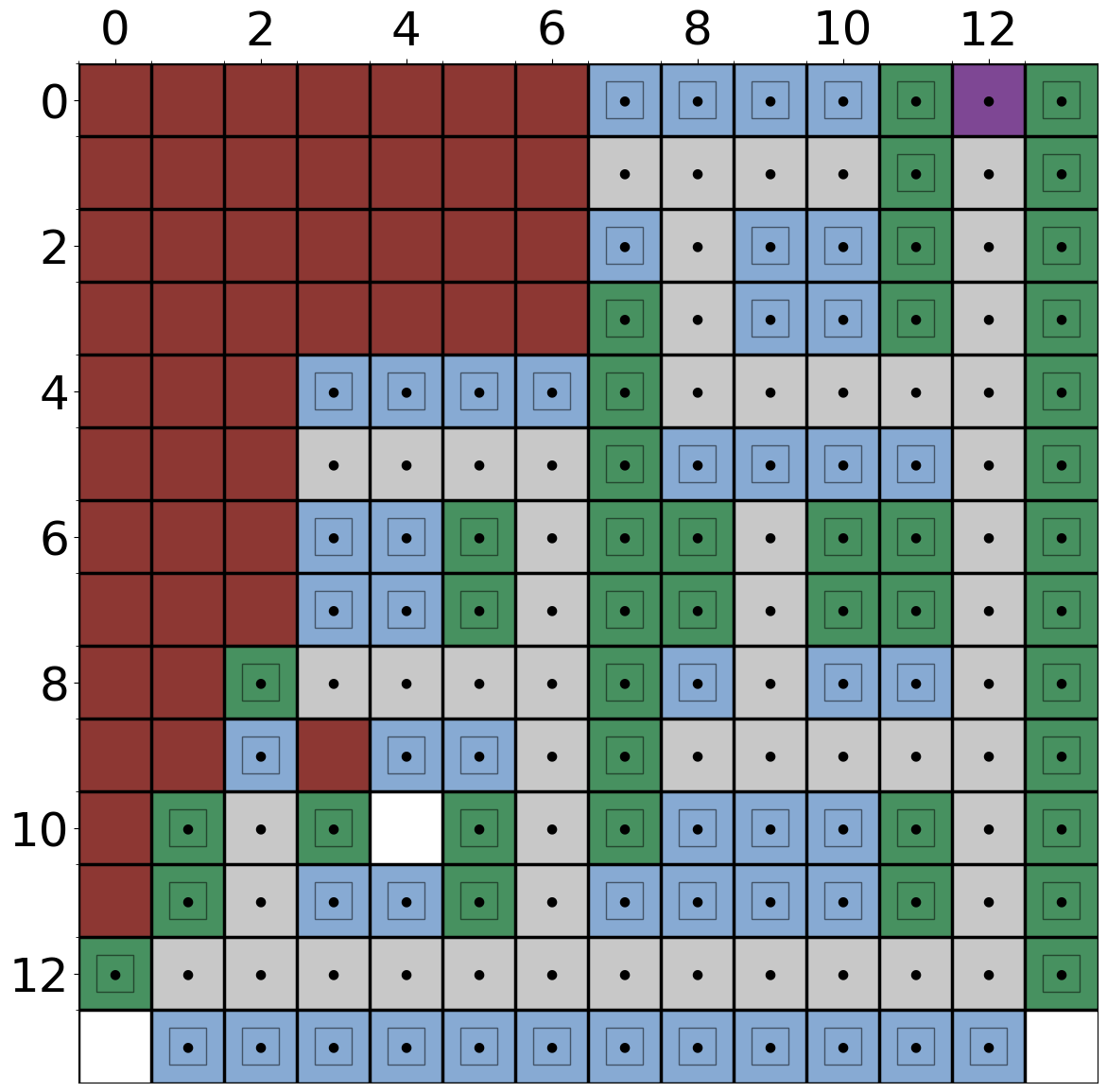}  
    \caption{Two-way layouts with different-sized parking and driving fields}
    \label{fig:two_way_different_sizes}
  \end{subfigure}
  \caption{Solutions for different parking and driving configurations for the layout used in \cite{stephan2021layout}}
    \label{fig:resolutions}
\end{figure}

Figures in \ref{fig:two_way_different_sizes} show examples of our model output for the same grid resolution as the middle panel of Figure \ref{fig:two_way_resolutions}, but with different parking/driving field sizes. The orientation is not relevant for square-shaped parking fields. The right-most panels of Figures \ref{fig:two_way_resolutions} and \ref{fig:two_way_different_sizes} show lower-bound solutions, as they were terminated after reaching a time limit. All these instances can be solved using the MIP formulation described in the next section. Throughout the paper, we present examples with different grid resolutions and parking and driving field sizes selected to best illustrate the proposed concepts.
\end{remark}

The remainder of this paper is organized as follows. Section \ref{sec:MIP} describes a generalized MIP formulation with two-way driving lanes. We discuss extensions to the one-way case in Section \ref{sec:one_way}. Section \ref{sec:alternate_twoway} suggests valid inequalities for the problem and a reformulation with a branch-and-cut algorithm that dynamically generates feasibility cuts. Section \ref{sec:mipresults} demonstrates the advantages of the proposed formulations using open data from New York City. Finally, Section \ref{sec:conclusion} summarizes the findings in this paper along with a discussion on possible extensions. Proofs establishing the equivalence of different formulations are provided in \ref{sec:equivalence}. Additional constraints addressing vehicle turning-radius limitations, multiple entrances and exits, and obstacle placement are discussed in  \ref{sec:practical_considerations}. Supplementary experimental results are included in \ref{sec:add_stats}.

\section{A generalized MIP formulation}
\label{sec:MIP}

\subsection{Problem setup}
\label{sec:perp problem setup}
This section describes a MIP formulation for the case where parking stalls are perpendicular to driveway lanes, and the driveway lanes are two-way streets. We first rasterize a given plot into a grid of squares or \textit{cells}. A \textit{field} is assumed to be a square or rectangular unit comprising multiple cells. Let the number of rows and columns in the grid be $\numrows$ and $\numcols$, respectively. The set of all cells in the parking lot is designated as $\lotset =\{(i,j): i \in \{0,\hdots,\numrows-1\}, j \in \{0,\hdots,\numcols-1\} \}$. We suppose that the lot may have obstacles such as buildings, trees, utility poles, shopping cart corrals, and structural columns (in the case of indoor parking). The set of such blocked cells is denoted as $\blockedset$. 

\begin{figure}[H]
  \centering
  \begin{subfigure}[b]{0.24\textwidth}
      \centering
    \includegraphics[width=\textwidth]{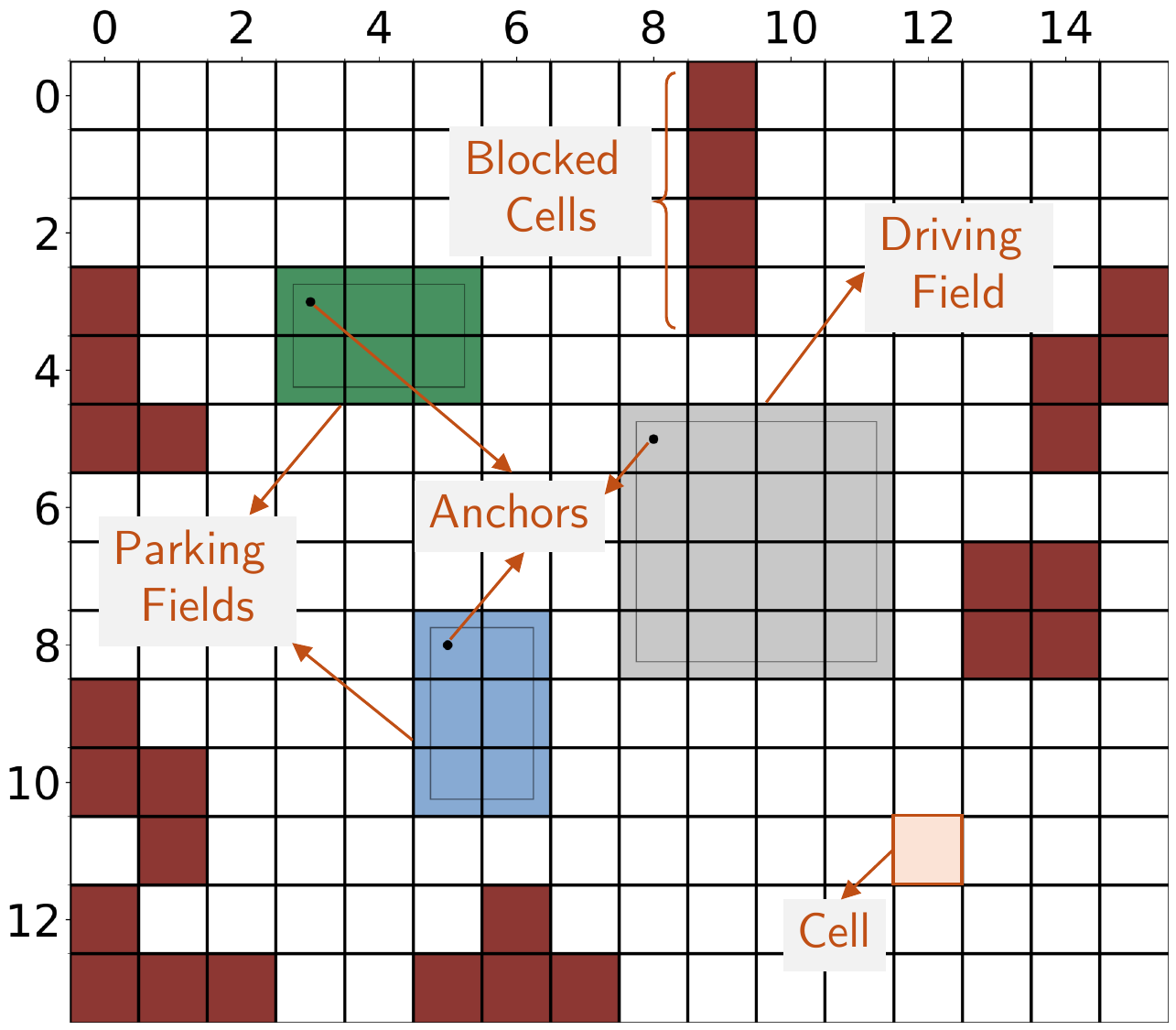}
    \caption{Sample components}
    \label{fig:sample_perpendicular_fields}
  \end{subfigure}
  \hfill
  \begin{subfigure}[b]{0.24\textwidth}
      \centering
    \includegraphics[width=\textwidth]{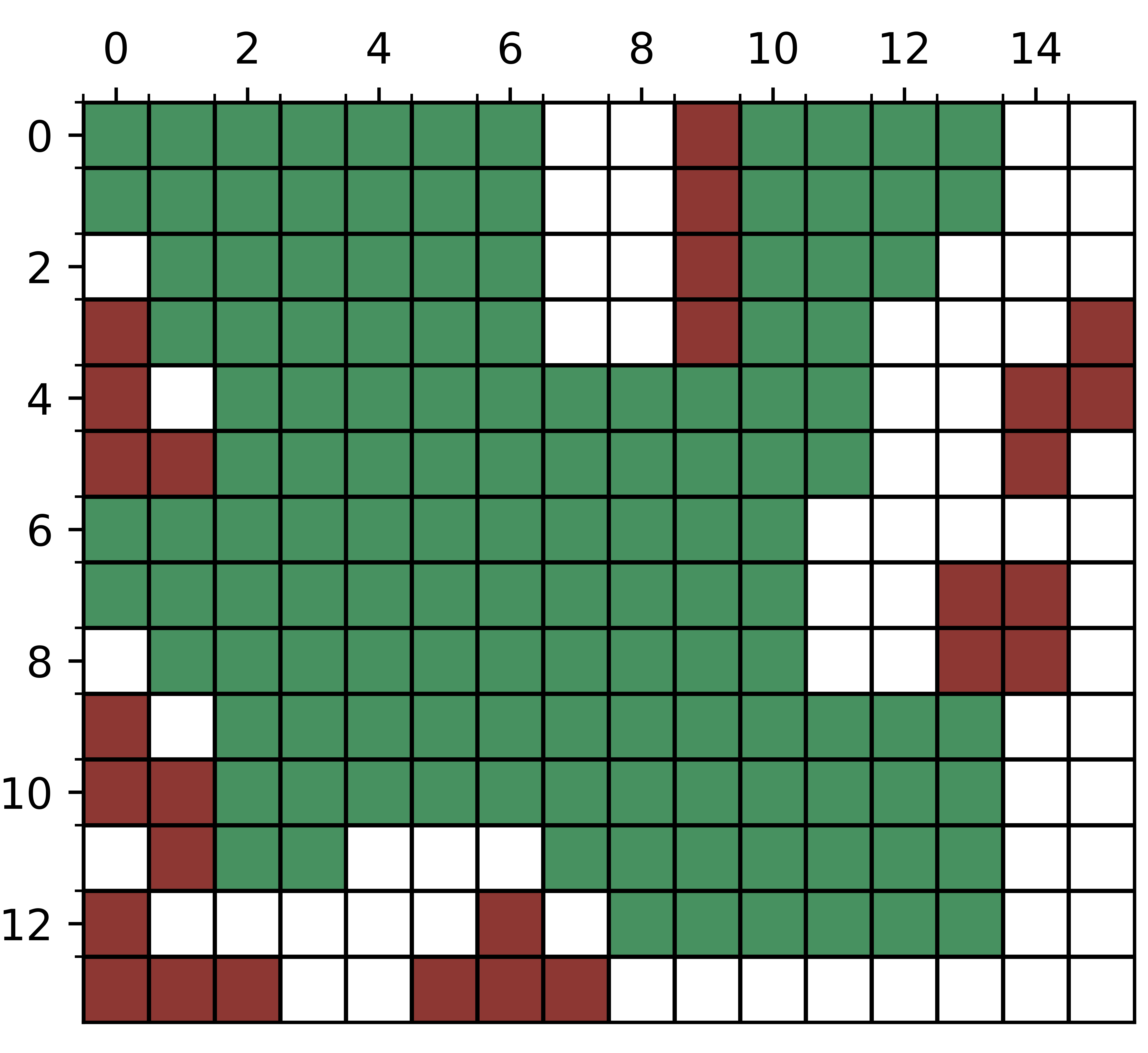}
    \caption{$\validzeroparkset$}
    \label{fig:valid_0deg_fields}
  \end{subfigure}
  \hfill
  \begin{subfigure}[b]{0.24\textwidth}
      \centering
    \includegraphics[width=\textwidth]{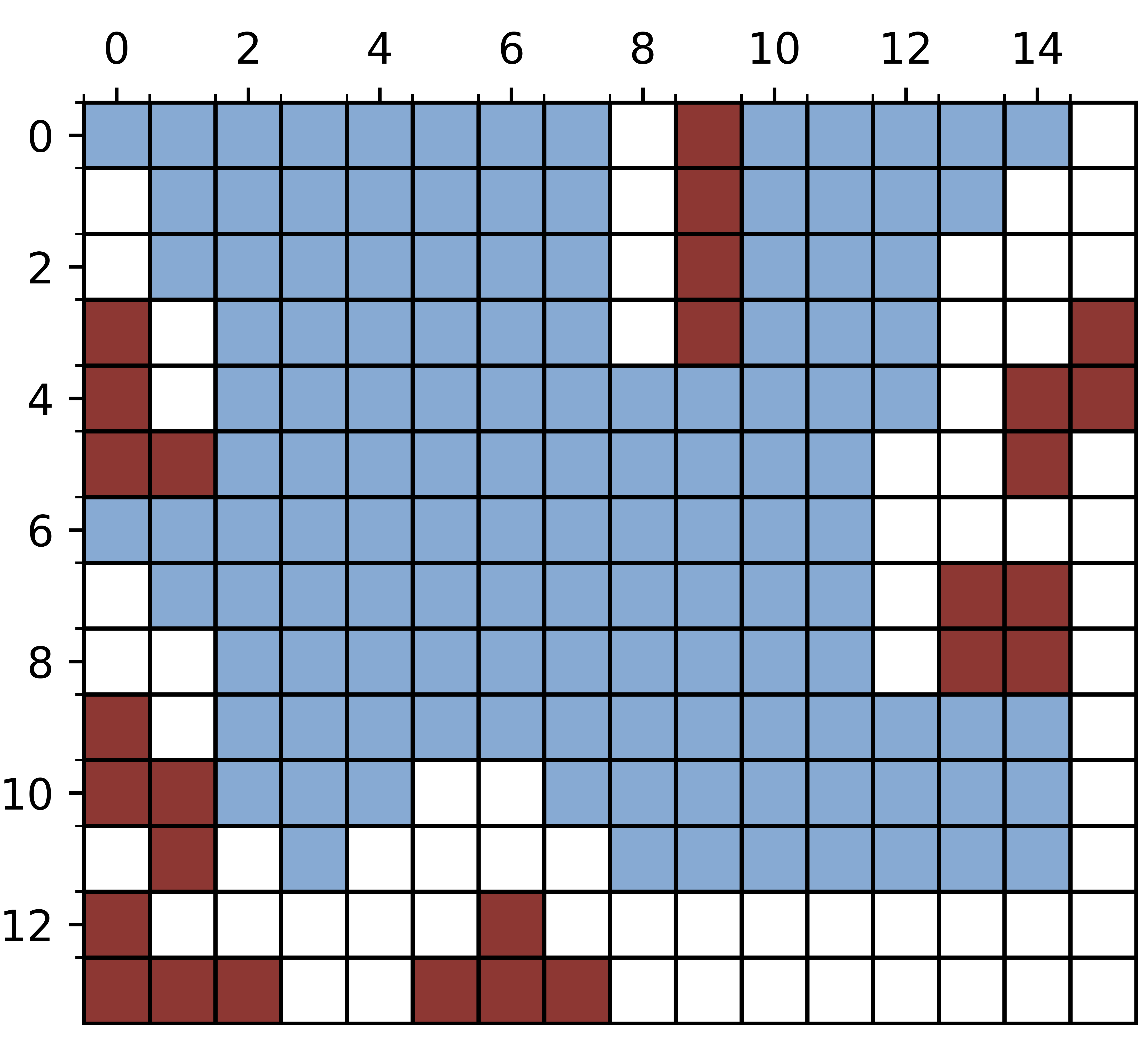}
    \caption{$\validninetyparkset$}
    \label{fig:valid_90deg_fields}
  \end{subfigure}
  \hfill
  \begin{subfigure}[b]{0.24\textwidth}
      \centering
    \includegraphics[width=\textwidth]{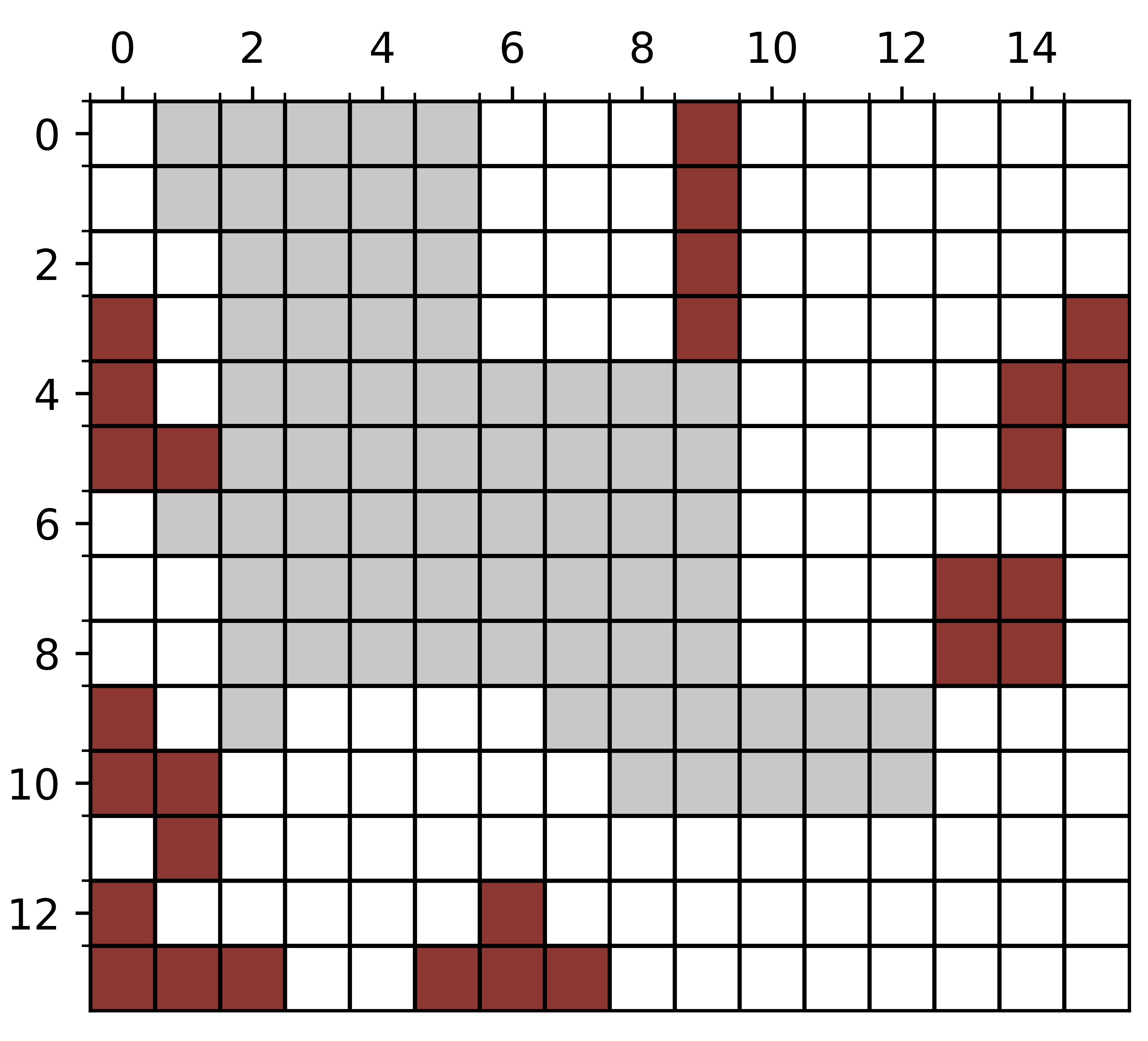}
    \caption{$\validdriveset$}
    \label{fig:valid_driving_fields}
  \end{subfigure}  
  \caption{Decision variables of the integer program. The example here assumes parking fields with $\parkwidth = 2$ and $\parklength=3$ and square driving fields with $\drivewidth=4$. Blocked cells are shown in red. The entrance field is anchored at $(0,5)$.}
  \label{fig:layout}
\end{figure}

A \textit{parking field} is assumed to have a width of $\parkwidth$ cells and a length of $\parklength$ cells and can be oriented along the lateral (east-west) direction or the longitudinal (north-south) direction. We use angles $0^{\circ}$ and $90^{\circ}$ to distinguish between these orientations. The east-west and $0^{\circ}$ nomenclature is used only to indicate the orientation with reference to the grid, but not the driving lanes or geographic directions. They are not to be confused with parallel parking. Parking is assumed to be orthogonal to the driving lanes throughout this work. The binary decision variables that indicate the presence/absence of $0^{\circ}$ and $90^{\circ}$ parking fields \textit{anchored} at cell $(i, j)$ are $\zeroparkvar{i}{j}$ and $\ninetyparkvar{i}{j}$, respectively. We simply refer to them as parking fields at $(i, j)$, although they may contain other cells in the lot depending on $\parkwidth$ and $\parklength$. The choice of the anchor within a parking field is simply a convention and can be set to any other cell in the field. Changing the anchor alters the model implementation but does not affect the formulation or the solution. The MIP formulations can also be written using all parking fields that can be anchored at a specific cell, without separate variables for the $0^{\circ}$ and $90^{\circ}$ parking fields. The variable $\zeroparkvar{i}{j}$ takes a value of one if an east-west parking field is anchored at $(i,j)$ and includes cells of the form $(i + k, j + l)$, where $k \in \{0, \ldots, \parkwidth-1 \}$ and $ l \in \{0, \ldots, \parklength-1\}$ (see the green rectangle in Figure \ref{fig:sample_perpendicular_fields}). Similarly, if $\ninetyparkvar{i}{j}$ is one, then a parking space anchored at $(i, j)$ is created in the north-south direction with cells of the form $(i + k, j + l)$, where $k \in \{0, \ldots, \parklength-1 \}$ and $l \in \{0, \ldots, \parkwidth-1\}$ (see the blue rectangle in Figure \ref{fig:sample_perpendicular_fields}). 

The \textit{driving fields}, on the other hand, are always squares of width $\drivewidth$ cells, where $\drivewidth \geq \parkwidth$ to provide sufficient space for the vehicles to drive and park. A driving field at $(i, j)$ is assumed to contain cells $(i + k, j + l)$, where $k, l \in \{0, \ldots, \drivewidth-1 \}$ and are determined using binary decision variables $y_{ij}$ (e.g., the gray square in Figure \ref{fig:sample_perpendicular_fields}). The decision variable $\drivevar{i}{j}$ is set to one for a driving field anchored at $(i, j)$. For two-way settings, assumed throughout this section, we require $ \drivewidth \geq 2\parkwidth$. Sometimes, parking lots may contain a street accessed by through traffic. We retain such existing driveways while optimizing the layout. The set of anchor points of the existing driving fields is indicated by $\existingdriveset$. 

We assume a single entrance/exit using a driving field anchored at $(\entrancex,\entrancey)$. In a two-way lane configuration, the same cell serves as the exit from the lot and is sometimes simply referred to as the entrance. The entrance cell to the parking lot is assumed to be known/predetermined, and the corresponding $y$ variable value is set to one. The entrance cell can also be inside the lot to model indoor parking garages with ramp access. We identify subsets of $\lotset$ in which $\zeroparkvar{i}{j}$, $\ninetyparkvar{i}{j}$, and $\drivevar{i}{j}$ are valid, i.e., they can potentially take a value of one at these cells, and name them $\validzeroparkset$, $\validninetyparkset$, and $\validdriveset$, respectively (see Figures \ref{fig:valid_0deg_fields}, \ref{fig:valid_90deg_fields}, and \ref{fig:valid_driving_fields}). Note that blocked cells can prohibit other nearby cells from being a parking or a driving field, or can restrict a cell from being connected to the entrance/exit.

\begin{longtable}{ll}
    \caption{Notation for flow-based parking MIP}\label{tab:glosssary}\\
    \hline
    \textbf{Symbol} & \textbf{Description}\\
    \hline
    \endfirsthead
    \hline
    Sets & \\
    \hline
    $\lotset$ & Rasterized cells in the parking lot\\
    $\blockedset$ & Blocked cells\\
    $\existingdriveset$ & Anchor cells of existing drive fields\\
    \vspace{-1.2mm} & \vspace{-1.2mm} \\
    $\validzeroparkset$ & Cells at which $0^\circ$ parking fields are valid\\
    $\validninetyparkset$ & Cells at which $90^{\circ}$ parking fields are valid\\
    $\validdriveset$ & Cells at which driving fields are valid\\
    \vspace{-1.2mm} & \vspace{-1.2mm} \\
    $\zeroneighborset{i}{j}$ & Driving fields that neighbor a $0^\circ$ parking field at $(i,j)$\\
    $\ninetyneighborset{i}{j}$ & Driving fields that neighbor a $90^{\circ}$ parking field at $(i,j)$\\
    \vspace{-1.2mm} & \vspace{-1.2mm} \\
    $\invzeroparkset{i}{j}$ & Anchors of  $0^\circ$ parking fields that contain cell  $(i,j)$\\
    $\invninetyparkset{i}{j}$ & Anchors of $90^{\circ}$ parking fields that contain cell $(i,j)$\\
    $\invdriveset{i}{j}$ & Driving fields that contain cell $(i,j)$\\
    \vspace{-1.2mm} & \vspace{-1.2mm} \\
    $\grid$ & A grid graph of cells in $\validdriveset$ as nodes and arcs connecting cells that share a side \\
    $\gridarcset$ & Arcs connecting the cells in the grid \\
    $\gridadjset{i}{j}$ & Cells at the head end of arcs emanating from $(i,j)$\\
    $\gridinvadjset{i}{j}$ & Cells at the tail end of arcs incoming to $(i,j)$\\
    \hline
    Parameters & \\
    \hline
    $\numrows$ & Number of rows in the rasterized parking lot\\
    $\numcols$ & Number of columns in the rasterized parking lot\\
    $\parkwidth$ & Width of the parking field in terms of the number of cells\\
    $\parklength$ & Length of the parking field in terms of the number of cells\\
    $\drivewidth$ & Width of the driving field in terms of the number of cells\\
    $(\entrancex, \entrancey)$ & Anchor cell of the entrance/exit driving field\\
    \hline
    Variables & \\
    \hline
    $\zeroparkvar{i}{j}$ & Binary $0^\circ$ park variable with cell $(i,j)$ as an anchor\\
    $\ninetyparkvar{i}{j}$ & Binary $90^\circ$ park variable with cell $(i,j)$ as an anchor\\
    $\drivevar{i}{j}$ & Binary drive variable with cell $(i,j)$ as an anchor \\
    $\flowvar{i}{j}{k}{l}$ & Flow from cell $(i,j)$ to cell $(k,l)$ \\
    \hline
\end{longtable}

\cite{stephan2021layout} formulated the problem of designing parking lots with two-way driveways as a MIP using three main types of constraints: (1) \textit{Single-purpose constraints} -- which require each cell to have at most one purpose, i.e., parking, driving, or being left empty, (2) \textit{Accessibility constraints} -- every parking field should be connected to a driving field, and (3) \textit{Driveway connectivity constraints} -- the driving fields must be connected to the entrances/exits. The objective is to maximize the number of parking stalls. They, however, create different models for different resolutions. In the following subsections, we discuss a generalized version.

\subsection{Single-purpose constraints}
\begin{figure}[t]
  \centering
  \begin{subfigure}[b]{0.32\textwidth}
    \centering
    \includegraphics[width=0.88\textwidth]{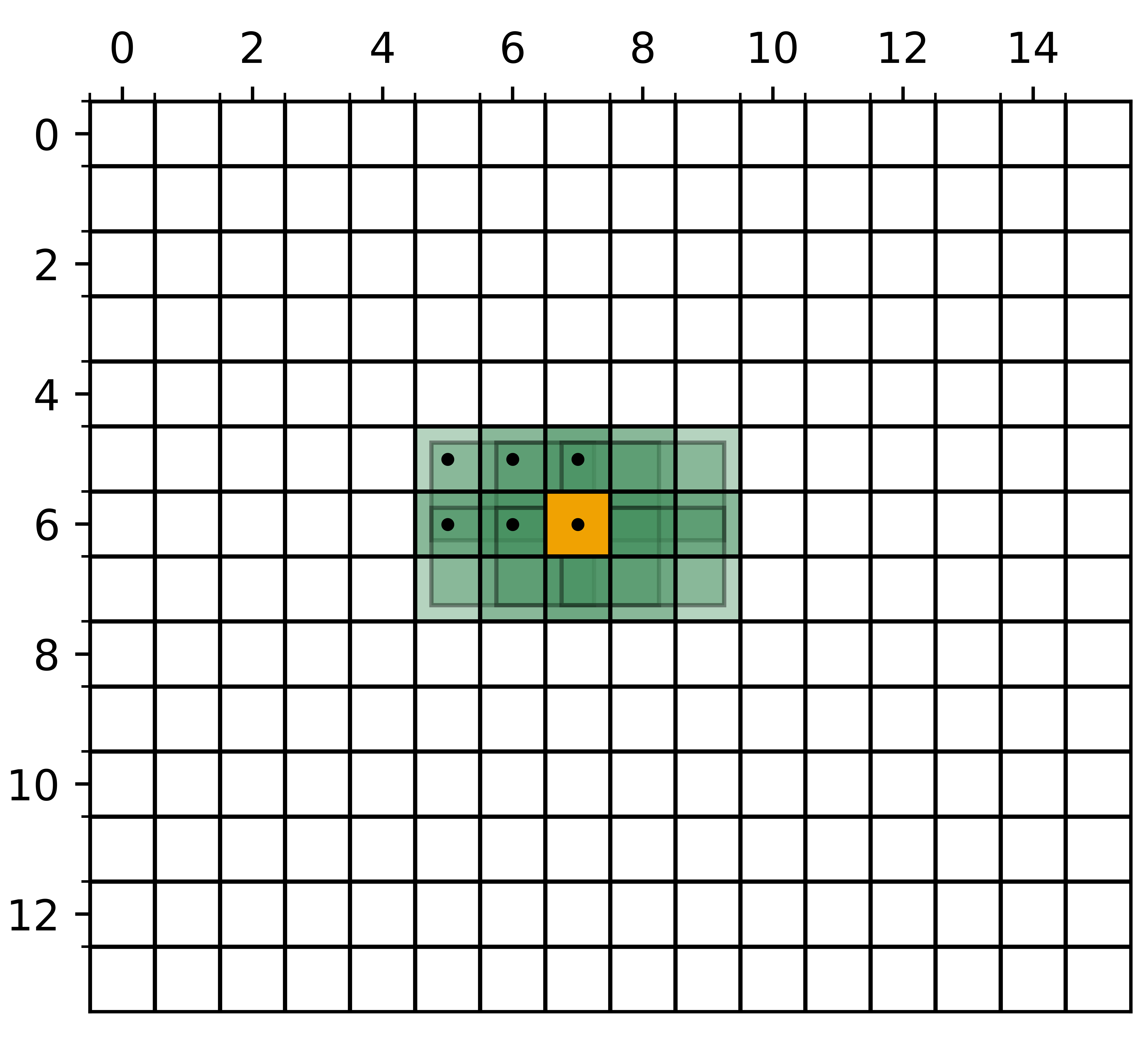}
    \caption{$\invzeroparkset{i}{j}$}
    \label{fig:figure2}
  \end{subfigure}
  \hfill
  \begin{subfigure}[b]{0.32\textwidth}
    \centering
    \includegraphics[width=0.88\textwidth]{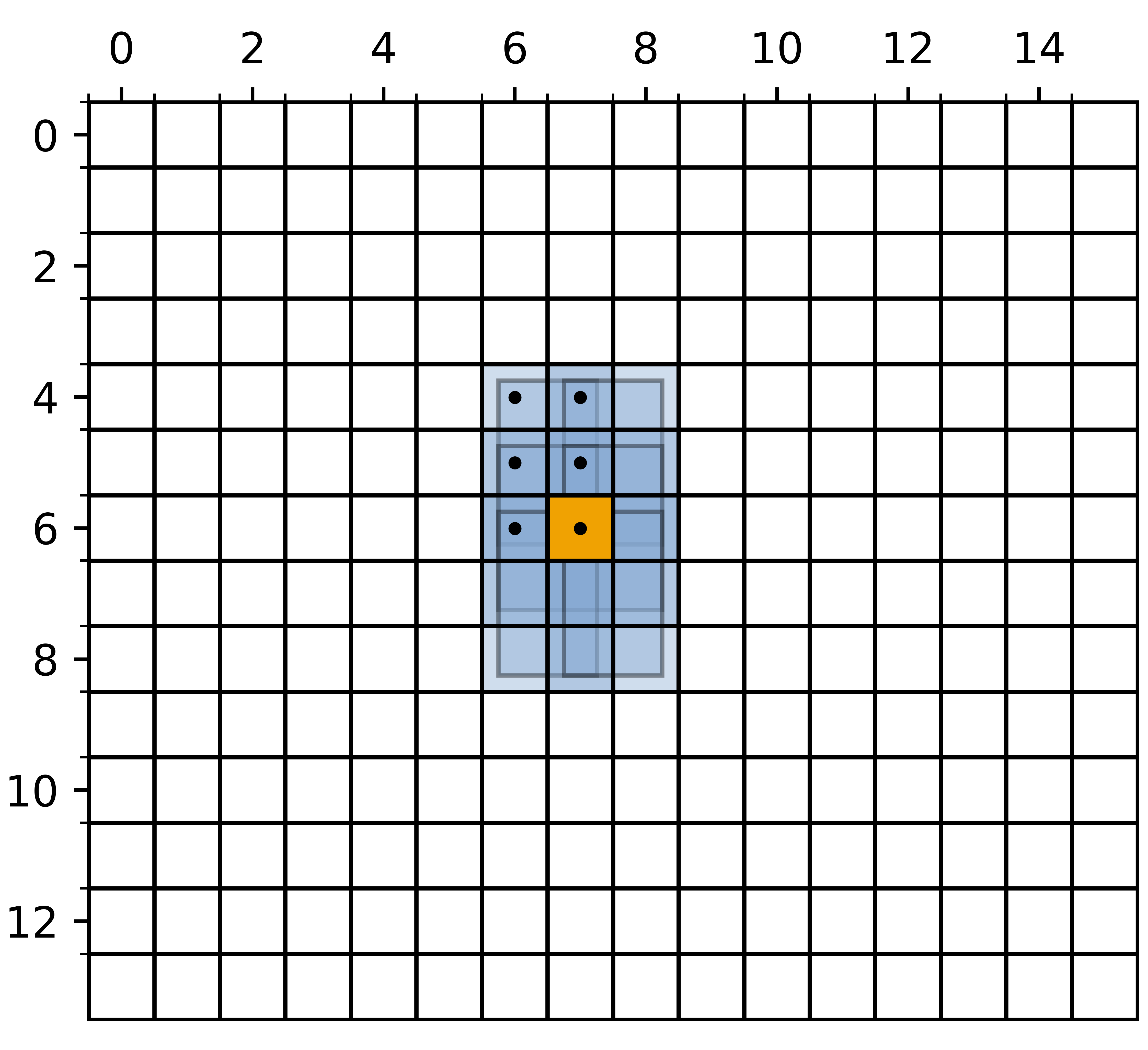}
    \caption{$\invninetyparkset{i}{j}$}
    \label{fig:figure3}
  \end{subfigure}
  \hfill
  \begin{subfigure}[b]{0.32\textwidth}
  \centering
    \includegraphics[width=0.88\textwidth]{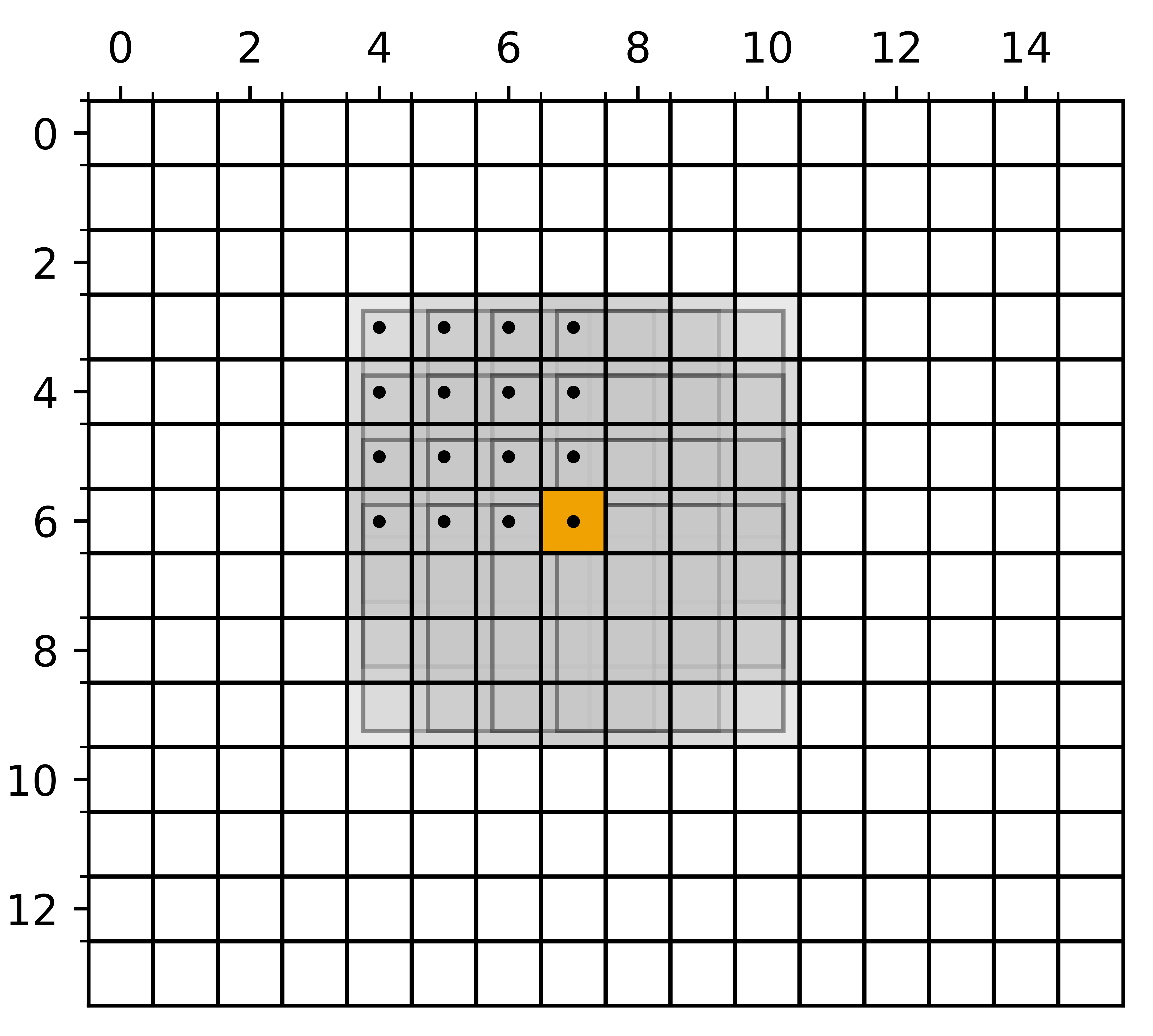}
    \caption{$\invdriveset{i}{j}$}
    \label{fig:figure4}
  \end{subfigure}  
  \caption{Sets of parking and driving fields containing the yellow cell $(6,7)$}
  \label{fig:invcells}
\end{figure}

Let the set of $0^\circ$ and $90^{\circ}$ parking fields that contain a cell $(i,j)$ be represented by $\invzeroparkset{i}{j}$ and $\invninetyparkset{i}{j}$, respectively. Likewise, let $\invdriveset{i}{j}$ denote the set of driving fields that contain cell $(i,j)$. Figure \ref{fig:invcells} illustrates different fields containing a cell $(i,j)$. For example, $\invzeroparkset{(6,}{7)} = \{(5, 5), (5, 6), (5, 7), (6, 5), (6,6), (6, 7) \}$. In the examples discussed in this paper, we sometimes write the subscripts $ij$ as $(i,j)$ to avoid ambiguity between row and column indices. We do not assume any blocked fields in this example. If there are blocked cells, only the valid fields that contain the cell $(i, j)$ are included in $\invzeroparkset{i}{j}$, $\invninetyparkset{i}{j}$, and $\invdriveset{i}{j}$.

Cells in the parking lot that are not blocked can be used for parking or driving or be left unused. Constraints \eqref{eq:singlepurpose_zero} and \eqref{eq:singlepurpose_nonzero} capture this feature. The term $\frac{1}{|\invdriveset{i}{j}|} \sum_{(k,l) \in  \invdriveset{i}{j}} \drivevar{k}{l}$ in the left-hand side allows for overlapping of driving fields as multiple driving fields can contain a cell $(i,j)$.
\begin{align}
& \sum_{(k,l) \in \invzeroparkset{i}{j}} \zeroparkvar{k}{l} + \sum_{(k,l) \in \invninetyparkset{i}{j}} \ninetyparkvar{k}{l} \leq 1  \qquad && \forall\,  (i,j) \in \lotset \setminus \blockedset,\,|\invdriveset{i}{j}| = 0  \label{eq:singlepurpose_zero} \\
& \sum_{(k,l) \in \invzeroparkset{i}{j}} \zeroparkvar{k}{l} + \sum_{(k,l) \in \invninetyparkset{i}{j}} \ninetyparkvar{k}{l} + \frac{1}{|\invdriveset{i}{j}|} \sum_{(k,l) \in  \invdriveset{i}{j}} \drivevar{k}{l} \leq 1  \qquad && \forall\,  (i,j) \in \lotset \setminus \blockedset,\,|\invdriveset{i}{j}| > 0 \label{eq:singlepurpose_nonzero}
\end{align}
Constraint \eqref{eq:singlepurpose_nonzero} can be written in a disaggregate form, which can tighten the LP relaxations of the problem. This new version \eqref{eq:singlepurpose_zero_disaggregate} also captures \eqref{eq:singlepurpose_zero} since the first term on the left-hand side does not contain the $y$ variable if $|\invdriveset{i}{j}| = 0$.
\begin{align}
& \drivevar{m}{n} + \sum_{(k,l) \in \invzeroparkset{i}{j}} \zeroparkvar{k}{l} + \sum_{(k,l) \in \invninetyparkset{i}{j}} \ninetyparkvar{k}{l} \leq 1  \qquad && \forall\,  (m, n) \in  \invdriveset{i}{j}, (i,j) \in \lotset \setminus \blockedset \label{eq:singlepurpose_zero_disaggregate} 
\end{align}

\subsection{Accessibility constraints}
\label{sec:Perpendicular Accessibility constraints}
Each parking field must be accessible from at least one side by a driving field, allowing vehicles to maneuver into the parking stall. Hence, we impose constraints \eqref{eq:zeroparkdriveconnector} and \eqref{eq:ninetyparkdriveconnector}. These constraints ensure that the park field variables can be set to one only if it has at least one active adjacent driving field. We assume that neighbor sets $\zeroneighborset{i}{j}$ and $\ninetyneighborset{i}{j}$ are constructed using fields that entirely overlap with the shorter edges of a parking field at $(i,j)$. Recall that, as previously assumed, the width of the driving fields is equal to or greater than twice the width of the parking fields.
\begin{figure}[H]
  \centering
  \begin{subfigure}[b]{0.32\textwidth}
  \centering
    \includegraphics[width=0.88\textwidth]{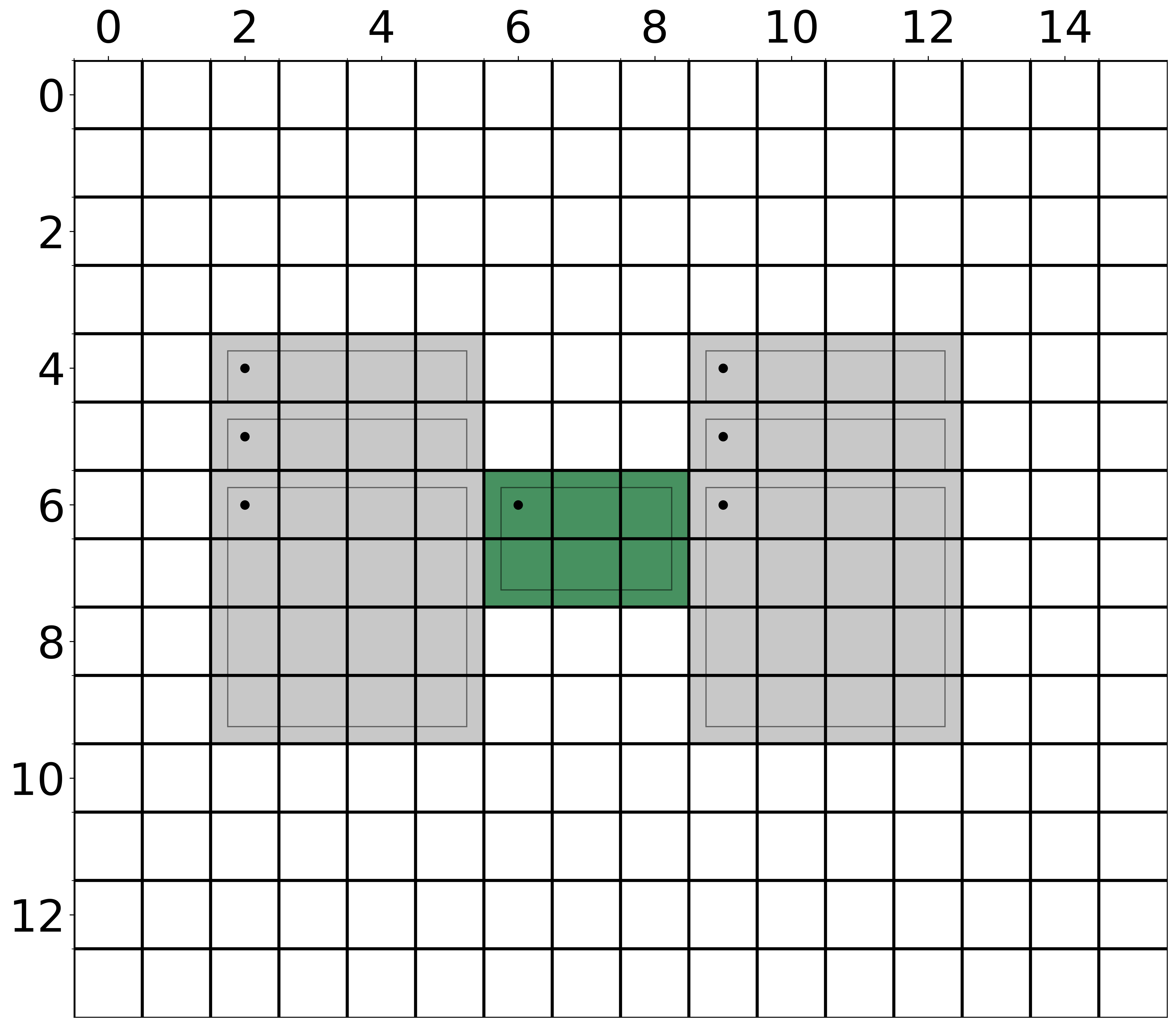}
    \caption{$\zeroneighborset{i}{j}$}
    \label{fig:parkdriveconstraints_a}
  \end{subfigure}
  \hfill
  \begin{subfigure}[b]{0.32\textwidth}
    \centering
    \includegraphics[width=0.88\textwidth]{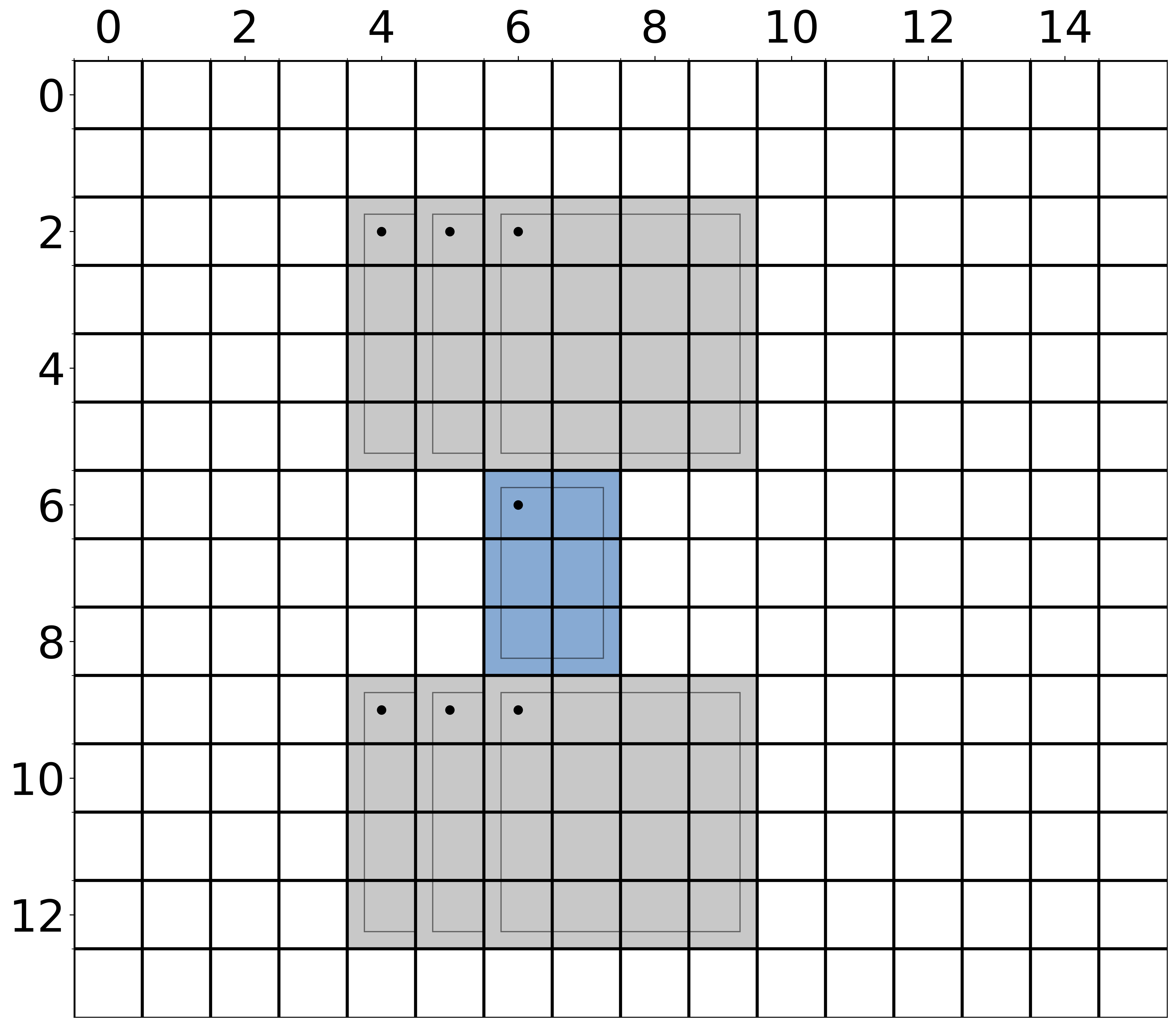}
    \caption{$\ninetyneighborset{i}{j}$}
    \label{fig:parkdriveconstraints_b}
  \end{subfigure}
  \hfill
  \begin{subfigure}[b]{0.32\textwidth}  \centering
    \includegraphics[width=0.88\textwidth]{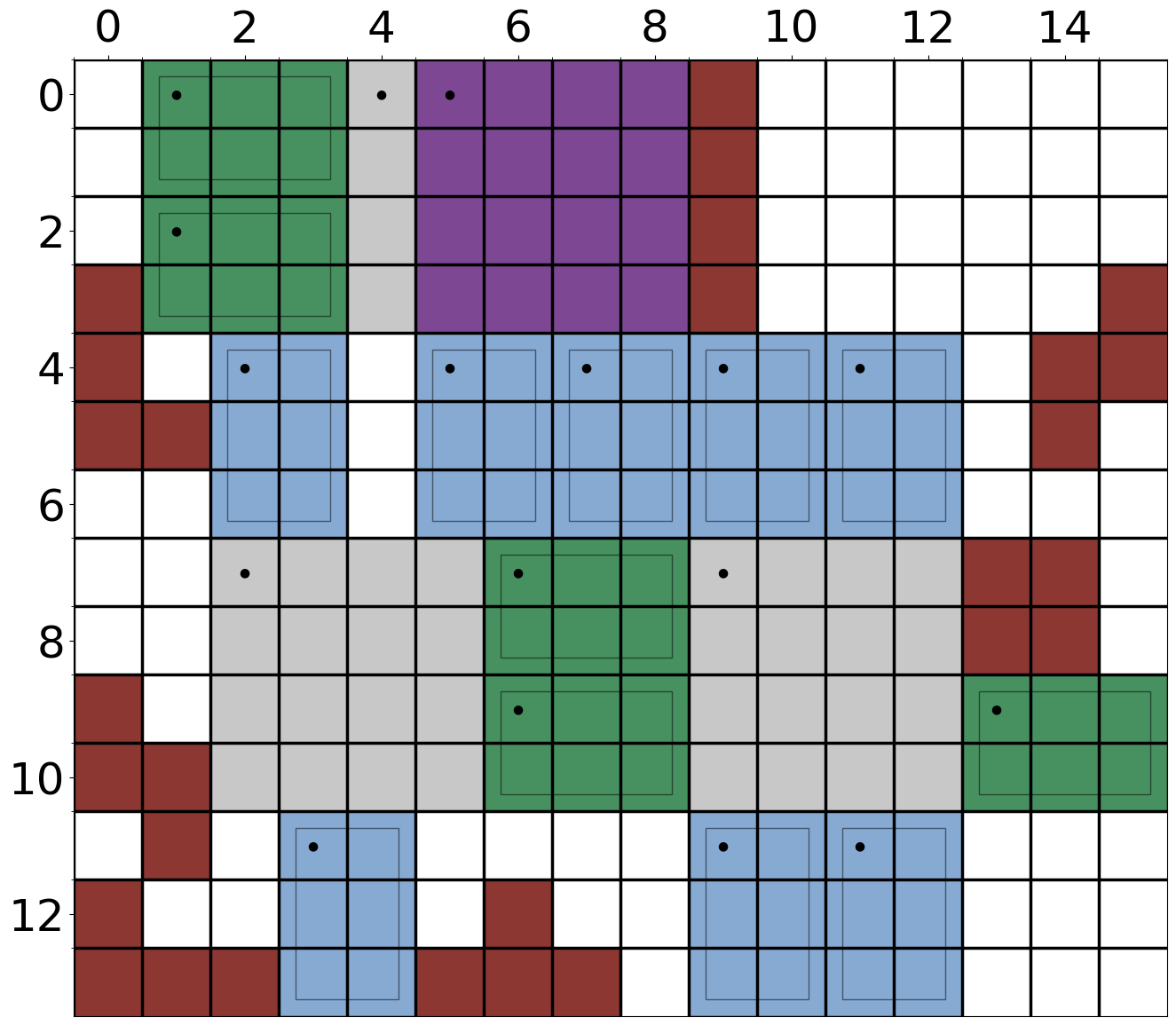}
    \caption{Lot with accessibility constraints}
    \label{fig:parkdriveconstraints_c}
  \end{subfigure}  
\caption{Driving field neighbors of parking field at $(6,6)$ (in left and center panels) and the need for additional constraints on the driving fields (right)}
\label{fig:parkdriveconstraints}
\end{figure}
Figures \ref{fig:parkdriveconstraints_a} and \ref{fig:parkdriveconstraints_b} show the neighbors of two sample parking fields. The driving anchors that allow access to the $0^\circ$ parking field at $(6,6)$ is shown using gray cells $\zeroneighborset{(6,}{6)}=\{ (4,2),(5,2),(6,2),(4,9),(5,9),(6,9) \}$. Similarly, $\ninetyneighborset{(6,}{6)}=\{ (2,4), (2,5), (2,6), (9,4), (9,5), (9,6) \}$ contains driving anchors that help reach the $90^\circ$ parking field at $(6,6)$. 
\begin{align}
    & \zeroparkvar{i}{j} \leq \sum_{(k,l) \in \zeroneighborset{i}{j}} \drivevar{k}{l} && \forall\,(i,j) \in \validzeroparkset \label{eq:zeroparkdriveconnector} \\ 
    & \ninetyparkvar{i}{j} \leq \sum_{(k,l) \in \ninetyneighborset{i}{j}} \drivevar{k}{l} && \forall\,(i,j) \in \validninetyparkset \label{eq:ninetyparkdriveconnector}
\end{align}
However, the above constraints are insufficient since they do not require connecting the driving fields. For example, consider Figure \ref{fig:parkdriveconstraints_c}, which shows a feasible solution that satisfies constraints \eqref{eq:singlepurpose_zero_disaggregate}--\eqref{eq:ninetyparkdriveconnector}. In this example, each parking field is accessible by a driving field, but the driving fields are not connected to the entrance.

\subsection{Driveway connectivity constraints}
\label{sec:Driveway connectivity constraints}
\begin{figure}[H]
  \centering
  \begin{subfigure}[b]{0.48\textwidth}
      \centering
    \includegraphics[width=0.63\textwidth]{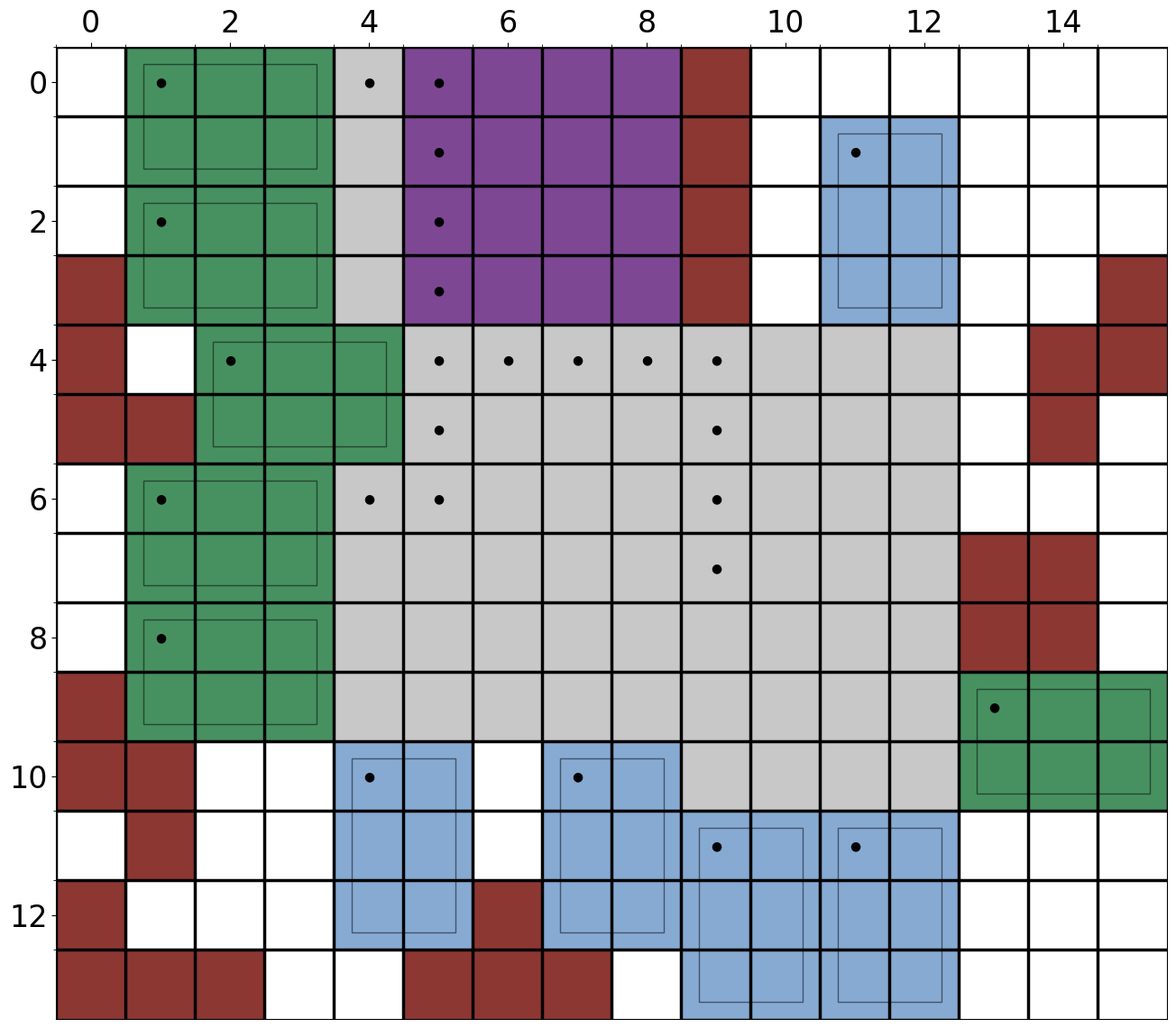}
    \caption{Optimal layout with flow constraints}
    \label{fig:figure3}
  \end{subfigure}
  \hfill
  \begin{subfigure}[b]{0.48\textwidth}
      \centering
    \includegraphics[width=0.58\textwidth]{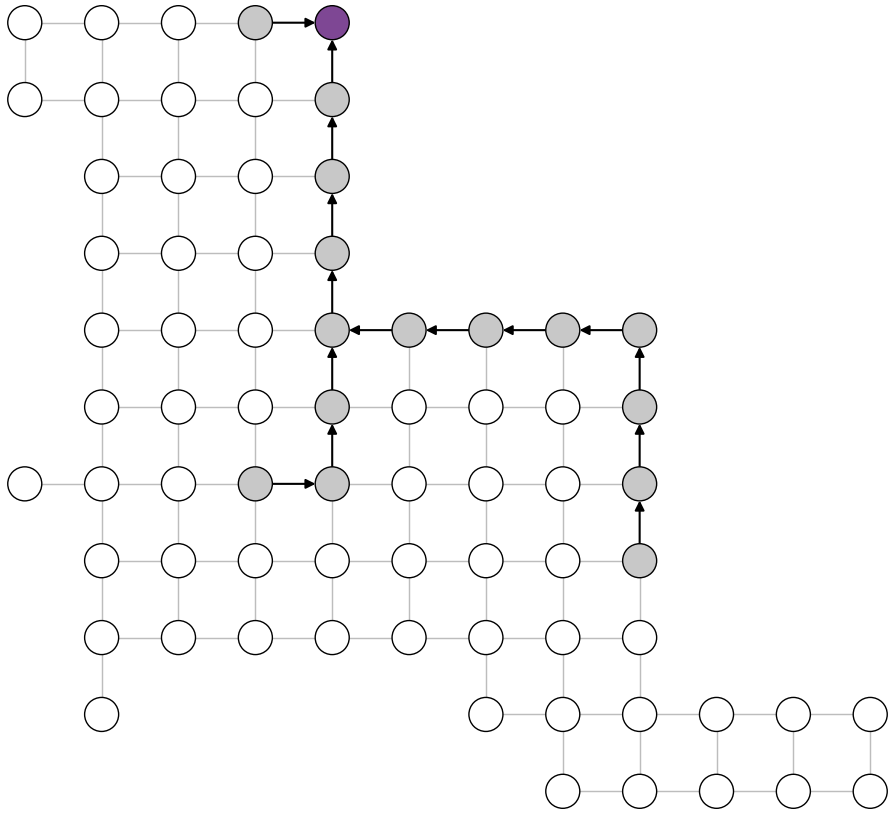}
    \caption{Optimal flow variables}
    \label{fig:figure4}
  \end{subfigure}
\caption{Driving fields must be connected to the entrance}
\label{fig:flowconnectivity}
\end{figure}

To connect the driving fields, we generate a grid-like strongly connected directed graph $\grid = (\validdriveset, \gridarcset)$, using the cells in $\validdriveset$ as nodes and the arcs $\gridarcset$ to connect adjacent cells (see Figure \ref{fig:figure4}). We use $\gridadjset{i}{j}$ to denote the set of cells that are directly reachable from $(i,j)$, and $\gridinvadjset{i}{j}$ to denote the set of cells from which $(i,j)$ is directly reachable. For example, in Figure \ref{fig:figure4}, $\gridadjset{(3,}{4)} = \gridinvadjset{(3,}{4)} = \{(2,4), (3,3), (4,4)\}$. Constraint \eqref{eq:floworigin} generates a unit flow at cells where a driving field is anchored. Flows are routed to the entrance node using constraint \eqref{eq:flowdest}. These flow-conservation constraints ensure that all active driving cells are connected to the entrance location. Additionally, we set the upper bounds on the flow variables using \eqref{eq:flowtail} and \eqref{eq:flowhead} so that positive flows are not generated in and out of a non-driving cell. A safe choice for $M$ is $|\validdriveset|-1$.
\begin{align}
    & \sum_{(k,l) \in \gridadjset{i}{j}} \flowvar{i}{j}{k}{l} - \sum_{(k,l) \in \gridinvadjset{i}{j}} \flowvar{k}{l}{i}{j} = \drivevar{i}{j} \qquad && \forall \, (i,j) \in \validdriveset \setminus \{(\entrancex, \entrancey)\} \label{eq:floworigin}\\
    & \sum_{(i,j) \in \gridinvadjset{\entrancex}{\entrancey}} \flowvar{i}{j}{\entrancex}{\entrancey} = \sum_{(i,j) \in \validdriveset} \drivevar{i}{j} -1 \qquad &&  \label{eq:flowdest}\\
    & \flowvar{i}{j}{k}{l} \leq M \drivevar{i}{j} \qquad && \forall \, ((i,j), (k,l)) \in \gridarcset \label{eq:flowtail}\\
    & \flowvar{i}{j}{k}{l} \leq M \drivevar{k}{l} \qquad && \forall \, ((i,j), (k,l)) \in \gridarcset \label{eq:flowhead}
\end{align}
Figure \ref{fig:flowconnectivity} shows an optimal solution with driveway connectivity constraints \eqref{eq:floworigin}--\eqref{eq:flowhead}. The left panel contains the lot configuration for this instance, and the links with positive flow are shown in the graph on the right.

We can thus write the complete MIP formulation for the parking lot optimization problem as
\begin{align}
   \max & \sum_{(i,j) \in \validzeroparkset} \zeroparkvar{i}{j} + \sum_{(i,j) \in \validninetyparkset} \ninetyparkvar{i}{j}  && [\text{Formulation } \mathcal{F}^{\textsc{flow}}_{\textsc{2W}}] \\
     \text{s.t. } & \eqref{eq:singlepurpose_zero_disaggregate} \text{--}\eqref{eq:flowhead}&& \notag \\
     &  \drivevar{i}{j} = 1 \qquad && \forall \, (i,j) \in   \existingdriveset \cup \entranceset \label{eq:existing drive vars} \\
    & \zeroparkvar{i}{j} \in \{0,1\} \qquad && \forall \, (i,j) \in \validzeroparkset \label{eq: x_zero is binary}\\
    & \ninetyparkvar{i}{j} \in \{0,1\} \qquad && \forall \, (i,j) \in \validninetyparkset \label{eq: x_ninety is binary}\\
    & \drivevar{i}{j} \in \{0,1\} \qquad && \forall \, (i,j) \in \validdriveset \label{eq: y var is binary}\\   
    & \flowvar{i}{j}{k}{l} \geq 0 \qquad && \forall \, ((i,j), (k,l)) \label{eq:flow var is continuous non negative}\in \gridarcset
\end{align}
The flow variables are continuous in \eqref{eq:flow var is continuous non negative}, but they take integral values at the optimal solution. This property follows from the total unimodularity of the coefficient matrix for any integer feasible combination of $x$ and $y$ variables. Note that the LP relaxation of the overall model may still produce fractional values for the $x$ and $y$ variables. All arc flows emerging from the entrance cell are zeros, i.e., $\flowvar{\entrancex}{\entrancey}{i}{j} = 0 \, \forall \, (i,j) \in \gridadjset{\entrancex}{\entrancey} \label{eq:outflow from entry cell is zero}$. One can verify this by adding \eqref{eq:floworigin} for all cells $(i,j) \in \validdriveset \setminus \{(\entrancex, \entrancey)\}$, and subtracting \eqref{eq:flowdest} from it. We also formulated and tested a multi-commodity version in which arc flows are disaggregated by valid driving cells. However, this approach significantly increased the number of variables and constraints, making it considerably slower to solve. We omit the subscripts and superscripts and write $(x, y, f)$ to denote the solution vector.

\section{One-way parking lot configuration}
\label{sec:one_way}
Many parking lots are designed with one-way driving lanes to help accommodate more vehicles. This section considers such a variant and describes the model for a single entrance $(\entrancex,\entrancey)$ and a single exit, denoted as $(\exitx,\exity)$. The driveways should allow vehicles entering to access all parking fields and enable them to return to the exit. Since driveways are one-way, we require $\drivewidth \geq \parkwidth$. 

\begin{figure}[H]
\centering
\begin{subfigure}[b]{0.49\textwidth}
    \centering
    \includegraphics[width=0.49\textwidth]{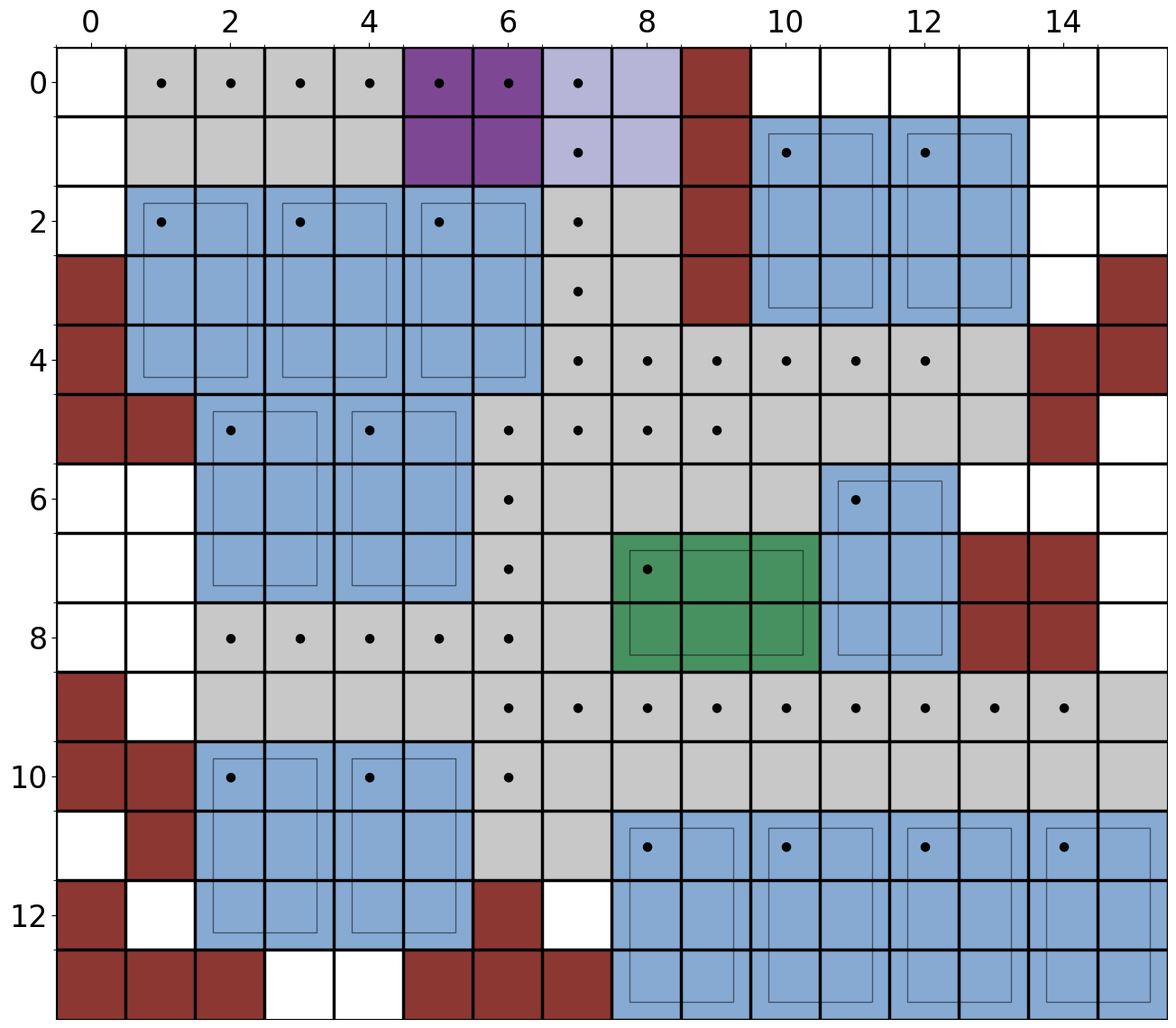}
    \includegraphics[width=0.49\textwidth]{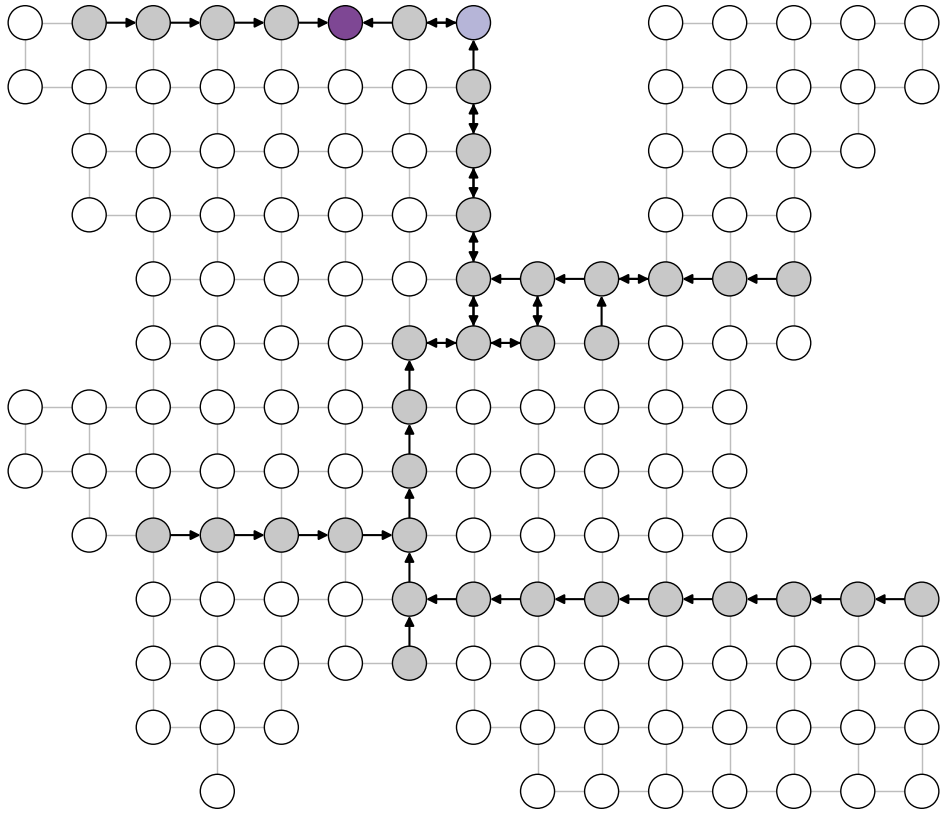}
    \caption{Lot and flow solution with only $f$ variables}
    \label{fig:f_only_constraints}
\end{subfigure} 
\begin{subfigure}[b]{0.49\textwidth}
\centering
\includegraphics[width=0.49\linewidth]{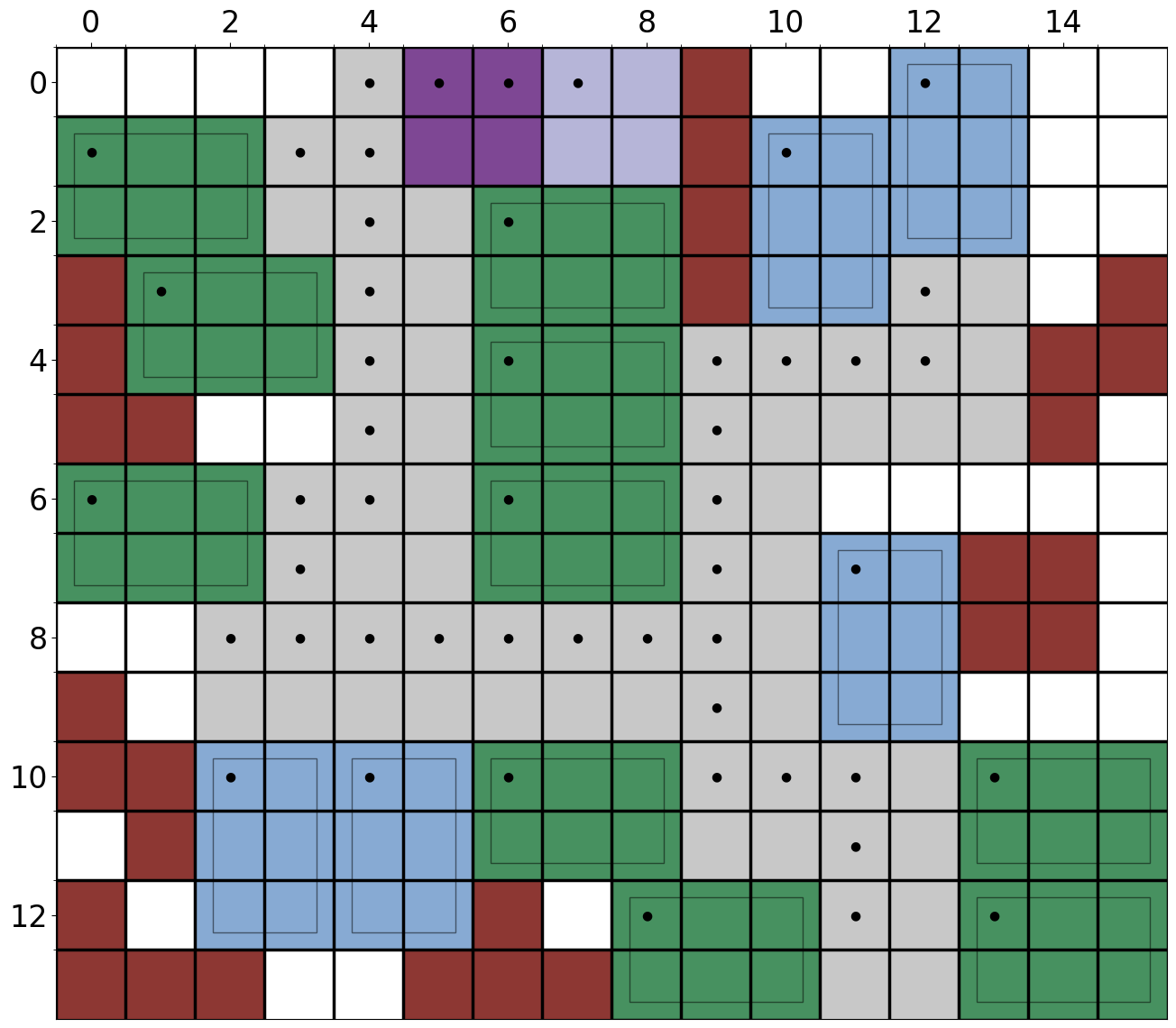}
\includegraphics[width=0.49\textwidth]{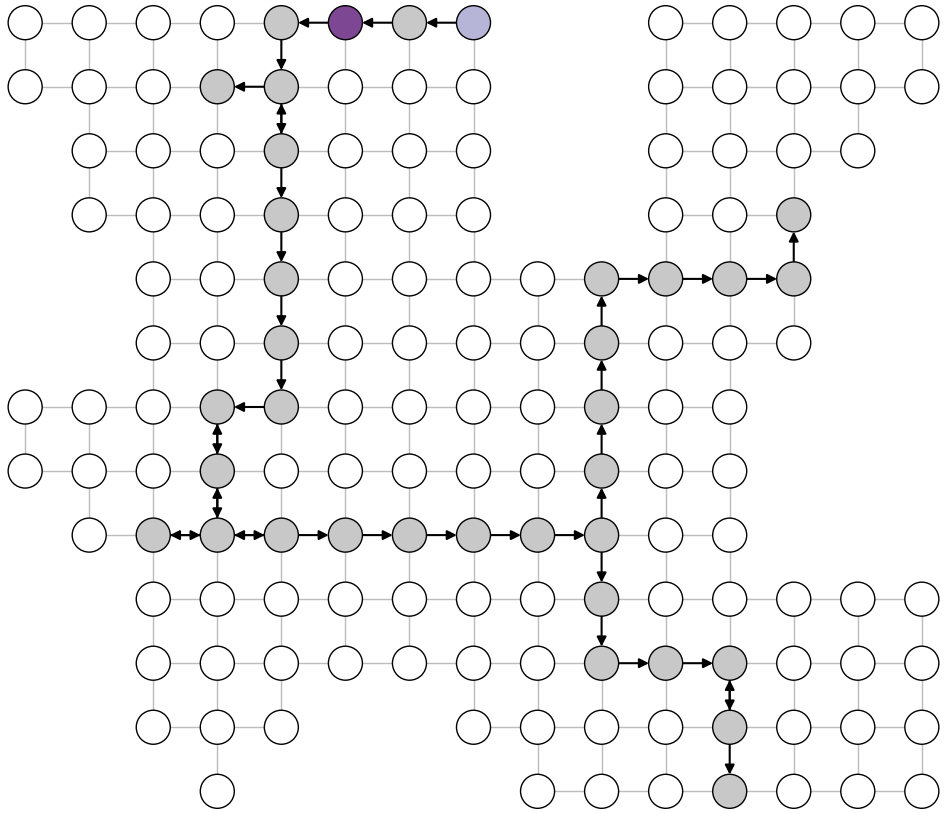}
\caption{Lot and flow solution with only $g$ variables}
\label{fig:g_only_constraints}
\end{subfigure} 
   \caption{Modifications to driveway connectivity constraints for the single-lane one-way configuration. The entrance and exit cells are located at (0, 5) and (0, 7), respectively, and are represented by dark and light purple nodes.}
\label{fig:one_flow_one_way}
\end{figure}

In Section \ref{sec:Driveway connectivity constraints}, every driveway anchor cell except the entrance generates one unit of flow in \eqref{eq:floworigin}, which is guided to the entrance cell using \eqref{eq:flowdest}. The entrance and exit locations are adjacent, and the driveways are two-way streets. However, imposing these constraints for the one-way case, by only changing the width of the driveway, generates layouts that reach the entrance and exit but are not truly one-way streets, as shown in Figure \ref{fig:f_only_constraints}. The $f$ variables ensure that all active driving cells are connected to the entrance, but not all are reachable from the exit. Note that the grid of valid driving cells differs from the previous grid since the driving fields are smaller.  

\subsection{Exit connectivity constraints}
The previously described flow formulation can be viewed as a \textit{multi-origin single-destination} problem, where flows originate at the driving cells and are directed toward the entrance cell, which serves as the destination. To connect parking fields to the exit, we additionally use \textit{single-origin multi-destination} flow variables $\flowvardash{i}{j}{k}{l}$ for every arc $((i,j),(k,l)) \in \gridarcset$, representing the flow of another type/commodity. The exit cell sends out a flow equal to one unit less than the total number of active driving field anchors in the grid, and using the $g$ variables, one unit is delivered at each driving cell as shown in \eqref{eq:flowdashdest_singlelane}. Every active driving cell, except the exit cell, draws one unit of this type of flow according to \eqref{eq:flow_dash_origin_singlelane}. In the two-way scenario, we follow the convention that flows originate from the parking fields and are directed toward the entrance/exit, which is consistent with \cite{stephan2021layout}. However, this necessitates reversing the interpretation of the $g$ variables, which are treated as flows from the exit rather than to the exit. Note that in one-way layouts, the roles of the entrance and exit, as well as the driving lanes, may be interchanged without affecting the positions of parking fields in the solution. Constraints \eqref{eq:flow_dash_tail_singlelane} and \eqref{eq:flow_dash_head_singlelane} set bounds on flow variables similar to \eqref{eq:flowtail} and \eqref{eq:flowhead}, ensuring that only active driving fields are connected to the exit. Finally, \eqref{eq:g_non_negative} restricts the flows to be non-negative.
\begin{align}
& \sum_{(i,j) \in \gridadjset{\exitx}{\exity}} \flowvardash{\exitx}{\exity}{i}{j} = \sum_{(i,j) \in \validdriveset} \drivevar{i}{j}  - 1 \qquad && \label{eq:flowdashdest_singlelane}\\
    & \sum_{(k,l) \in \gridinvadjset{i}{j}} \flowvardash{k}{l}{i}{j} - \sum_{(k,l) \in \gridadjset{i}{j}} \flowvardash{i}{j}{k}{l} = \drivevar{i}{j}  \qquad && \forall \, (i,j) \in \validdriveset \setminus \exitset \label{eq:flow_dash_origin_singlelane}\\
    & \flowvardash{i}{j}{k}{l} \leq M \drivevar{i}{j}  \qquad && \forall \, ((i,j), (k,l)) \in \gridarcset \label{eq:flow_dash_tail_singlelane}\\
    & \flowvardash{i}{j}{k}{l} \leq M \drivevar{k}{l}  \qquad && \forall \, ((i,j), (k,l)) \in \gridarcset \label{eq:flow_dash_head_singlelane}\\
    & \flowvardash{i}{j}{k}{l} \geq 0 \qquad && \forall \, ((i,j), (k,l)) \label{eq:g_non_negative}\in \gridarcset
\end{align}

\begin{figure}[H]
    \centering
    \begin{subfigure}[b]{0.32\textwidth}
        \centering
        \includegraphics[width=0.88\textwidth]{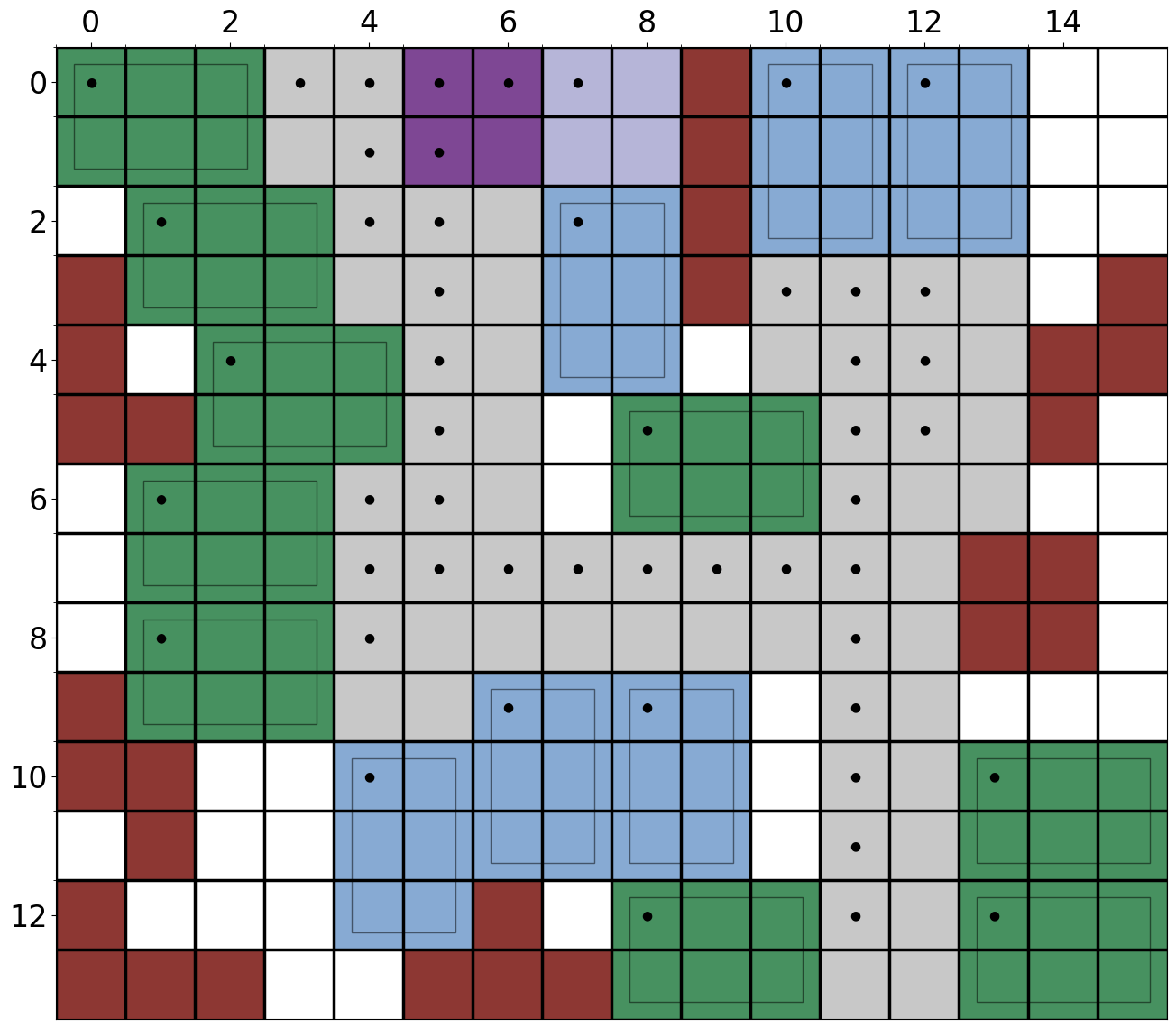}
        \caption{Lot with $f$ and $g$ variables}
        \label{fig:lot_joint}
    \end{subfigure} 
    \hfill
    \begin{subfigure}[b]{0.32\textwidth}
        \centering
        \includegraphics[width=0.88\textwidth]{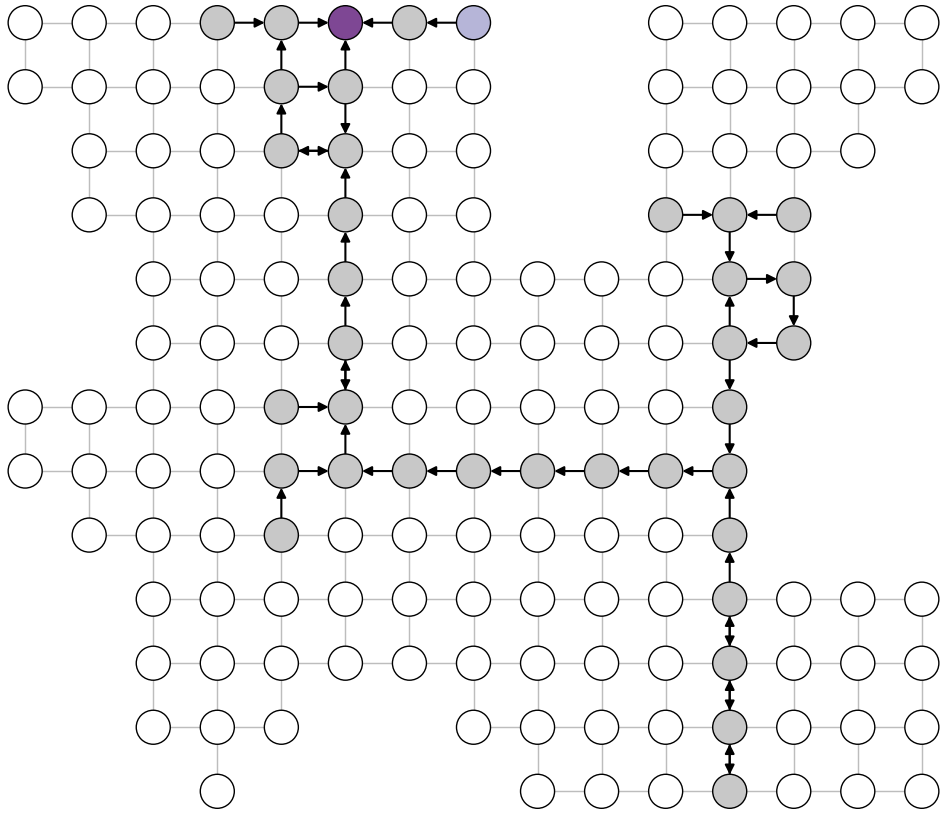}
        \caption{Flow solution $f$}
        \label{fig:f_only_joint}
    \end{subfigure}     
    \begin{subfigure}[b]{0.32\textwidth}
    \centering
    \includegraphics[width=0.88\linewidth]{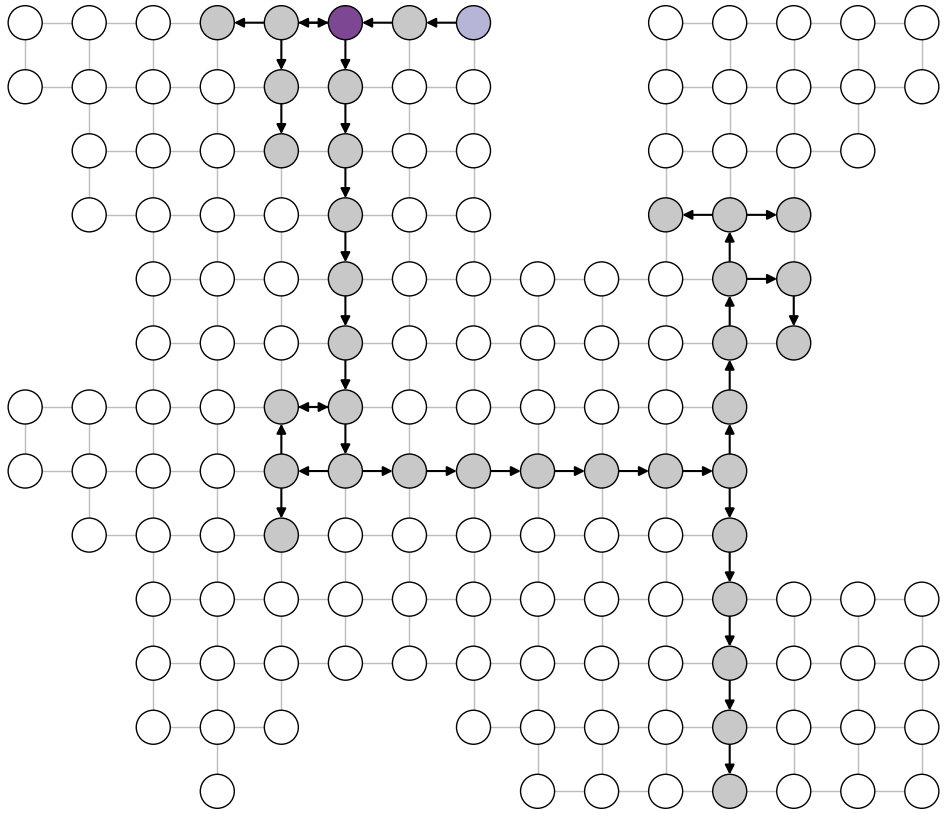}
    \caption{Flow solution $g$}  
    \label{fig:g_only_joint}
    \end{subfigure}  
       \caption{Parking lot and flow solutions with both $f$ and $g$ variables}
    \label{fig:joint_f_g}
\end{figure}

A formulation with only the $g$ variables will connect the entrances and exits but may not produce one-way layouts, as seen from the parking lot and the flow solution in Figure \ref{fig:g_only_constraints}. Since we need connectivity from the entrance to the parking stalls and back to the exit, we include both types of driveway connectivity constraints involving $f$ and $g$ variables. However, this modification does not still solve the problem, as seen in the outcome shown in Figure \ref{fig:joint_f_g}. 

Note that the total outgoing flows from a valid driving cell $(i, j)$ that is not the entrance or exit equals the sum of the incoming flows. Mathematically, 
\begin{align}
    \sum_{(k,l) \in \gridadjset{i}{j}} \big (\flowvar{i}{j}{k}{l}+\flowvardash{i}{j}{k}{l}\big) = \sum_{(k,l) \in \gridinvadjset{i}{j}} \big (\flowvar{k}{l}{i}{j}+\flowvardash{k}{l}{i}{j} \big) \qquad \forall \, (i,j) \in \validdriveset \setminus \{(\entrancex, \entrancey), (\exitx, \exity)\} \label{eq:inflow_outflow_oneway}
\end{align}
This is a direct consequence of \eqref{eq:floworigin} and \eqref{eq:flow_dash_origin_singlelane}. Apart from this equality, no constraints link $f$ and $g$, particularly those that induce one-way restrictions.

\subsection{Directionality constraints and dead ends}
\label{sec:dead_end}
In the two-way case, one can take a U-turn and move in the reverse direction at dead ends. However, this feature must be restricted in the case of one-way driveways. In the flow solutions in Figures \ref{fig:one_flow_one_way} and \ref{fig:joint_f_g}, many dead ends exist, and between several pairs of adjacent driving nodes, there are arc flows in both directions to ensure connectivity to the entrance or exit. Having arc flows in both directions poses ambiguities in fixing driveway directions and leads to infeasible one-way layouts. To navigate such a layout, vehicles must use the same driveway link in both directions, requiring U-turns, which is impossible with one-way driveways. To prevent such solutions, we define new binary \textit{driveway direction} decision variables $\countvar{i}{j}{k}{l}$ for each arc $((i,j),(k,l)) \in \gridarcset$ that keep track whether the arc carries non-zero flows and sets the direction of vehicular movement.
Bounds on $\countvar{i}{j}{k}{l}$ with respect to commodity flows can be set as shown in \eqref{eq:count_var_bounds} and \eqref{eq:count_dash_var_bounds}. 
\begin{align} 
    & \countvar{i}{j}{k}{l} \leq \flowvar{i}{j}{k}{l} \leq  M \countvar{i}{j}{k}{l} \qquad && \forall \, ((i,j), (k,l)) \in \gridarcset  \label{eq:count_var_bounds}\\
    & \countvar{i}{j}{k}{l} \leq \flowvardash{i}{j}{k}{l} \leq M \countvar{i}{j}{k}{l} \qquad && \forall \, ((i,j), (k,l)) \in \gridarcset   \label{eq:count_dash_var_bounds}     
\end{align}
Constraint \eqref{eq:count_var_bounds} ensures that $\countvar{i}{j}{k}{l}$ is zero if $\flowvar{i}{j}{k}{l}$ is zero, and $\countvar{i}{j}{k}{l}$ is one if $\flowvar{i}{j}{k}{l} > 0$. Similar relations hold for the $g$ variables. We can observe from \eqref{eq:count_var_bounds} and \eqref{eq:count_dash_var_bounds} that it is not possible to have arcs with positive $f$ values and zero $g$ values, and vice versa. Thus, both $f$ and $g$ variables on an arc are positive (and consequently, the $z$ variable is one), or all three variables are zero. The $\countvar{i}{j}{k}{l}$ variable cannot be one if it is not connected to active driving fields at either end. Hence, we can relate the $z$ and $y$ variables using \eqref{eq:count_drive1} and \eqref{eq:count_drive2}. As before, we set $M$ to $|\validdriveset|-1$. Note that these constraints along with \eqref{eq:count_var_bounds} and \eqref{eq:count_dash_var_bounds} render bounds on $f$ and $g$ variables \eqref{eq:flowtail}, \eqref{eq:flowhead}, \eqref{eq:flow_dash_tail_singlelane}, and \eqref{eq:flow_dash_head_singlelane} redundant. 
\begin{align}
& \countvar{i}{j}{k}{l} \leq \drivevar{i}{j}  \qquad && \forall \, ((i,j), (k,l)) \in \gridarcset \label{eq:count_drive1}\\
    & \countvar{i}{j}{k}{l} \leq \drivevar{k}{l}  \qquad && \forall \, ((i,j), (k,l)) \in \gridarcset \label{eq:count_drive2}
\end{align}
As seen in \eqref{eq:inflow_outflow_oneway}, if a driving cell that is not the entrance or the exit has inflow (with respect to both $f$ and $g$ variables), it must also have some outflow. To enforce one-way restrictions and avoid inflows and outflows between a given set of adjacent driving cells, we ensure that flow can occur only along one of the two possible arcs connecting them. 
\begin{align}
& \countvar{i}{j}{k}{l} + \countvar{k}{l}{i}{j} \leq 1 \qquad && \forall \, ((i,j), (k,l)) \in \gridarcset \label{eq:flow_one_direction}
\end{align}
Due to the above constraint, the driving layouts will also not contain dead ends or cul-de-sacs. In case there exists a driving cell that is a leaf node or a dead end, connectivity to both the entrance and exit requires an outgoing and an incoming arc, which would violate \eqref{eq:flow_one_direction}.

\begin{figure}[H]
    \centering
    \begin{subfigure}[b]{0.49\textwidth}
        \centering
        \includegraphics[width=0.63\textwidth]{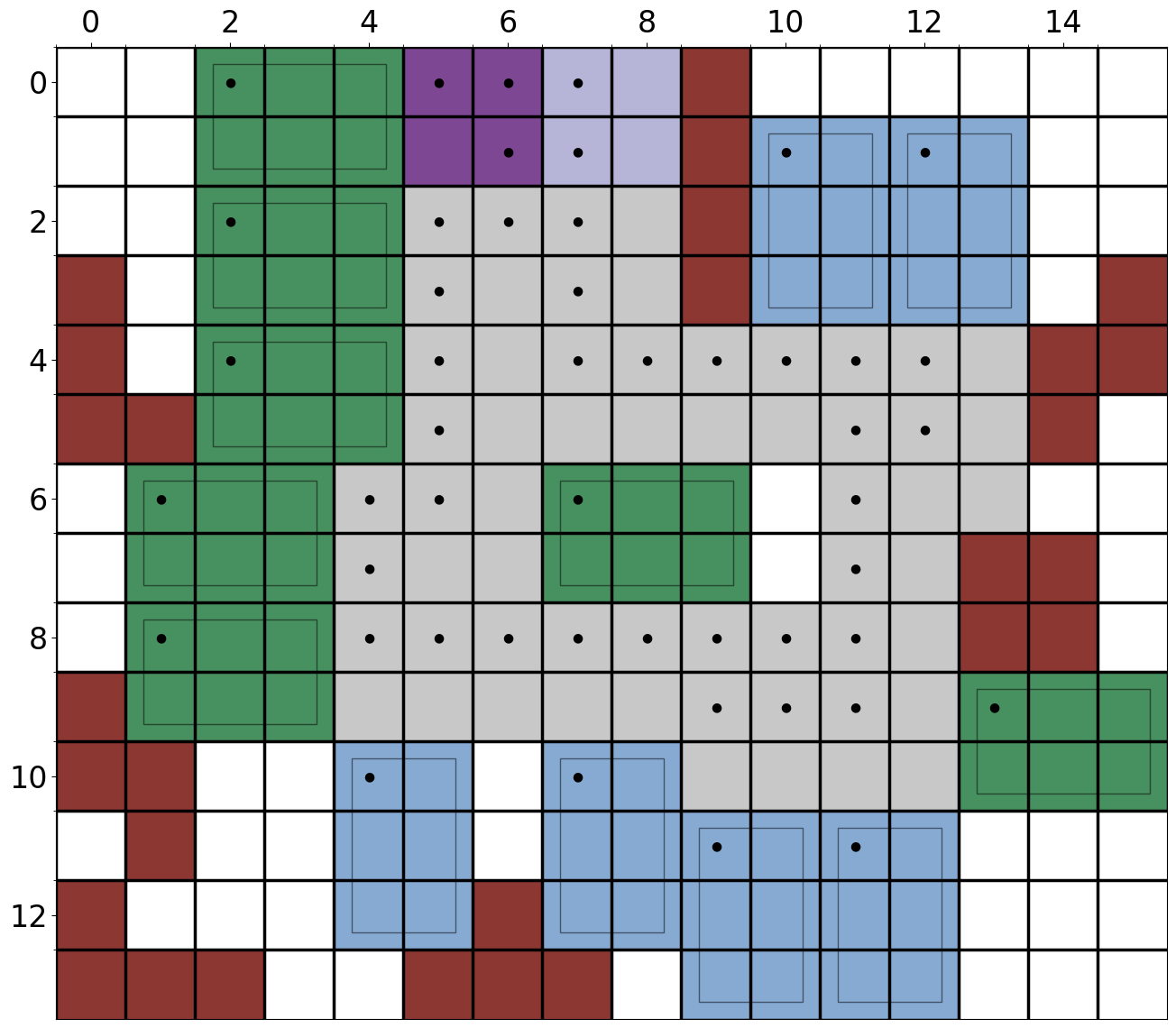}
        \caption{Parking lot layout}
    \end{subfigure} 
    \hfill
    \begin{subfigure}[b]{0.49\textwidth}
        \centering
        \includegraphics[width=0.68\textwidth]{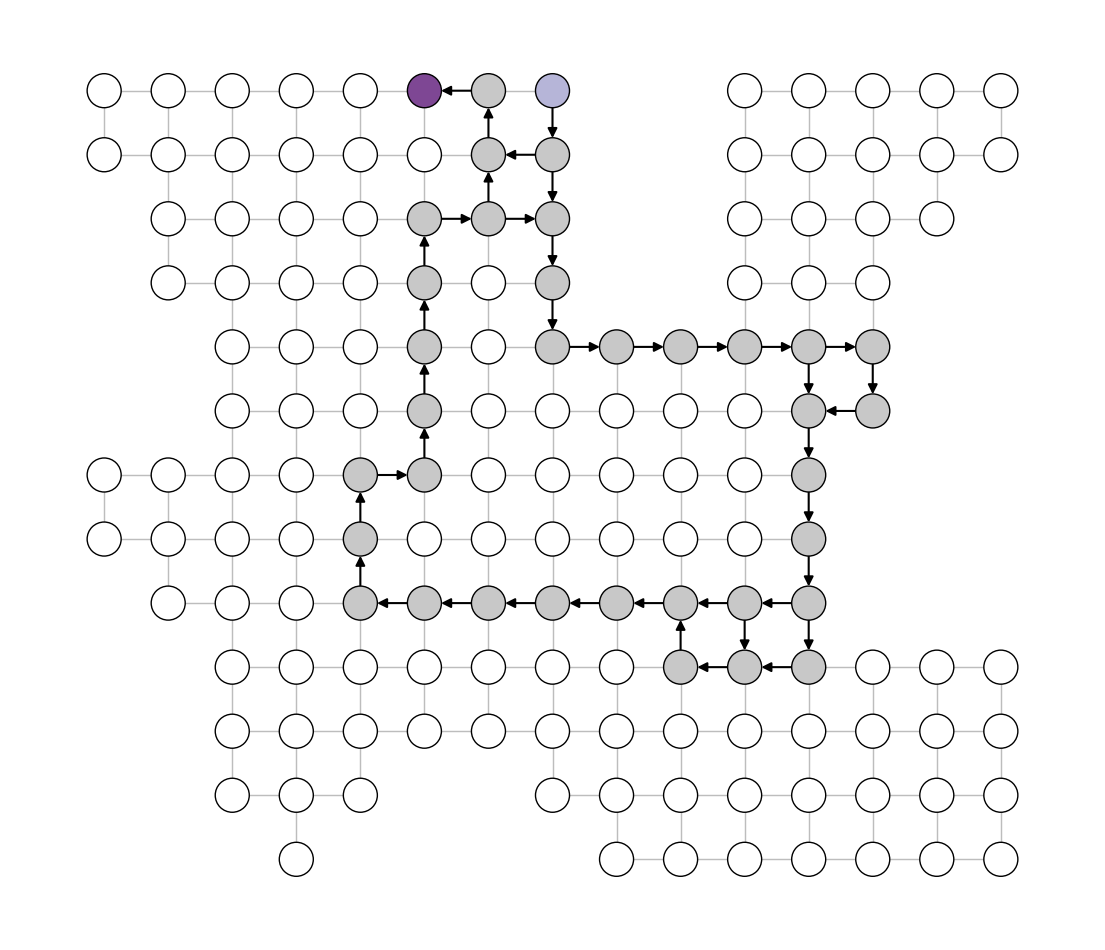}
        \caption{Driveway direction variables}
    \end{subfigure}     
       \caption{Optimal solution to a one-way lane configuration}
    \label{fig:optimal singlelane}
\end{figure}

Similar to the two-way case, adding \eqref{eq:flow_dash_origin_singlelane} for all $(i,j) \in \validdriveset \setminus \exitset$ and comparing it with \eqref{eq:flowdashdest_singlelane} implies 
$\flowvardash{i}{j}{\exitx}{\exity} = 0~ \forall~(i,j) \in \gridinvadjset{\exitx}{\exity}$. Since the two-commodity flows can be either zero or positive on any arc, $\flowvar{i}{j}{\exitx}{\exity} = 0~\forall ~(i,j) \in \gridinvadjset{\exitx}{\exity}$ and $\flowvardash{\entrancex}{\entrancey}{i}{j} = 0~\forall~(i,j) \in \gridadjset{\entrancex}{\entrancey}$. In other words, there are no outgoing flows from the entrance and no incoming flows to the exit. The complete formulation for the one-way driveway case can thus be written as follows. Figure \ref{fig:optimal singlelane} illustrates an optimal solution for this formulation.
\begin{align}
    \max &  \sum_{(i,j) \in \validzeroparkset} \zeroparkvar{i}{j} + \sum_{(i,j) \in \validninetyparkset} \ninetyparkvar{i}{j}  && [\text{Formulation } \mathcal{F}^{\textsc{flow}}_{\textsc{1W}}] \notag\\
     \text{s.t. } &  \eqref{eq:singlepurpose_zero_disaggregate}\text{--}\eqref{eq:flowdest}, 
     \eqref{eq: x_zero is binary}\text{--} \eqref{eq:flow_dash_origin_singlelane}, \eqref{eq:g_non_negative}, 
     \eqref{eq:count_var_bounds} \text{--}\eqref{eq:flow_one_direction}&& \notag \\
     &  \drivevar{i}{j} = 1 \qquad && \forall\,(i,j) \in  \existingdriveset \cup \{(\entrancex, \entrancey), (\exitx, \exity)\} \label{eq:entrance_exit_driveways}\\ 
    & \countvar{i}{j}{k}{l} \in \{0,1\} \qquad && \forall \, ((i,j), (k,l)) \in \gridarcset   \label{eq:binary_countvar}  
\end{align}

Notice that while the lot appears navigable by vehicles of dimensions $\parkwidth =2$ and $\parklength=3$, the orientation of the direction variables at few locations, e.g., at cells (8,9), (8,10), (9,9), and (9,10), would force sharp turns that are impossible to execute (unless the problem setting involves multiple smaller vehicles packed within each field or robots in a warehouse that can slide laterally). Turning can also be an issue in two-way parking layouts generated by the optimization models in Section \ref{sec:MIP}, as the lanes are not clearly delineated. We try to address these concerns in a general manner through additional constraints in \ref{sec:practical_considerations} for both two-way and one-way lanes. However, for the rest of the paper, we ignore the issues associated with turning since we focus on the R2 resolution \citep{stephan2021layout} with $\parkwidth = 1$, $\parklength=2$, and $\drivewidth=2$ for two-way lanes and $\drivewidth=1$ for the one-way case, where such concerns are less pronounced.

\begin{remark}
So far, we have assumed that the entry and exit locations are known. More generally, the entrance and exit locations could also be optimized along the lot's perimeter (or in the interior if designing a level of a multi-storied parking facility). The number of feasible entrance and exit options may be limited by the geometry of adjoining streets and their proximity to intersections (or building plans in the case of a multi-storied parking facility). One possibility is to run these models with different entry/exit combinations and select the lot with the best configuration. Alternatively, the MIP formulations can be modified to allow simultaneous optimization of the entrances and exits. However, this may introduce several additional variables and constraints, slowing the computation of optimal solutions.  

In the two-way case, we can create a dummy entrance/exit or a super sink $(\entrancex, \entrancey)$ to which all valid driving cells on the boundary are connected through dummy arcs. See the purple node in the right panel of Figure \ref{fig:two_way_entrance}. In this instance, we assume that all four sides of the lot are accessible from adjoining streets; in practice, however, one can restrict attention to a subset of boundary cells that provide feasible access. To optimize the entrance/exit, we require an additional constraint ensuring that only one of these dummy arcs from cells in $\gridinvadjset{p}{q}$ to $(\entrancex, \entrancey)$ carries positive flow. This necessitates introducing driveway direction-like $z$ variables for each of these arcs and constraints of the form $\flowvar{i}{j}{p}{q} \leq M \countvar{i}{j}{p}{q}~\forall~(i, j) \in \gridinvadjset{p}{q}$ and $\sum_{(i, j) \in \gridinvadjset{p}{q}} \countvar{i}{j}{p}{q} = 1$. The resulting optimal layout is shown in the left panel of Figure \ref{fig:two_way_entrance}. Any one of the active driving fields anchored at the boundary, i.e., the driving field at (0, 4) or (0, 6), can serve as the entrance/exit. In this example, the optimal number of parking spaces remains the same as before, but the stall configuration is slightly different.

For the one-way case, we proceed similarly by connecting valid driving cells on the boundary to a dummy entrance as well as a dummy exit $(r,s)$, shown by the light purple node in the right panel of Figure \ref{fig:one_way_entrance}. In addition to the constraints added for the two-way case, we include constraints of the form $\flowvardash{r}{s}{i}{j} \leq M \countvar{r}{s}{i}{j}~\forall~(i,j) \in \gridadjset{r}{s}$ and $\sum_{(i, j) \in \gridadjset{r}{s}} \countvar{r}{s}{i}{j} = 1$. Additionally, we include constraints on the $z$ variables that prevent the entrance and exit fields from overlapping. For example, adding $\countvar{(0,}{4)}{p}{q} + \countvar{r}{s}{(0,}{3)} + \countvar{r}{s}{(0,}{4)} + \countvar{r}{s}{(0,}{5)} \leq 1$ avoids an exit field from overlapping with an entrance field at (0, 4). The resulting optimal layout accommodates 16 parking spaces, three more than in Figure 10, where the entrance and exit were predetermined. The driving fields at (0, 4) and (12, 11) can serve as the entrance and exit, respectively.

\begin{figure}[H]
\centering
\begin{subfigure}[b]{0.49\textwidth}
    \centering
    \includegraphics[width=0.49\textwidth]{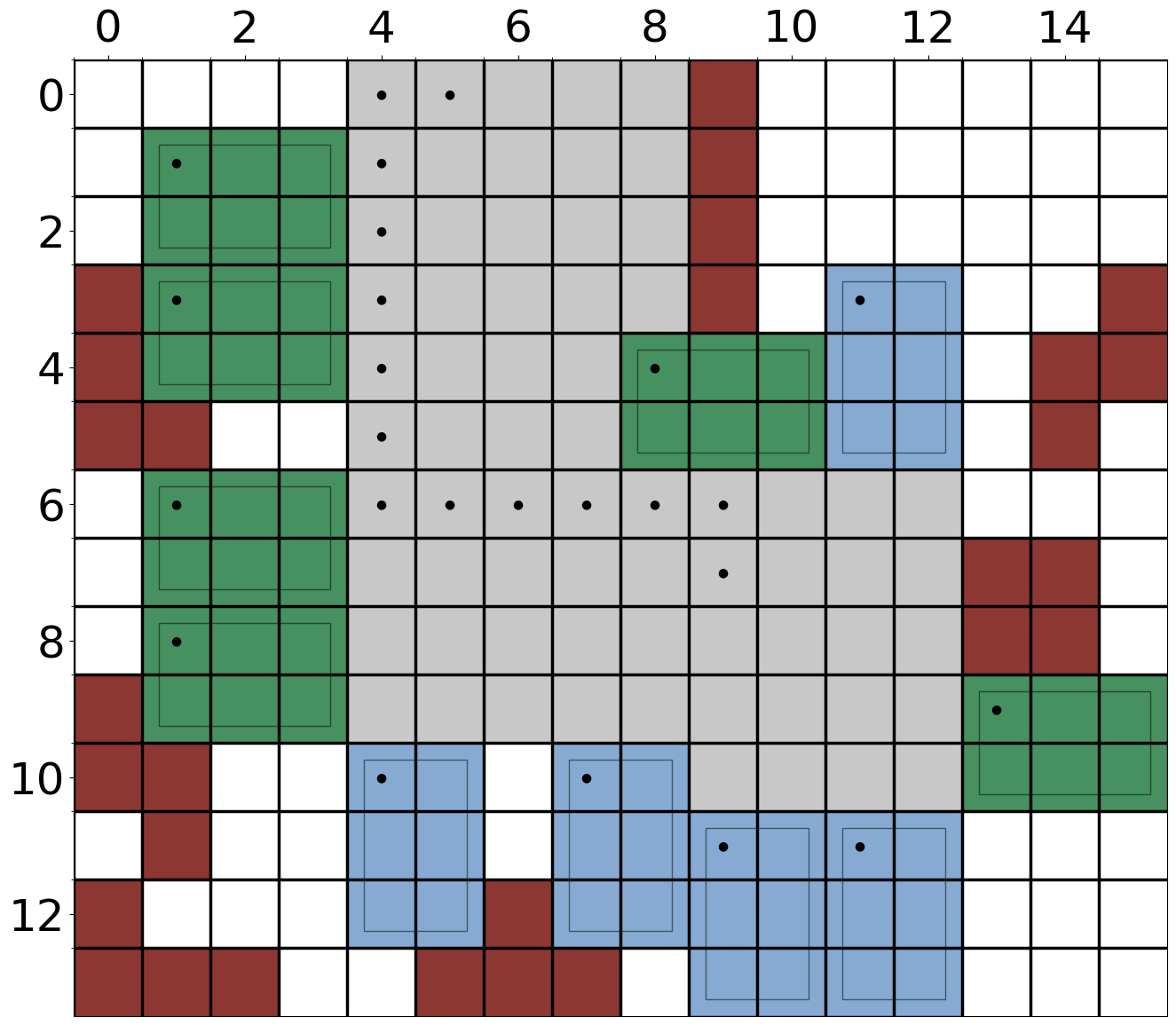}
    \includegraphics[width=0.49\textwidth]{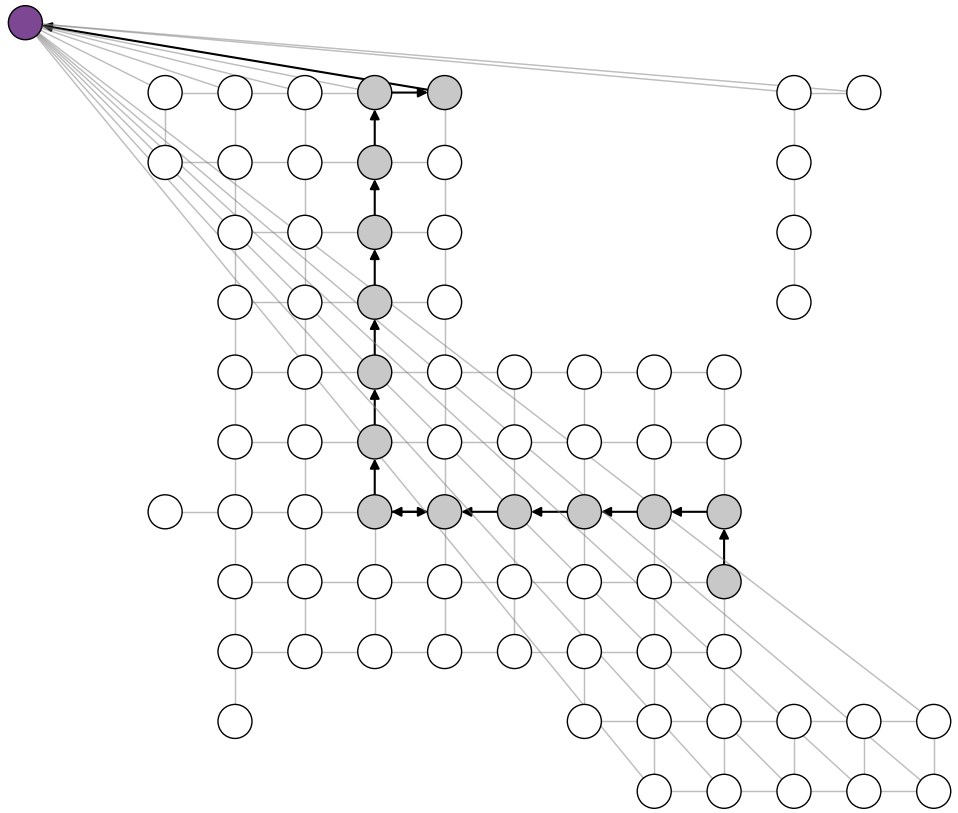}
    \caption{Two-way scenario}
    \label{fig:two_way_entrance}
\end{subfigure}
\begin{subfigure}[b]{0.49\textwidth}
    \centering
    \includegraphics[width=0.49\textwidth]{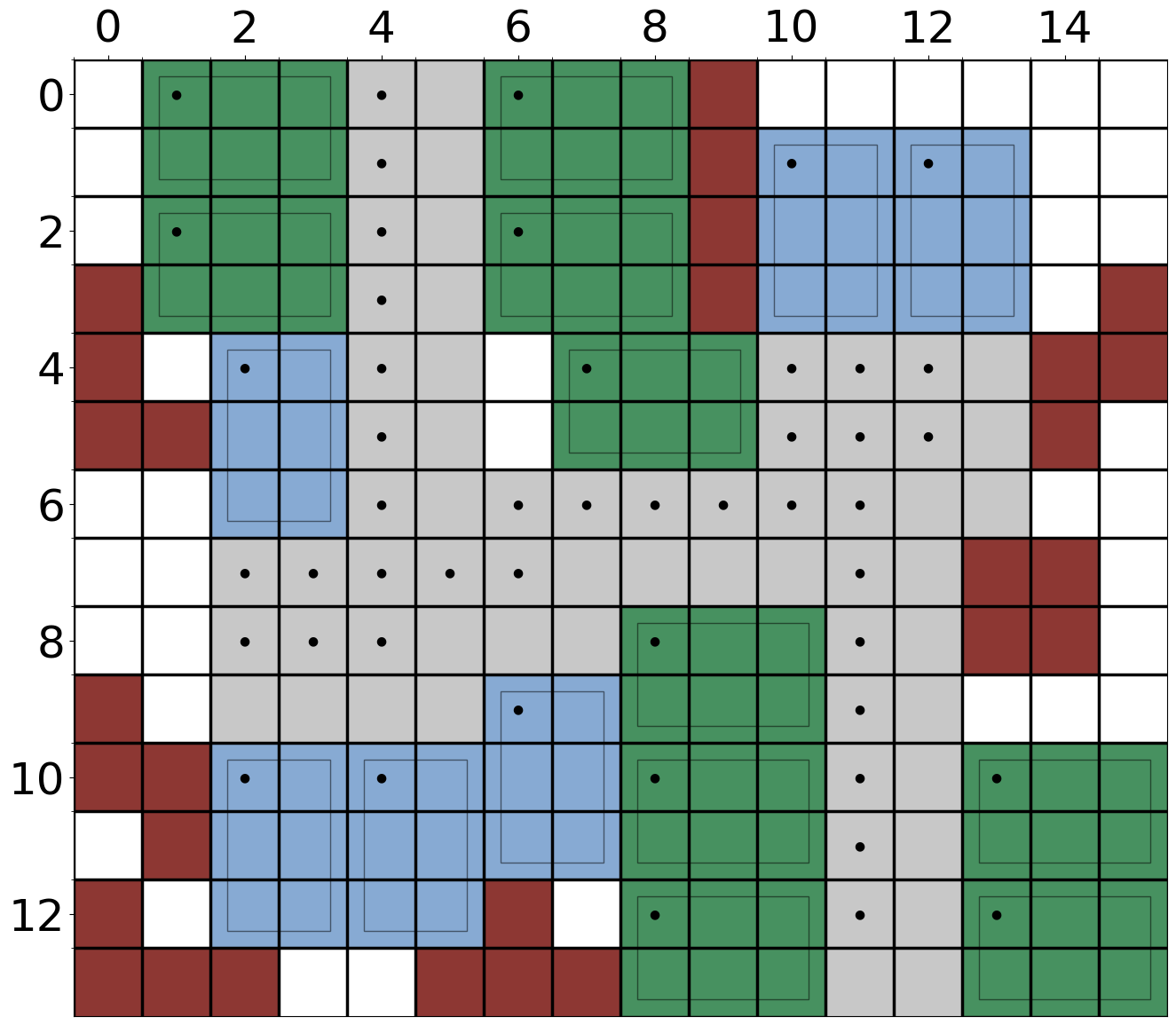}
    \includegraphics[width=0.49\textwidth]{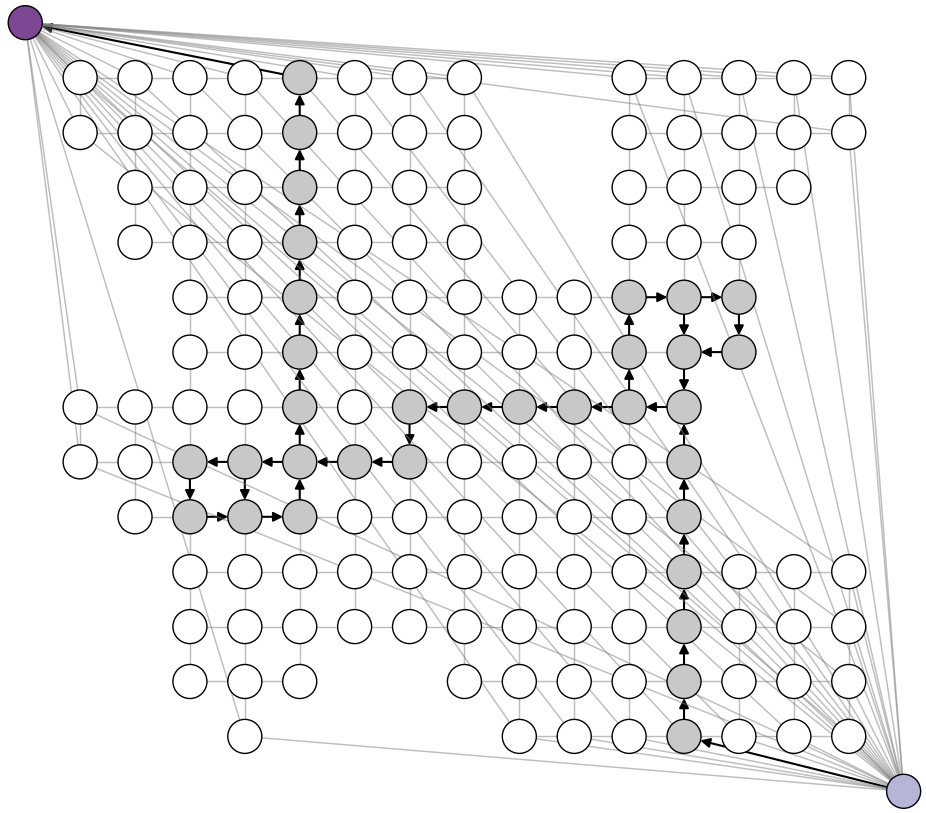}
    \caption{One-way scenario}
    \label{fig:one_way_entrance}
\end{subfigure}
\caption{Layouts and flow solutions with optimized entrances and exits}
\end{figure}
\end{remark}

\section{A branch-and-cut algorithm}
\label{sec:alternate_twoway}
The formulations discussed so far work well for small instances but are slow to converge for medium-scale problems. One reason is that the forcing constraints, which relate the flow variables to the driving or direction variables, lead to poor LP relaxations. \cite{TRISTAN2022website} proposed an alternative formulation that avoids flow variables and, consequently, the use of a constant big-M. We build on this idea by first reformulating the problem using only binary variables for parking and driving, and strengthening the corresponding inequalities, which we refer to as \textit{Connectivity Cut Constraints} (CCCs). We then devise a method in which a limited number of these cuts are added upfront, and the rest are generated dynamically within a branch-and-cut framework. Although the facet-defining properties are not explored in this study, the empirical performance of this approach is encouraging.

\begin{figure}[H]
    \centering
    \includegraphics[width=\linewidth]{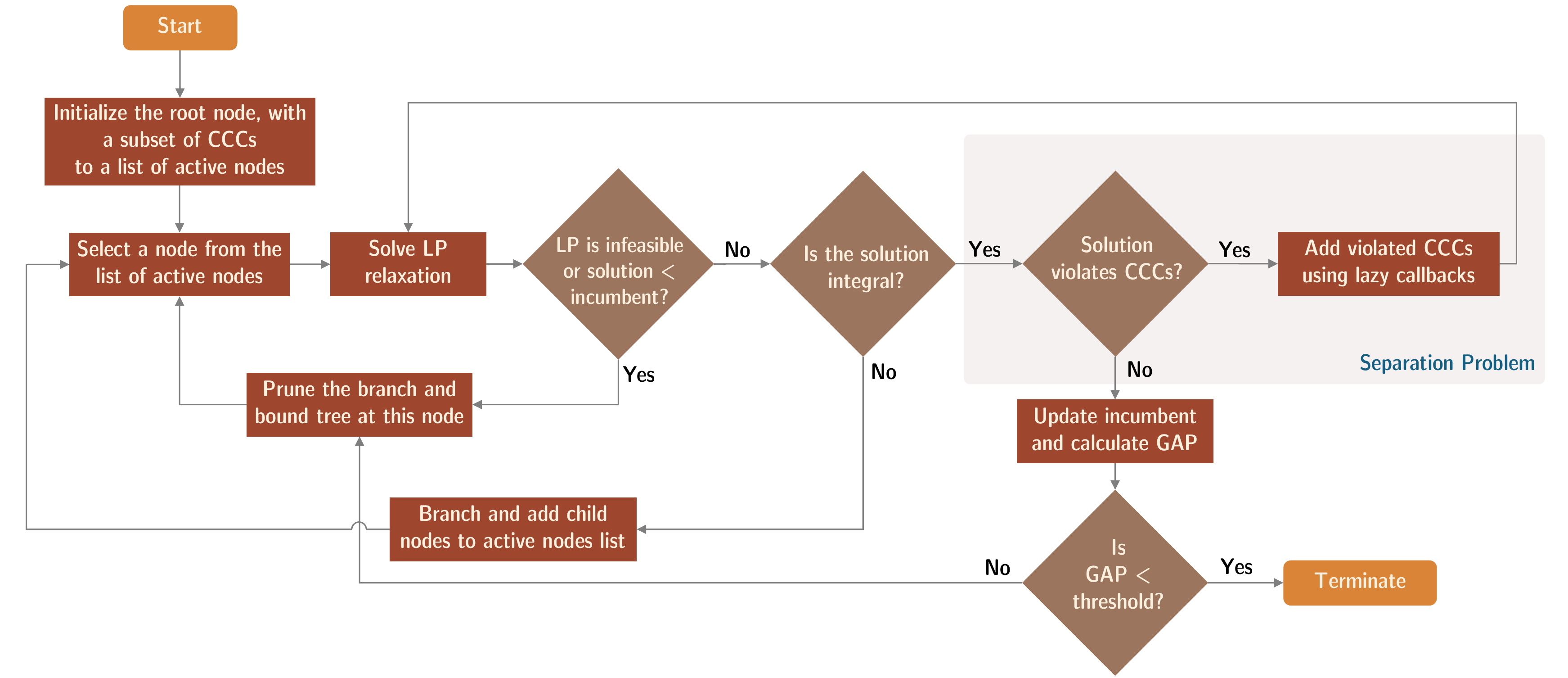}
    \caption{Flowchart of the branch-and-cut algorithm}
    \label{fig:bnc_flowchart}
\end{figure}

The reformulation involves an exponential number of CCCs. To address the resulting computational challenges, we employ a branch-and-cut algorithm that combines branch-and-bound with dynamic cut generation. This approach is particularly well-suited to large-scale MIPs, as it avoids the explicit enumeration of all constraints and instead introduces them only when they are violated. The branch-and-cut method has been successfully applied to several classical logistics problems, including the traveling salesperson problem, vehicle routing problems, and their variants \citep{padberg1991, lysgaard2004, ascheuer2001, tamke2021}, as well as vehicle scheduling problems and their variants \citep{hadjar2006, chai2024}.

Figure \ref{fig:bnc_flowchart} summarizes the steps involved in the branch-and-cut algorithm used. The procedure begins with a limited set of constraints, and LP relaxations are solved at different nodes of the branch-and-bound tree. When an integer solution is obtained, a callback is triggered to verify whether the solution satisfies all omitted constraints. If violations are detected, new CCCs are generated by solving a separation problem to exclude the current solution, and the procedure is repeated. All other branch-and-bound components, such as pruning nodes based on bounds or infeasibility and branching on fractional solutions, are handled using the solver’s default settings.

The proposed method differs slightly for the one-way and two-way scenarios. We first explain the method for the two-way driving configuration in Section  \ref{sec:two-way-branch-and-cut} and discuss the one-way case in Section \ref{sec:one-way-branch-and-cut}.

\subsection{Two-way driveways}
\label{sec:two-way-branch-and-cut}
\begin{figure}[H]
\centering
\begin{minipage}{0.48\textwidth}
    \flushright
    \includegraphics[width=0.63\textwidth]{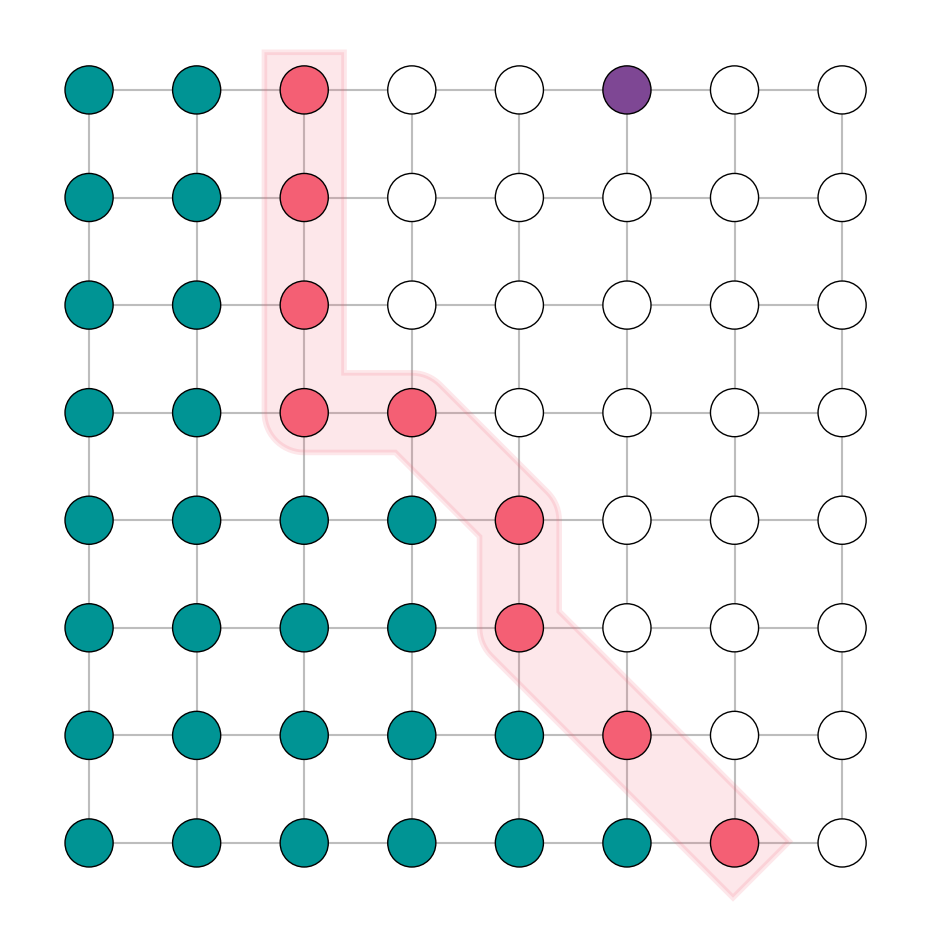}
\end{minipage}
\hfill
\begin{minipage}{0.48\textwidth}
    \flushleft
    \begin{tikzpicture}
        \definecolor{pink}{RGB}{244,95,116}
        \definecolor{teal}{RGB}{0,148,148}
        \definecolor{purple}{RGB}{126,71,148}
        
        \filldraw[fill=pink, draw=black] (0,1cm) circle (0.14cm);
        \filldraw[fill=teal, draw=black] (0,0) circle (0.14cm);
        \filldraw[fill=purple, draw=black] (0,-1cm) circle (0.14cm);
        \filldraw[fill=white, draw=black] (0,-2cm) circle (0.14cm);

        \node[right=5pt] at (0, 1cm)  {\footnotesize{Drive cells in $\vertexcutset$}};
        \node[right=5pt] at (0, 0)    {\footnotesize{Drive cells in $\twowaypartitioncomp{\entrancex \entrancey}$}};
        \node[right=5pt] at (0,-1cm)  {\footnotesize{Entrance $(\entrancex, \entrancey)$ which is part of $\twowaypartition{\entrancex \entrancey}$}};
        \node[right=5pt] at (0,-2cm)  {\footnotesize{Drive cells in $\twowaypartition{\entrancex \entrancey}$}};
    \end{tikzpicture}
\end{minipage}
\caption{Illustration of a vertex cut set}
\label{fig:vertex_cut_set}
\end{figure}

An alternative formulation can be developed for the parking lot design problem by replacing the flow variables with modified driveway connectivity constraints or CCCs of the form \eqref{eq:cutconstraints}. This constraint ensures that a driving field can take a value of one only if at least one of the driving fields within specific vertex cut sets is active. A \textit{vertex cut set} or a \textit{vertex separator} is a subset of nodes in $\grid$ which, when removed, partitions the graph into exactly two disconnected components with the entrance in one of them. If more than two components are generated, smaller vertex cut sets can be created to separate the graph into two components, with the entrance/exit in one of them. Let $\vertexpartitions{pq}$ be a collection of all vertex cut sets that separate the entrance/exit $(\entrancex, \entrancey)$. We denote the partitions created by a vertex separator as a tuple $(\twowaypartition{\entrancex \entrancey}, \vertexcutset, \twowaypartitioncomp{\entrancex \entrancey})$, where $\vertexcutset \in \vertexpartitions{pq}$, $\twowaypartition{\entrancex \entrancey}$ is the set of nodes in the component that contains the entrance/exit $(\entrancex, \entrancey)$, and $\twowaypartitioncomp{\entrancex \entrancey} = \validdriveset \setminus (\twowaypartition{\entrancex \entrancey} \cup \vertexcutset) \neq \emptyset$. To generate tighter cuts, we require every cell in $V$ to be connected to at least one element of $\twowaypartitioncomp{\entrancex \entrancey}$. Figure~\ref{fig:vertex_cut_set} illustrates an example driving grid with a vertex cut set $\vertexcutset$ (shown in pink) that separates the purple entrance cell $(p, q)$ from the teal-colored cells $\twowaypartitioncomp{\entrancex \entrancey}$.
If the drive cell variables corresponding to a cut set $\vertexcutset$ are all zero, all driving fields in $\twowaypartitioncomp{\entrancex \entrancey}$ are forced to be inactive. This relationship can be expressed using \eqref{eq:cutconstraints}. 
\begin{align}
\drivevar{i}{j} \leq \sum_{(k,l) \in \vertexcutset} \drivevar{k}{l} \qquad \forall \, \vertexcutset \in \vertexpartitions{pq}, (i, j) \in \twowaypartitioncomp{\entrancex \entrancey}  \label{eq:cutconstraints}
\end{align}

We first show that these inequalities are valid for the flow-based formulation using Propositions \ref{prop:2w_aggregate} and \ref{prop:2w_disaggregate}. We then strengthen it by adding extra terms to the left-hand side and prove its validity in Proposition \ref{prop:2w_strengthened}. 

\begin{proposition}
\label{prop:2w_aggregate}
Let $\vertexcutset \in \vertexpartitions{pq}$ be a vertex separator and let the partitions created in $\grid$ after removing it be $\twowaypartition{\entrancex \entrancey}$ and $\twowaypartitioncomp{\entrancex \entrancey}$. Then the inequality 
\begin{align}
\sum_{(i, j) \in  \twowaypartitioncomp{\entrancex \entrancey}} \drivevar{i}{j} \leq \sum_{(k, l) \in \vertexcutset} M c_{kl} \drivevar{k}{l} \label{eq:aggregate_two_way}
\end{align}
where $c_{kl} = |\gridinvadjset{k}{l} \cap \twowaypartitioncomp{\entrancex \entrancey}|$ is valid for the feasible region of the formulation $\mathcal{F}^{\textsc{flow}}_{\textsc{2W}}$.
\end{proposition}
\begin{proof}
Summing the flow-conservation constraints \eqref{eq:floworigin} for all cells in $\twowaypartitioncomp{\entrancex \entrancey}$, 
\begin{align}
\sum_{(i,j) \in \twowaypartitioncomp{\entrancex \entrancey}} \drivevar{i}{j} & = \sum_{\substack{((i,j), (k,l)) \in \gridarcset :\\
(i,j) \in \twowaypartitioncomp{\entrancex \entrancey},
(k,l) \in \vertexcutset
}} f_{ij, kl} -
\sum_{\substack{((k,l), (i,j)) \in \gridarcset :\\
(k,l) \in \vertexcutset,
(i,j) \in \twowaypartitioncomp{\entrancex \entrancey}
}} f_{kl, ij} \notag \\
& \leq \sum_{\substack{((i,j), (k,l)) \in \gridarcset :\\
(i,j) \in \twowaypartitioncomp{\entrancex \entrancey},
(k,l) \in \vertexcutset
}} f_{ij, kl} \notag \\
& \leq \sum_{\substack{((i,j), (k,l)) \in \gridarcset :\\
(i,j) \in \twowaypartitioncomp{\entrancex \entrancey},
(k,l) \in \vertexcutset
}} M \drivevar{k}{l}  \qquad \,\,\,\,  \text{[using \eqref{eq:flowhead}]} \notag\\
& = \sum_{(k, l) \in \vertexcutset} M c_{kl} \drivevar{k}{l} \notag
\end{align}
\end{proof}

\begin{proposition}
\label{prop:2w_disaggregate}
Constraints \eqref{eq:cutconstraints} are valid for the feasible region of the formulation $\mathcal{F}^{\textsc{flow}}_{\textsc{2W}}$.
\end{proposition}
\begin{proof}
We prove the contrapositive by showing that a solution violating \eqref{eq:cutconstraints} does not belong to the feasible region of the formulation $\mathcal{F}^{\textsc{flow}}_{\textsc{2W}}$. For a given vertex separator $\vertexcutset$ and $(i,j) \in \twowaypartitioncomp{\entrancex \entrancey}$, suppose that $\drivevar{i}{j} > \sum_{(k,l) \in \vertexcutset} \drivevar{k}{l}$ for some solution $y$. Since the $y$ variables are binary, it follows that $\drivevar{i}{j} = 1$ and $\sum_{(k,l) \in \vertexcutset} \drivevar{k}{l} = 0$, which implies $\drivevar{k}{l} = 0~\forall~ (k,l) \in \vertexcutset$, violating \eqref{eq:aggregate_two_way}. Hence, this solution is not feasible to $\mathcal{F}^{\textsc{flow}}_{\textsc{2W}}$.    
\end{proof}

\subsubsection{Strengthened valid inequalities}
\label{sec:vis_two_way}
If the drive variables in a cut set are zero, it constrains the drive variables in the partition without the entrance to be zero due to \eqref{eq:cutconstraints}. Additionally, some parking variables can also be made inactive. Specifically, all parking fields whose drive field neighbors lie exclusively in $\twowaypartitioncomp{\entrancex \entrancey} \cup \vertexcutset$ would become inaccessible. To strengthen the connectivity cut constraints based on this observation, consider a cut set $\vertexcutset$ and define the following sets. For a given cell $(i, j) \in \validdriveset$ and a set $S \subset \validdriveset$, let $\invzeroparksetcut{i}{j}{S}$ be the set of $0^{\circ}$ parking fields containing $(i, j)$ and whose driving field neighbors lie completely within $S$. Mathematically, $\invzeroparksetcut{i}{j}{S} = \{(k, l) \in \invzeroparkset{i}{j} : \zeroneighborset{k}{l} \subseteq S \}$. Likewise, define $\invninetyparksetcut{i}{j}{S} = \{(k, l) \in \invninetyparkset{i}{j} : \ninetyneighborset{k}{l} \subseteq S\}$ and $\invdrivesetcut{i}{j}{S} = \invdriveset{i}{j} \cap S $. With these new sets, we can strengthen \eqref{eq:cutconstraints} to \eqref{eq:strongcutconstraints}. Notice that the left-hand side of this strengthened inequality is similar to that of the single-purpose constraint \eqref{eq:singlepurpose_zero_disaggregate}. Replacing all the flow-conservation constraints with these inequalities yields an equivalent model $\mathcal{F}^{\textsc{CCC}}_{\textsc{2W}}$ shown below with identical feasible $(x, y)$ solutions because the new constraints can be viewed as feasibility cuts for a Benders' reformulation (see \ref{sec:equivalence} for a proof). However, the main drawback of these constraints is that their number grows exponentially with the size of $\validdriveset$, since the number of possible vertex cut sets can be extremely huge. We tackle this problem using a branch-and-cut algorithm.
\begin{align}
    \max & \sum_{(i,j) \in \validzeroparkset} \zeroparkvar{i}{j} + \sum_{(i,j) \in \validninetyparkset} \ninetyparkvar{i}{j}  && [\text{Formulation } \mathcal{F}^{\textsc{CCC}}_{\textsc{2W}}] \notag\\
      \text{s.t. }  & \eqref{eq:singlepurpose_zero_disaggregate}\text{--}\eqref{eq:ninetyparkdriveconnector}, \eqref{eq:existing drive vars}\text{--}\eqref{eq: y var is binary} && \notag \\
    & \drivevar{m}{n} + \sum_{(k, l) \in \invzeroparksetcut{i}{j}{{\twowaypartitioncomp{\entrancex \entrancey} \cup \vertexcutset}}} \zeroparkvar{k}{l} + \sum_{(k, l) \in \invninetyparksetcut{i}{j}{\twowaypartitioncomp{\entrancex \entrancey} \cup \vertexcutset}} \ninetyparkvar{k}{l} \leq \sum_{(k,l) \in \vertexcutset} \drivevar{k}{l} \qquad && \forall \, \vertexcutset \in \vertexpartitions{pq}, (i, j) \in \twowaypartitioncomp{\entrancex \entrancey}, (m, n) \in \invdrivesetcut{i}{j}{\twowaypartitioncomp{\entrancex \entrancey}}   \label{eq:strongcutconstraints}
\end{align}

\begin{proposition}
\label{prop:2w_strengthened}
The strengthened inequalities \eqref{eq:strongcutconstraints} are valid for the feasible region of the formulation $\mathcal{F}^{\textsc{flow}}_{\textsc{2W}}$.
\end{proposition}
\begin{proof}
We once again proceed by proving the contrapositive. Observe that \eqref{eq:singlepurpose_zero_disaggregate} ensures that the left-hand side of \eqref{eq:strongcutconstraints} never exceeds one; therefore, the constraint holds trivially whenever the right-hand side is at least one. Suppose \eqref{eq:strongcutconstraints} is violated by a binary solution $(x, y)$ for some vertex cut set $\vertexcutset$ and $\sum_{(k,l) \in \vertexcutset} \drivevar{k}{l} = 0$. Three cases arise. \textbf{Case (1):} $y_{mn} = 1$ for some $(m, n) \in \invdrivesetcut{i}{j}{\twowaypartitioncomp{\entrancex \entrancey}}$. Since $(m,n) \in \twowaypartitioncomp{\entrancex \entrancey}$, \eqref{eq:cutconstraints} is violated, i.e., $\drivevar{m}{n} > \sum_{(k,l) \in \vertexcutset} \drivevar{k}{l}$, implying that $(x, y)$ is infeasible for the formulation $\mathcal{F}^{\textsc{flow}}_{\textsc{2W}}$. \textbf{Case (2):} $\zeroparkvar{k}{l} = 1$ for some $(k, l) \in \invzeroparksetcut{i}{j}{{\twowaypartitioncomp{\entrancex \entrancey} \cup \vertexcutset}}$. From \eqref{eq:zeroparkdriveconnector}, $\drivevar{u}{v} = 1$ for some $(u,v) \in \zeroneighborset{k}{l}$. But by construction $(u, v) \in \twowaypartitioncomp{\entrancex \entrancey}$ and hence $\drivevar{u}{v} > \sum_{(k,l) \in \vertexcutset} \drivevar{k}{l}$, which implies that the solution is infeasible for $\mathcal{F}^{\textsc{flow}}_{\textsc{2W}}$. \textbf{Case (3):}  $\ninetyparkvar{k}{l} = 1$ for some $(k, l) \in \invzeroparksetcut{i}{j}{{\twowaypartitioncomp{\entrancex \entrancey} \cup \vertexcutset}}$ is similar to Case (2). 
\end{proof}

The exponential number of constraints \eqref{eq:strongcutconstraints} required to enforce connectivity is similar to the sub-tour elimination constraints of the traveling salesperson problem. To tackle this challenge, we employ three strategies. (1) We propose a family of valid inequalities added in advance. (2) To avoid slowing down the LP relaxations, we limit the number of these inequalities by only considering $(m, n) = (i, j)$ instead of all $(m, n) \in \invdrivesetcut{i}{j}{\twowaypartitioncomp{pq}}$.  (3) The model with these limited valid inequalities typically produces disconnected layouts, but the number of driving-field components is usually small. We propose a simple and effective procedure for generating feasibility cuts that reject such invalid integer solutions and iteratively add them to a branch-and-cut scheme when a violation is found, as shown in Figure \ref{fig:bnc_flowchart}. For the constraints added upfront, we describe a simple family of connectivity cut constraints, which we refer to as \textit{bidirectional hop inequalities}.  These are added to the model as regular or lazy constraints and provide adequate restrictions to obtain reasonable upper bounds. We will use the example layout from the previous section to illustrate these methods, but with the R2 resolution. 

\begin{figure}[H]
  \centering
  \begin{subfigure}[b]{0.32\textwidth}
    \includegraphics[width=0.9\textwidth]{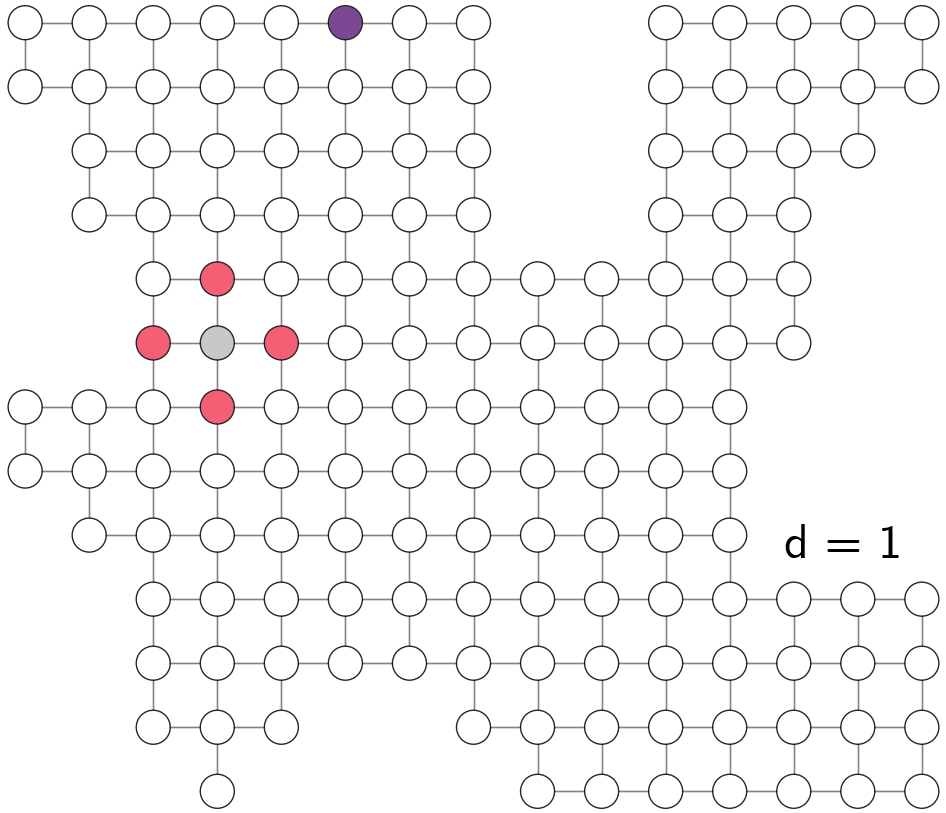}
  \end{subfigure}
    \begin{subfigure}[b]{0.32\textwidth}
    \includegraphics[width=0.9\textwidth]{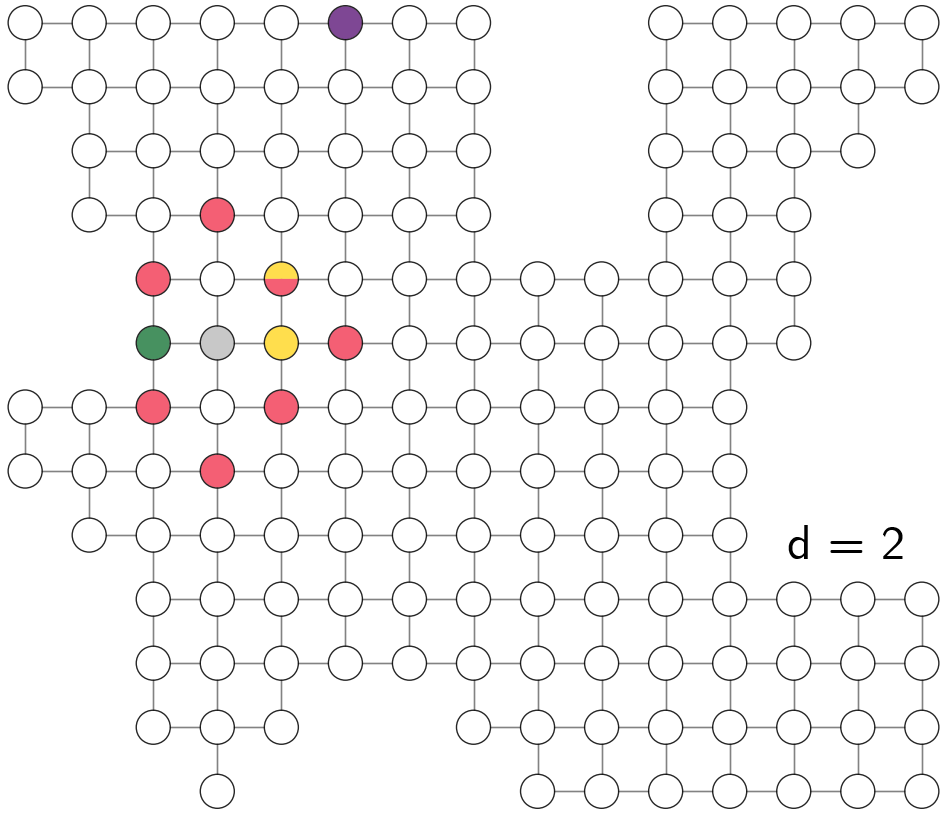}
  \end{subfigure}
      \begin{subfigure}[b]{0.32\textwidth}
    \includegraphics[width=0.9\textwidth]{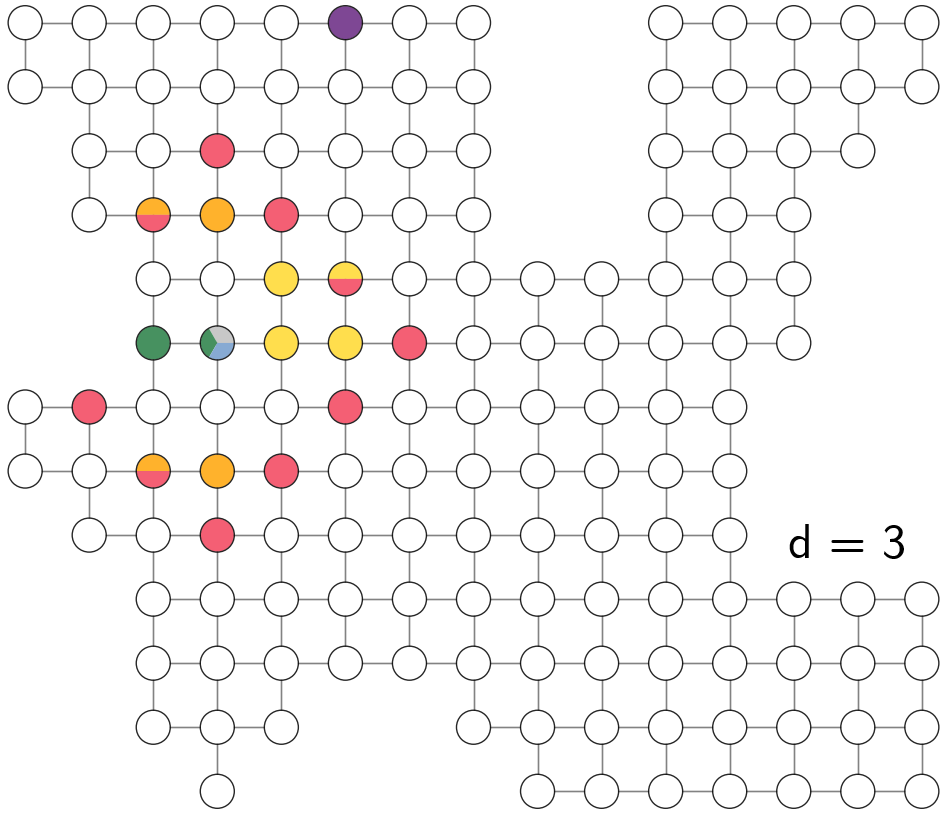}
  \end{subfigure}\\
    \vspace{3mm}
  \begin{subfigure}[b]{0.32\textwidth}
    \includegraphics[width=0.9\textwidth]{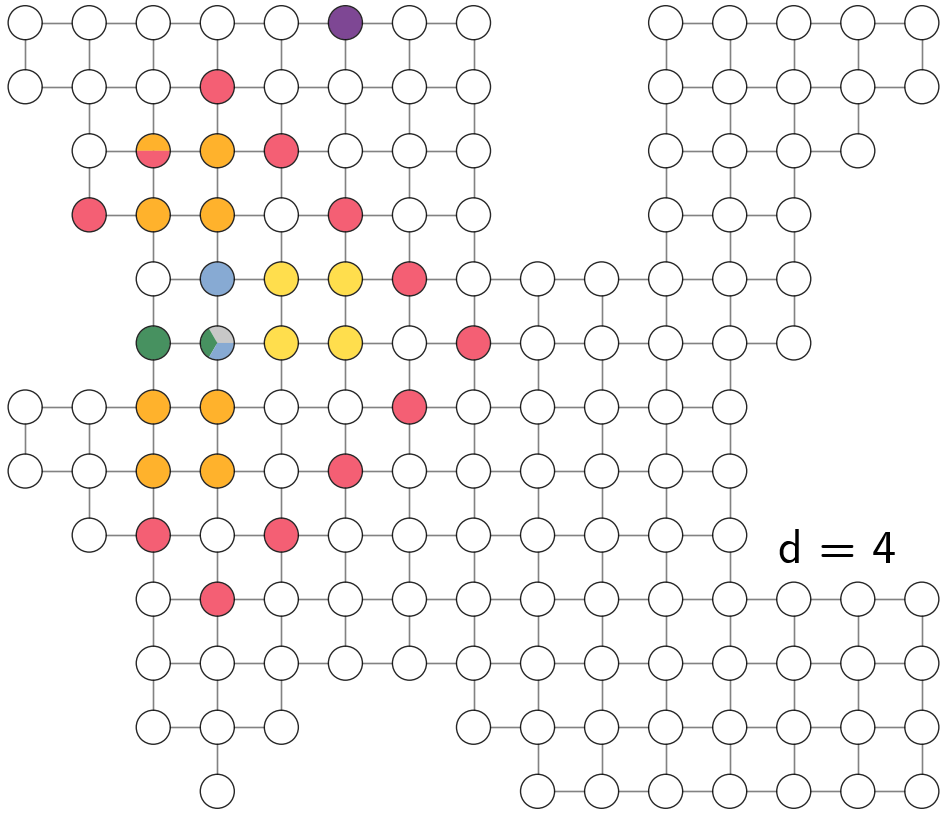}
  \end{subfigure}
    \begin{subfigure}[b]{0.32\textwidth}
    \includegraphics[width=0.9\textwidth]{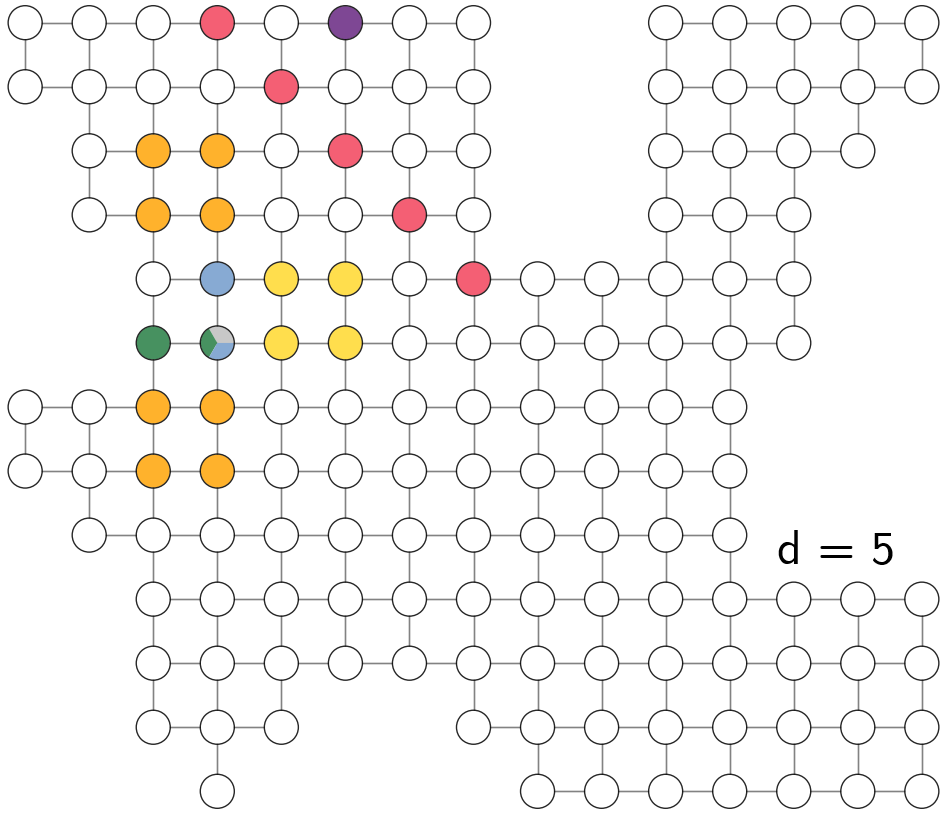}
  \end{subfigure}
      \begin{subfigure}[b]{0.32\textwidth}
    \includegraphics[width=0.9\textwidth]{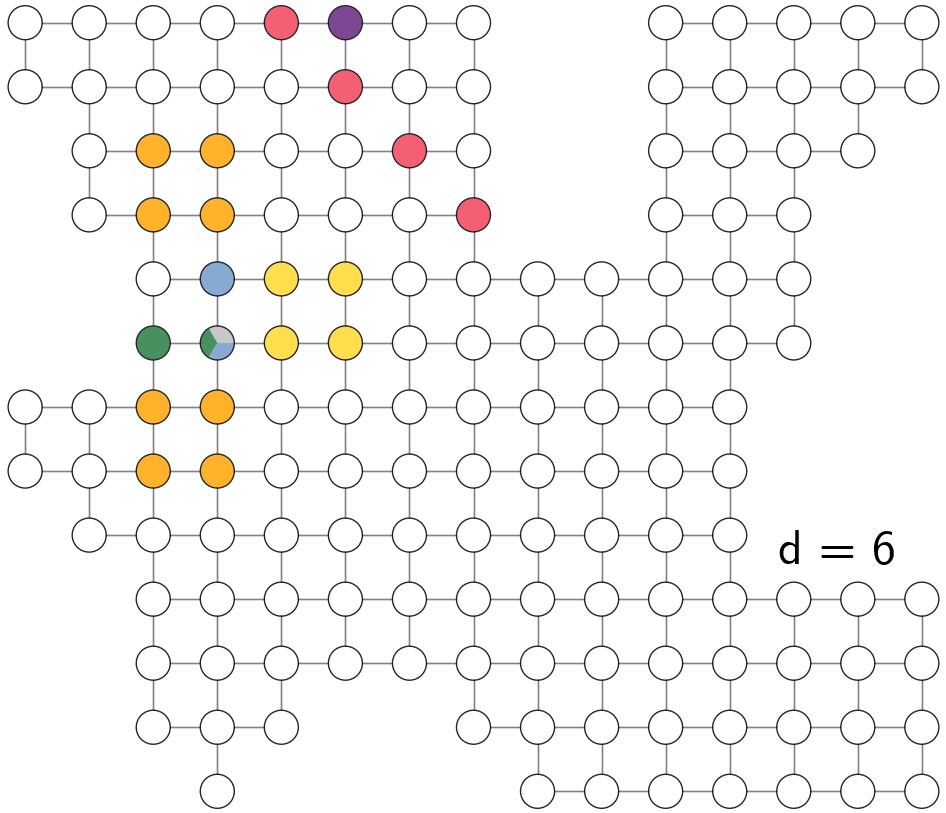}
  \end{subfigure}
  \begin{tikzpicture}
        \definecolor{pink}{RGB}{244,95,116}
        \definecolor{purple}{RGB}{126,71,148}
        \definecolor{blue}{RGB}{135,170,211}
        \definecolor{green}{RGB}{71,145,96}
        \definecolor{grey}{RGB}{200,200,200}
        \definecolor{yellow}{RGB}{255,222,77}
        \definecolor{gold}{RGB}{255,178,44}
        
        \filldraw[fill=pink, draw=black] (-8,0.25cm) circle (0.12cm);
        \filldraw[fill=green, draw=black] (-6.4,0.25cm) circle (0.12cm);
        \filldraw[fill=blue, draw=black] (-3.6,0.25cm) circle (0.12cm);
        \filldraw[fill=yellow, draw=black] (-0.8,0.25cm) circle (0.12cm);
        \filldraw[fill=gold, draw=black] (2.3,0.25cm) circle (0.12cm);
        \filldraw[fill=grey, draw=black] (5.2,0.25cm) circle (0.12cm);
        \filldraw[fill=purple, draw=black] (8,0.25cm) circle (0.12cm);

        \node[right=5pt] at (-8,0.23cm)  {\footnotesize{$H_{(5,3)}^d$}};
        \node[right=5pt] at (-6.4,0.23cm)    {\footnotesize{$\invzeroparksetcut{(5,}{3)}{{\twowaypartitioncomp{\entrancex \entrancey} \cup \vertexcutset}}$}};
        \node[right=5pt] at (-3.6,0.23cm)  {\footnotesize{$\invninetyparksetcut{(5,}{3)}{\twowaypartitioncomp{\entrancex \entrancey} \cup \vertexcutset}$}};
        \node[right=5pt] at (-0.8,0.23cm)  {{\footnotesize $\zeroneighborset{i}{j}$ of green cells}};
        \node[right=5pt] at (2.3,0.23cm)  {{\footnotesize $\ninetyneighborset{i}{j}$ of blue cells}};
        \node[right=5pt] at (5.2,0.23cm)  {{\footnotesize Drive cell $(5,3)$}};
        \node[right=5pt] at (8,0.23cm)  {{\footnotesize $(\entrancex, \entrancey)$}};
    \end{tikzpicture}
    \vspace{-2mm}
\caption{Illustration of forward hop inequalities. The gray driving field (5, 3) is active only if accessible through its $d$-hop neighbors  (pink). Inequalities can be strengthened by adding additional parking variables (green and blue). Yellow and orange nodes are neighbors of parking fields in the inequality, and the purple node is the entrance.}
\label{fig:hop_connectivity}
\end{figure}

\textbf{Forward hop inequalities:} Finding vertex cuts for a given driving cell $(i, j)$ that are a certain hop distance away is one way to generate connectivity cuts. Figure \ref{fig:hop_connectivity} illustrates vertex cut sets with pink nodes for the valid driving cell (5, 3), marked in gray, at distances ranging from $d=1$ to 6 hops for the instance from Figure \ref{fig:sample_perpendicular_fields} but with R2 resolution. The central node (5, 3) and its neighbors are colored based on the variables included in \eqref{eq:strongcutconstraints}: gray for drive variables, green for $0^{\circ}$ variables, and blue for $90^{\circ}$ variables. Nodes with multiple colors represent all the corresponding variables that can be added to the valid inequality. The yellow and orange nodes denote the neighbors of the park variables on the left-hand side of \eqref{eq:strongcutconstraints}. Initially, for a hop distance of one, only the drive variable is included on the left-hand side (see the $d=1$ case). If all immediate neighbors of (5, 3) are inactive, connectivity to the entrance or exit from the gray cell cannot be ensured. In other words, at least one pink node has to be active to access the central driving field. Hence, we can add a cut $y_{(5, 3)} \leq y_{(5,2)} + y_{(4, 3)} + y_{(5,4)} + y_{(6,3)}$. Additional park field variables may be added to the left-hand side for vertex cut sets with higher hop distances. For example, for a hop distance of 3, it is not possible to have a $0^{\circ}$ park field at (5, 2) and (5, 3) if all pink cell drive variables are inactive. Similarly, the $90^{\circ}$ park variable at (5, 3) can be included on the valid inequality's left-hand side. This is because their neighboring driving anchors, from which they can be accessed, are either part of the vertex separator or entirely within the component that does not contain the entrance/exit. If a neighbor of a parking field is in the component containing the entrance/exit, the parking field remains accessible even if all driving cells in the vertex cut set are inactive. For a hop limit of 5 and 6, note that the vertex cut set creates more than one component. In these cases, as shown in the figure, the vertex separator is updated to include only essential drive field variables that separate the entrance/exit component from the cell (5, 3), resulting in a tighter cut.

Let $H_{ij}^d$ be the set of nodes $d$ hops from a driving field $(i, j)$, and let $\pi_{ij}$ be the shortest hop distance from $(i, j)$ to the entrance $(p, q)$. For a driving field $(i, j)$, we write valid inequalities for hop distances up to $\pi_{ij}$. We set $\vertexcutset = H_{ij}^d$ and add a single constraint for cell $(i, j)$, from which $H_{ij}^d$ is derived, instead of each cell in $(H_{ij}^d)'_{pq}$. This minimizes the number of lazy constraints, preventing LP relaxation slowdowns, and because other hop vertex separators are anyway generated for cells in $(H_{ij}^d)'_{pq}$. Specifically, we include the following constraint in our formulation.  
\begin{align}
    & \drivevar{i}{j} + \sum_{(k, l) \in \invzeroparksetcut{i}{j}{{\twowaypartitioncomp{\entrancex \entrancey} \cup \vertexcutset}}} \zeroparkvar{k}{l} + \sum_{(k, l) \in \invninetyparksetcut{i}{j}{\twowaypartitioncomp{\entrancex \entrancey} \cup \vertexcutset}} \ninetyparkvar{k}{l} \leq \sum_{(k,l) \in \vertexcutset} \drivevar{k}{l} \,  && \forall \, (i, j) \in \validdriveset \setminus \entranceset,  d \in \{1, \ldots, \pi_{ij}\}, \vertexcutset = H_{ij}^d
     \label{eq:hop1}
\end{align}
Intuitively, using these inequalities for different hop distances increases the likelihood of maintaining connectivity all the way to the entrance of the parking lot. Although for larger hop distances, there is a greater chance of including more decision variables on the left-hand side, it may also increase the number of drive variables on the right-hand side, potentially reducing the constraint's effectiveness. Therefore, in our experiments, we add these inequalities as regular constraints for lower hop distances or when the right-hand side has fewer variables. The remaining hop inequalities are incorporated as lazy constraints.

\begin{figure}[H]
  \centering
  \begin{subfigure}[b]{0.32\textwidth}
    \includegraphics[width=0.9\textwidth]{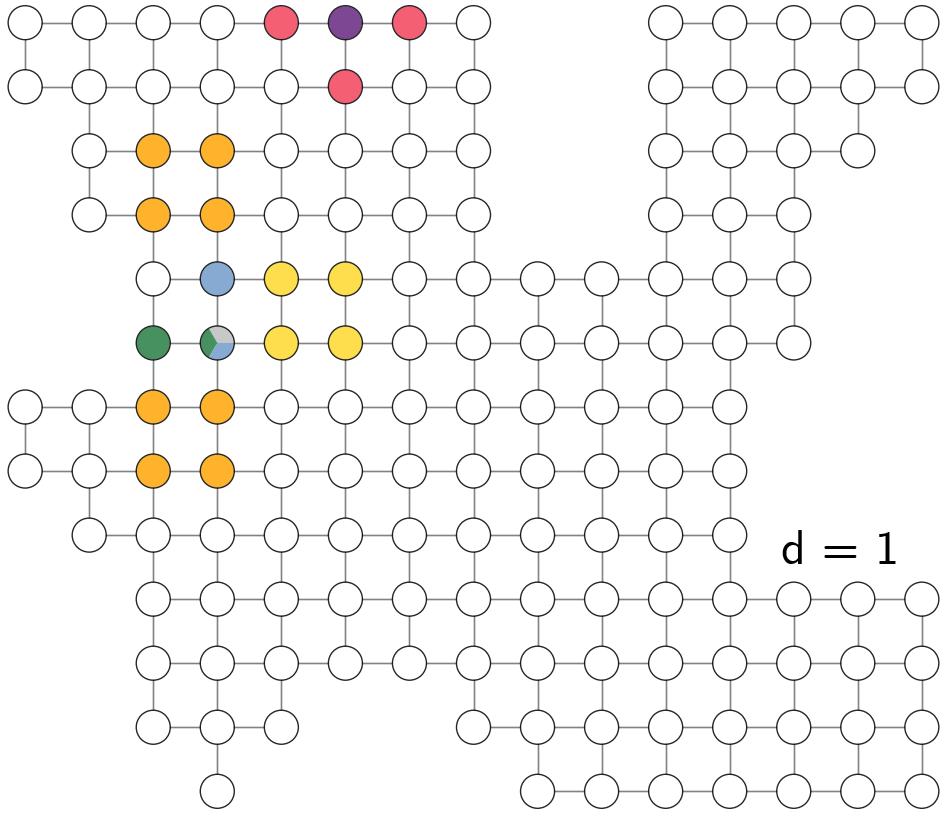}
  \end{subfigure}
    \begin{subfigure}[b]{0.32\textwidth}
    \includegraphics[width=0.9\textwidth]{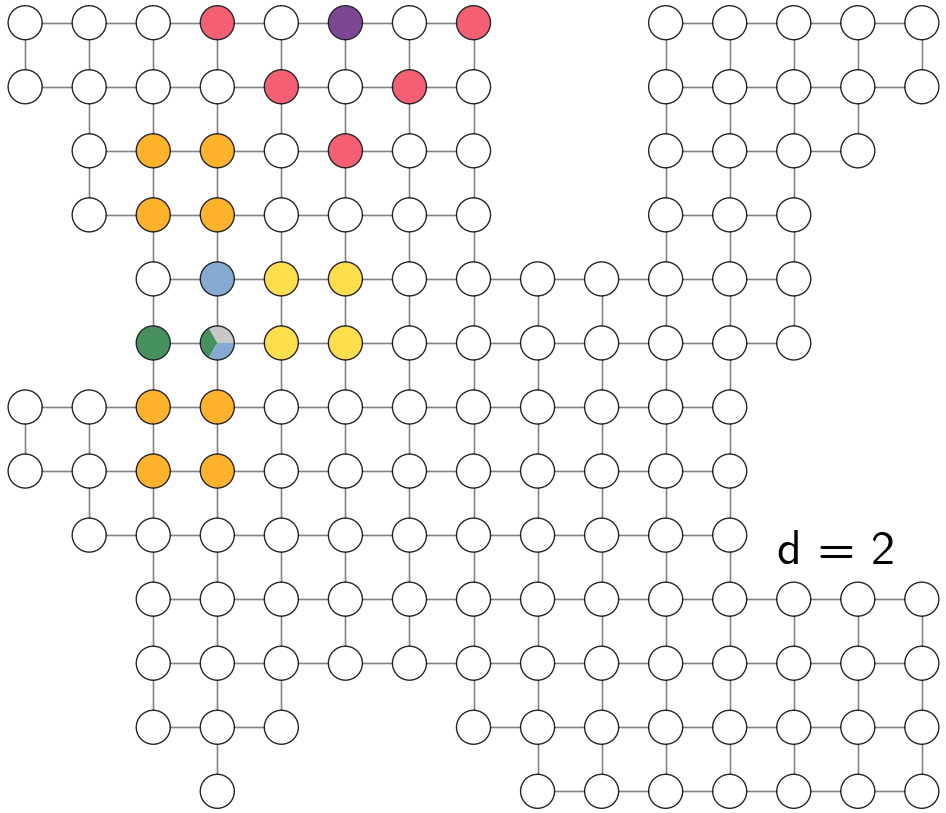}
  \end{subfigure}
      \begin{subfigure}[b]{0.32\textwidth}
    \includegraphics[width=0.9\textwidth]{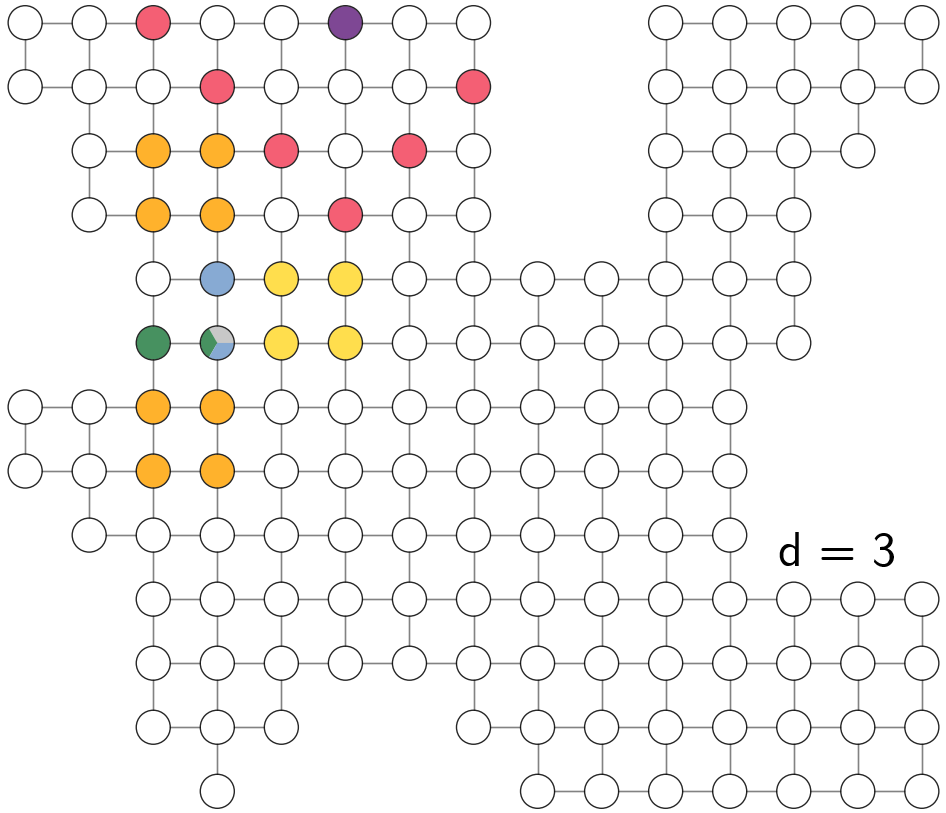}
  \end{subfigure}\\
  \vspace{3mm}
  \begin{subfigure}[b]{0.32\textwidth}
    \includegraphics[width=0.9\textwidth]{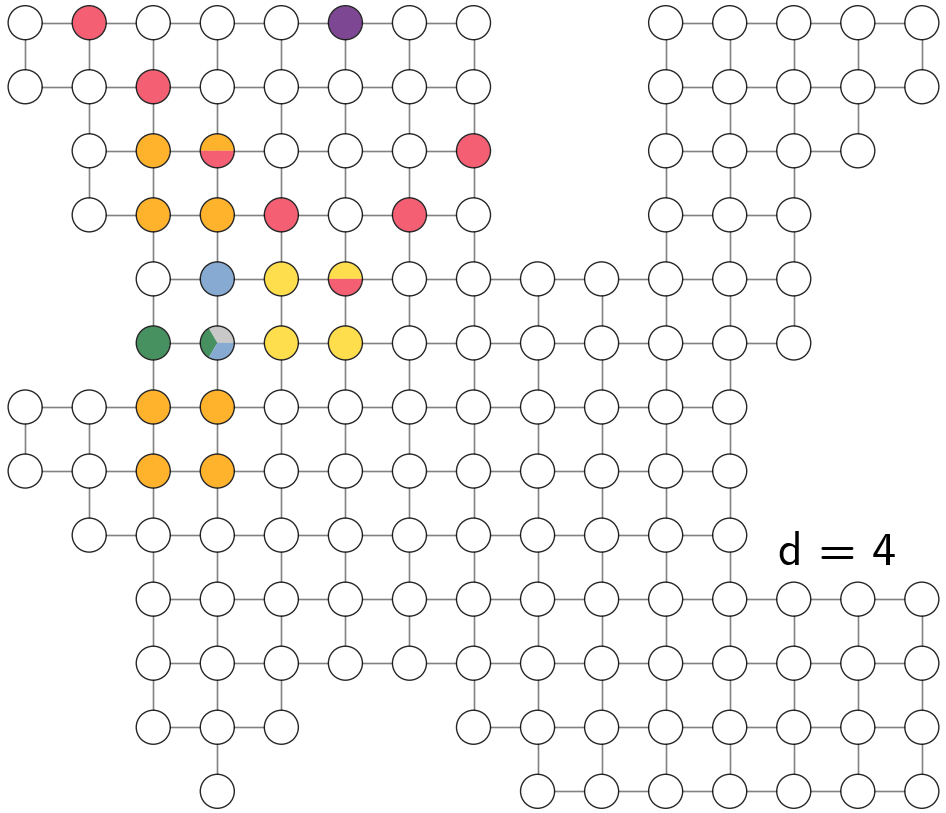}
  \end{subfigure}
    \begin{subfigure}[b]{0.32\textwidth}
    \includegraphics[width=0.9\textwidth]{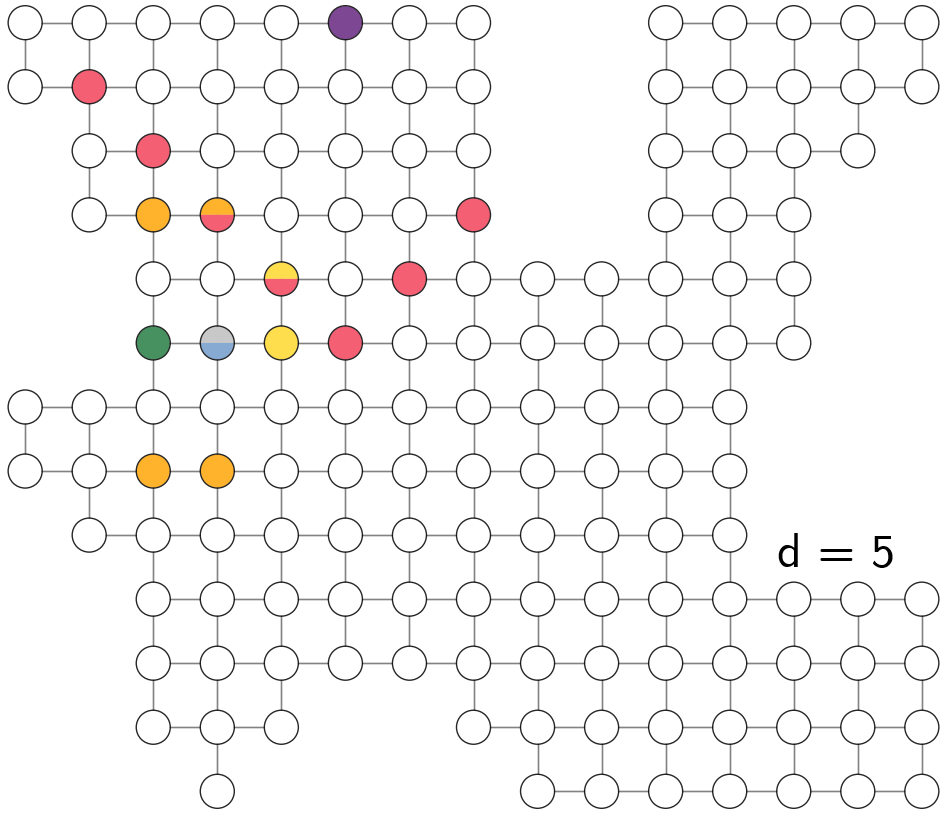}
  \end{subfigure}
      \begin{subfigure}[b]{0.32\textwidth}
    \includegraphics[width=0.9\textwidth]{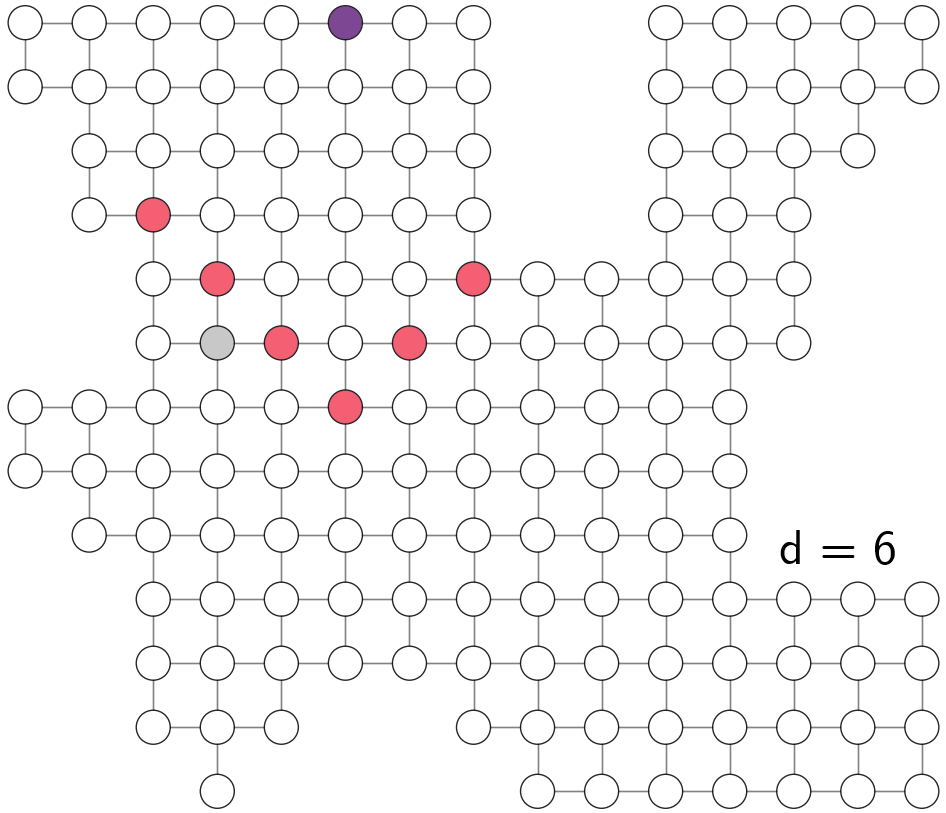}
  \end{subfigure}
    \begin{tikzpicture}
        \definecolor{pink}{RGB}{244,95,116}
        \definecolor{purple}{RGB}{126,71,148}
        \definecolor{blue}{RGB}{135,170,211}
        \definecolor{green}{RGB}{71,145,96}
        \definecolor{grey}{RGB}{200,200,200}
        \definecolor{yellow}{RGB}{255,222,77}
        \definecolor{gold}{RGB}{255,178,44}
        
        \filldraw[fill=pink, draw=black] (-8,0.25cm) circle (0.12cm);
        \filldraw[fill=green, draw=black] (-6.5,0.25cm) circle (0.12cm);
        \filldraw[fill=blue, draw=black] (-3.6,0.25cm) circle (0.12cm);
        \filldraw[fill=yellow, draw=black] (-0.7,0.25cm) circle (0.12cm);
        \filldraw[fill=gold, draw=black] (2.4,0.25cm) circle (0.12cm);
        \filldraw[fill=grey, draw=black] (5.3,0.25cm) circle (0.12cm);
        \filldraw[fill=purple, draw=black] (8,0.25cm) circle (0.12cm);

        \node[right=5pt] at (-8,0.23cm)  {\footnotesize{$H_{pq}^d$}};
        \node[right=5pt] at (-6.5,0.23cm)    {\footnotesize{$\invzeroparksetcut{(5,}{3)}{{\twowaypartitioncomp{\entrancex \entrancey} \cup \vertexcutset}}$}};
        \node[right=5pt] at (-3.6,0.23cm)  {\footnotesize{$\invninetyparksetcut{(5,}{3)}{\twowaypartitioncomp{\entrancex \entrancey} \cup \vertexcutset}$}};
        \node[right=5pt] at (-0.7,0.23cm)  {{\footnotesize$\zeroneighborset{i}{j}$ of green cells}};
        \node[right=5pt] at (2.4,0.23cm)  {{\footnotesize$\ninetyneighborset{i}{j}$ of blue cells}};
        \node[right=5pt] at (5.3,0.23cm)  {{\footnotesize Drive cell $(5,3)$}};
        \node[right=5pt] at (8,0.23cm)  {{\footnotesize$(\entrancex, \entrancey)$}};
    \end{tikzpicture}
    \vspace{-2mm}
\caption{Pink nodes show vertex separators generated at different hop distances from the entrance (purple). Driving (gray) and parking fields (green and blue) can be included in the valid inequality if the neighbors (yellow and orange) are not in the component containing the entrance.}
\label{fig:reverse_hop_connectivity}
\end{figure}
\textbf{Reverse hop inequalities:} To help with connectivity closer to the entrance, we also use a reverse version of hop inequalities by considering vertex cut sets that are a certain hop distance away from the entrance cell. Specifically, we set $\vertexcutset = H_{\entrancex \entrancey}^d$ in \eqref{eq:strongcutconstraints}, which consists of cells that are $d$ hops away from the entrance $(\entrancex, \entrancey)$, where $d = \{1, 2, \ldots, \pi_{max}\}$ and $\pi_{max}$ is the maximum hop distance from the entrance to other cells in $\validdriveset$. As illustrated in Figure \ref{fig:reverse_hop_connectivity}, with increased hop distance, the vertex cut is closer to the central node, thereby reducing the odds of including additional decision variables on the left-hand side. If removing a vertex cut set creates more than one disconnected component, we create multiple vertex cut sets from it, ensuring that each vertex cut set produces exactly two components of the graph. Then, we enforce the reverse hop inequality. An example instance that satisfies these bidirectional hop inequalities is shown in Figure \ref{fig:after_hops_two_way}. The driveways are not connected, but there are only a few disconnected components. We also generated other families of vertex cut sets using rows and columns of the grid, or by successive applications of minimum-cardinality vertex separators. However, they were found to be less effective than bidirectional hop inequalities.

\subsubsection{Cut separation procedure}
\label{sec:disc_infeasible_two_way}
The integer feasible solutions with the earlier inequalities may be disconnected. Hence, we discard them using a series of vertex contractions and minimum-cardinality vertex-cut problems to dynamically generate additional connectivity cut constraints. In our implementation, they are added using CPLEX's lazy constraint callback. These separation problems can be solved in polynomial time using the Edmonds-Karp algorithm \citep{edmonds1972theoretical}. Consider an integer feasible solution with disconnected driveways and active driving fields represented by $D^{int}$. Let $G(D^{int})$ denote the subgraph induced by $D^{int}$, consisting of the node set $D^{int}$ and all arcs in $\grid$ whose endpoints both lie in $D^{int}$. First, we identify the sets of nodes in the connected components of $\grid(D^{int})$, denoted by $\mathcal{C}(D^{int})$. Let the nodes in the component containing the entrance be represented as $D^{int}_{pq}$. We then pick a component $C \in \mathcal{C}(D^{int}) \setminus D^{int}_{pq}$ and perform two sets of node contractions. (1) We contract the nodes in the set $\{(i, j) : (i, j) \in D^{int}_{pq} \cup (\mathcal{C}(D^{int}) \setminus C) \}$ to create a super-source. (2) We also contract the nodes in $C$ to create a super-sink.

Finally, we solve a minimum cardinality vertex-cut problem between the super-source and the super-sink to obtain a vertex cut set, say, $V$. We can then write connectivity cut constraints to avoid the current solution using \eqref{eq:strongcutconstraints_feasibility} for all cells in $C$. This process is then repeated for every component $C \in \mathcal{C}(D^{int}) \setminus D^{int}_{pq}$. 
\begin{align}
     & \drivevar{i}{j} + \sum_{(k, l) \in \invzeroparksetcut{i}{j}{{\twowaypartitioncomp{pq} \cup \vertexcutset}}} \zeroparkvar{k}{l} + \sum_{(k, l) \in \invninetyparksetcut{i}{j}{\twowaypartitioncomp{pq} \cup \vertexcutset}} \ninetyparkvar{k}{l} \leq \sum_{(k,l) \in \vertexcutset} \drivevar{k}{l} \qquad && \forall \, (i, j) \in C   \label{eq:strongcutconstraints_feasibility}
\end{align}

\begin{figure}[H]
  \centering
  \begin{subfigure}[b]{0.24\textwidth}
      \centering
    \includegraphics[width=\textwidth]{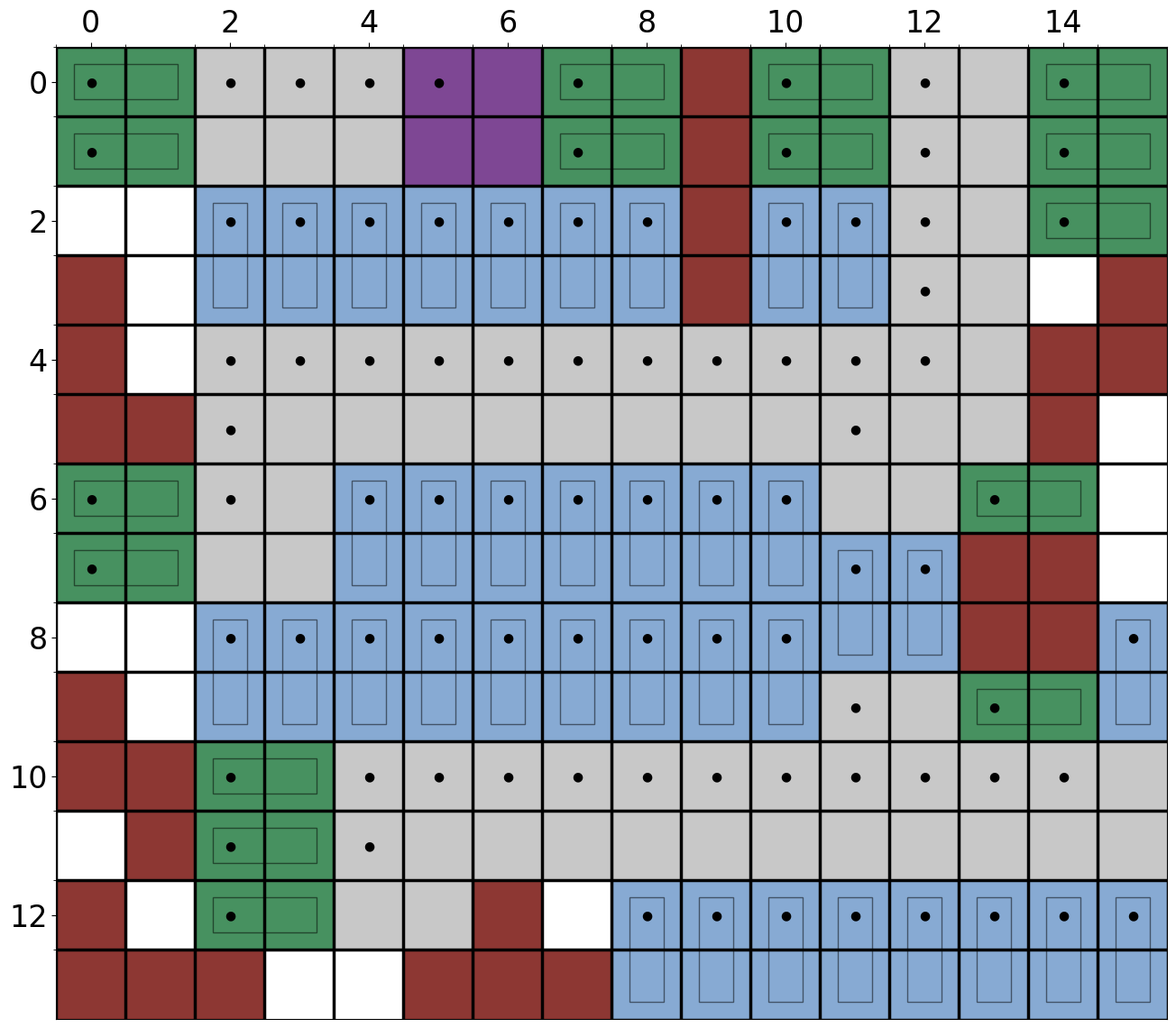}
    \caption{Infeasible layout}
    \label{fig:after_hops_two_way}
  \end{subfigure}
  \hfill
    \begin{subfigure}[b]{0.24\textwidth}
      \centering
    \includegraphics[width=\textwidth]{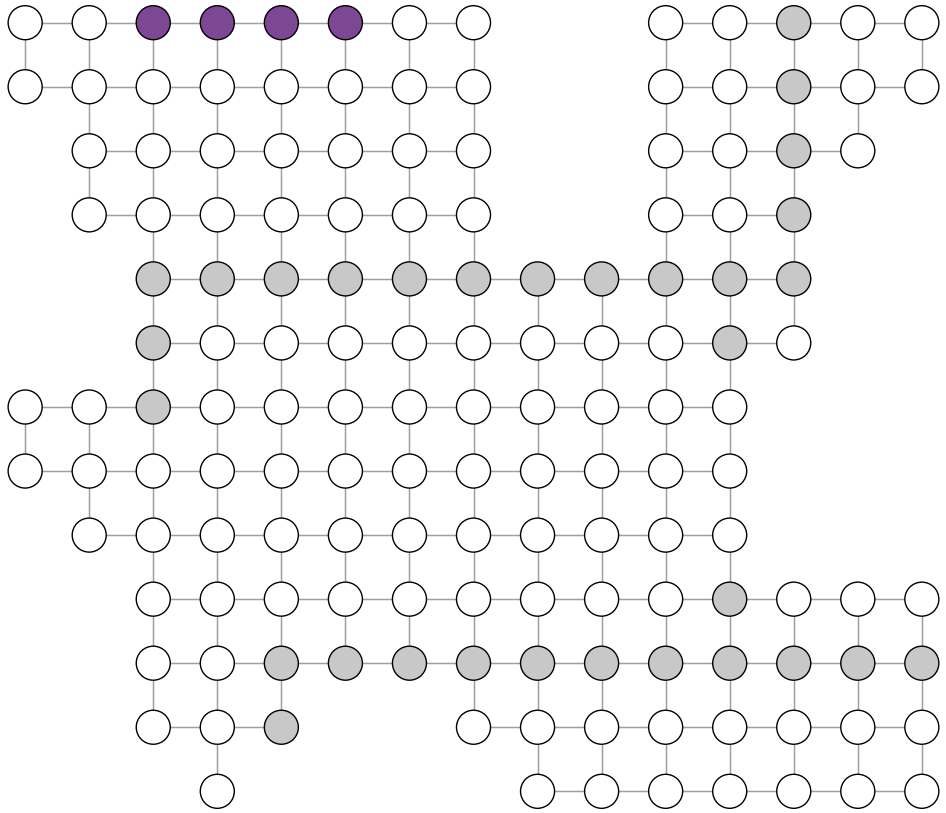}
    \caption{Connected components}
    \label{fig:connected_components_two_way}
  \end{subfigure}
  \hfill
  \begin{subfigure}[b]{0.49\textwidth}
      \centering
    \includegraphics[width=0.49\textwidth]{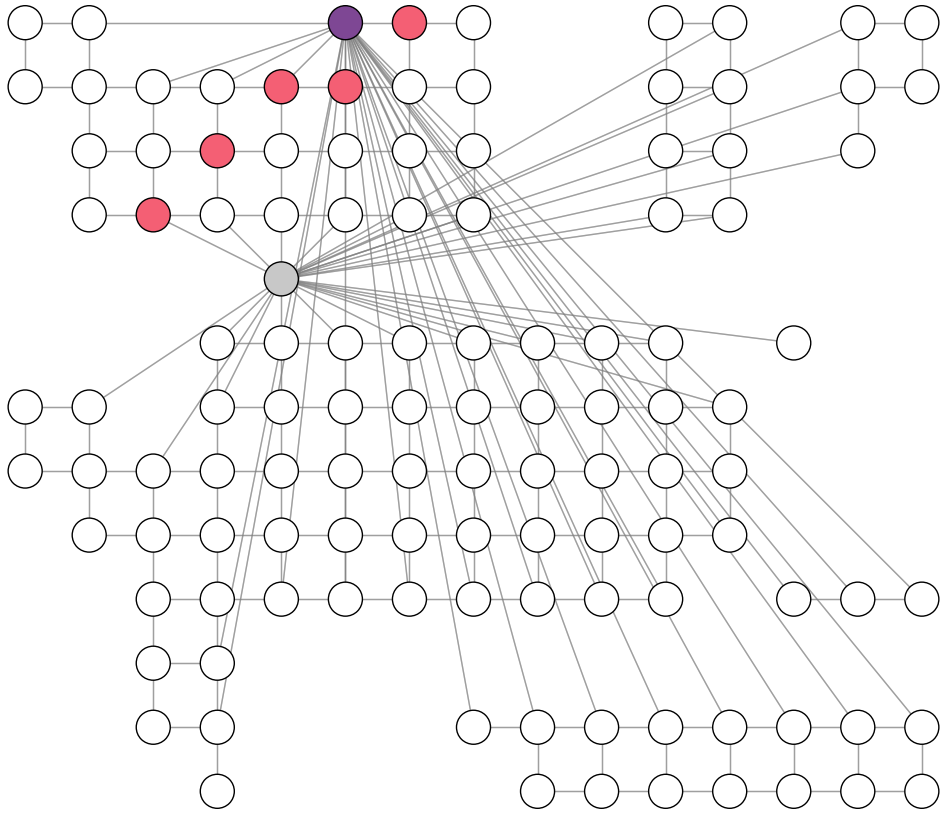}
    \includegraphics[width=0.49\textwidth]{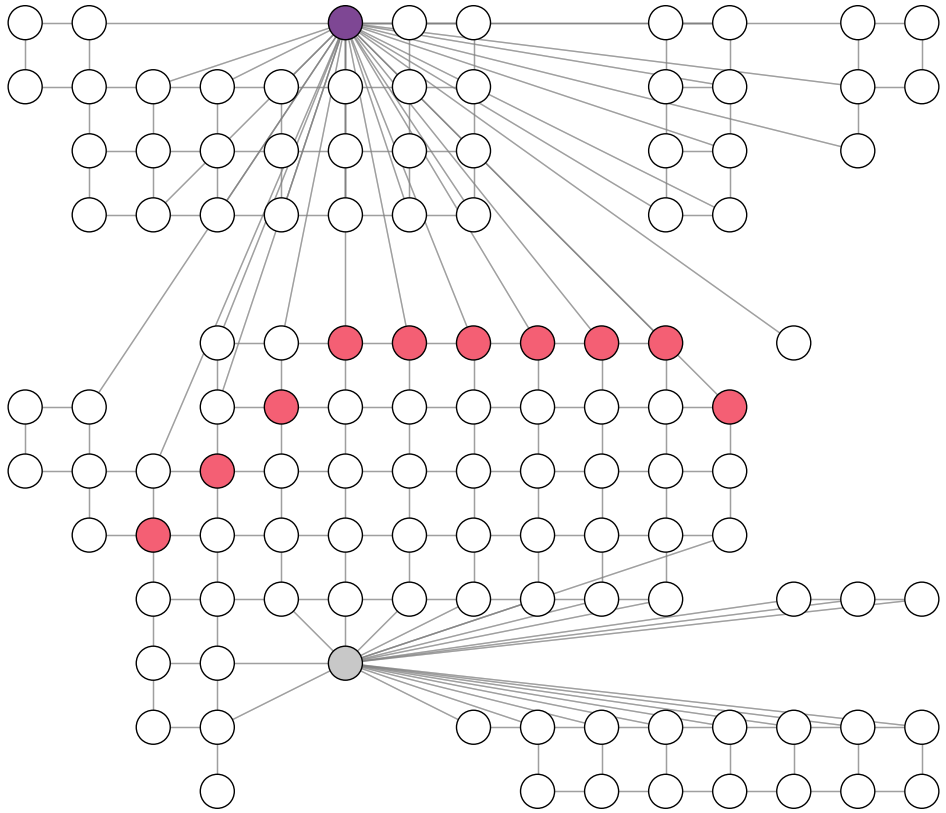}
    \caption{Vertex separators}
    \label{fig:vertex_separators}
  \end{subfigure}
\caption{Generating vertex cut sets within the branch-and-cut algorithm}
\label{fig:disc_infeasible}
\end{figure}
Note that node contractions are necessary before finding the vertex cut set. Otherwise, we might find a cut set that contains active driving nodes from one of the other components. Consequently, we would not discard the current solution because the right-hand side of \eqref{eq:strongcutconstraints_feasibility} would be greater than or equal to one. Figure \ref{fig:disc_infeasible} shows the steps in discovering these feasibility cuts. The infeasible layout has three connected components. The nodes in $D^{int}_{pq}$ are indicated in purple in Figure \ref{fig:connected_components_two_way}. The other driveway components in $\mathcal{C}(D^{int})$ are colored gray. We contract the union of all purple nodes and the bottom component of the driveway graph to create a super-source and contract the middle component to generate a super-sink. The resulting vertex cut set between them is shown in the left panel of Figure \ref{fig:vertex_separators}. 

A similar process is repeated by contracting the purple nodes and the nodes in the middle component. The super-source and the contracted super-sink are then used to find a minimum vertex cut set, which is shown using the pink nodes in Figure \ref{fig:vertex_separators}. As seen in the case of hop inequalities, it is sometimes possible to get a vertex separator that splits the contracted graphs into more than two components. In such cases, to make the cuts tighter, we shrink the vertex separator by only including the nodes that separate the entrance from the elements in $C$. Note that just one of the above-mentioned cuts is sufficient to avoid the current solution. However, adding multiple constraints and letting the solver decide where else to apply them can tighten the LP relaxation along the branch-and-bound tree.

\subsection{One-way driveways}
\label{sec:one-way-branch-and-cut}
The branch-and-cut algorithm can also be designed for the problem of parking with one-way streets using only binary variables. The connectivity cut inequalities discussed earlier are also valid for this variant and ensure that each driving field is connected to the entrance $(\entrancex, \entrancey)$ and the exit $(\exitx, \exity)$. For example, a solution in Figure \ref{fig:one_way_undirected} with $\parkwidth=\drivewidth=1$, and $ \parklength=2$, satisfies the bidirectional hop constraints and restrictions on the $z$ variables that determine the directions of the driveways. However, this is not a feasible one-way configuration, as it requires us to backtrack along the path from the entrance to reach the exit. Note that there is always an optimal solution that does not use a direct connection between the entrance and exit cells, e.g., (0,5) to (0,6); therefore, these arcs can be removed from the graph $\grid$.

\begin{figure}[H]
    \centering
    \begin{subfigure}[b]{0.49\textwidth}
        \centering
        \includegraphics[width=0.63\textwidth]{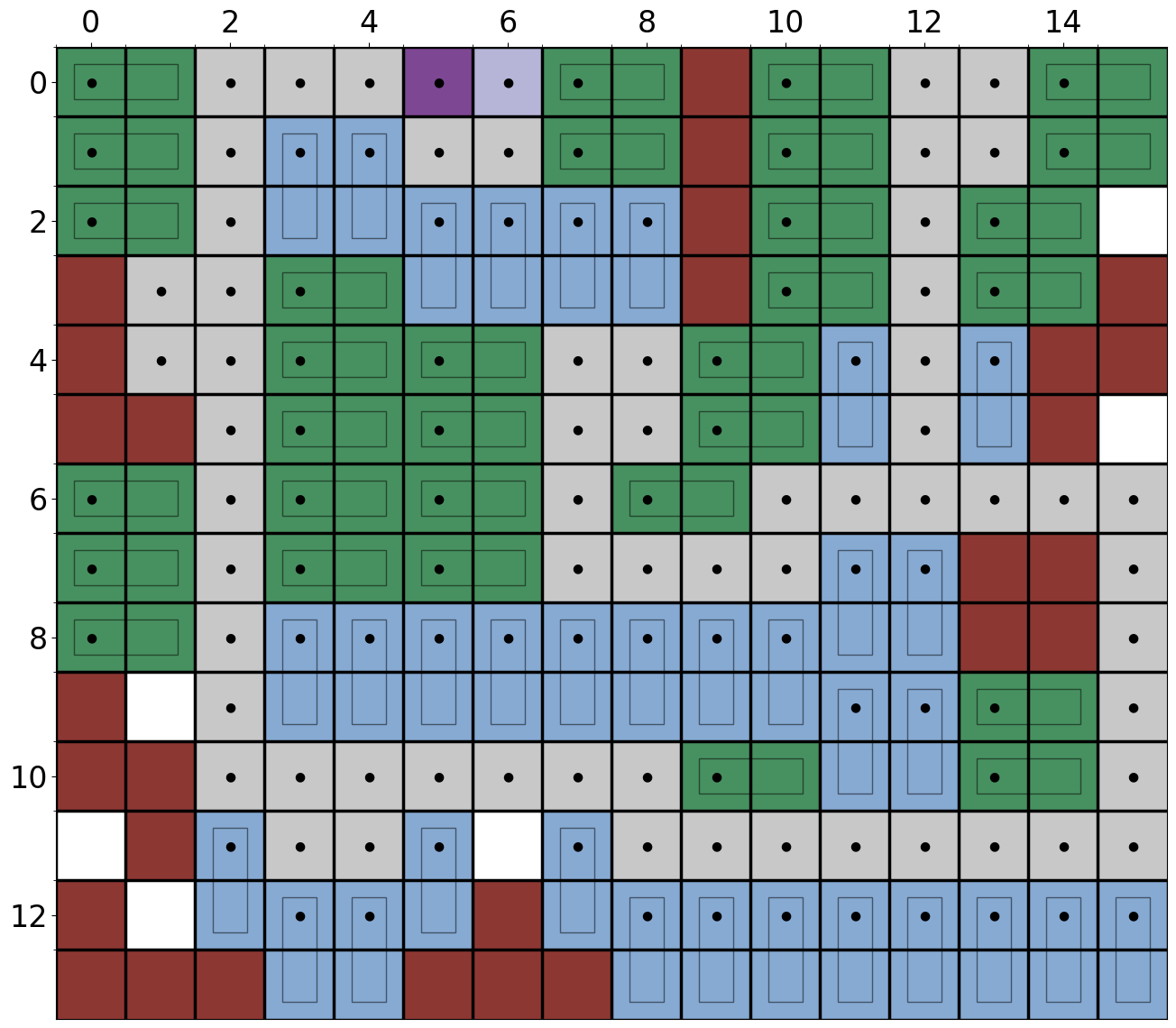}
        \caption{Layout with bidirectional hop inequalities}
    \end{subfigure} 
    \hfill
    \begin{subfigure}[b]{0.49\textwidth}
        \centering
        \includegraphics[width=0.63\textwidth]{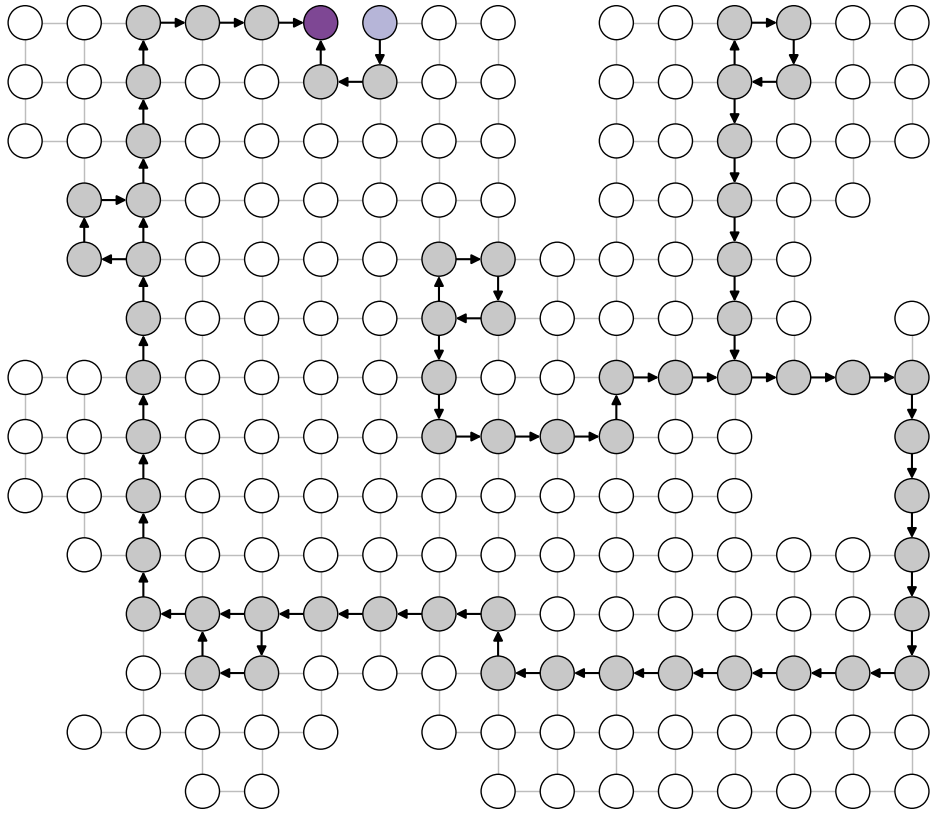}
        \caption{Direction variables $z$}
    \end{subfigure}     
       \caption{Infeasible solutions from the two-way branch-and-cut approach for the one-way scenario}
    \label{fig:one_way_undirected}
\end{figure}

The primary issue is that the previously defined connectivity cut constraints operate in an undirected sense. To fix this limitation, we redefine them using \textit{edge separators} or \textit{edge cut sets} that partition the graph of driving fields into two strongly connected components with node sets $\onewaypartition \in \onewaypartitionsets{ij}$, which contain a valid driving cell $(i, j)$ and its complement $ \onewaypartitioncomp \neq \emptyset$. The set $\onewaypartitionsets{ij}$ represents all possible subsets of valid drive cells that contain $(i, j)$, each of which has nodes that are strongly connected. Figure~\ref{fig:one_way_undirected} shows an example grid with two edge cut sets (highlighted in pink) that separate the entrance and exit cells (in purple) from the teal-colored drive cells. The teal drive cells cannot be active if the $z$ variables associated with either of the two edge cut sets are zero.

\begin{figure}[H]
    \centering
    \begin{subfigure}[b]{0.32\textwidth}
        \centering
        \includegraphics[width=0.98\textwidth]{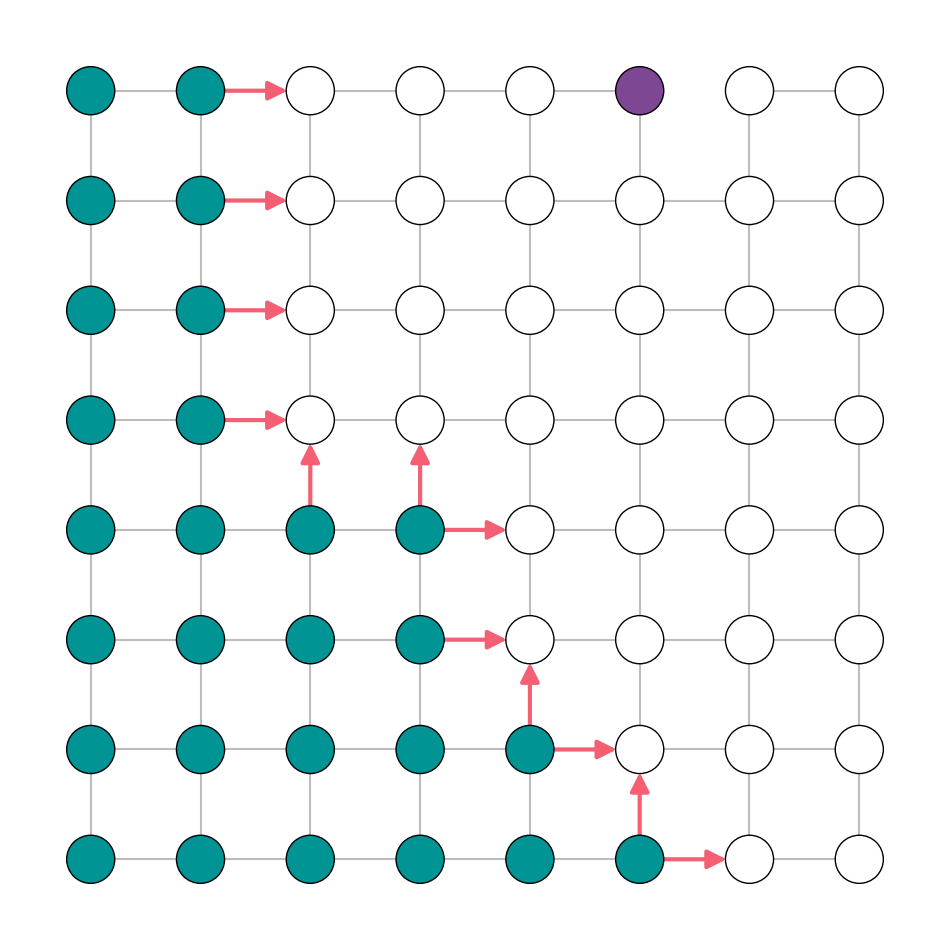}
    \end{subfigure}
    \hfill
    \begin{subfigure}[b]{0.32\textwidth}
        \centering
        \includegraphics[width=0.98\textwidth]{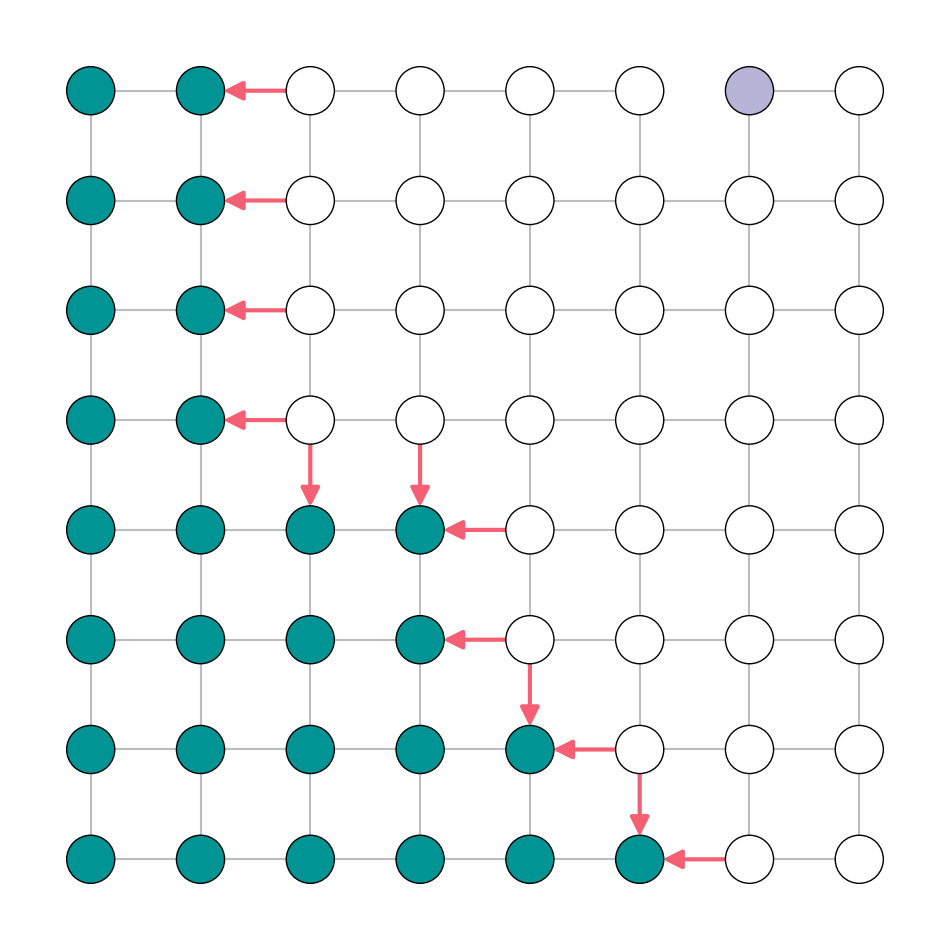}
    \end{subfigure} 
    \hfill
    \begin{subfigure}[b]{0.32\textwidth}
        \begin{tikzpicture}[baseline=-80pt]
        \definecolor{pink}{RGB}{244,95,116}
        \definecolor{teal}{RGB}{0,148,148}
        \definecolor{purple}{RGB}{126,71,148}
        \definecolor{lavender}{RGB}{182,181,216}
        
\draw[->, thick, pink, >=stealth] (-0.25cm,1.6cm) -- (0.25cm,1.6cm);
\filldraw[fill=teal, draw=black] (0,0.8cm) circle (0.14cm);
\filldraw[fill=purple, draw=black] (0,0cm) circle (0.14cm);
\filldraw[fill=lavender, draw=black] (0,-0.8cm) circle (0.14cm);
\filldraw[fill=white, draw=black] (0,-1.6cm) circle (0.14cm);

\node[right=5pt] at (0,1.6cm) {\footnotesize{Arcs in the edge cut sets}};
\node[right=5pt] at (0,0.8cm) {\footnotesize{Drive cells in $\onewaypartitioncomp$}};
\node[right=5pt] at (0,0cm) {\footnotesize{Entrance $(\entrancex, \entrancey)$}};
\node[right=5pt] at (0,-0.8cm) {\footnotesize{Exit $(\exitx, \exity)$}};
\node[right=5pt] at (0,-1.6cm) {\footnotesize{Drive cells in $\onewaypartition \in \onewaypartitionsets{pq}, \onewaypartitionsets{rs}$}};

        \end{tikzpicture}
    \end{subfigure}
    \caption{Illustration of edge cut sets}
    \label{fig:one_way_undirected}
\end{figure}

Constraints \eqref{eq:cutconstraint_entrance} and \eqref{eq:cutconstraint_exit} force the driving field variables in $\onewaypartitioncomp$ to zero if they cannot access the entrance or cannot be reached by the exit. This situation can arise when the $z$ variables for specific edge sets between the two components $\onewaypartition$ and $\onewaypartitioncomp$ are inactive. 
\begin{align}
    & \drivevar{m}{n} \leq \sum_{\substack{((i,j),(k,l)) \in \gridarcset: \\ (i,j) \in \onewaypartitioncomp, (k, l) \in \onewaypartition}} \countvar{i}{j}{k}{l} \qquad && \forall \,  \onewaypartition \in \onewaypartitionsets{\entrancex \entrancey}, (m, n) \in \onewaypartitioncomp \label{eq:cutconstraint_entrance}\\
    & \drivevar{m}{n} \leq  \sum_{\substack{((i,j),(k,l)) \in \gridarcset: \\(i,j) \in \onewaypartition, (k, l) \in \onewaypartitioncomp}} \countvar{i}{j}{k}{l} \qquad && \forall \, \onewaypartition \in \onewaypartitionsets{\exitx \exity}, (m, n) \in \onewaypartitioncomp \label{eq:cutconstraint_exit}
\end{align}

\begin{proposition}
\label{prop:1w_aggregate}
Constraints \eqref{eq:cutconstraint_entrance} and \eqref{eq:cutconstraint_exit} are valid for the feasible region of the formulation $\mathcal{F}^{\textsc{flow}}_{\textsc{1W}}$.
\end{proposition}
\begin{proof}
Suppose $\onewaypartition \in \onewaypartitionsets{\entrancex \entrancey}$. Adding the flow-conservation constraints \eqref{eq:floworigin} for all cells in $(m, n) \in \onewaypartitioncomp$, 
\begin{align}
\sum_{(m,n) \in \onewaypartitioncomp} \drivevar{m}{n} & = \sum_{\substack{((i,j),(k,l)) \in \gridarcset: \\(i,j) \in \onewaypartitioncomp, (k, l) \in \onewaypartition}} f_{ij, kl} -
\sum_{\substack{((i,j),(k,l)) \in \gridarcset: \\(k, l) \in \onewaypartition, (i,j) \in \onewaypartitioncomp}} f_{kl, ij} \notag \\
& \leq \sum_{\substack{((i,j),(k,l)) \in \gridarcset: \\(i,j) \in \onewaypartitioncomp, (k, l) \in \onewaypartition}} f_{ij, kl} \notag \\
& \leq \sum_{\substack{((i,j),(k,l)) \in \gridarcset: \\(i,j) \in \onewaypartitioncomp, (k, l) \in \onewaypartition}} M \countvar{i}{j}{k}{l} \qquad \text{[using \eqref{eq:count_var_bounds}]} \label{eq:aggregate_vi_1w_a}
\end{align}
Since $y$ and $z$ variables are binary, using a contrapositive argument, we can conclude that \eqref{eq:cutconstraint_entrance} must hold, in a manner similar to Proposition~\ref{prop:2w_disaggregate}.

Likewise, when $\onewaypartition \in \onewaypartitionsets{\exitx \exity}$, adding the flow-conservation constraints \eqref{eq:flow_dash_origin_singlelane} for all cells in $(m, n) \in \onewaypartitioncomp$, 
\begin{align}
\sum_{(m,n) \in \onewaypartitioncomp} \drivevar{m}{n} & = \sum_{\substack{((i,j),(k,l)) \in \gridarcset: \\(i, j) \in \onewaypartition, (k,l) \in \onewaypartitioncomp}} g_{ij, kl} - \sum_{\substack{((i,j),(k,l)) \in \gridarcset: \\(k,l) \in \onewaypartitioncomp, (i, j) \in \onewaypartition}} g_{kl, ij} \notag \\
& \leq \sum_{\substack{((i,j),(k,l)) \in \gridarcset: \\(i, j) \in \onewaypartition, (k,l) \in \onewaypartitioncomp}} g_{ij, kl}  \notag \\
& \leq \sum_{\substack{((i,j),(k,l)) \in \gridarcset: \\(i, j) \in \onewaypartition, (k,l) \in \onewaypartitioncomp}} M \countvar{i}{j}{k}{l} \qquad \text{[using \eqref{eq:count_dash_var_bounds}]} \label{eq:aggregate_vi_1w_b}
\end{align}
Again using a contrapositive argument, we can show that \eqref{eq:cutconstraint_exit} is valid for the feasible region of formulation  $\mathcal{F}^{\textsc{flow}}_{\textsc{1W}}$.
\end{proof}

The feasibility cuts can also be strengthened by including specific parking field variables on the left-hand side of the inequalities whenever possible, thereby deriving a new formulation $\mathcal{F}^{\textsc{CCC}}_{\textsc{1W}}$. The proof of validity for these strengthened valid inequalities \eqref{eq:cutconstraint_entrance_strengthened} and \eqref{eq:cutconstraint_exit_strengthened} parallels that of Proposition \ref{prop:2w_strengthened}. The equivalence of the integer-feasible solutions $(x, y, z)$ for the flow-based and cut-based formulations is proven in \ref{sec:equivalence}.
\begin{align}
    \max & \sum_{(i,j) \in \validzeroparkset} \zeroparkvar{i}{j} + \sum_{(i,j) \in \validninetyparkset} \ninetyparkvar{i}{j}  && [\text{Formulation } \mathcal{F}^{\textsc{CCC}}_{\textsc{1W}}] \notag\\
     \text{s.t. } &\eqref{eq:singlepurpose_zero_disaggregate}\text{--}\eqref{eq:ninetyparkdriveconnector}, 
     \eqref{eq: x_zero is binary}\text{--}\eqref{eq: y var is binary}, 
     \eqref{eq:count_drive1} \text{--} \eqref{eq:binary_countvar} && \notag \\
      & \drivevar{u}{v} + \sum_{(k, l) \in \invzeroparksetcut{m}{n}{{\onewaypartitioncomp}}} \zeroparkvar{k}{l} + \sum_{(k, l) \in \invninetyparksetcut{m}{n}{\onewaypartitioncomp}} \ninetyparkvar{k}{l} \leq \sum_{\substack{((i,j),(k,l)) \in \gridarcset: \\(i,j) \in \onewaypartitioncomp, (k, l) \in \onewaypartition}} \countvar{i}{j}{k}{l}  && \forall \,  \onewaypartition \in \onewaypartitionsets{\entrancex \entrancey}, (m, n) \in \onewaypartitioncomp, (u, v) \in  \invdriveset{m}{n}(\onewaypartitioncomp) \label{eq:cutconstraint_entrance_strengthened}\\
    & \drivevar{u}{v} + \sum_{(k, l) \in \invzeroparksetcut{m}{n}{\onewaypartitioncomp}} \zeroparkvar{k}{l} + \sum_{(k, l) \in \invninetyparksetcut{m}{n}{\onewaypartitioncomp}} \ninetyparkvar{k}{l} \leq  \sum_{\substack{((i,j),(k,l)) \in \gridarcset: \\(i,j) \in \onewaypartition, (k, l) \in \onewaypartitioncomp}} \countvar{i}{j}{k}{l}  && \forall \, \onewaypartition \in \onewaypartitionsets{\exitx \exity}, (m, n) \in \onewaypartitioncomp, (u, v) \in  \invdriveset{m}{n}(\onewaypartitioncomp)\label{eq:cutconstraint_exit_strengthened}
\end{align}

\subsubsection{Other valid inequalities}
\label{sec:vis_one_way}
As noted in the previous section, the number of feasibility cuts grows exponentially for the one-way case, and this is also true for the new connectivity cut constraints. We take a similar approach as before, generating a limited number of valid inequalities added as regular and lazy constraints.

\textbf{Bidirectional hop constraints:} The constraints from Section \ref{sec:vis_two_way} are also valid for the one-way scenario. One could introduce two types of valid inequalities, where hop distances are calculated with respect to both the entrance and exit. However, this led to an excessive number of constraints and marginal differences in runtime. Hence, we only included bidirectional hop inequalities with respect to the entrance cell. 

\textbf{Dead end prevention constraints:} In the flow formulation version discussed in Section \ref{sec:dead_end}, constraints involving the lane direction variables $z$ and the flow variables $f$ and $g$ prevented dead ends. Alternatively, it is possible to write valid inequalities that eliminate dead ends without using the flow variables. This can be done (1) purely using drive field variables $y$ or (2) using both $y$ and $z$ variables.

For constraints with only $y$ variables, a valid driving cell with a single neighbor cannot anchor an active driving field (see \eqref{eq:one_neigbhor}). In other situations, an active driving cell must have at least two active neighbors to prevent dead ends. The logic for this constraint is captured in \eqref{eq:at_least_two_neigbhors}.
\begin{align}
    &  \drivevar{i}{j} = 0 \qquad && \forall \, (i,j) \in \validdriveset, |\gridadjset{i}{j}| = 1 
     \label{eq:one_neigbhor}\\
     & \sum_{(k, l) \in \gridadjset{i}{j}} \drivevar{k}{l}  \geq 2\drivevar{i}{j} \qquad && \forall \, (i,j) \in  \validdriveset \setminus \{(\entrancex, \entrancey), (\exitx, \exity)\}, |\gridadjset{i}{j}| > 1
     \label{eq:at_least_two_neigbhors}     
\end{align}
Additionally, we can also enforce \eqref{eq:outflow_count_singlelane} and \eqref{eq:inflow_count_singlelane}. These constraints ensure that at least one active outgoing and one active incoming arc exist for all active drive cells, excluding the entrance and exit cells. Proofs of the validity of these inequalities are trivial and are hence omitted. 
\begin{align}
     & \sum_{(k,l) \in \gridadjset{i}{j}} \countvar{i}{j}{k}{l}  \geq \drivevar{i}{j}  \qquad && \forall \, (i,j) \in  \validdriveset \setminus \entranceset
     \label{eq:outflow_count_singlelane}\\
     & \sum_{(k,l) \in \gridinvadjset{i}{j}} \countvar{k}{l}{i}{j} \geq \drivevar{i}{j}  \qquad && \forall \, (i,j) \in  \validdriveset \setminus \exitset
     \label{eq:inflow_count_singlelane}     
\end{align}

\subsubsection{Cut separation procedure}
\label{sec:disc_infeasible_one_way}
\begin{figure}[H]
  \centering
  \begin{subfigure}[b]{0.32\textwidth}
      \centering
    \includegraphics[width=0.95\textwidth]{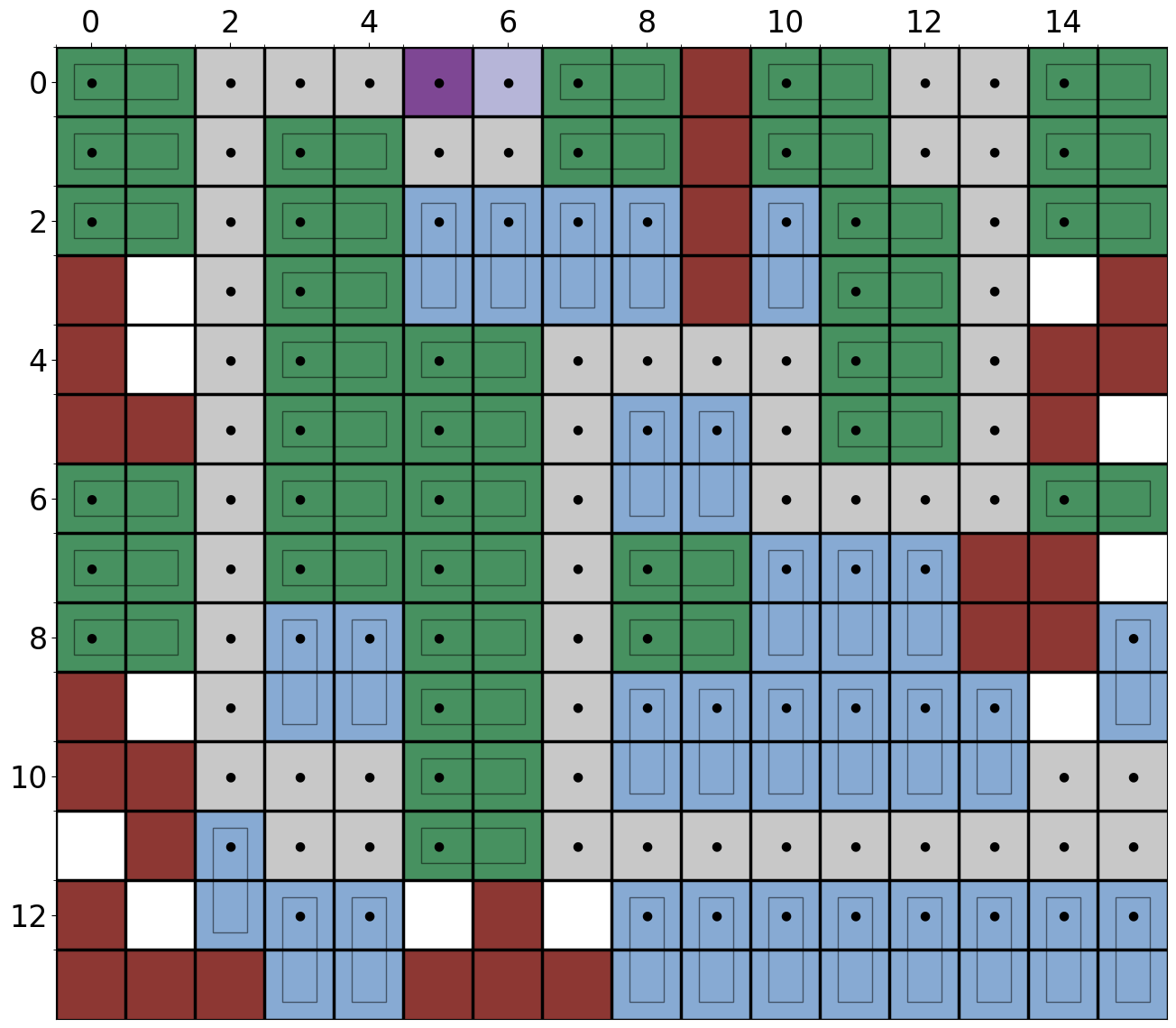}
    \caption{Infeasible layout}
    \label{fig:disc_infeasible2a}
  \end{subfigure}
  \hfill
  \begin{subfigure}[b]{0.32\textwidth}
      \centering
    \includegraphics[width=0.95\textwidth]{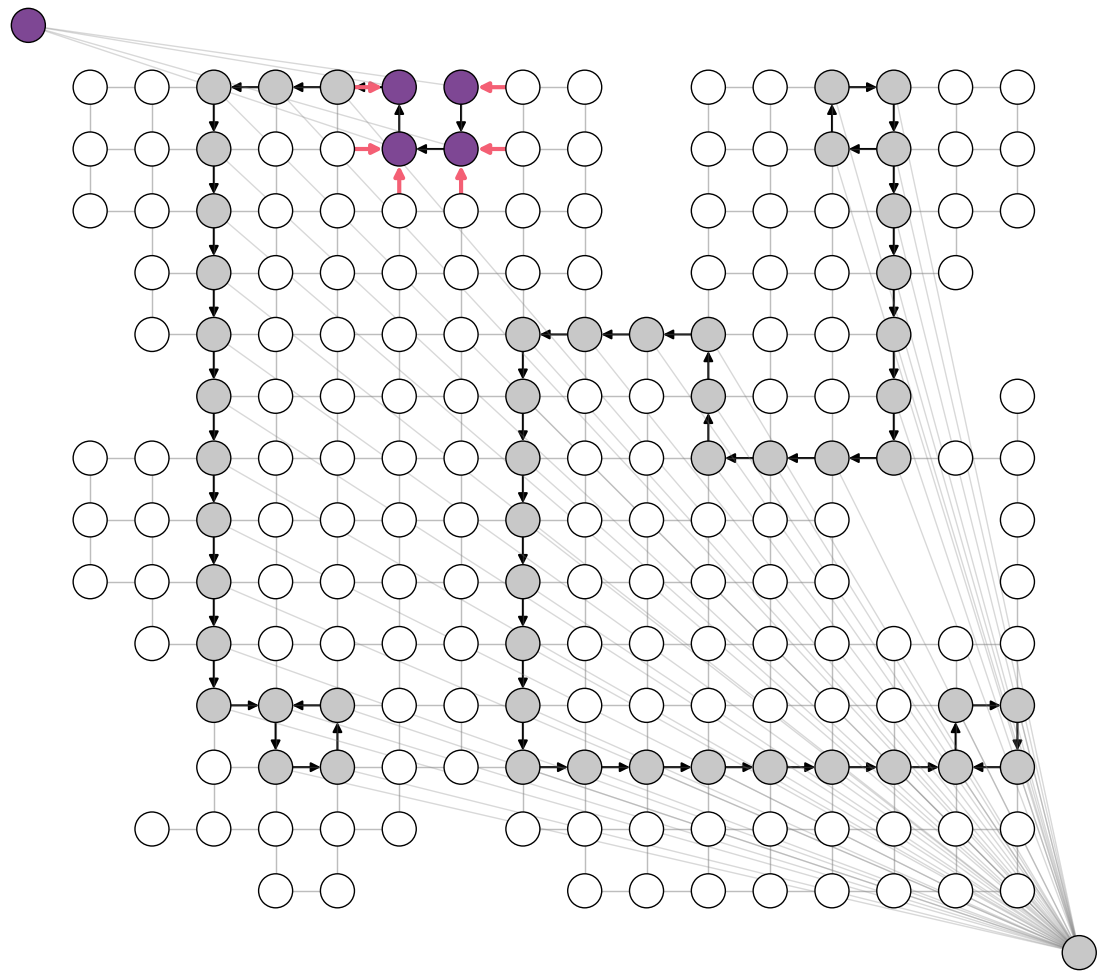}
    \caption{Connectivity cut for the entrance}
    \label{fig:disc_infeasible2b}
  \end{subfigure}
  \hfill
  \begin{subfigure}[b]{0.32\textwidth}
      \centering
    \includegraphics[width=0.95\textwidth]{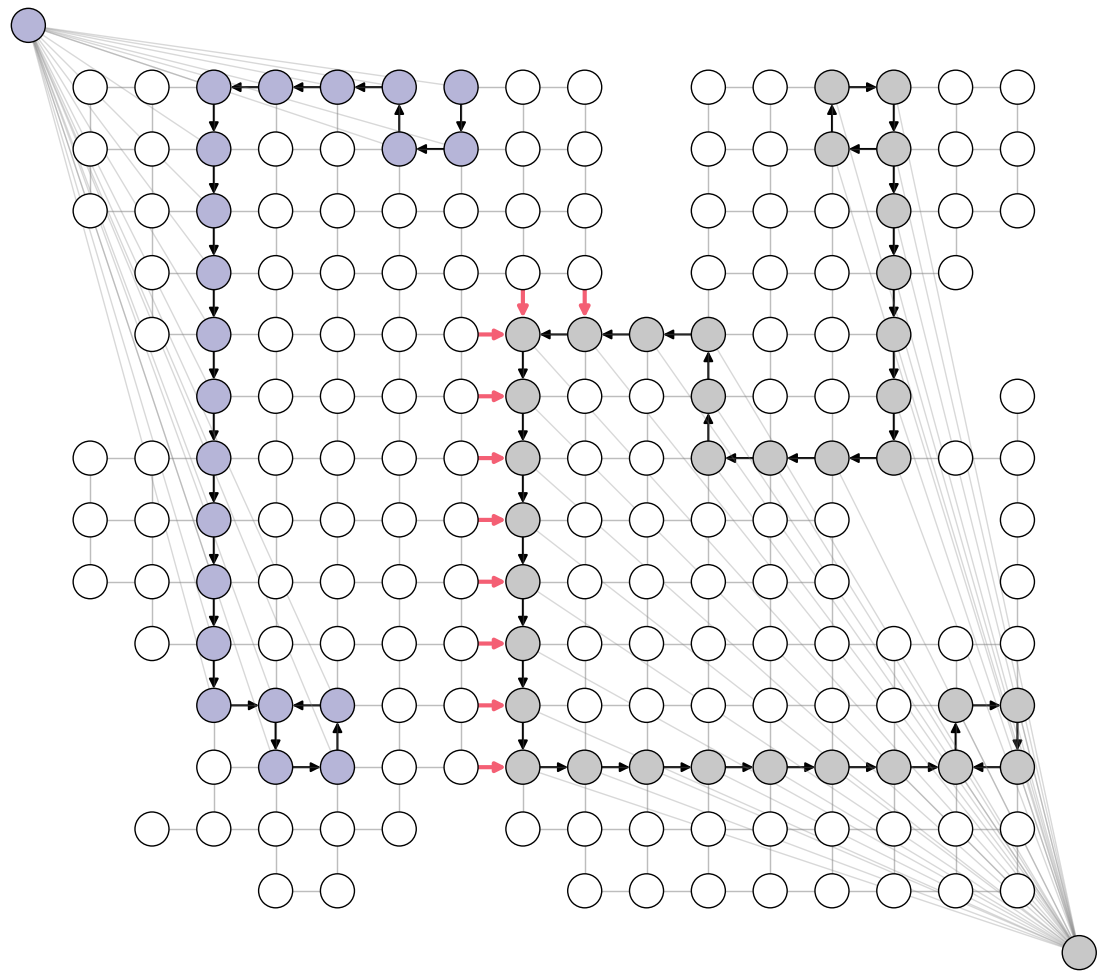}
    \caption{Connectivity cut for the exit}
        \label{fig:disc_infeasible2c}
  \end{subfigure}
\caption{Generating edge cut sets within the branch-and-cut algorithm}
\label{fig:disc_infeasible2}
\end{figure}

Adding the above inequalities can still result in layouts with few disconnected components (e.g., see Figure \ref{fig:disc_infeasible2a}). We generate violated constraints and add them to the branch-and-bound process as a remedy. This step is done separately for the entrance and exit. We again use CPLEX's lazy constraint callback feature for implementing this step. As before, let $D^{int}$ indicate the driving cells currently active in the infeasible solution. We let $D^{int}_{\entrancex \entrancey}$ denote the driving cells that have a directed path to the entrance cell and $D^{int}_{\exitx \exity}$ be the set of active driving cells that can be reached via directed paths from the exit. These sets are shown using the dark and light purple nodes in Figures \ref{fig:disc_infeasible2b} and \ref{fig:disc_infeasible2c}, respectively.

To ensure the connectivity of driveways to the entrance, if $|D^{int}_{\entrancex \entrancey}| < |D^{int}|$, we first assign weights to the arcs in $\grid$. Arcs connecting active drive cells, i.e., the ones for which the $z$ variables are one, are assigned a weight $M$. Every other arc is assigned a weight of one. We then create a new super-source and connect arcs from it to all nodes in $D^{int}_{\entrancex \entrancey}$. Similarly, we connect all other active driving cells $D^{int} \setminus D^{int}_{\entrancex \entrancey}$ to a newly created super-sink. The arcs incident to the super-source and super-sink are also assigned a weight of $M$. In Figure \ref{fig:disc_infeasible2b}, the purple and gray nodes at the top left and bottom right represent the super-source and super-sink, respectively. We then solve a minimum weighted edge cut problem between these nodes to find partitions $\onewaypartitioncomp_{\entrancex \entrancey}$ and $\onewaypartition_{\entrancex \entrancey}$ that separate $D^{int} \setminus D^{int}_{\entrancex \entrancey}$ and $D^{int}_{\entrancex \entrancey}$, respectively. By assigning a weight of one to arcs with inactive $z$ variables and large weights to the rest, this edge cut set has minimum cardinality (see the pink arcs in the example). Consequently, the resulting feasibility cut will feature fewer terms on the right-hand side, enhancing its effectiveness. The number of arcs in $\gridarcset$ is a good choice for $M$. Note that $G(D^{int} \setminus D^{int}_{\entrancex \entrancey})$ can have more than one connected component, as seen in Figure \ref{fig:disc_infeasible2b}. In such situations, multiple edge cut sets can be generated to separate each component from the rest, but this would require more computational effort. The partitions and edge cut set generate a feasibility cut, which can be strengthened with parking variables as shown in \eqref{eq:strongedgecutconstraints1}.
\begin{align}
     & \drivevar{m}{n} + \sum_{(k, l) \in \invzeroparksetcut{m}{n}{{\onewaypartitioncomp_{\entrancex \entrancey}}}} \zeroparkvar{k}{l} + \sum_{(k, l) \in \invninetyparksetcut{m}{n}{\onewaypartitioncomp_{\entrancex \entrancey}}} \ninetyparkvar{k}{l} \leq \sum_{\substack{((i,j),(k,l)) \in \gridarcset: \\(i,j) \in \onewaypartitioncomp_{\entrancex \entrancey}, (k, l) \in \onewaypartition_{\entrancex \entrancey}}} \countvar{i}{j}{k}{l}\qquad && \forall \, (m, n) \in D^{int} \setminus D^{int}_{\entrancex \entrancey} \label{eq:strongedgecutconstraints1}
\end{align}

We follow a similar approach to connect the exit to valid driving cells. We first check if $|D^{int}_{\exitx \exity}| < |D^{int}|$. If so, a new super-source is created and is linked to all active driving cells in $D^{int} \setminus D^{int}_{\exitx \exity}$. We then link all driving cells that can be accessed from the exit $D^{int}_{\exitx \exity}$ to a new super-sink.  As before, arc weights are set to one when the $z$ variables are inactive and to $M$ for others to find a minimum weighted edge cut. Figure \ref{fig:disc_infeasible2c} shows an example of this procedure. The gray node at the bottom right and the purple node at the top left depict the super-source and super-sink, respectively. The resulting edge cut set is shown in pink. Let the partitions that separate $D^{int}_{\exitx \exity}$ and $D^{int} \setminus D^{int}_{\exitx \exity}$ be denoted as $\onewaypartition_{\exitx \exity}$ and $\onewaypartitioncomp_{\exitx \exity}$, respectively. Using this notation, adding \eqref{eq:strongedgecutconstraints2} to the branch-and-bound algorithm will help discard the current infeasible solution.
\begin{align}
     & \drivevar{m}{n} + \sum_{(k, l) \in \invzeroparksetcut{m}{n}{\onewaypartitioncomp_{\exitx \exity}}} \zeroparkvar{k}{l} + \sum_{(k, l) \in \invninetyparksetcut{m}{n}{\onewaypartitioncomp_{\exitx \exity}}} \ninetyparkvar{k}{l} \leq \sum_{\substack{((i,j),(k,l)) \in \gridarcset: \\(i,j) \in \onewaypartition_{\exitx \exity}, (k, l) \in \onewaypartitioncomp_{\exitx \exity}}} \countvar{i}{j}{k}{l} \qquad && \forall \, (m, n) \in D^{int} \setminus D^{int}_{\exitx \exity}  \label{eq:strongedgecutconstraints2}
\end{align}

\section{Experiments}
\label{sec:mipresults}

\subsection{Dataset}
\begin{figure}[t]
  \centering
  \hfill
    \begin{subfigure}{0.36\textwidth}
    \centering
    \includegraphics[width=0.9\textwidth]{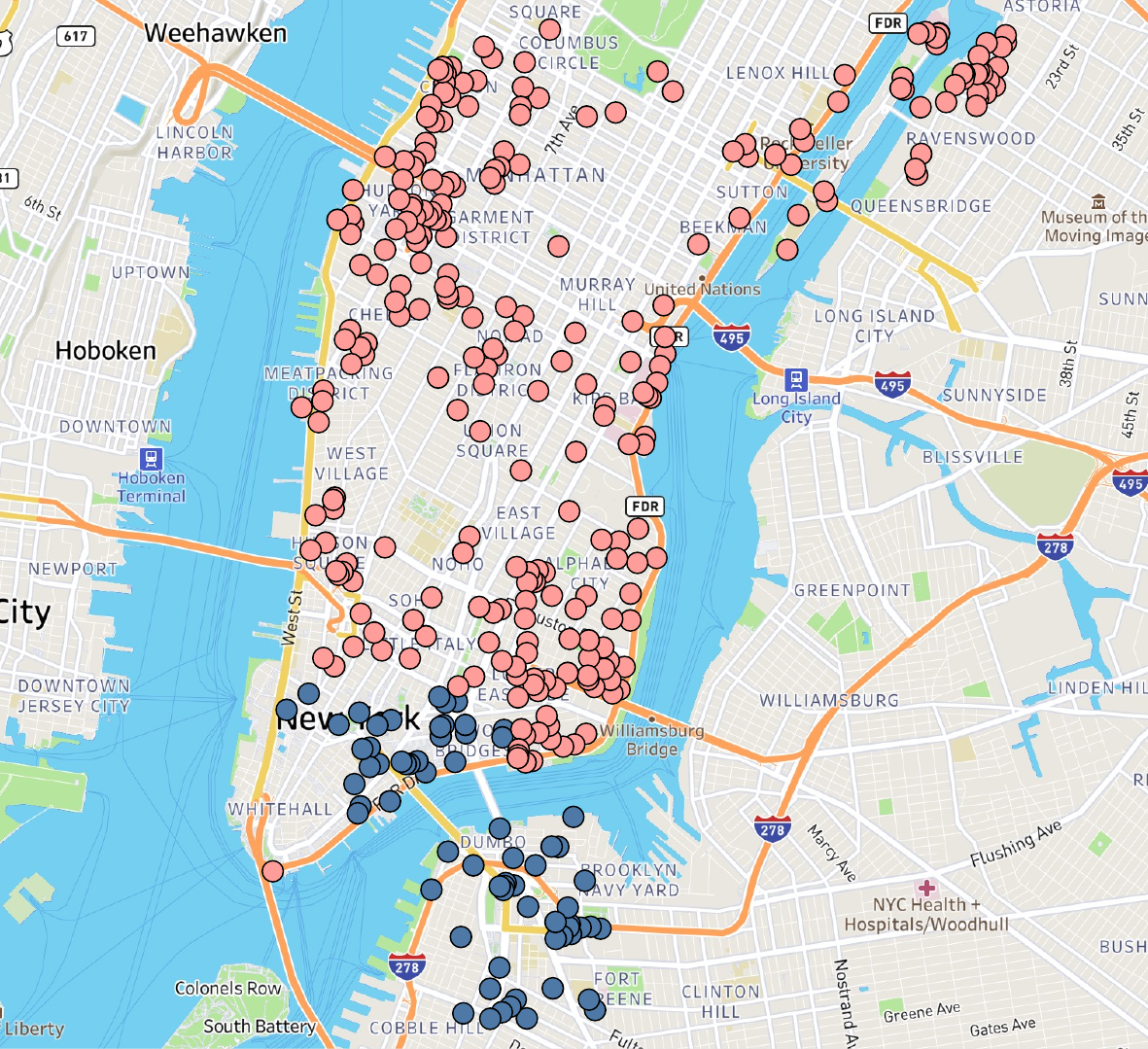}
    \caption{NYC lot locations}
        \label{fig:nyc_downtown}
  \end{subfigure}  
  \begin{subfigure}{0.62\textwidth}
    \centering
    \includegraphics[width=0.86\textwidth]{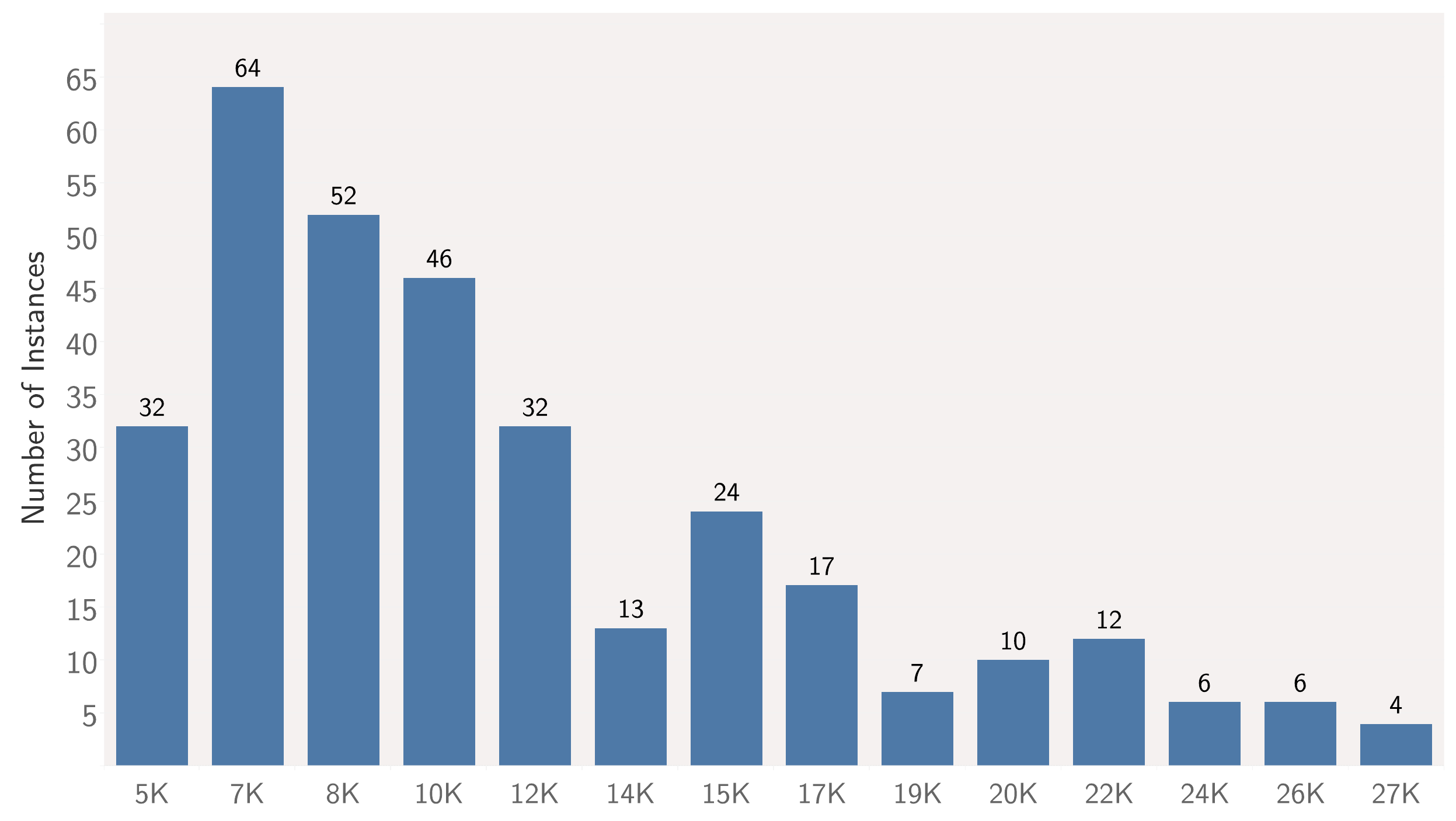}
    \caption{Area (in sq ft) distribution of parking lots}
        \label{fig:nyc_area_distribution}
  \end{subfigure}
\caption{Real-world instances of parking lots from New York City}
\end{figure}

We tested the performance of the MIP models for 325 real-world parking lots from \href{https://data.cityofnewyork.us/City-Government/Parking-Lot/h7zy-iq3d}{New York City (NYC) Open Data}. These lots correspond to open-air parking lots in lower and midtown Manhattan, as well as parts of Brooklyn and Queens, and have areas ranging from 6011 sq ft to 28436 sq ft (roughly up to half the size of a football field). The chosen lots have IDs in the 125 and 215 series, and their locations are illustrated in different colors in Figure \ref{fig:nyc_downtown}. The area distribution is shown in Figure \ref{fig:nyc_area_distribution}. Although data on the exact number of stalls in these parking lots is unavailable, we study the utility of models from a computational standpoint and analyze the differences between two-way and one-way layouts. Experiments on exact methods in \citep{stephan2021layout} were run on slightly smaller layouts in the range 6458--19375 sq ft, although heuristics were tested on larger instances in their paper. 

We designate the parking lot's entrance and exit points to be next to each other. For rasterization, we considered 3 m wide driveways. The size of the parking stalls was set to $6$ m $\times$ $3$ m, i.e., a parking field occupies two cells. At the polygon's boundaries, we can apply a threshold-based rule to include a cell in the parking lot, depending on the fraction of its area that is contained within the polygon. For the experiments in this section, any cell that overlaps with the polygon is assumed to be part of the parking lot, and the other cells are marked as blocked cells.

All experiments were conducted on a Dell Precision workstation with a 32-core AMD Ryzen Threadripper Pro 7975WX CPU @ 4.0 GHz and 256 GB of RAM. The algorithms were coded in Python 3.10 with CPLEX 22.1.1, utilizing the callback functionality for the branch-and-cut algorithm and NetworkX for finding vertex and edge cut sets. Branching was performed using CPLEX’s default strategies. A time limit of 900 s is set for each instance. Experiments were also conducted on a subset of instances with a time budget of two hours. Additional details supporting the main results are provided in \ref{sec:add_stats}. The source codes are available at \href{https://github.com/transnetlab/parking-lot-design.git}{github.com/transnetlab/parking-lot-design}.

\subsection{Computational performance}
\label{sec:comp_performance}
We compare the computational performance of the following three formulations and specific solution methods using the instances described earlier.
\begin{enumerate}
    \item \textbf{Flow-based formulation}: We implement the formulations $\mathcal{F}^{\textsc{flow}}_{\textsc{2W}}$ and $\mathcal{F}^{\textsc{flow}}_{\textsc{1W}}$, discussed in Sections~\ref{sec:MIP} and~\ref{sec:one_way}, as our baseline. These formulations are similar to those introduced by \cite{stephan2021layout} and are solved directly using CPLEX.
    \item \textbf{Flow-based models with valid inequalities}: For the two-way case, this formulation augments $\mathcal{F}^{\textsc{flow}}_{\textsc{2W}}$ with the bidirectional hop inequalities developed in Section~\ref{sec:vis_two_way}. For the one-way case, it extends $\mathcal{F}^{\textsc{flow}}_{\textsc{1W}}$ by incorporating both the bidirectional hop inequalities and the dead-end constraints described in Section~\ref{sec:vis_one_way}. These valid inequalities are added to the models upfront, and the formulations are solved directly using CPLEX.  
    \item \textbf{Connectivity-cut formulation}: These models solve $\mathcal{F}^{\textsc{CCC}}_{\textsc{2W}}$ and $\mathcal{F}^{\textsc{CCC}}_{\textsc{1W}}$, adding a limited number of valid inequalities (e.g., bidirectional hop inequalities, dead-end constraints) upfront and generating the remainder using the polynomial-time separation procedures outlined in Sections~\ref{sec:disc_infeasible_two_way} and~\ref{sec:disc_infeasible_one_way}. They are solved using the branch-and-cut method through CPLEX callbacks.
\end{enumerate}
\begin{figure}[H]
  \centering
  \begin{subfigure}[b]{0.49\textwidth}
          \centering
    \includegraphics[width=0.98\linewidth]{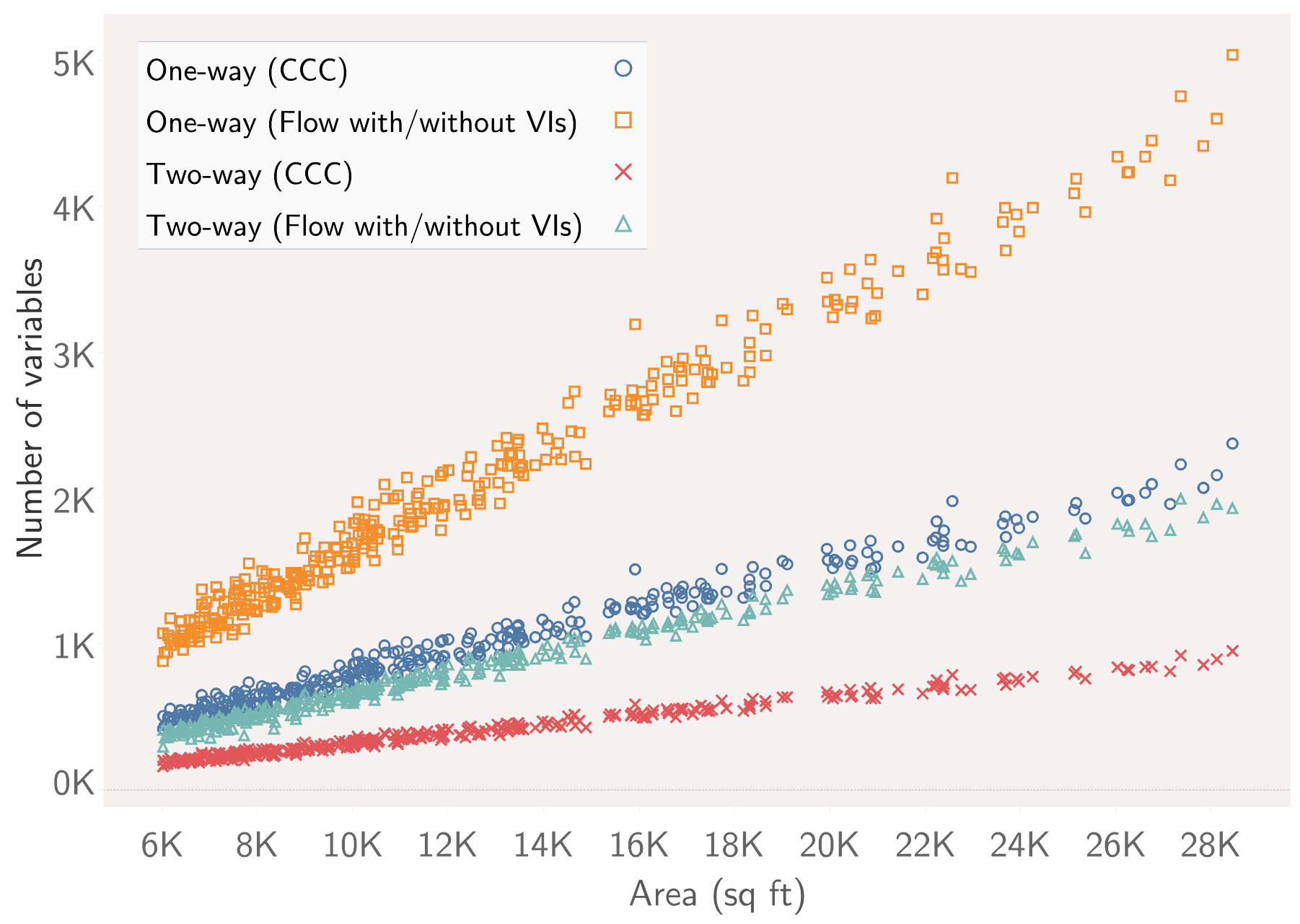}
    \caption{Comparison of the number of variables}
    \label{fig:num_variables}
  \end{subfigure}
  \hfill
  \begin{subfigure}[b]{0.49\textwidth}
          \centering
    \includegraphics[width=0.98\linewidth]{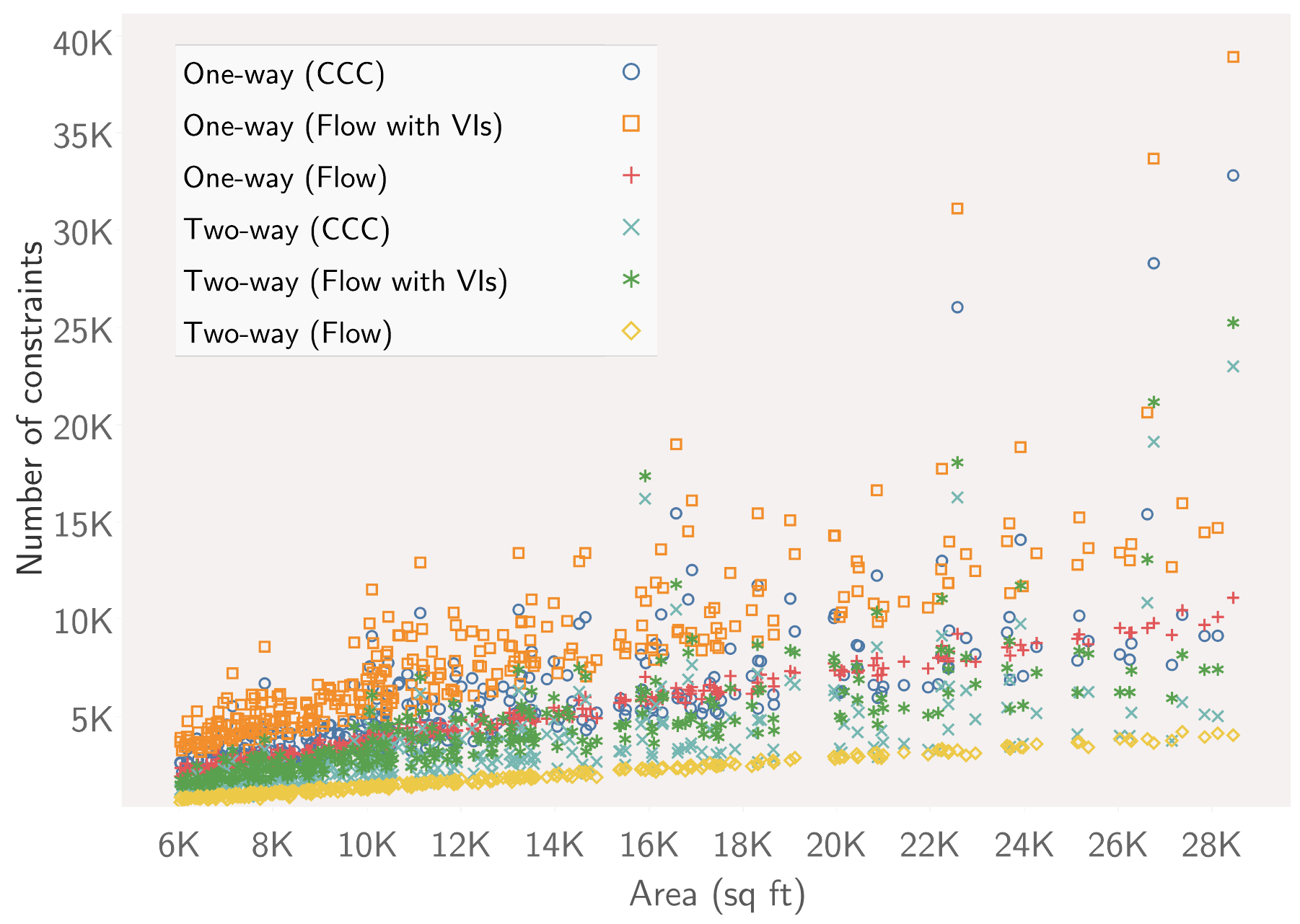}
    \caption{Comparison of the number of constraints}
    \label{fig:num_constraints}
  \end{subfigure}
  \hfill
\caption{Problem sizes for parking lots for different instances}
\label{fig:problem_size}
\end{figure}

\textbf{Problem size:}
Figure \ref{fig:num_variables} illustrates how the number of variables changes with the area for the 325 parking lots used in this research. The relationship is more or less linear for all formulations. In the two-way (and one-way) case, flow-based formulations have a lot more variables than the CCC formulation because of the $f$ (and $g$) variables. The valid inequalities add about 135\% more constraints on average for the two-way case and 74\% for the one-way case (see  Figure \ref{fig:num_constraints}). The CCC formulation generally has fewer constraints. Although valid inequalities are present in this formulation, the constraint count usually falls between those of flow-based models with and without valid inequalities due to the absence of flow constraints. However, in many cases, particularly in one-way configurations, the number of constraints in the CCC formulation may be lower than those in the flow-based formulation. These statistics only indicate the structural constraints and do not include the feasibility cuts added during the branching process. The constraints can increase sharply for larger lots, especially with fewer blocked fields. In these cases, the bidirectional hop inequalities can be significantly large in number. 

\begin{figure}[H]
  \centering
  \begin{subfigure}[b]{0.49\textwidth}
          \centering
    \includegraphics[width=0.98\linewidth]{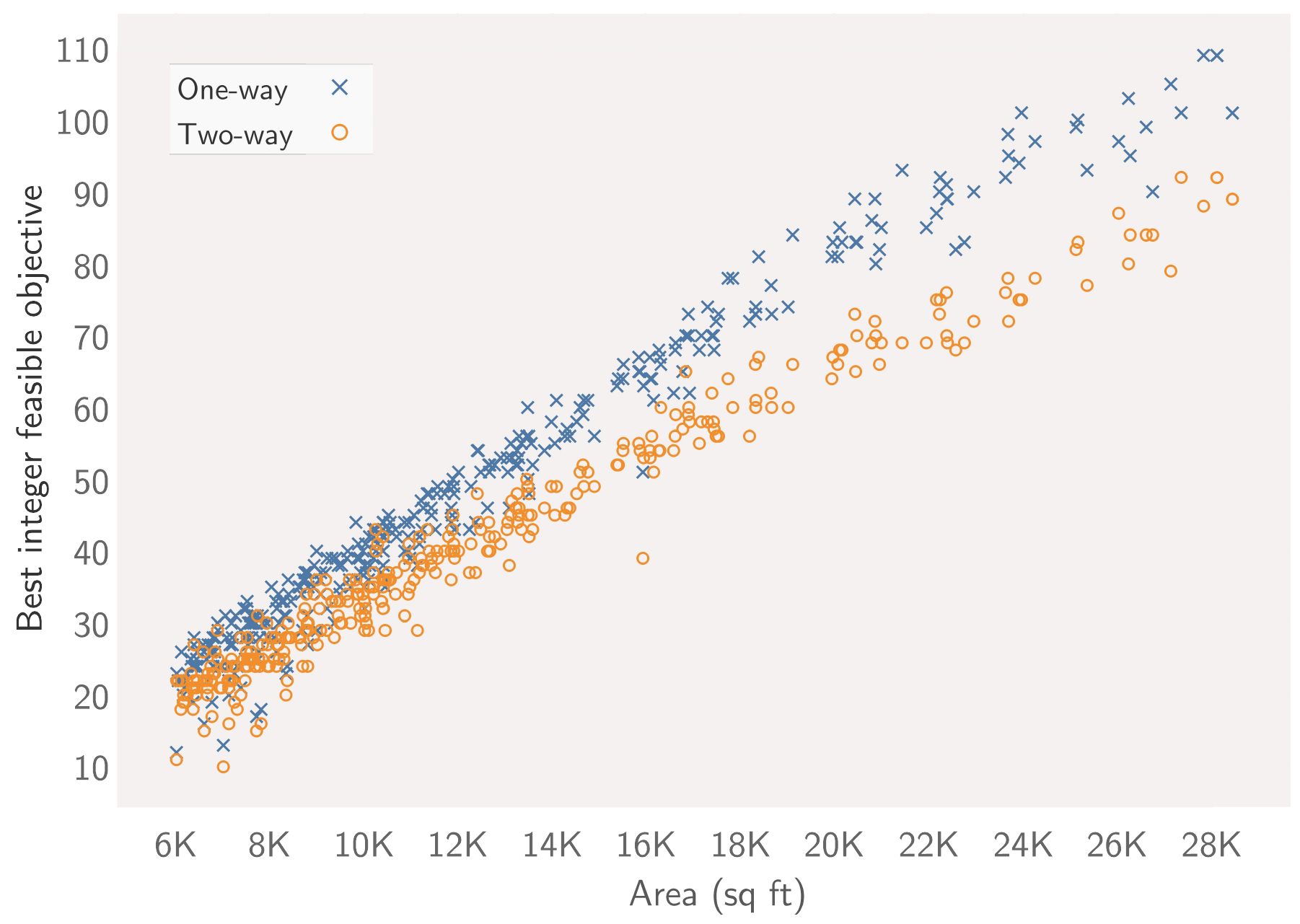}
    \caption{Lower bounds (One-way vs. Two-way)}
    \label{fig:upper_bounds}
\end{subfigure}
  \begin{subfigure}[b]{0.49\textwidth}
          \centering
    \includegraphics[width=0.98\linewidth]{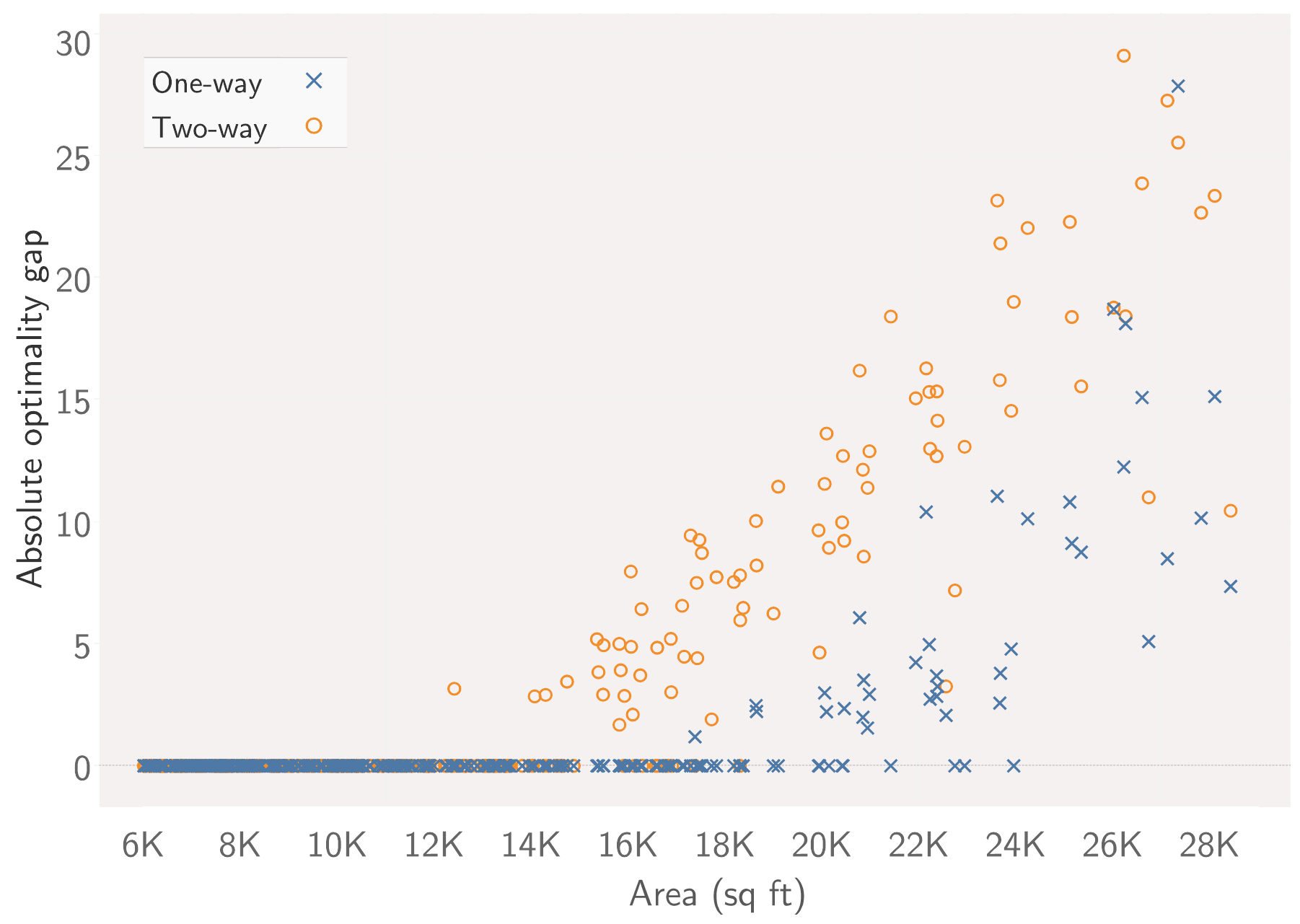}
    \caption{Optimality gap (One-way vs. Two-way)}
    \label{fig:optimality_gap}
  \end{subfigure}
\caption{Comparison between one-way and two-way solutions}
\label{fig:best_solutions}
\end{figure}

\textbf{Analysis of one-way and two-way driving layout solutions:} Figure \ref{fig:upper_bounds} shows the best lower bounds of the three formulations---flow-based model discussed in Sections \ref{sec:MIP} and \ref{sec:one_way} (which is the baseline model), flow-based model with valid inequalities (bidirectional hop and dead-end constraints), and the CCC version---for the one-way and two-way configurations. While many of these instances are solved optimally, others are terminated due to the time limit with an optimality gap. A total of 13291 stalls are discovered for the two-way version. An additional 2476 stalls can be accommodated if the drive layouts are one-way, which is an 18.63\% increase. This improvement can be attributed to the reduced lane width in the one-way configuration (3 m versus a total carriageway width of 6 m in the two-way configuration), enabling the reclaimed space to accommodate additional vehicles. In fact, when the entry and exit fields are adjacent, models with a one-way driving configuration can be regarded as a relaxed version of the two-way scenario and thus achieve a better objective value. The difference in the number of stalls between the one-way and two-way configurations is relatively less for smaller lots. However, it grows significantly with the parking lot area. The absolute difference between the upper and lower bounds is shown in Figure \ref{fig:optimality_gap}. Out of 325 instances, the gap is zero for 245 two-way and 288 one-way instances. The absolute optimality gap in the two-way models was generally higher than in the one-way version when the time limits were reached. 

\begin{table}[htbp]
\small
  \centering
  \caption{Summary statistics for the two-way and one-way experiments. All times are in seconds. LB and UB indicate lower bounds and upper bounds, respectively.}
    \begin{tabular}{l|l|l|rr|l|l|l|rr}
    \hline
    \multicolumn{5}{c|}{\textbf{Two-Way}} & \multicolumn{5}{c}{\textbf{One-Way}} \\
    \hline
    \textbf{Group} & \multicolumn{2}{l|}{\textbf{Statistics}} & \textbf{Average} & \textbf{Median} & \textbf{Group} & \multicolumn{2}{l|}{\textbf{Statistics}} & \textbf{Average} & \textbf{Median} \\
    \hline
    \multirow{4}[2]{*}{I (187)} & \multicolumn{2}{l|}{Objective value} & 28.68 & 28    & \multirow{4}[2]{*}{I (149)} & \multicolumn{2}{l|}{Objective value} & 35.71 & 35 \\
          & \multicolumn{2}{l|}{Time (Flow)} & 100.97 & 12.76 &       & \multicolumn{2}{l|}{Time (Flow)} & 111.44 & 43.23 \\
          & \multicolumn{2}{l|}{Time (Flow \& VIs)} & 7.12  & 3.23  &       & \multicolumn{2}{l|}{Time (Flow \& VIs)} & 143.12 & 47.31 \\
          & \multicolumn{2}{l|}{Time (CCC)} & 5.34  & 2.63  &       & \multicolumn{2}{l|}{Time in s (CCC)} & 12.19 & 5.43 \\
    \hline
    \multirow{9}[6]{*}{II (138)} & \multirow{3}[2]{*}{Flow} & LB    & 55.97 & 55    & \multirow{9}[6]{*}{II (176)} & \multirow{3}[2]{*}{Flow} & LB    & 55.20 & 55 \\
          &       & UB    & 69.22 & 66.03 &       &       & UB    & 65.13 & 62.88 \\
          &       & Gap   & 0.22  & 0.20  &       &       & Gap   & 0.19  & 0.12 \\
\cline{2-5}\cline{7-10}          & \multirow{3}[2]{*}{Flow \& VIs} & LB    & 57.17 & 56    &       & \multirow{3}[2]{*}{Flow \& VIs} & LB    & 55.85 & 55 \\
          &       & UB    & 65.16 & 61.90 &       &       & UB    & 65.81 & 63.99 \\
          &       & Gap   & 0.12  & 0.12  &       &       & Gap   & 0.32  & 0.12 \\
\cline{2-5}\cline{7-10}          & \multirow{3}[2]{*}{CCC} & LB    & 56.60 & 55    &       & \multirow{3}[2]{*}{CCC} & LB    & 59.28 & 58 \\
          &       & UB    & 63.77 & 60.00 &       &       & UB    & 60.85 & 58 \\
          &       & Gap   & 0.10  & 0.08  &       &       & Gap   & 0.02  & 0 \\
    \hline
    \end{tabular}%
  \label{tab:summary_stats}%
\end{table}%

\textbf{Two-way layouts:} We analyzed the computational performance of the three formulations for the two-way case by grouping all instances into two categories: (I) 187 instances in which all methods discover optimal solutions within 15 minutes, and (II) 138 instances that are timed out for at least one method (which in this case was always the flow formulation). See Table \ref{tab:summary_stats} for a summary. Note that the average of 100.97 s for the flow formulation is significantly higher than the median since some instances take longer to solve, but the optimal solution is found within 15 min. Hence, we will make comparisons based on the median values. For the first group, the bidirectional hop inequalities in the flow-based formulation significantly reduce the median runtimes by 74.70\%. A further reduction of 18.40\% is noticed when using the branch-and-cut algorithm. For the slightly harder Group II instances, we compare the integer lower bounds, best LP relaxation upper bounds, and the relative gap (expressed as a ratio). The branch-and-cut method could solve 58 of these instances optimally. Overall, the flow model with valid inequalities performed better, and across all instances in Group II, the flow model discovered 7724 stalls, the flow version with VIs found 7889 stalls, and the branch-and-cut method on the CCC formulation yielded 7811 parking stalls.

\textbf{One-way layouts: }For the one-way case, we again classified the instances into two groups: (I) 149 instances in which all three methods terminate within 15 min, and (II) 176 instances in which at least one of the three methods has a non-zero gap. Unlike the two-way case, some of the instances in group II terminated with zero gaps when using the flow formulation. Surprisingly, the flow version with valid inequalities did not improve the runtimes and the gaps. Compared to the flow formulation in Group I, the branch-and-cut algorithm was 87.43\% faster in terms of median times. Across all instances in Group II, the flow formulations with and without valid inequalities discovered 9715 and 9830 parking stalls, respectively, whereas the branch-and-cut found 10434 stalls. Of the 176 instances in Group II, the branch-and-cut method found an optimal solution in 139 instances, as reflected in the average and median gap values.  

\textbf{Effect of bidirectional hop and dead-end constraints on the branch-and-cut algorithm:} To evaluate the impact of the valid inequalities added at the start (bidirectional hop and dead-end constraints), we tested a version of the branch-and-cut algorithm without them. For the two-way case, with a time limit of 15 min, this version did not solve six instances optimally among those solved by the earlier CCC formulation. The median and average runtimes for instances where both methods produced optimal solutions were 42.51\% and 71.64\% slower than the previous branch-and-cut version, respectively. For the more difficult instances in which either method failed to find an optimal solution in 15 min, the median and average gaps without valid inequalities were 0.19 and 0.22, respectively, compared to 0.15 and 0.17 for the version with valid inequalities. 

For the one-way case, the model performed very poorly without valid inequalities. Out of the 288 instances solved optimally in 15 min using the previous branch-and-cut algorithm, an optimal solution was found in only 132 instances when run without valid inequalities. For instances that could be solved optimally by both methods, the median and average runtimes without valid inequalities worsened to 198.58 s and 363.92 s (compared to 6.08 s and 8.77 s when using valid inequalities). Among the remaining problems, in 169 instances, no integer feasible solution was found after 15 min. Without valid inequalities, both models rely on dynamically generated feasibility cuts to achieve connectivity, which is greater in number, particularly in the one-way case. This drastically slows down the LP relaxations at intermediate nodes of the branch-and-bound tree. 

\begin{table}[H]
\small
  \centering
  \caption{Runtimes (in seconds) and gaps of the two-way parking lot formulations}
    \begin{tabular}{l|rrrr|rrrr|rrrrr}
    \hline
    \multicolumn{1}{r|}{\multirow{2}[4]{*}{\textbf{Instance ID}}} & \multicolumn{4}{c|}{\textbf{Flow }} & \multicolumn{4}{c|}{\textbf{Flow  with VIs}} & \multicolumn{5}{c}{\textbf{CCC}} \\
\cline{2-14}          & \textbf{Time} & \textbf{LB} & \textbf{UB} & \textbf{Gap} & \textbf{Time} & \textbf{LB} & \textbf{UB} & \textbf{Gap} & \textbf{Time} & \textbf{LB} & \textbf{UB} & \textbf{Gap} & \textbf{Cuts} \\
        \hline
    125-068 & 1.6   & \textbf{20} & 20    & 0     & 0.7   & \textbf{20} & 20    & 0     & 0.5   & \textbf{20} & 20    & 0     & 12 \\
    125-140 & 4.5   & \textbf{26} & 26    & 0     & 2.4   & \textbf{26} & 26    & 0     & 1.4   & \textbf{26} & 26    & 0     & 85 \\
    125-299 & 51.3  & \textbf{33} & 33    & 0     & 12.9  & \textbf{33} & 33    & 0     & 10.5  & \textbf{33} & 33    & 0     & 500 \\
    125-425 & 785.8 & \textbf{41} & 41    & 0     & 9.4   & \textbf{41} & 41    & 0     & 7.6   & \textbf{41} & 41    & 0     & 442 \\
    125-116 & 2524.6 & \textbf{38} & 38    & 0     & 69.2  & \textbf{38} & 38    & 0     & 60.4  & \textbf{38} & 38    & 0     & 421 \\
    125-026 & 7200.1 & \textbf{46} & 48.54 & 0.06  & 63.4  & \textbf{46} & 46    & 0     & 43.7  & \textbf{46} & 46    & 0     & 2518 \\
    125-269 & 7200.8 & \textbf{56} & 60.24 & 0.08  & 4460.5 & \textbf{56} & 56    & 0     & 2850.3 & \textbf{56} & 56    & 0     & 8462 \\
    125-158 & 7200.4 & \textbf{64} & 73.24 & 0.14  & 6562.6 & \textbf{64} & 64    & 0     & 1265.5 & \textbf{64} & 64    & 0     & 2578 \\
    125-355 & 7200.6 & \textbf{68} & 78.44 & 0.15  & 7207.9 & \textbf{68} & 71.75 & 0.06  & 3851.9 & \textbf{68} & 68    & 0     & 6780 \\
    125-372 & 7200.4 & 70    & 88.17 & 0.26  & 7213.5 & \textbf{73} & 84.04 & 0.15  & 7215.0 & 72    & 81.28 & 0.13  & 14826 \\
    215-145 & 7200.4 & 70    & 88.21 & 0.26  & 7211.9 & 70    & 83.05 & 0.19  & 7216.2 & \textbf{71} & 80.69 & 0.14  & 21274 \\
    125-094 & 7200.4 & 84    & 104.58 & 0.24  & 7214.0 & \textbf{85} & 101.65 & 0.20  & 7217.0 & 80    & 101.48 & 0.27  & 47929 \\
    125-407 & 7200.3 & 82    & 106.48 & 0.30  & 7212.4 & \textbf{87} & 102.85 & 0.18  & 7215.0 & 80    & 100.57 & 0.26  & 23286 \\
    125-112 & 7200.4 & 92    & 119.28 & 0.30  & 7215.9 & \textbf{93} & 115.30 & 0.24  & 7221.3 & 87    & 113.96 & 0.31  & 52962 \\
    \hline
    \end{tabular}%
  \label{tab:longer_two_way}%
\end{table}%

\textbf{Extended analysis: }We also tested the performance of all three formulations for both driving layout configurations with a longer time limit of two hours. These experiments demonstrated that the trends in the computational times and gaps hold when run for longer periods. 

We choose lots with areas closest to the median area values of the instances in the 14 bins shown in Figure \ref{fig:nyc_area_distribution}. Additional details of these instances are provided in \ref{sec:add_stats}. The results of these trials are shown in Tables \ref{tab:longer_two_way} and \ref{tab:longer_one_way} for the two-way and one-way cases, respectively. The bold LB values represent the best solution among the three formulations. The number of user cuts dynamically generated within the branch-and-cut method is indicated in the last column. The restrictions imposed by the added inequalities help reduce the upper bounds to a relatively greater degree. Needless to say, in the other cases, the branch-and-cut algorithm reached optimality more quickly. Even with longer runtime thresholds, the flow-based formulation did not ``catch up" with respect to the bounds, and including the valid inequalities helped discover better lower bounds. The branch-and-cut method yielded the tightest upper bounds among all three formulations.

\begin{table}[H]
\small
  \centering
  \caption{Runtimes (in seconds) and gaps of the one-way parking lot formulations}
    \begin{tabular}{l|rrrr|rrrr|rrrrr}
    \hline
    \multirow{2}[4]{*}{\textbf{Instance ID}} & \multicolumn{4}{c|}{\textbf{Flow }} & \multicolumn{4}{c|}{\textbf{Flow  with VIs}} & \multicolumn{5}{c}{\textbf{CCC}} \\
\cline{2-14}          & \textbf{Time} & \textbf{LB} & \textbf{UB} & \textbf{Gap} & \textbf{Time} & \textbf{LB} & \textbf{UB} & \textbf{Gap} & \textbf{Time} & \textbf{LB} & \textbf{UB} & \textbf{Gap} & \textbf{Cuts} \\
        \hline
    125-068 & 9.0   & \textbf{24} & 24    & 0     & 4.4   & \textbf{24} & 24    & 0     & 2.2   & \textbf{24} & 24    & 0     & 938  \\
    125-140 & 77.2  & \textbf{30} & 30    & 0     & 44.0  & \textbf{30} & 30    & 0     & 4.7   & \textbf{30} & 30    & 0     & 1532 \\
    125-299 & 126.4 & \textbf{39} & 39    & 0     & 170.3 & \textbf{39} & 39    & 0     & 28.4  & \textbf{39} & 39    & 0     & 2880 \\
    125-425 & 7200.1 & 42    & 47.41 & 0.13  & 7206.5 & 42    & 48.77 & 0.16  & 33.1  & \textbf{44} & 44    & 0     & 3860 \\
    125-116 & 1340.0 & \textbf{46} & 46    & 0     & 3592.6 & \textbf{46} & 46    & 0     & 51.7  & \textbf{46} & 46    & 0     & 6374 \\
    125-026 & 1560.6 & \textbf{56} & 56    & 0     & 6614.3 & \textbf{56} & 56    & 0     & 38.3  & \textbf{56} & 56    & 0     & 6720 \\
    125-269 & 3223.6 & \textbf{64} & 64    & 0     & 7203.9 & \textbf{64} & 66    & 0     & 369.5 & \textbf{64} & 64    & 0     & 9055 \\
    125-158 & 7200.2 & \textbf{78} & 80.18 & 0.03  & 7209.9 & \textbf{78} & 80.58 & 0.03  & 159.0 & \textbf{78} & 78    & 0     & 16143 \\
    125-355 & 7200.2 & \textbf{83} & 84.83 & 0.02  & 7214.6 & \textbf{83} & 84.88 & 0.02  & 146.0 & \textbf{83} & 83    & 0     & 10198 \\
    125-372 & 7200.2 & 88    & 94.29 & 0.07  & 7221.2 & 88    & 95.99 & 0.09  & 2146.6 & \textbf{89} & 89    & 0     & 21389 \\
    215-145 & 7200.4 & 82    & 96.40 & 0.18  & 7216.8 & 88    & 96.44 & 0.10  & 5166.2 & \textbf{89} & 89    & 0     & 50973 \\
    125-094 & 7200.4 & \textbf{104} & 108.06 & 0.04  & 7218.3 & 101   & 110.08 & 0.09  & 7217.2 & 103   & 107.60 & 0.04  & 22696 \\
    125-407 & 7200.5 & 95    & 112.03 & 0.18  & 7216.4 & 101   & 112.13 & 0.11  & 7217.9 & \textbf{103} & 111.28 & 0.08  & 64819 \\
    125-112 & 7204.6 & 103   & 123.46 & 0.20  & 7219.8 & 112   & 124.34 & 0.11  & 7222.2 & \textbf{113} & 122.68 & 0.09  & 35412 \\
    \hline
    \end{tabular}%
  \label{tab:longer_one_way}%
\end{table}%

\begin{figure}[H]
    \centering
    \includegraphics[width=0.88\linewidth]{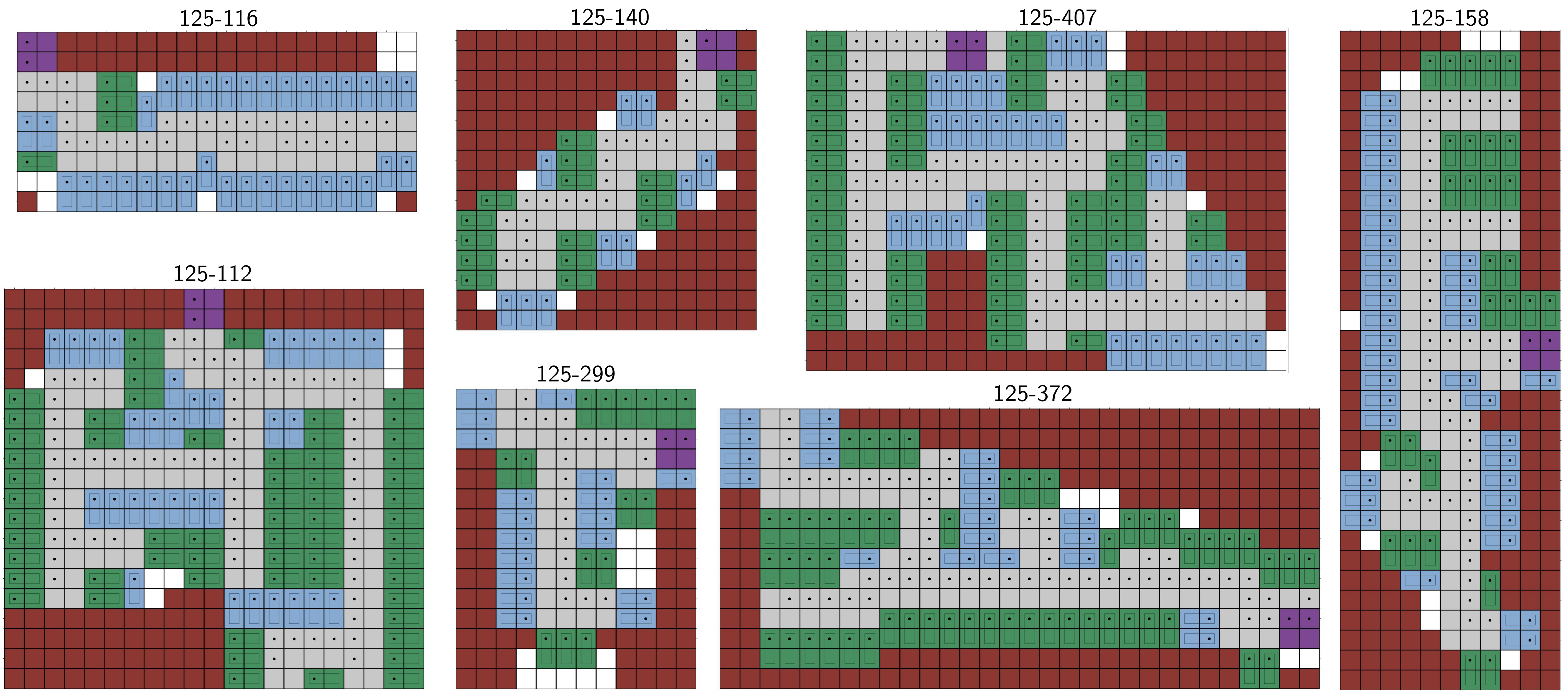}
    \caption{Sample parking lots for the two-way driveway configuration}
    \label{fig:two_way_collage}
\end{figure}

\begin{figure}[t]
    \centering
    \includegraphics[width=0.88\linewidth]{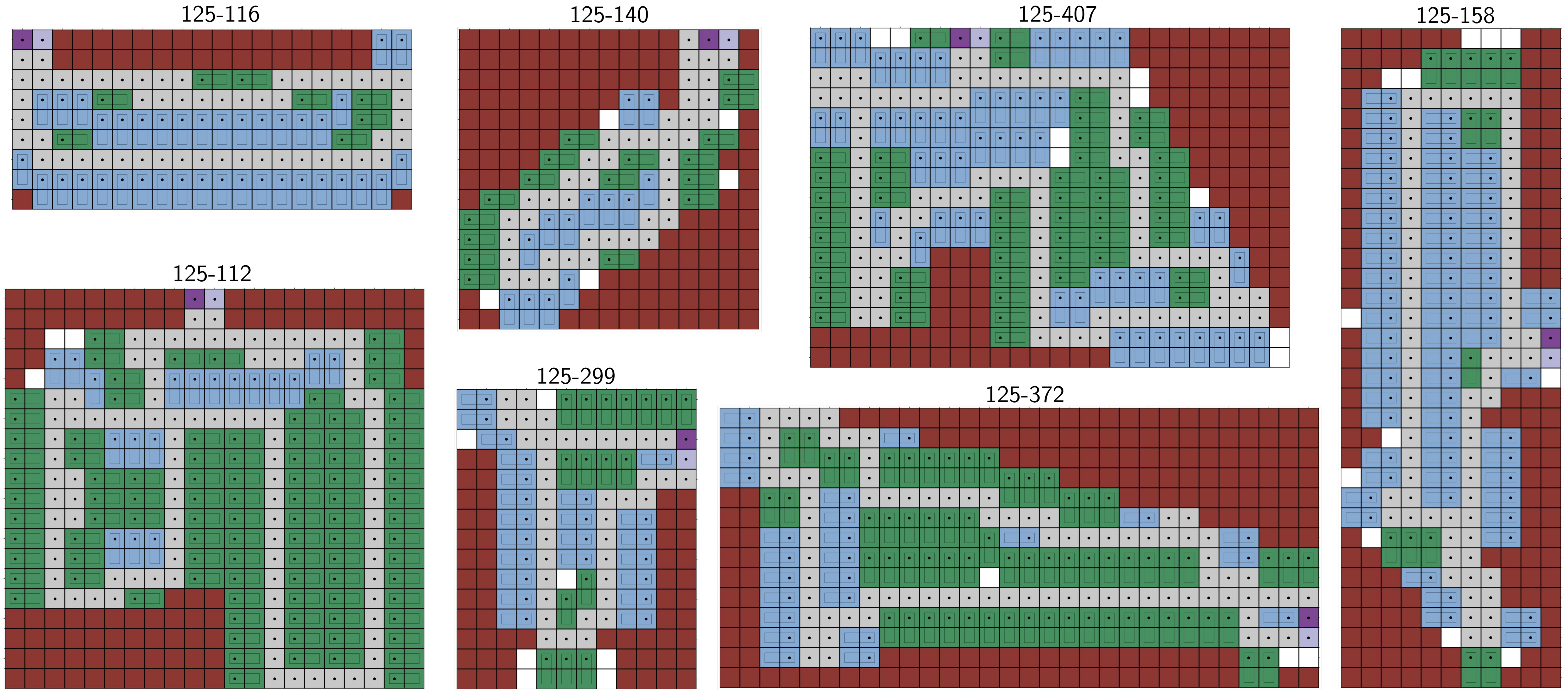}
    \caption{Sample parking lots for the one-way driveway configuration}
    \label{fig:one_way_collage}
\end{figure}

In the one-way case, the branch-and-cut algorithm discovered the same or more parking lots and terminated with a smaller gap when the time limit was reached, except for one case. As noted earlier, the branch-and-cut method proved more effective for the one-way scenario than the two-way case and was exceedingly fast in generating optimal solutions. In some instances, optimal solutions were generated in minutes, whereas flow-based formulations could not find or certify optimality even after two hours due to poor upper bounds. Figures \ref{fig:two_way_collage} and \ref{fig:one_way_collage} show sample parking lot solutions for some of the instances from these experiments.

\section{Conclusion}
\label{sec:conclusion}
This paper proposed alternative integer programming formulations to optimize the number of stalls in a parking lot. We divided the lot into a grid where each cell is associated with variables that determine whether it is designated for parking, driving, or left empty. The constraints ensure that each cell is used for a single purpose and that parking fields are connected using driveways. For one-way driveways, we included additional flow and directionality variables to help decide the flow of traffic on the grid and prevent dead ends. In the one-way and two-way cases, the flow variables appear in forcing constraints with big-M terms, leading to poorer runtimes when using an off-the-shelf solver. We reformulated these problems by defining connectivity cut constraints to link driveways to entrances and exits. A simple family of valid inequalities was proposed using neighbors at different hop distances from a central cell and from the entrance/exit. These constraints resulted in layouts with few disconnected driveway components. We repaired such solutions by generating additional constraints using minimum-cardinality vertex or edge cut sets during the branch-and-bound process. Our experiments on 325 NYC parking layouts highlighted the advantages of the valid inequalities and the branch-and-cut method. The branch-and-cut algorithm's median runtimes were significantly faster than the baseline flow-based formulation. 

There are several potential directions for future research. While we focused on maximizing the number of parking stalls, evaluating different layouts from a traffic perspective can provide a comprehensive approach to parking lot design. Additionally, parking lots are often heterogeneous, accommodating various vehicle sizes with both parallel and angled parking options. Requirements for differently-abled parking and electric vehicle charging can add practical value. From a computational perspective, while LP relaxations improve with added cuts, these cuts are added only for integer infeasible solutions and not at other tree nodes. Adding cuts at other nodes can quickly prune less promising search space regions, but would require a fast maximum separation subroutine. One could also explore branching rules that fix a block of parking field variables or partitioning techniques that divide a lot into smaller regions. Discovering other valid inequalities and analyzing their mathematical properties can also significantly improve model formulations. 

\vspace{5mm}
{\large \textbf{Acknowledgments}}\\
\\
The authors appreciate the constructive comments provided by the reviewers and the associate editor, which have significantly improved the quality of the paper.

\newpage
\appendix
\label{sec:appendix}
\renewcommand{\thesection}{Appendix \Alph{section}}
\renewcommand{\thesubsection}{\Alph{section}.\arabic{subsection}}

\section{Equivalence of formulations}
\label{sec:equivalence}
In this section, we establish that the formulations with connectivity cut constraints ($\mathcal{F}_{\textsc{2W}}^{\textsc{CCC}}$ and  $\mathcal{F}_{\textsc{1W}}^{\textsc{CCC}}$) yield the same set of integer components as the flow-based formulations ($\mathcal{F}_{\textsc{2W}}^{\textsc{flow}}$ and $\mathcal{F}_{\textsc{1W}}^{\textsc{flow}}$). We will demonstrate that \eqref{eq:cutconstraints} can be derived from Benders feasibility cuts, in which we project the continuous variables to obtain a pure 0-1 integer program. The proof is adapted from the One-commodity Pickup-and-Delivery Traveling Salesman Problem (1-PDTSP) literature \citep{hernandez2003one}. However, unlike the 1-PDTSP, the constraints in our model are written for connected subsets that do not contain the entrance/exit, rather than all subsets of the valid drive cells; the reason for this will become apparent in the proof.

\begin{theorem}
\label{theorem:two_way}
Suppose $(x, y)$ denotes an integer feasible parking and driving-field assignment. The formulations $\mathcal{F}_{\textsc{2W}}^{\textsc{flow}}$ and $\mathcal{F}_{\textsc{2W}}^{\textsc{CCC}}$ induce the same set of integer feasible values for $(x,y)$.
\end{theorem} 
\begin{proof}
Consider the Benders reformulation of the problem $\mathcal{F}^{\textsc{flow}}_{\textsc{2W}}$ in which the binary variables $x$ and $y$ are treated as master problem variables and the flow variables $f$ are part of the subproblem. 
\begin{align}
    \max & \sum_{(i,j) \in \validzeroparkset} \zeroparkvar{i}{j} + \sum_{(i,j) \in \validninetyparkset} \ninetyparkvar{i}{j} + R(y)  && \\
     \text{s.t. } &  \eqref{eq:singlepurpose_zero_disaggregate} \text{--}\eqref{eq:ninetyparkdriveconnector}, \eqref{eq:existing drive vars} \text{--} \eqref{eq: y var is binary}  \notag 
\end{align}
where the recourse function $R(y)$ is the solution to the following subproblem 
\begin{align}
    \max & \, \, \, \, 0  && \notag \\
    \text{s.t. } & \sum_{(k,l) \in \gridadjset{i}{j}} \flowvar{i}{j}{k}{l} - \sum_{(k,l) \in \gridinvadjset{i}{j}} \flowvar{k}{l}{i}{j} \geq \drivevar{i}{j} \qquad && \forall \, (i,j) \in \validdriveset \setminus \{(\entrancex, \entrancey)\} \label{eq:floworigin_sp}\\
     & \sum_{(i,j) \in \gridinvadjset{\entrancex}{\entrancey}} \flowvar{i}{j}{\entrancex}{\entrancey} \leq \sum_{(i,j) \in \validdriveset} \drivevar{i}{j} -1 \qquad &&  \label{eq:flowdest_sp}\\
    & \flowvar{i}{j}{k}{l} \leq M \drivevar{k}{l} \qquad && \forall \, ((i,j), (k,l)) \in \gridarcset \label{eq:flowhead_sp} \\
    & \flowvar{i}{j}{k}{l} \geq 0 \qquad && \forall \, ((i,j), (k,l)) \label{eq:flow var is continuous non negative_sp}\in \gridarcset
\end{align}
The recourse function does not contribute to the overall objective as it does not contain any $f$ variables. Hence, we write its objective as 0. Thus, the Benders reformulation would only include feasibility cuts and no optimality cuts. Note that the $y$ variables are fixed in the above problem, and the constraints differ slightly from the flow constraints used earlier. We use inequalities for \eqref{eq:floworigin_sp} and \eqref{eq:flowdest_sp} instead of an equality because the dual feasible space becomes a pointed cone, which simplifies the proof. The modified signs ensure that the feasible region remains the same. We also remove a redundant coupling constraint between the flow variables of an arc and the drive variables of its tail end. Let the dual variables associated with \eqref{eq:floworigin_sp} and \eqref{eq:flowdest_sp} be denoted by $\alpha_{ij}$, and those associated with \eqref{eq:flowhead_sp} be denoted by $\beta_{ij, kl}$. Consider the dual of the above formulation.
\begin{align}
     \min & \sum_{(i,j) \in \validdriveset \setminus \{(\entrancex, \entrancey)\}} \drivevar{i}{j} \alpha_{ij} + \Bigg(\sum_{(i,j) \in \validdriveset} \drivevar{i}{j} -1 \Bigg) \alpha_{pq} + \sum_{((i,j), (k,l)) \in \gridarcset}  M y_{kl} \beta_{ij, kl}  \label{eq:sp_dual}
\end{align}
\begin{align}
    \text{s.t. } & \alpha_{ij} - \alpha_{kl} + \beta_{ij, kl} \geq 0 && \forall \, ((i,j), (k,l)) \in \gridarcset, (i, j), (k, l) \neq (p, q)  \label{eq:sp_dual_1} \\
     & - \alpha_{kl} + \beta_{pq, kl} \geq 0 && \forall \,  (k, l) \in \gridadjset{p}{q}  \label{eq:sp_dual_2} \\
    & \alpha_{pq} + \beta_{ij, pq} \geq 0 && \forall \,  (i, j) \in \gridinvadjset{p}{q}  \label{eq:sp_dual_3} \\
    & \alpha_{ij} \leq 0 && \forall \, (i,j) \in \validdriveset \setminus \{(\entrancex, \entrancey)\} \label{eq:sp_dual_4} \\
    & \alpha_{pq} \geq 0 &&  \label{eq:sp_dual_4.5} \\
    & \beta_{ij, kl} \geq 0 && \forall \, ((i,j), (k,l)) \in \gridarcset \label{eq:sp_dual_5} 
\end{align}
Using the dual, we can reformulate the master problem by adding constraints that avoid $y$ variables that lead to an infeasible subproblem, or equivalently, an unbounded dual. Let $E$ be the set of all extreme directions of the polyhedral cone defined by \eqref{eq:sp_dual_1}--\eqref{eq:sp_dual_5}. To prevent an unbounded dual, we need the objective in \eqref{eq:sp_dual} to be greater than or equal to 0 for all $(\alpha, \beta) \in E$. Hence, the Benders reformulation takes the form
\begin{align}
    &\max \sum_{(i,j) \in \validzeroparkset} \zeroparkvar{i}{j} + \sum_{(i,j) \in \validninetyparkset} \ninetyparkvar{i}{j}   && \notag\\
     & \text{s.t. }  \eqref{eq:singlepurpose_zero_disaggregate} \text{--}\eqref{eq:ninetyparkdriveconnector}, \eqref{eq:existing drive vars} \text{--} \eqref{eq: y var is binary}  && \notag \\
     & \sum_{(i,j) \in \validdriveset \setminus \{(\entrancex, \entrancey)\}} \drivevar{i}{j} \alpha_{ij} + \Bigg(\sum_{(i,j) \in \validdriveset} \drivevar{i}{j} -1 \Bigg) \alpha_{pq} + \sum_{((i,j), (k,l)) \in \gridarcset}  M y_{kl} \beta_{ij, kl} \geq 0  && \forall \, (\alpha, \beta) \in E \label{eq:two_way_benders_feasibility_cut}
\end{align}

We start by analyzing the following three classes of directions that satisfy the constraints of the dual problem \eqref{eq:sp_dual_1}--\eqref{eq:sp_dual_5} and establish that they are extreme directions.

\begin{enumerate}
    \item A solution of the form $\alpha_{pq} = 1$ with all the remaining variables set to zero is clearly an extreme direction. Let this direction belong to the singleton $E_1$. 
    \item Another class of trivial extreme directions $E_2$, one for each arc $((i, j), (k, l)) \in \gridarcset$, of the form $\beta_{ij, kl} = 1$, with all remaining variables set to zero.
    \item For the third class of directions consider a set $S \subseteq \validdriveset \setminus \{(p,q)\}$, such that its induced subgraph $\grid(S)$ is connected. Elements of $E_3$ are of the form $\alpha_{ij} = -1 ~\forall~ (i, j) \in S$ and 0 otherwise (e.g., the teal cells in Figures~\ref{fig:proof1a} and \ref{fig:proof1b}). Define $\beta_{ij,kl} = 1 ~\forall~ ((i,j), (k, l)) \in \gridarcset : (i, j) \in S, (k, l) \in \validdriveset \setminus S$ and zero otherwise. That is, they are ones for all arcs that connect nodes in $S$ to its complement (see teal edges with pink ends in Figures~\ref{fig:proof1a} and \ref{fig:proof1b}).
\end{enumerate}

\begin{figure}[H]
  \centering
  \begin{subfigure}[b]{0.32\textwidth}
          \centering
    \includegraphics[width=\linewidth]{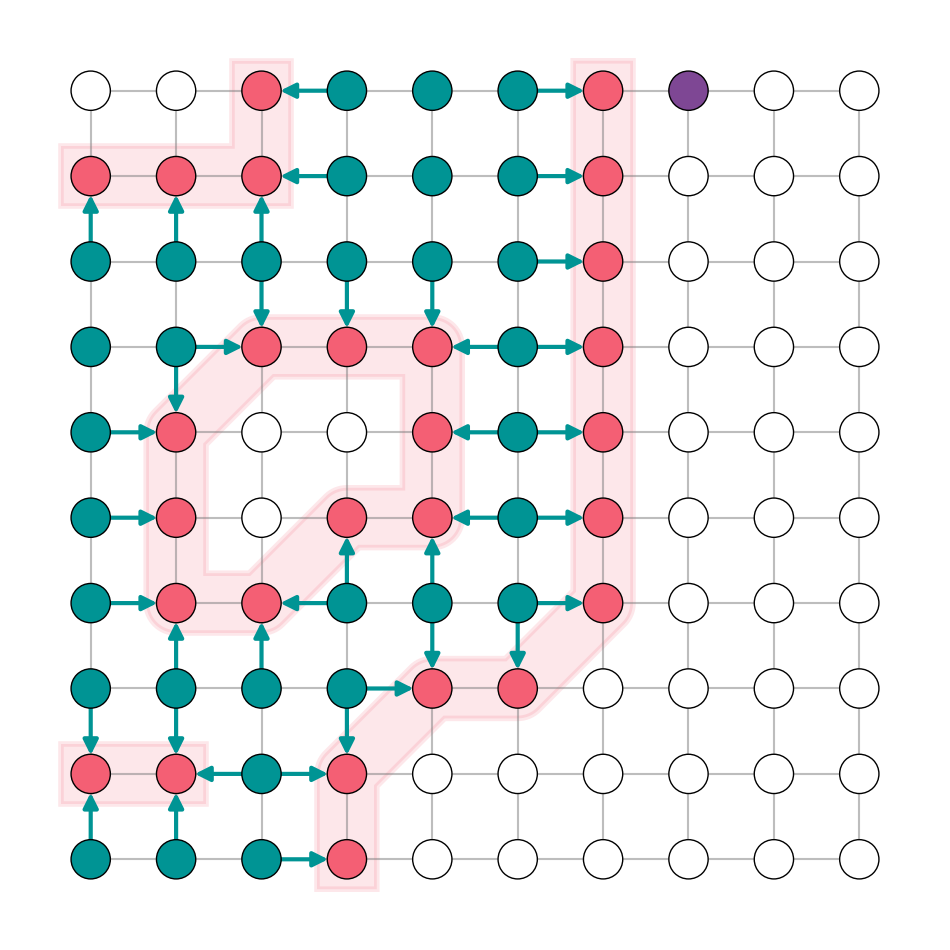}
    \caption{An extreme direction}
    \label{fig:proof1a}
\end{subfigure}
  \begin{subfigure}[b]{0.32\textwidth}
          \centering
    \includegraphics[width=\linewidth]{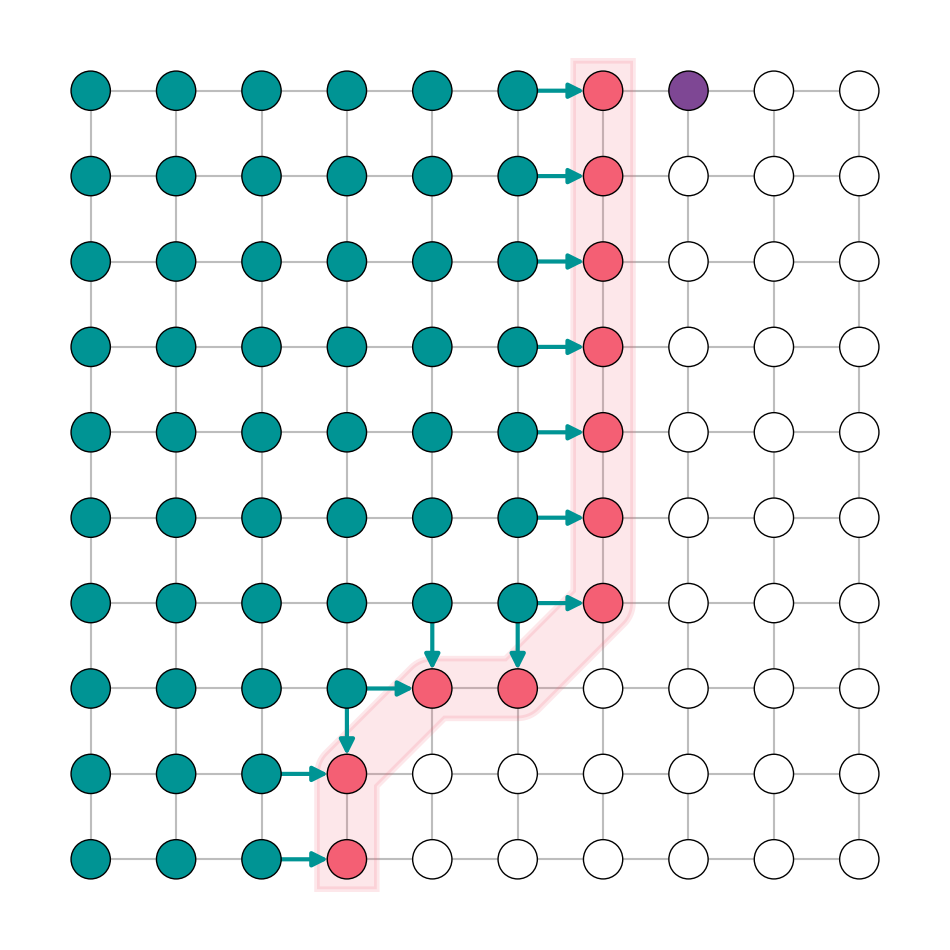}
    \caption{Extreme direction from a vertex cut}
    \label{fig:proof1b}
  \end{subfigure}
  \begin{subfigure}[b]{0.32\textwidth}
          \centering
    \includegraphics[width=\linewidth]{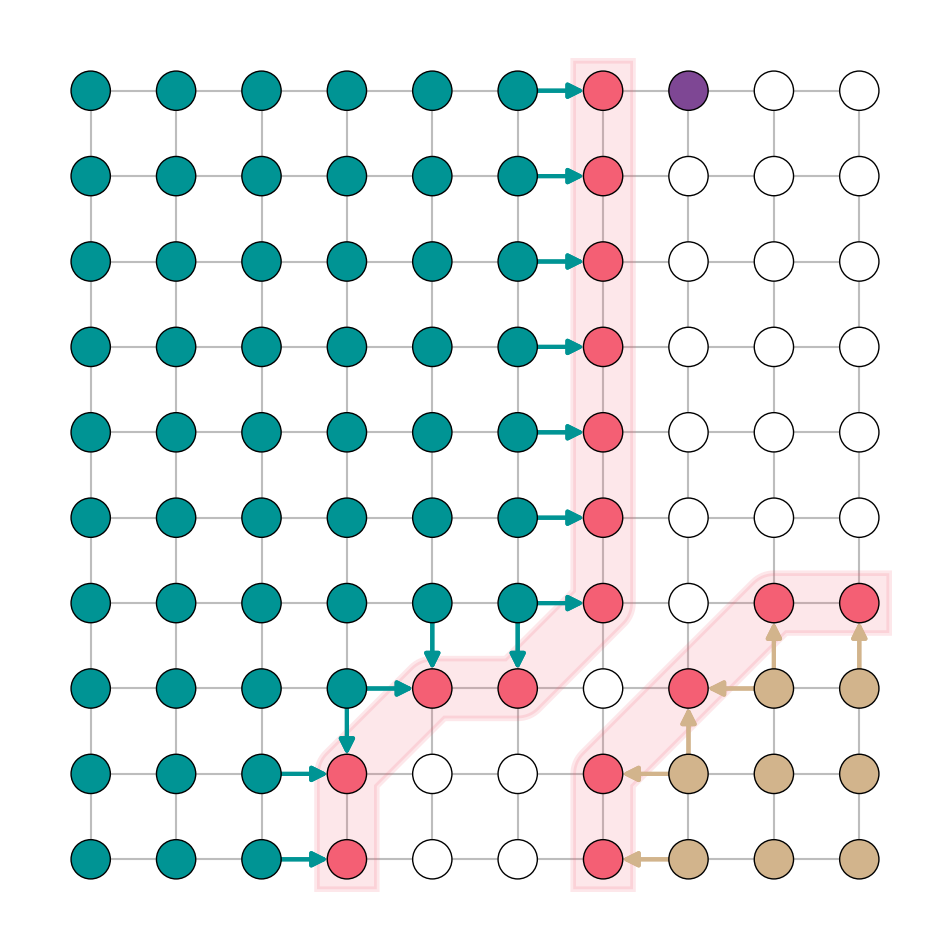}
    \caption{A non-extreme direction}
    \label{fig:proof1c}
  \end{subfigure}
\caption{Examples of direction vectors}
\label{fig:proof1}
\end{figure}

We first show that a direction $(\alpha, \beta)$ in $E_3$ generated from a subset $S$ also belongs to $E$. Assume, for the sake of contradiction, that it is not. Then, we can write it as a positive linear combination two distinct directions of \eqref{eq:sp_dual_1}--\eqref{eq:sp_dual_5}, i.e., $(\alpha, \beta) = \lambda_1 (\alpha^1, \beta^1) + \lambda_2 (\alpha^2, \beta^2)$, where $\lambda_1, \lambda_2 > 0$. It follows that the zero elements of $(\alpha^1, \beta^1)$ and $(\alpha^2, \beta^2)$ match with that of $(\alpha, \beta)$. Consider the arc connecting cells $(i, j)$ and $(k, l)$ in $S$. Inequality \eqref{eq:sp_dual_1}, implies that $\alpha_{ij}^1 \geq \alpha_{kl}^1$ and $\alpha_{ij}^2 \geq \alpha_{kl}^2$. Using a similar inequality for the reverse arc, we conclude that $\alpha_{ij}^1=c_1$ and $\alpha_{ij}^2=c_2$ for all $(i,j) \in S$, where $c_1$ and $c_2$ are constants such that $\lambda_1 c_1 + \lambda_2 c_2 = -1$. Now consider an arc $((i,j), (k, l)) \in \gridarcset : (i, j) \in S, (k, l) \in \validdriveset \setminus S$. Again from \eqref{eq:sp_dual_1}, $\beta_{ij, kl}^1 \geq -c_1$ and $\beta_{ij, kl}^2 \geq -c_2$. However, $\lambda_1 \beta_{ij, kl}^1 + \lambda_2 \beta_{ij, kl}^2 = \beta_{ij, kl} = 1 = - \lambda_1 c_1 - \lambda_2 c_2$, which implies $\beta_{ij, kl}^1 = -c_1$ and $\beta_{ij, kl}^2 = -c_2$. This contradicts the assumption that $(\alpha^1, \beta^1)$ and $(\alpha^2, \beta^2)$ are distinct, since one of them can be written as a positive scalar multiple of the other. 

The requirement that $\grid(S)$ be connected becomes evident now. Without this assumption, setting $\alpha$ values to 1 for all elements in $S$ would prevent us from concluding that $\alpha_{ij}^1=c_1$ and $\alpha_{ij}^2=c_2$ for all $(i,j) \in S$ in the above argument.  In such cases, the directions can be represented as positive linear combinations of directions in $E_3$, one for each connected component. See Figure \ref{fig:proof1c} for example. If $S$ included the teal and brown cells, combining two extreme directions, each of which corresponds to connected cells of the same color, can generate a solution where the $\alpha_{ij}=-1$ for all $(i,j) \in S$. 

Next, we claim that there are no other extreme directions beyond those described earlier by the three classes, i.e., $E = E_1 \cup E_2 \cup E_3$. To this end, we start with any feasible $(\alpha, \beta) \equiv (\bar{\alpha}^1, \bar{\beta}^1) \in E$ and demonstrate a procedure to express it as a positive linear combination of the extreme directions in $E_1 \cup E_2 \cup E_3$. To illustrate the procedure, consider the example shown in Figure \ref{fig:proof_1ed_a}. The shaded drive cells correspond to those with a negative $\alpha$ value. Darker shades correspond to higher values. The $\beta$ values for the edges and $\alpha_{pq}$ are not shown. Let $S_1 = \arg \max_{(i, j) \in \validdriveset} \{\bar{\alpha}^1_{ij} : \bar{\alpha}^1_{ij} < 0\}$. If $S_1 = \emptyset$, we can easily express $(\bar{\alpha}^1, \bar{\beta}^1)$ as $\sum_{u=1}^{\varphi_1} \lambda_{1u} (\alpha^{1u}, \beta^{1u})$, where $(\alpha^{1u}, \beta^{1u}) \in E_1 \cup E_2, \lambda_{1u} >0~\forall~u \in \{1, \ldots, \varphi_1\}$. Else, let $\Lambda_1 = -\bar{\alpha}^1_{ij}$, where $(i, j) \in S_1$. Let $\vartheta_1$ denote the number of connected components in $G(S_1)$. We represent the cells in the individual connected components as $S_{1v}$, where $v \in \{1, \ldots, \vartheta_1\}$. Select a collection of extreme directions from $E_3$,  $\{(\alpha^{1v}, \beta^{1v})\}_{v=1}^{\vartheta_1}$, where $\alpha^{1v}_{ij} = -1~\forall~(i,j) \in S_{1v}$ and $\beta^{1v}_{ij, kl} = 1$, where $((i,j), (k, l)) \in \gridarcset: (i,j) \in S_{1v}, (k, l) \in \validdriveset \setminus S_{1v}$. Define $(\bar{\alpha}^2, \bar{\beta}^2) = (\bar{\alpha}^1, \bar{\beta}^1) - \Lambda_1 \sum_{v=1}^{\vartheta_1} (\alpha^{1v}, \beta^{1v})$. In Figure~\ref{fig:proof_1ed_a}, $G(S_1)$ consists of two connected components, shown within the green polygons. By construction, among the cells with negative $\bar{\alpha}^1$ values, those with the largest $\bar{\alpha}^1$ values have $\bar{\alpha}^2 = 0$, and the remaining $\bar{\alpha}^2$ values are increased accordingly, as illustrated in Figure~\ref{fig:proof_1ed_b}.

\begin{figure}[H]
  \centering
  \begin{subfigure}[b]{0.32\textwidth}
          \centering
    \includegraphics[width=\linewidth]{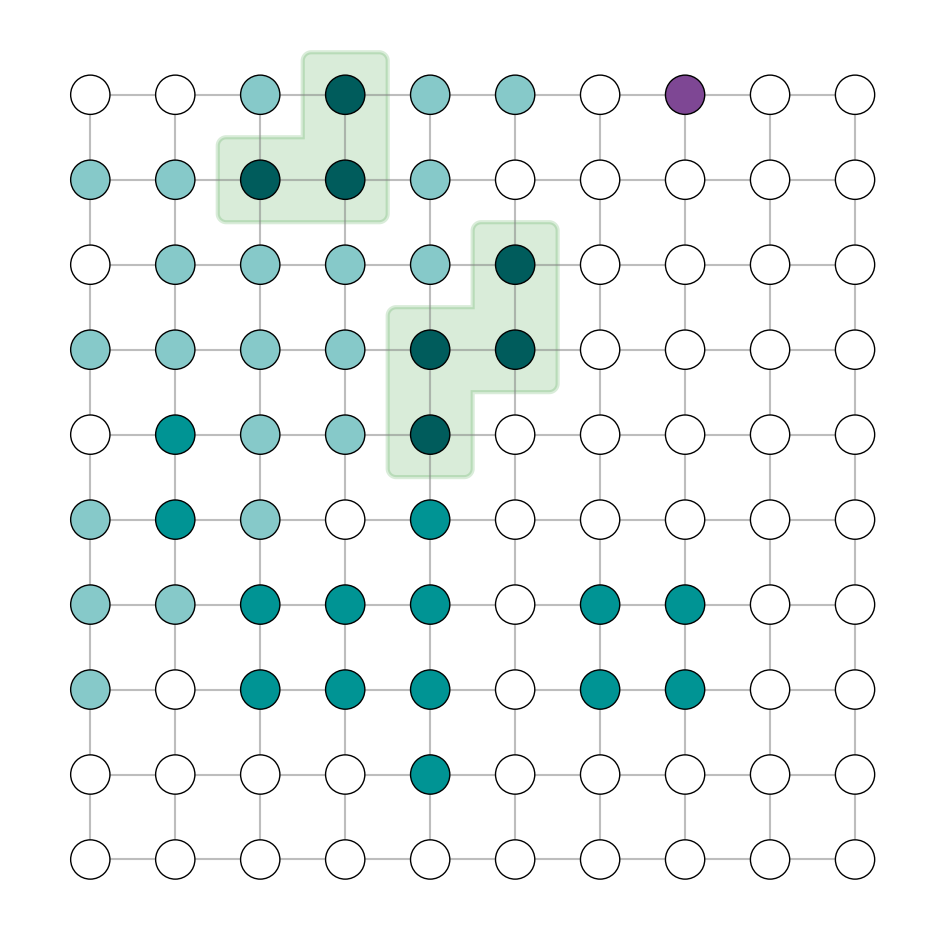}
    \caption{$G(S_1)$}
    \label{fig:proof_1ed_a}
\end{subfigure}
  \begin{subfigure}[b]{0.32\textwidth}
          \centering
    \includegraphics[width=\linewidth]{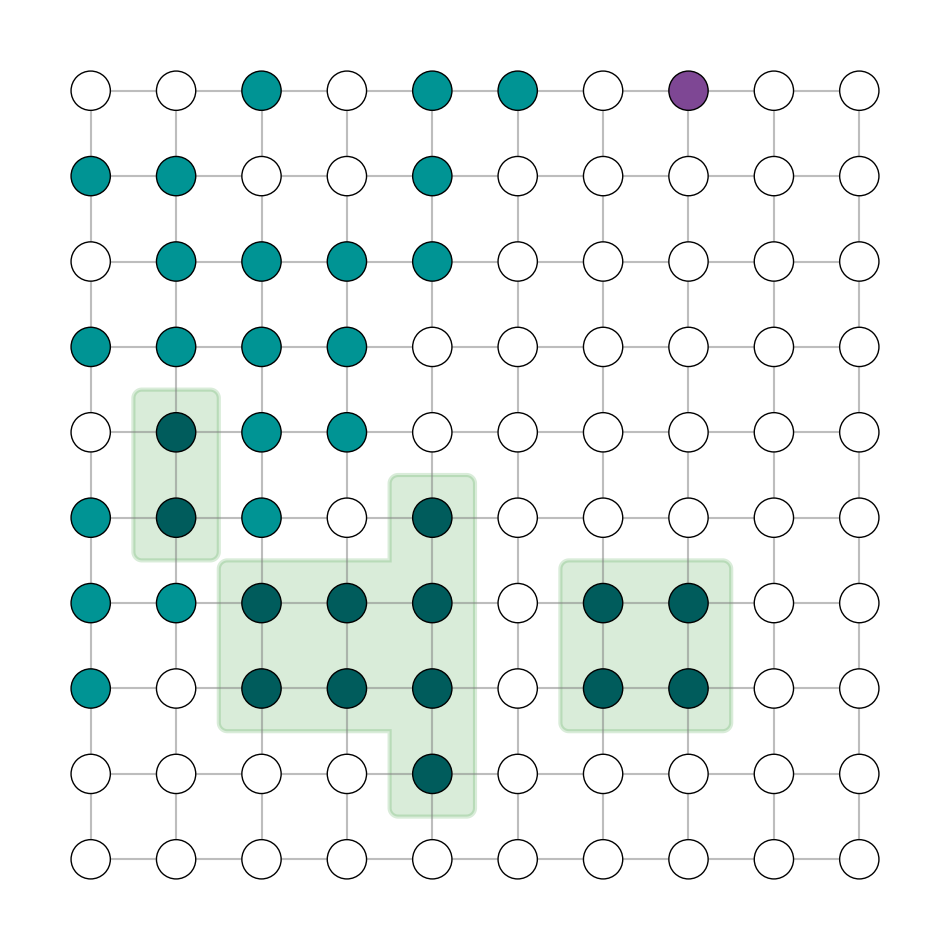}
    \caption{$G(S_2)$}
    \label{fig:proof_1ed_b}
  \end{subfigure}
  \begin{subfigure}[b]{0.32\textwidth}
          \centering
    \includegraphics[width=\linewidth]{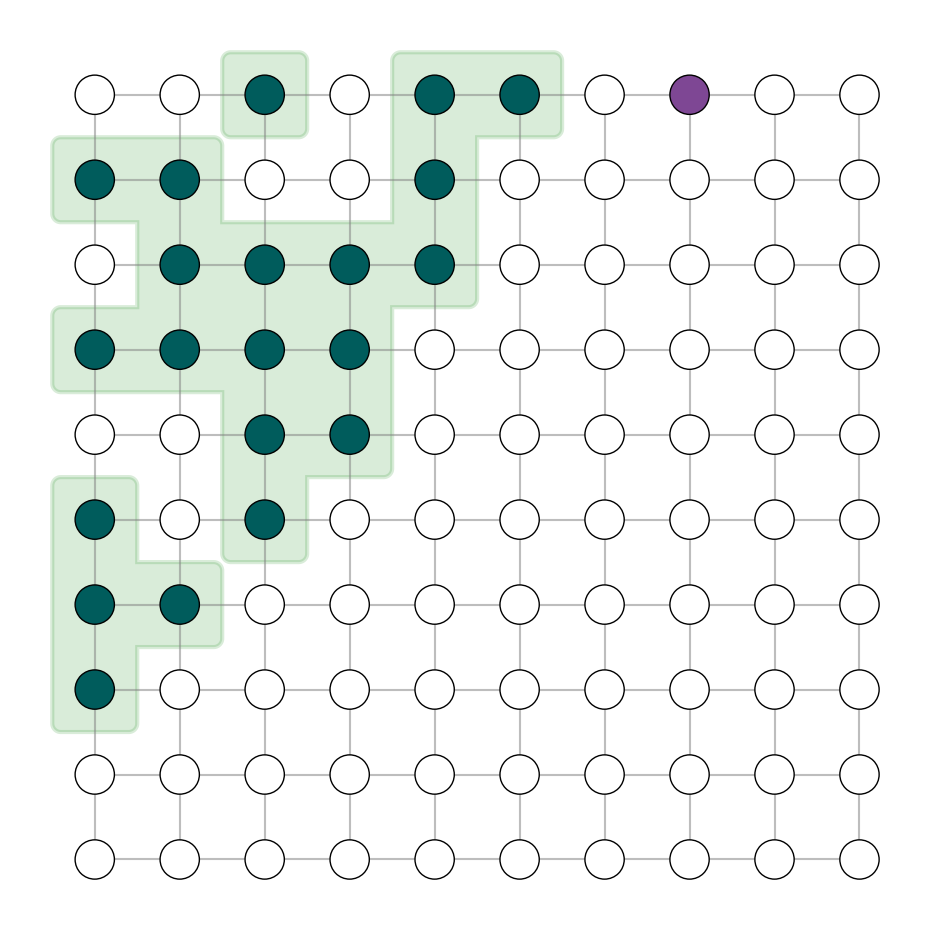}
    \caption{$G(S_3)$}
    \label{fig:proof_1ed_c}
  \end{subfigure}
\caption{Decomposing a feasible solution using extreme directions}
\label{fig:proof_1ed}
\end{figure}

Clearly $(\bar{\alpha}^2, \bar{\beta}^2)$ is feasible, and hence, we can apply the same procedure again by defining $S_2 = \arg \max_{(i,j) \in \validdriveset} \{ \bar{\alpha}^2_{ij}: \bar{\alpha}^2_{ij} < 0\}$ and if $S_2 = \emptyset$, we express $(\bar{\alpha}^2, \bar{\beta}^2)$ as a linear combination of extreme directions in $E_1 \cup E_2$. Else, we define $\Lambda_2 = -\bar{\alpha}^2_{ij}$, where $(i, j) \in S_2$, and proceed similarly. In Figure~\ref{fig:proof_1ed_b}, $G(S_2)$ has three components, and hence we construct the solution shown in Figure~\ref{fig:proof_1ed_c} using three extreme directions from $E_3$. After $\eta$ iterations, suppose $S_{\eta} = \emptyset$, then we can write 
\begin{align}
    (\bar{\alpha}^{\eta}, \bar{\beta}^{\eta}) = (\bar{\alpha}^{(\eta - 1)}, \bar{\beta}^{(\eta - 1)}) - \Lambda_{(\eta-1)} \sum_{v=1}^{\vartheta_{(\eta-1)}} (\alpha^{(\eta-1)v}, \beta^{(\eta-1)v}) = \sum_{u=1}^{\varphi_{\eta}} \lambda_{\eta u}(\alpha^{\eta u}, \beta^{\eta u})
\end{align}
Using the expansion for $(\bar{\alpha}^{(\eta - 1)}, \bar{\beta}^{(\eta - 1)})$, we can proceed backwards to conclude 
\begin{align}
    (\alpha, \beta) = (\bar{\alpha}^1, \bar{\beta}^1) = \sum_{n=1}^{\eta-1} \Lambda_n \sum_{v=1}^{\vartheta_n} (\alpha^{nv}, \beta^{nv}) + \sum_{u=1}^{\varphi_{\eta}} \lambda_{\eta u}(\alpha^{\eta u}, \beta^{\eta u})
\end{align}
where each $(\alpha^{\eta u}, \beta^{\eta u})$ in the above expression belongs to $E_1 \cup E_2$ and each $(\alpha^{nv}, \beta^{nv})$ is in $E_3$. Thus, $E = E_1 \cup E_2 \cup E_3$. 

Finally, we revisit the feasibility cuts for the extreme directions in $E$. Inequality \eqref{eq:two_way_benders_feasibility_cut} is trivially satisfied for $(\alpha, \beta) \in E_1 \cup E_2$. For directions in $E_3$, the third class of extreme directions, the feasibility cut \eqref{eq:two_way_benders_feasibility_cut} simplifies to 
\begin{align}
\sum_{(i,j) \in S} \drivevar{i}{j} \leq  \sum_{\substack{
  ((i,j),(k,l)) \in \gridarcset:\\
  (i,j) \in S,\,
  (k,l) \in D \setminus S
}} M \drivevar{k}{l}  \label{eq:simplified_cut}
\end{align}
However, notice that several of these constraints are dominated. For instance, consider the example in Figure~\ref{fig:proof1a}. The associated feasibility cut has fewer terms on the left-hand side and more terms on the right-hand side compared to the corresponding constraint for Figure \ref{fig:proof1b}. Thus, it is sufficient to consider only the extreme directions for subsets $S$ whose removal results in a single connected graph. Additionally, if the neighboring cells of $S$ contain the entrance/exit $(p, q)$, the right-hand side of \eqref{eq:simplified_cut} has at least one active $y$ variable and hence the inequality is always satisfied. In other words, we can restrict our attention to subsets $S$ that are vertex separators $V$ which exclude the entrance/exit and induce exactly two components, as described in Section \ref{sec:two-way-branch-and-cut}. Based on this observation, we derive the following cut.
\begin{align}
 \sum_{(i,j) \in \twowaypartitioncomp{\entrancex \entrancey}} \drivevar{i}{j} \leq \sum_{(k,l) \in \vertexcutset} M |\gridinvadjset{k}{l} \cap \twowaypartitioncomp{\entrancex \entrancey}| \drivevar{k}{l} \label{eq:benders_tw_vi}
\end{align}
Inequality \eqref{eq:benders_tw_vi} is the same as the valid inequality derived in Proposition \ref{prop:2w_aggregate}. Using the results from Propositions \ref{prop:2w_disaggregate} and \ref{prop:2w_strengthened}, we conclude that the integer feasible solutions $(x, y)$ from $\mathcal{F}_{\textsc{2W}}^{\textsc{flow}}$ and $\mathcal{F}_{\textsc{2W}}^{\textsc{CCC}}$ are the same.
\end{proof}

\begin{theorem}
Suppose $(x, y, z)$ denotes an integer feasible parking, driving-field, and direction assignment. The formulations $\mathcal{F}_{\textsc{1W}}^{\textsc{flow}}$ and $\mathcal{F}_{\textsc{1W}}^{\textsc{CCC}}$ induce the same set of integer feasible values for $(x,y,z)$.
\end{theorem}
\begin{proof}
We again write the Benders reformulation of the problem $\mathcal{F}_{\textsc{1W}}^{\textsc{flow}}$ using two subproblems, since the constraints with $f$ and $g$ variables have a block diagonal structure for a fixed $(x, y, z)$. 
\begin{align}
    \max & \, \, \, \, 0  && \notag \\
      \text{s.t. }  & \sum_{(k,l) \in \gridadjset{i}{j}} \flowvar{i}{j}{k}{l} - \sum_{(k,l) \in \gridinvadjset{i}{j}} \flowvar{k}{l}{i}{j} \geq \drivevar{i}{j} \qquad && \forall \, (i,j) \in \validdriveset \setminus \{(\entrancex, \entrancey)\} \label{eq:floworigin_1w_sp}\\
     & \sum_{(i,j) \in \gridinvadjset{\entrancex}{\entrancey}} \flowvar{i}{j}{\entrancex}{\entrancey} \leq \sum_{(i,j) \in \validdriveset} \drivevar{i}{j} -1 \qquad &&  \label{eq:flowdest_1w_sp}\\
    & \flowvar{i}{j}{k}{l} \leq M \countvar{i}{j}{k}{l} \qquad && \forall \, ((i,j), (k,l)) \in \gridarcset \label{eq:flowhead_1w_sp} \\
    & \flowvar{i}{j}{k}{l} \geq 0 \qquad && \forall \, ((i,j), (k,l)) \in \gridarcset \label{eq:flowvar_1w_nonneg}
\end{align}
Let $\alpha$ be the dual variables for \eqref{eq:floworigin_1w_sp}--\eqref{eq:flowdest_1w_sp}, and $\beta$ denote the duals for \eqref{eq:flowvar_1w_nonneg}. The dual of the first subproblem is
\begin{align}
    \min & \sum_{(i,j) \in \validdriveset \setminus \{(\entrancex, \entrancey)\}} \drivevar{i}{j} \alpha_{ij} + \Bigg(\sum_{(i,j) \in \validdriveset} \drivevar{i}{j} -1 \Bigg) \alpha_{pq} + \sum_{((i,j), (k,l)) \in \gridarcset}  M \countvar{i}{j}{k}{l} \beta_{ij, kl}  \\
    \text{s.t. } & \eqref{eq:sp_dual_1} \text{--} \eqref{eq:sp_dual_5}
\end{align}
Note that the feasible region is identical to that of the dual subproblem in Theorem \ref{theorem:two_way}. Hence, the non-trivial extreme directions in $E$ yield feasibility cuts of the form
\begin{align}
\sum_{(i,j) \in S} \drivevar{i}{j} \leq  \sum_{\substack{
  ((i,j),(k,l)) \in \gridarcset:\\
  (i,j) \in S,\,
  (k,l) \in D \setminus S
}} M \countvar{i}{j}{k}{l}  \label{eq:simplified_cut_1w_a}
\end{align}
Similarly, the second subproblem can be written as 
\begin{align}
    \max & \, \, \, \, 0  && \notag \\
     \text{s.t. } & \sum_{(k,l) \in \gridinvadjset{i}{j}} \flowvardash{k}{l}{i}{j} - \sum_{(k,l) \in \gridadjset{i}{j}} \flowvardash{i}{j}{k}{l} \geq \drivevar{i}{j}  \qquad && \forall \, (i,j) \in \validdriveset \setminus \exitset \label{eq:flow_dash_origin_singlelane_sp}\\
    & \sum_{(i,j) \in \gridadjset{\exitx}{\exity}} \flowvardash{\exitx}{\exity}{i}{j} \leq \sum_{(i,j) \in \validdriveset} \drivevar{i}{j}  - 1 \qquad && \label{eq:flowdashdest_singlelane_sp}\\
    & \flowvardash{i}{j}{k}{l} \leq M \countvar{i}{j}{k}{l}  \qquad && \forall \, ((i,j), (k,l)) \in \gridarcset \label{eq:flow_dash_tail_singlelane_sp}\\
    & \flowvardash{i}{j}{k}{l} \geq 0 \qquad && \forall \, ((i,j), (k,l)) \in \gridarcset \label{eq:g_non_negative_sp}
\end{align}
Let $\pi$ denote the dual variables for \eqref{eq:flow_dash_origin_singlelane_sp}--\eqref{eq:flowdashdest_singlelane_sp} and $\xi$ represent the dual variables for \eqref{eq:flow_dash_tail_singlelane_sp}. 
The dual subproblem for the above formulation is given by 
\begin{align}
    \min &\sum_{(i,j) \in \validdriveset \setminus \exitset} \drivevar{i}{j}\pi_{ij} + \Bigg(\sum_{(i,j) \in \validdriveset} \drivevar{i}{j}  - 1 \Bigg) \pi_{rs} + \sum_{((i,j), (k,l)) \in \gridarcset}  M \countvar{i}{j}{k}{l} \xi_{ij, kl} \label{eq:sp_dual_1w}
\end{align}
\begin{align}
     \text{s.t. } & - \pi_{ij} + \pi_{kl} + \xi_{ij, kl} \geq 0 && \forall \, ((i,j), (k,l)) \in \gridarcset, (i, j), (k, l) \neq (r, s)  \label{eq:sp_dual_1w_6} \\
     & - \pi_{ij} + \xi_{ij, rs} \geq 0 && \forall \,  (i,j) \in \gridinvadjset{r}{s}  \label{eq:sp_dual_1w_7} \\
    & \pi_{rs} + \xi_{rs, kl} \geq 0 && \forall \,  (k,l) \in \gridadjset{r}{s}  \label{eq:sp_dual_1w_8} \\
    & \pi_{ij} \leq 0 && \forall \, (i,j) \in \validdriveset \setminus \{(r, s)\} \label{eq:sp_dual_1w_9} \\
    & \pi_{rs} \geq 0 && \label{eq:sp_dual_1w_9.5} \\
     & \xi_{ij, kl} \geq 0 && \forall \, ((i,j), (k,l)) \in \gridarcset \label{eq:sp_dual_1w_10} 
\end{align}
The extreme directions of the pointed cone induced by \eqref{eq:sp_dual_1w_6}--\eqref{eq:sp_dual_1w_10} are given below. The justification is identical to that in Theorem \ref{theorem:two_way} and is therefore omitted.
\begin{enumerate}
\item A singleton $F_1$ consisting of a vector with $\pi_{rs} = 1$, and all remaining variables equal to zero.
\item A set of directions $F_2$, one for each arc $((i,j),(k,l)) \in \gridarcset$, with $\xi_{ij,kl} = 1$ and all other entries equal to zero.
\item A collection of vectors $F_3$, where each vector is defined by a subset 
$S \subseteq \validdriveset \setminus \{(r,s)\}$ such that the induced subgraph 
$\grid(S)$ is connected, and has components 
$\pi_{ij} = -1$ for all $(i,j) \in S$, 
$\xi_{ij,kl} = 1$ for all $((i,j),(k,l)) \in \gridarcset$ with $(k,l) \in S$ and $(i,j) \in \validdriveset \setminus S$, 
and all remaining components equal to zero.
\end{enumerate}
The extreme directions from $F_1 \cup F_2$ produce trivial feasibility cuts. The cuts from those in $F_3$ can be written as
\begin{align}
\sum_{(i,j) \in S} \drivevar{i}{j} \leq  \sum_{\substack{
  ((i,j),(k,l)) \in \gridarcset:\\
  (i,j) \in D \setminus S,\,
  (k,l) \in S
}} M \countvar{i}{j}{k}{l}  \label{eq:simplified_cut_1w_b}
\end{align}
Constraints \eqref{eq:simplified_cut_1w_a} and \eqref{eq:simplified_cut_1w_b} are the same as the valid inequalities \eqref{eq:aggregate_vi_1w_a} and \eqref{eq:aggregate_vi_1w_b}, respectively. Together with Proposition \ref{prop:1w_aggregate}, this implies that the integer feasible regions of $\mathcal{F}_{\textsc{1W}}^{\textsc{flow}}$ and $\mathcal{F}_{\textsc{1W}}^{\textsc{CCC}}$ coincide.
\end{proof}

\section{Practical considerations}
\label{sec:practical_considerations}

\textbf{Turn restrictions for two-ways:} The two-way formulations discussed so far ensure the connectivity and accessibility of parking lots and driving fields. However, the driving fields are not translated into the actual lane design. This is challenging to accomplish purely with the drive field or flow variables. The lack of designated lanes makes it tricky to anticipate vehicle paths to parking stalls and identify the issues they might encounter while turning. \cite{stephan2021layout} proposed additional restrictions using the drive and flow variables, which provide sufficient space for vehicles to turn. We can incorporate similar constraints using only the drive variables to make the method compatible with the branch-and-cut framework. Developing a fully generalized model that works in all cases seems challenging, as it largely depends on the required level of maneuverability. For example, tight turns might be sufficient for long-term parking or when using robotic units. This section provides some ways to `polish' the driveway layouts. These modifications mainly attempt to (1) keep the angle of turns at $90^{\circ}$ and the width of the driveways uniform and (2) avoid sharp turns. These adjustments become particularly crucial when the rasterization is coarse. For this section, we denote $y_{ij}$ and $z_{ij, kl}$ more elaborately as $y_{(i,j)}$ and $z_{(i,j),(k,l)}$ to avoid confusion over the indices in the examples.

\begin{figure}[H]
    \centering
    \begin{subfigure}[b]{0.49\textwidth}
        \centering
        \includegraphics[width=0.44\textwidth]{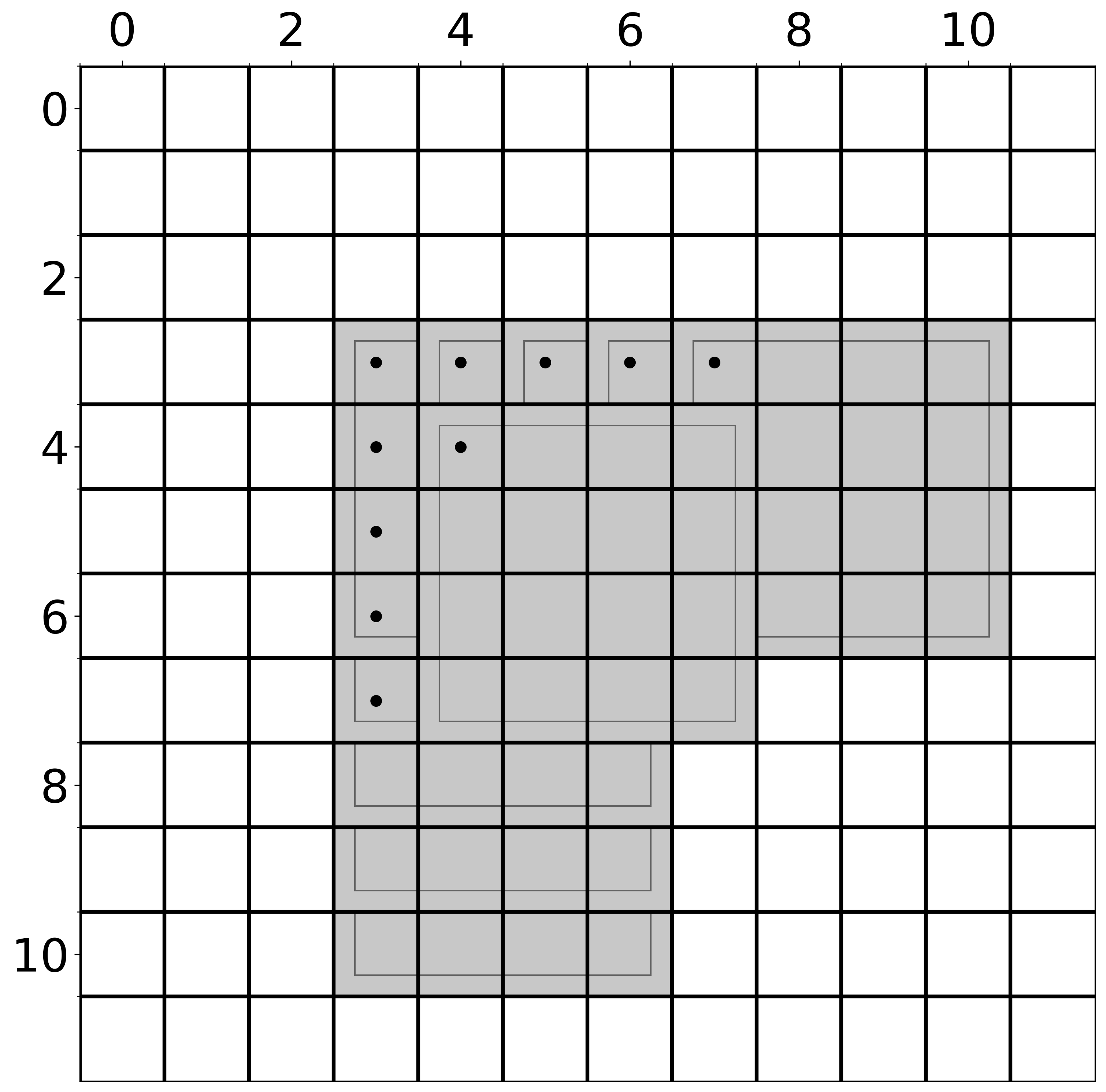}
    \includegraphics[width=0.44\textwidth]{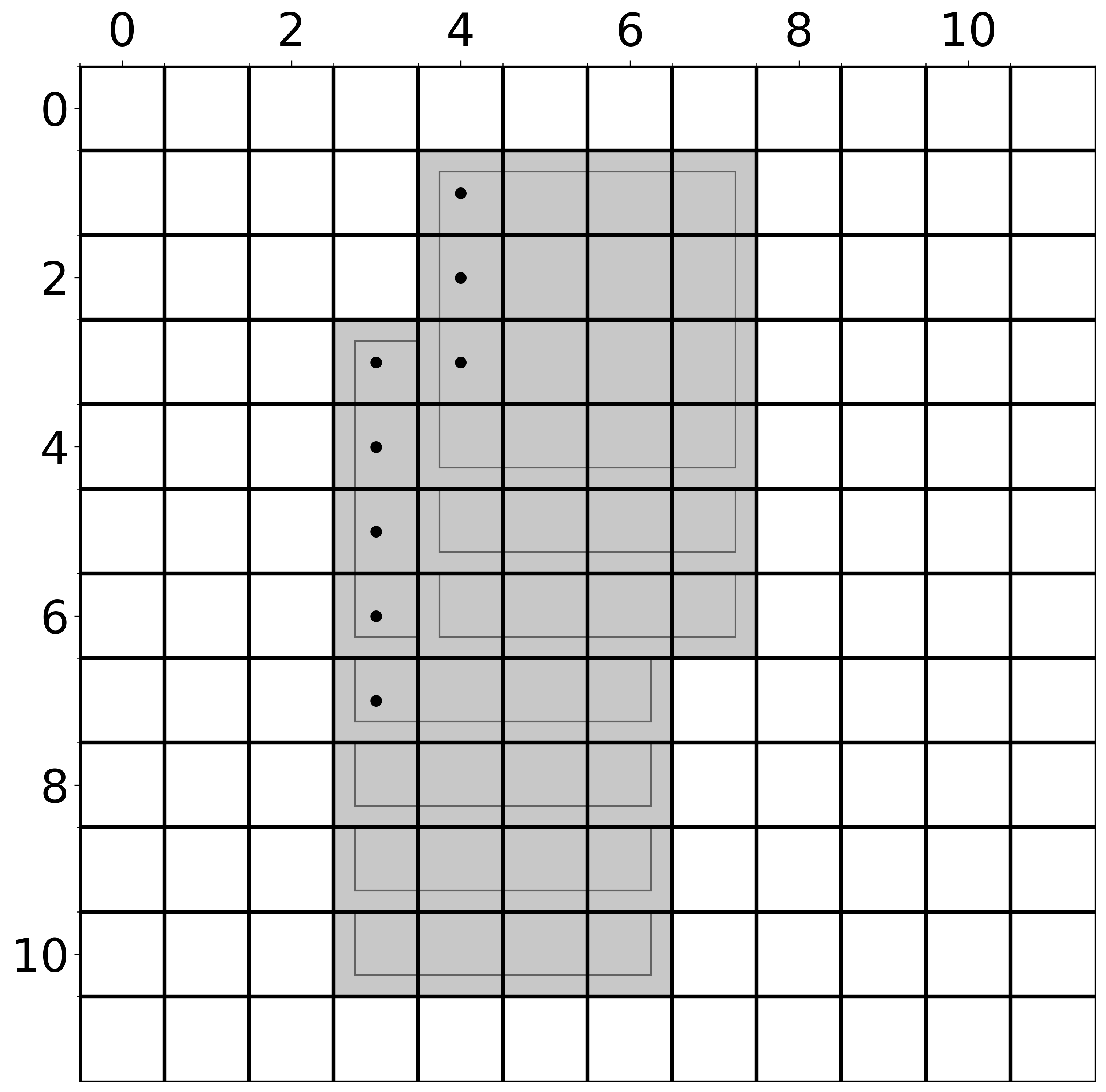}
        \caption{An indistinct perpendicular turn and a sharp turn}
        \label{fig:sharp_turn_two_way}   
    \end{subfigure}
    \begin{subfigure}[b]{0.49\textwidth}
        \centering
    \includegraphics[width=0.49\textwidth]{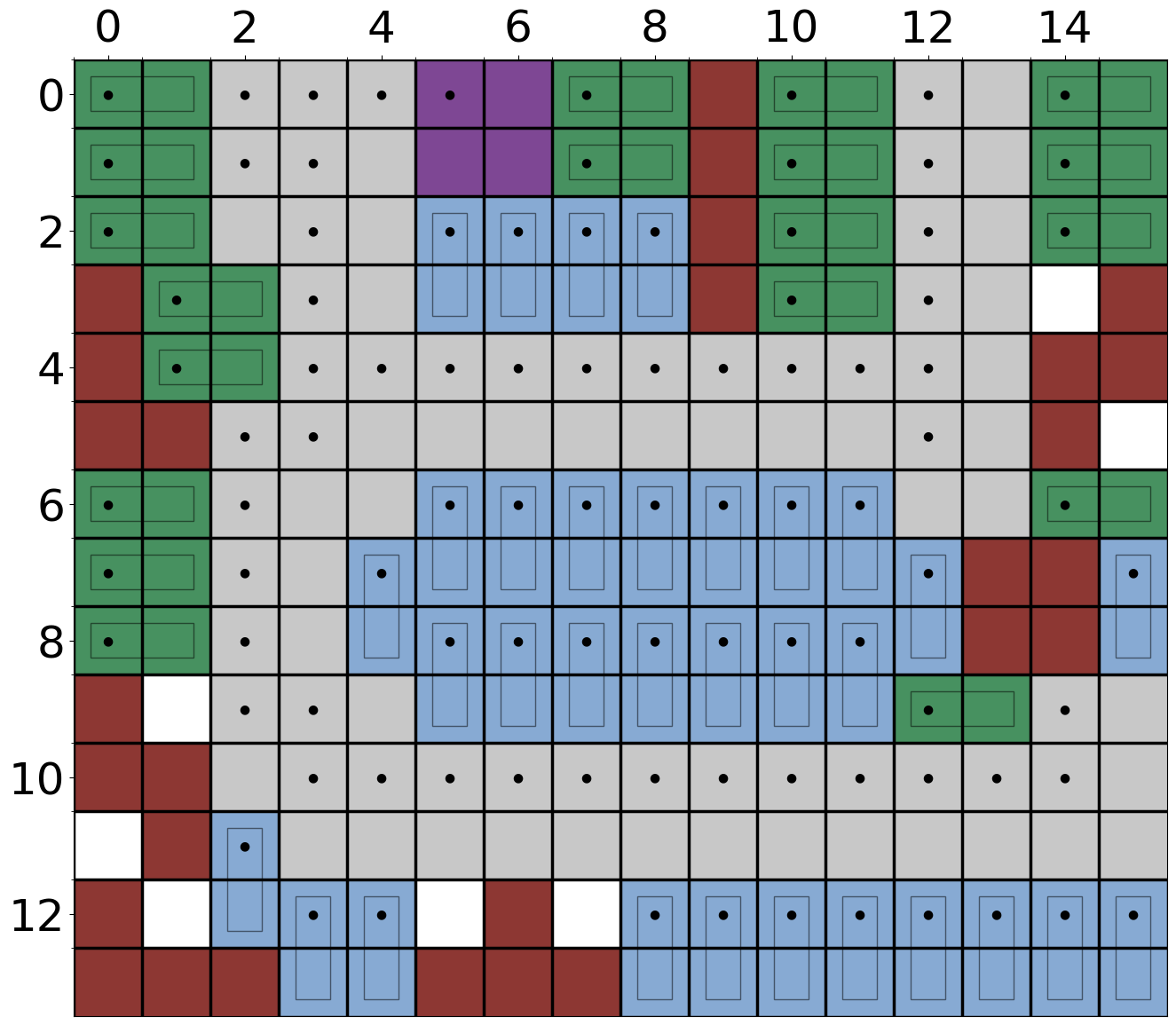}
    \includegraphics[width=0.49\textwidth]{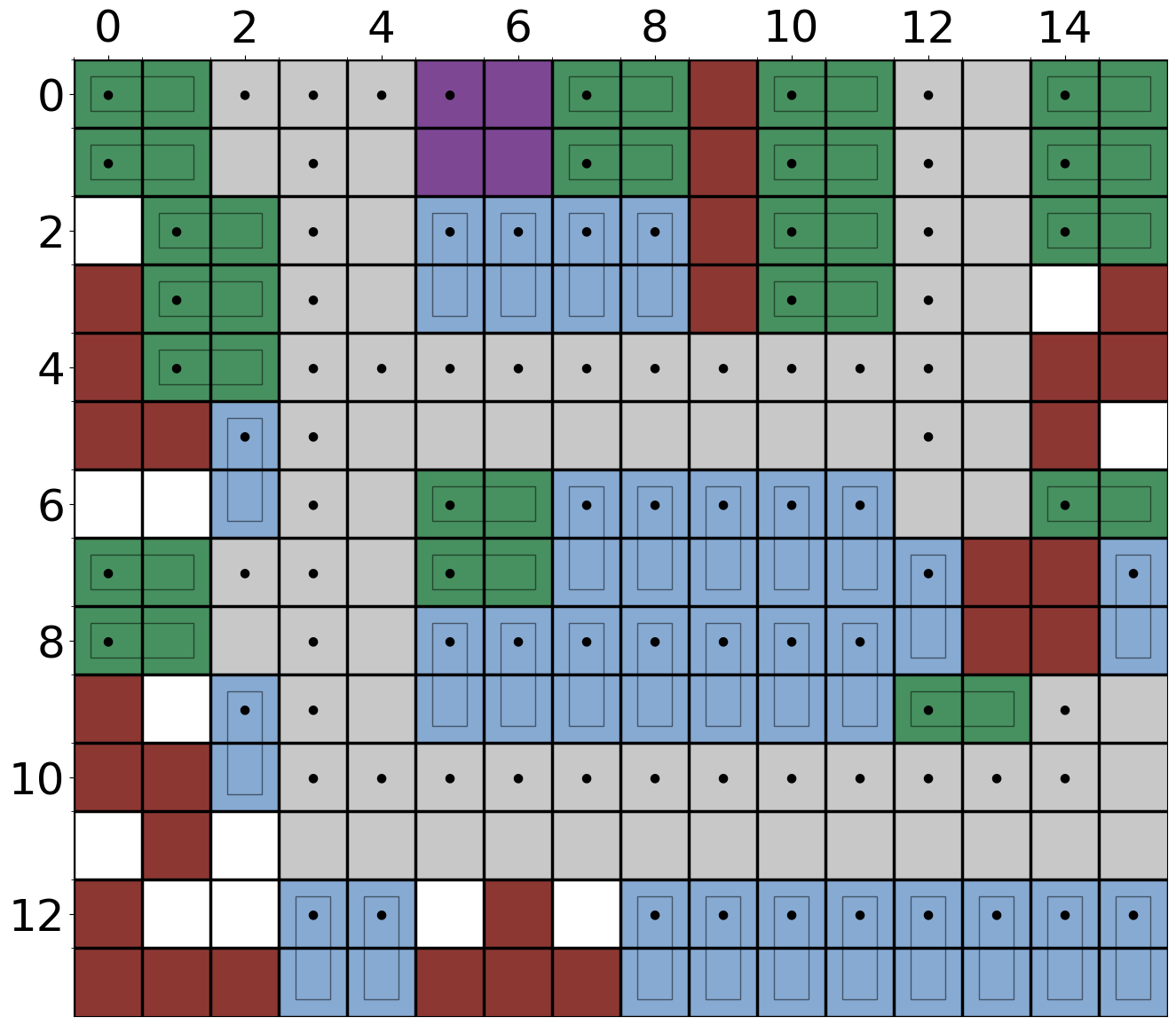}
        \caption{Layouts without and with turn restrictions}
        \label{fig:turn_restrictions}     
    \end{subfigure}
    \caption{Imposing turn restrictions for two-way configurations}
    \label{fig:imposing_restrictions} 
\end{figure}

Since driving fields can overlap, instances shown in Figure \ref{fig:sharp_turn_two_way} can occur (here $\drivewidth = 4$), where the perpendicular turns are indistinct, leading to non-uniform driving lanes, or the layouts may feature sharp turns, where drivers take consecutive right and left turns, complicating maneuvering. (We may want to permit the latter if we want to approximate curved driving sections using finer rasterization.) With additional restrictions, the number of stalls is reduced by one, but the driveways have uniform widths and well-defined turns. One solution to circumvent such issues is to prohibit overlapping driving fields. This can be modeled by defining arcs in $\grid$ that connect every node to other nodes that are $\drivewidth$ hops away in longitudinal (north-south) and lateral (east-west) directions. However, such an approach would significantly reduce the number of parking stalls, as driving lanes can only be created in lengths that are multiples of $\drivewidth$. Instead, we could restrict combinations of driving fields from being active. The example in Figure \ref{fig:turn_restrictions} depicts a layout without and with these restrictions. This parking lot solution is based on the R2 resolution and has a drive field width $\drivewidth = 2$.

 \begin{figure}[H]
    \centering
    \begin{subfigure}[b]{0.24\textwidth}
        \centering
        \includegraphics[width=\textwidth]{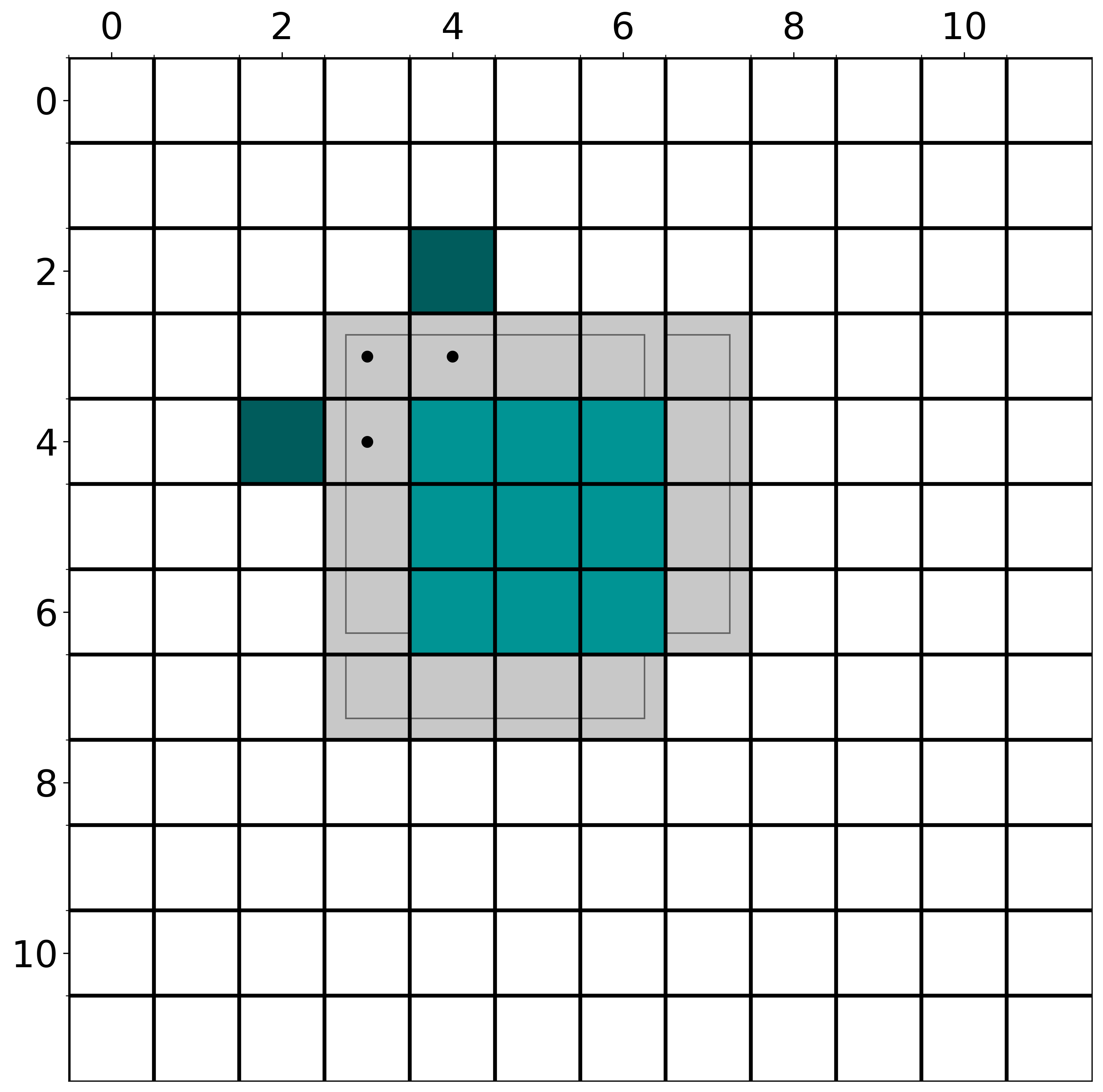}
    \end{subfigure}
    \hfill
    \begin{subfigure}[b]{0.24\textwidth}
        \centering
        \includegraphics[width=\textwidth]{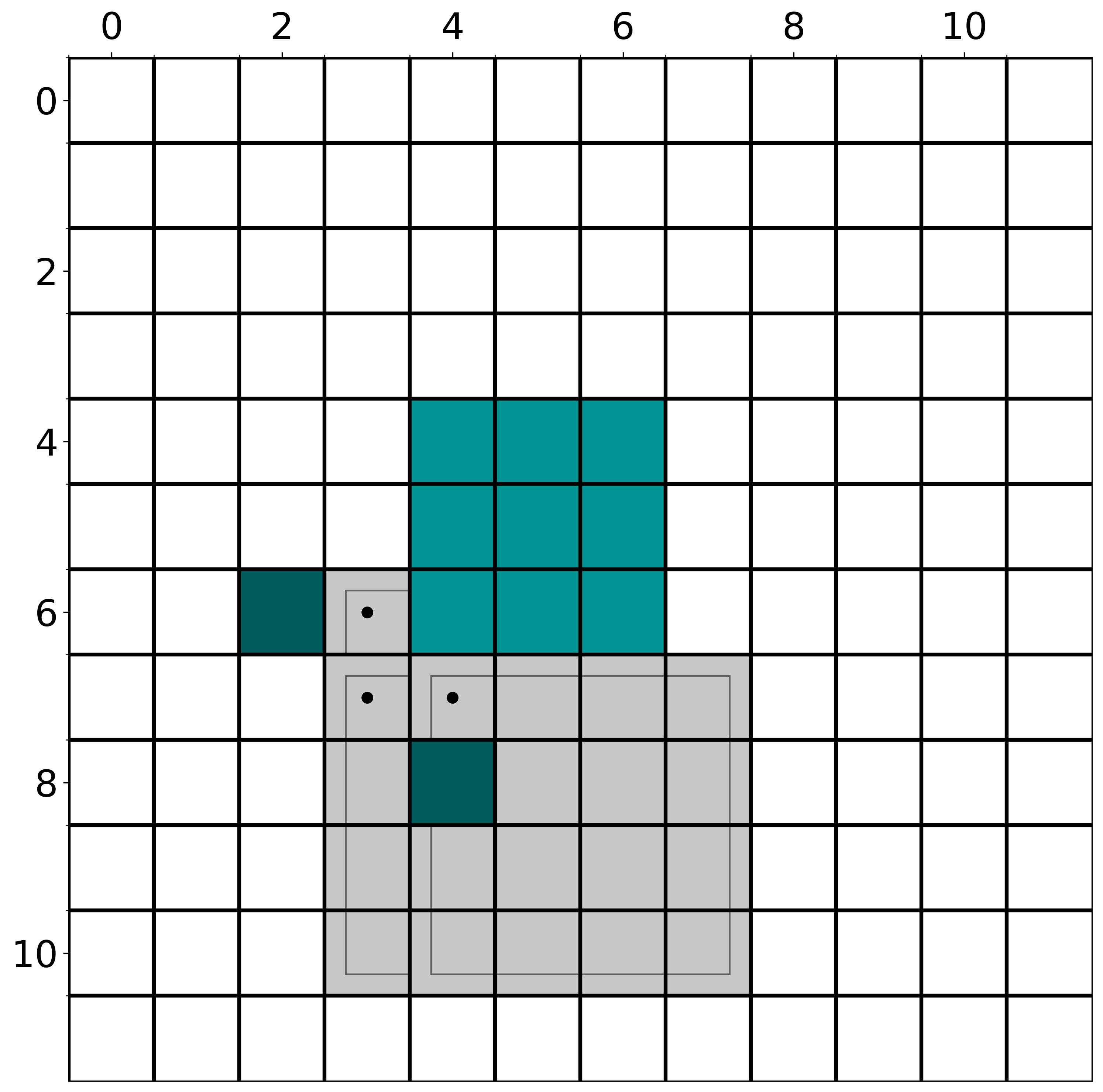}
    \end{subfigure}
    \hfill
    \begin{subfigure}[b]{0.24\textwidth}
        \centering
        \includegraphics[width=\textwidth]{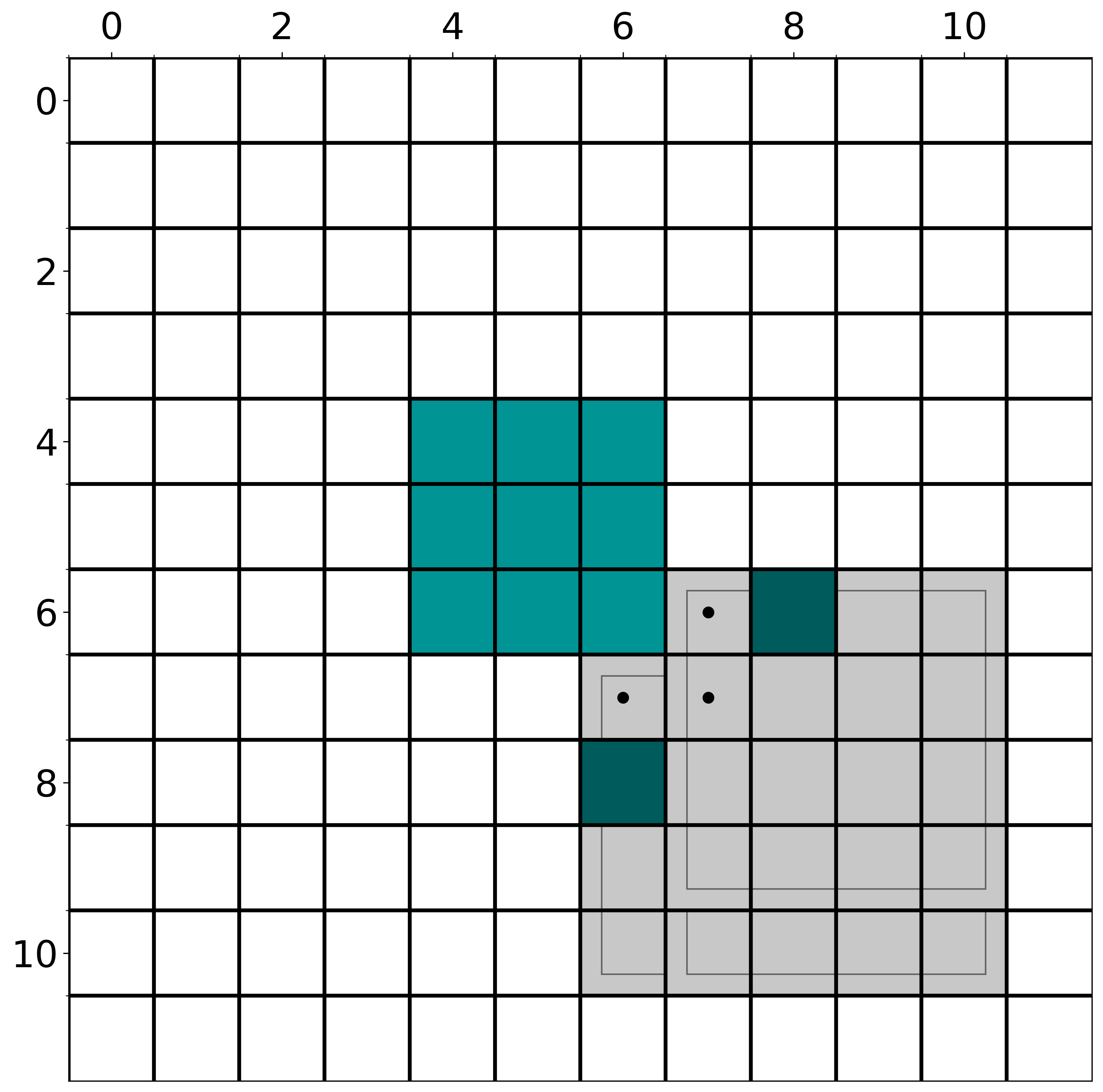}
    \end{subfigure}
    \hfill
    \begin{subfigure}[b]{0.24\textwidth}
        \centering
        \includegraphics[width=\textwidth]{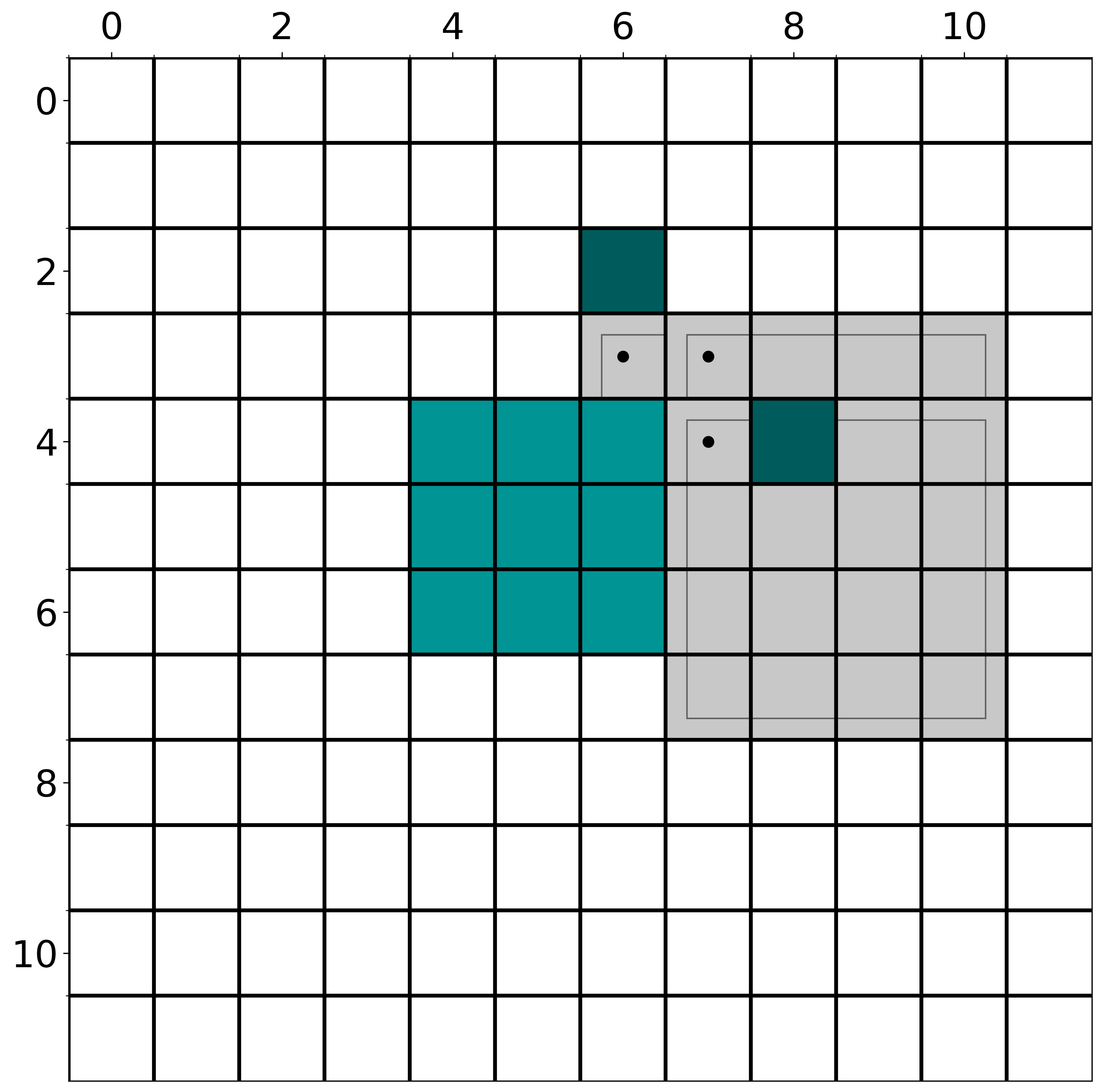}
    \end{subfigure}
    \caption{Restrictions for handling turn-related issues}
    \label{fig:turntypes}    
\end{figure}

In this approach, a turn is characterized by three active driving cells in a 2 $\times$ 2 field. Four possibilities can occur, as shown using the anchors in the example in Figure \ref{fig:turntypes}, where $\drivewidth = 4$. To avoid indistinct perpendicular turns and keep the width of the driveway uniform, we can prevent cells in the teal $3 \times 3$ (which in the general case would be $(\delta - 1)\times (\delta - 1)$) block from being active if the three cells around its corner have active driving fields. For example, in the panel on the left, we set constraints of the form $y_{(i + k, j + k)} + y_{(i+1,j)} + y_{(i, j)} + y_{(i,j+1)} \leq 3$ where $k \in \{1, \ldots, \delta - 1\}$, which forces at most three of these $y$ variables to be active. Similar constraints can be written in the other cases shown in the figure. Additionally, to prevent sharp turns, one could add such constraints using the cells shown in dark teal. They would prevent zig-zag driveways and can also be extended to include more cells, controlling exactly where the next turn is allowed. In the example from the left-most panel, these constraints would be of the form $y_{(i + 1, j - 1)} + y_{(i+1,j)} + y_{(i, j)} + y_{(i,j+1)} \leq 3$ and $y_{(i -1, j + 1)} + y_{(i+1,j)} + y_{(i, j)} + y_{(i,j+1)} \leq 3$.

\begin{figure}[H]
  \centering
  \begin{subfigure}[b]{0.3\textwidth}
      \centering
    \includegraphics[width=0.89\textwidth]{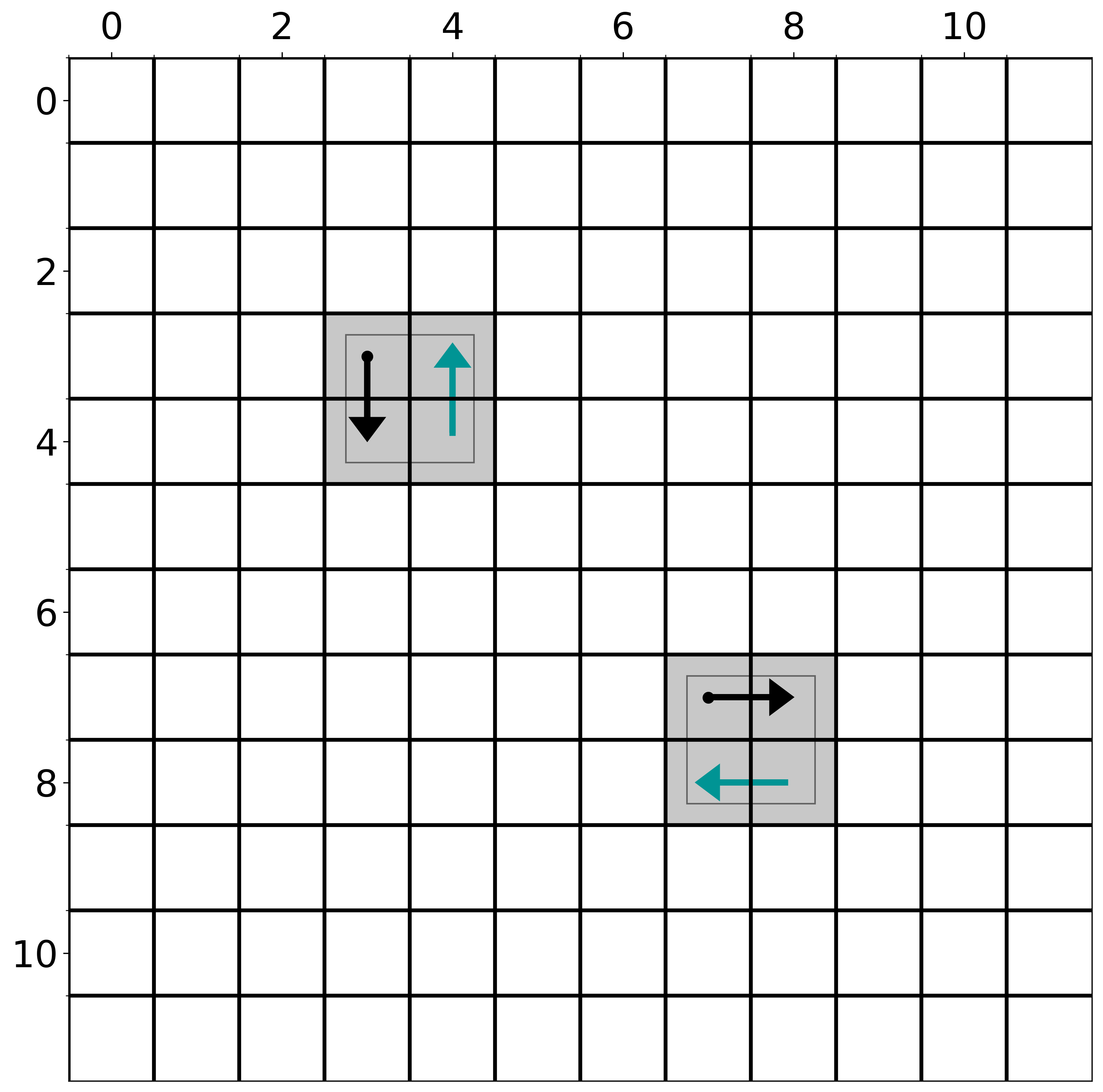}
    \caption{Preventing infeasible overlaps}
    \label{fig:one_way_opposite}
  \end{subfigure}
  \hfill
  \begin{subfigure}[b]{0.68\textwidth}
      \centering
    \includegraphics[width=0.395\textwidth]{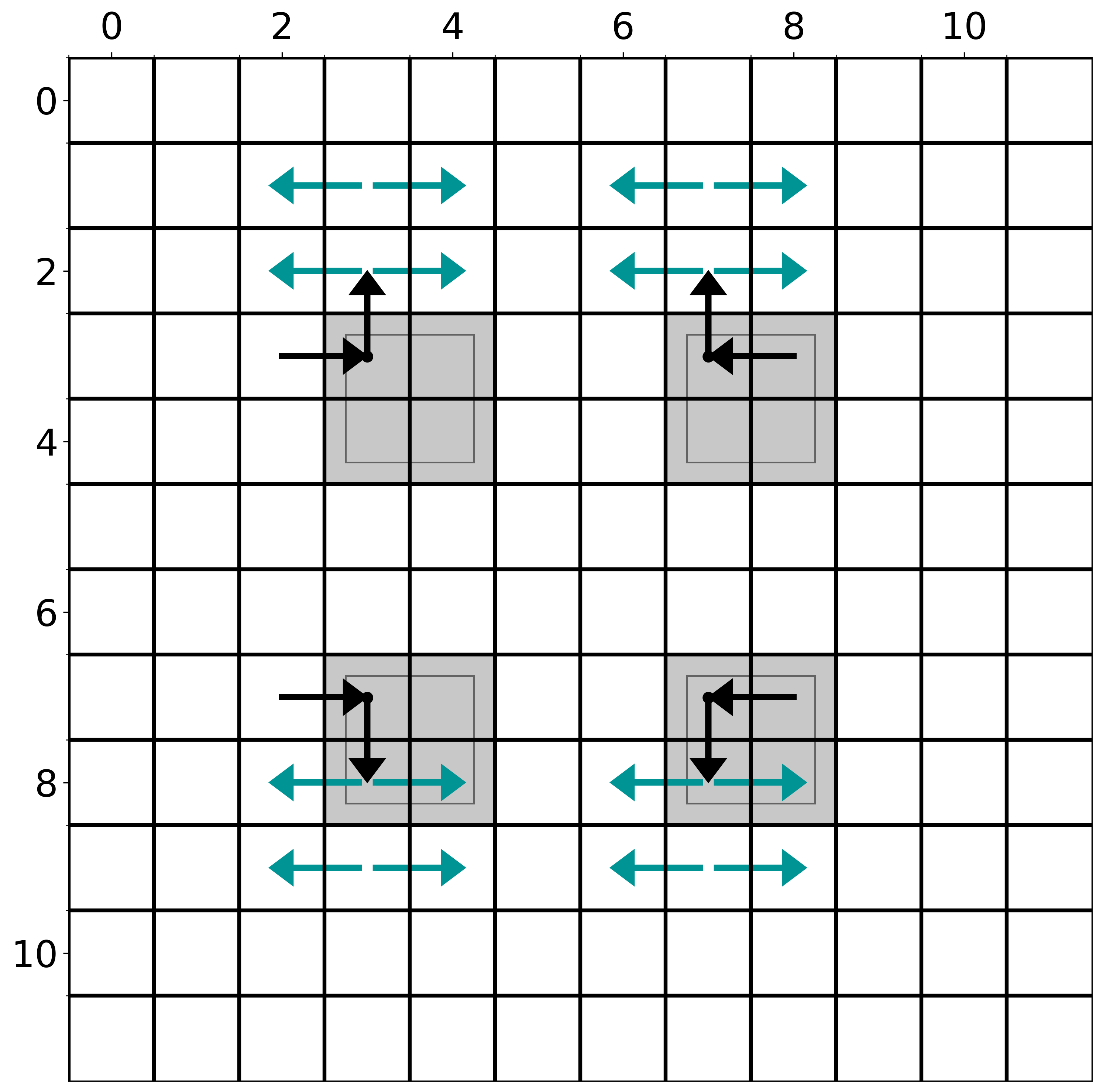}
    \hspace{3mm}
    \includegraphics[width=0.395\textwidth]{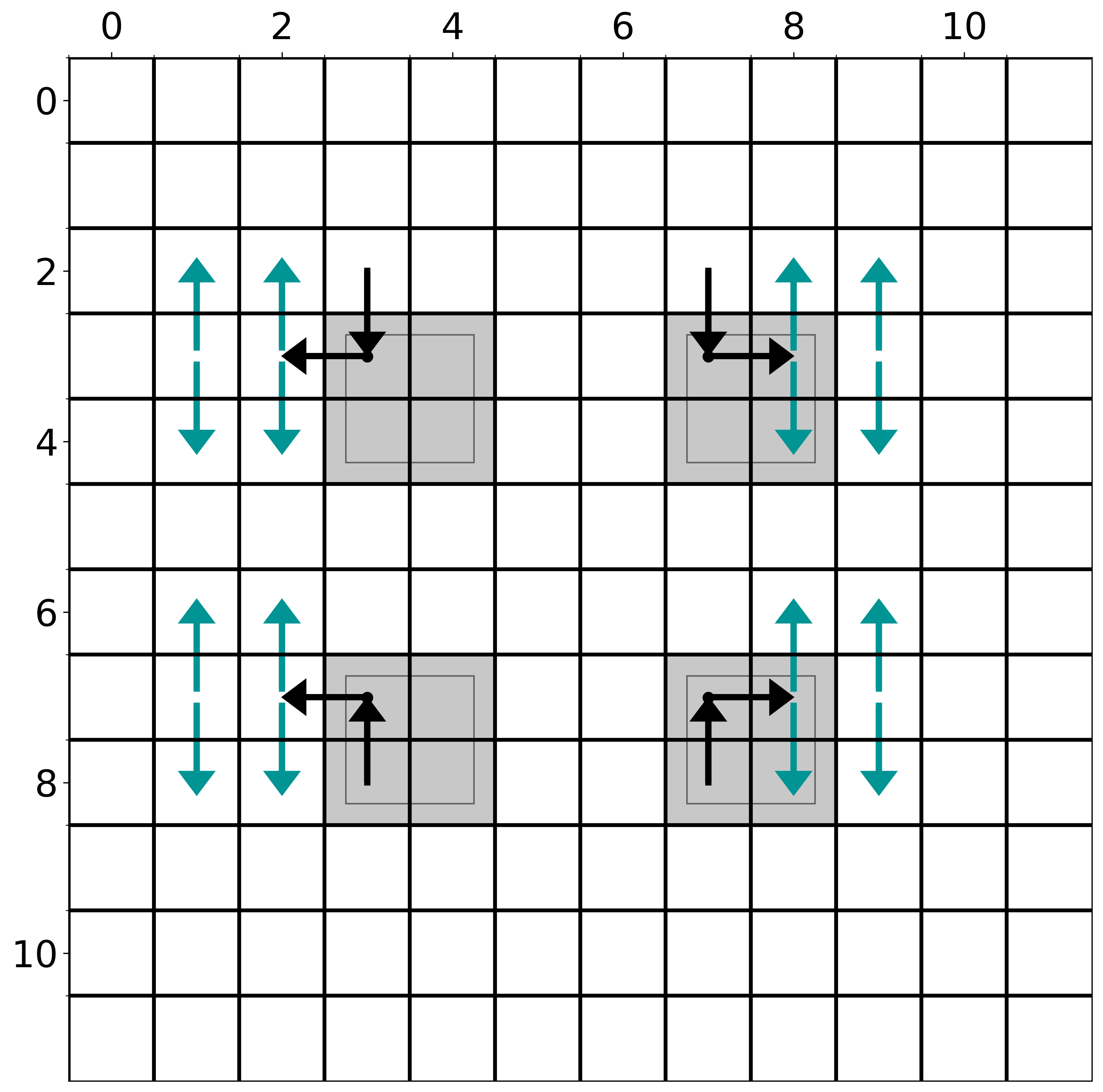}
    \caption{Avoiding sharp turns}
        \label{fig:one_way_avoid_sharp_turns}
  \end{subfigure}
\caption{Combinations of direction variables that are avoided for $\delta =2$. Black arcs represent active $z$ variables, and teal arcs are restricted to be inactive.}
\label{fig:one_way_turns}
\end{figure}

\textbf{Turn restrictions for one-ways:}
For one-way configurations, the presence of the direction variables $z$ gives us greater control over the turning restrictions. To avoid overlapping driving fields with opposite-moving traffic, we first restrict combinations of direction variables as shown in Figure \ref{fig:one_way_opposite}, where $\delta = 2$. More generally, if a direction variable $z_{(i, j), (i+1, j)}$ is active, we require variables in the opposite direction $z_{(i+1, j+k), (i, j+k)}$ to be inactive for $k \in \{1, \ldots, \delta - 1\}$. Likewise, if $z_{(i, j), (i, j+1)}$ is active, then $z_{(i+k, j+1), (i+k, j)}$ must be inactive for the same set of $k$ values. These constraints can be simply enforced using an upper bound of one on the sum of pairs of combinations that are restricted. For example, in Figure \ref{fig:one_way_opposite}, we set $z_{(3,3),(4,3)} + z_{(4,4),(3,4)} \leq 1$ and $z_{(7,7),(7,8)} + z_{(8, 8),(8, 7)} \leq 1$.

\begin{figure}[H]
\centering
\begin{subfigure}[b]{0.49\textwidth}
    \centering
    \includegraphics[width=0.49\textwidth]{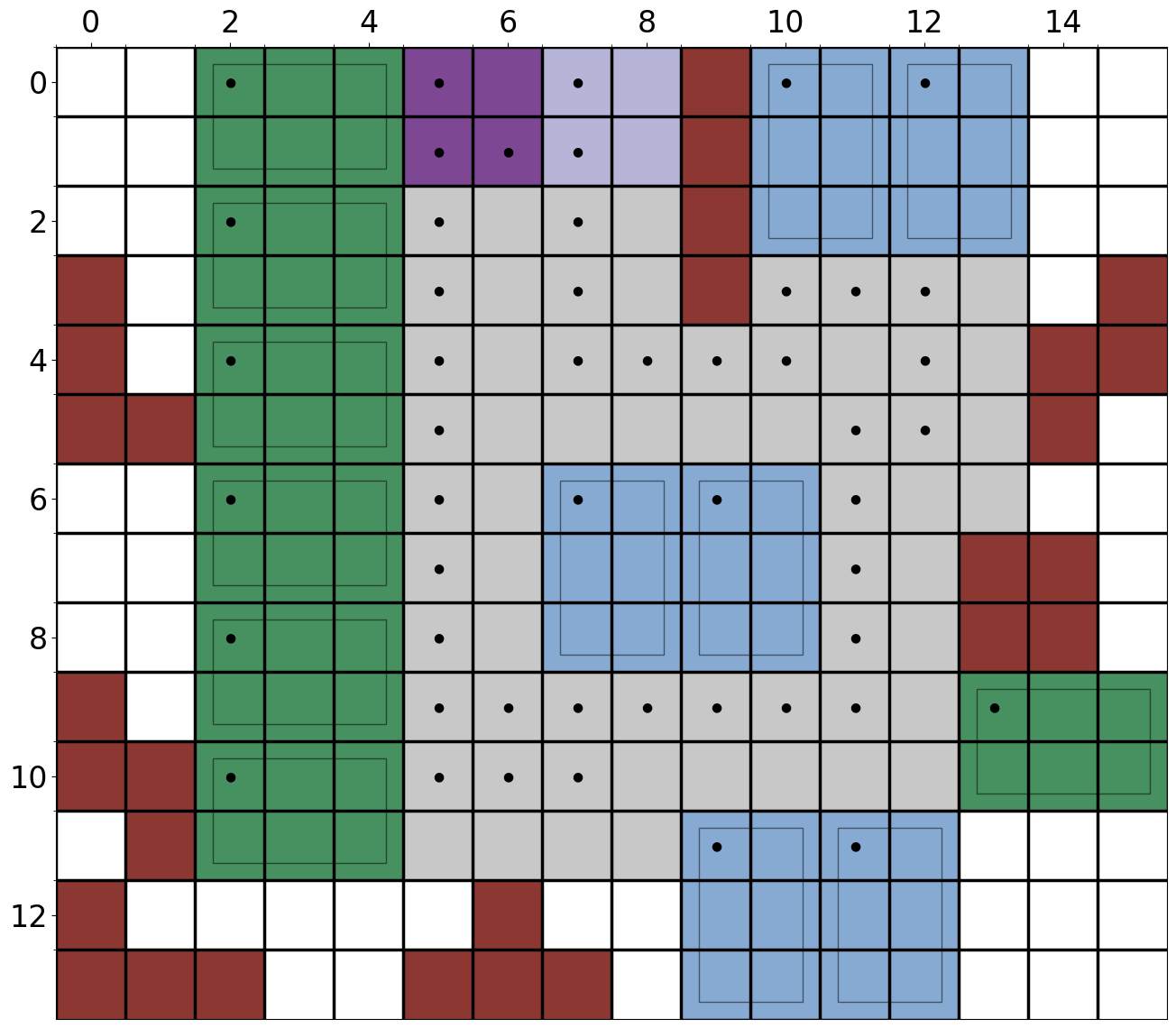}
    \includegraphics[width=0.49\textwidth]{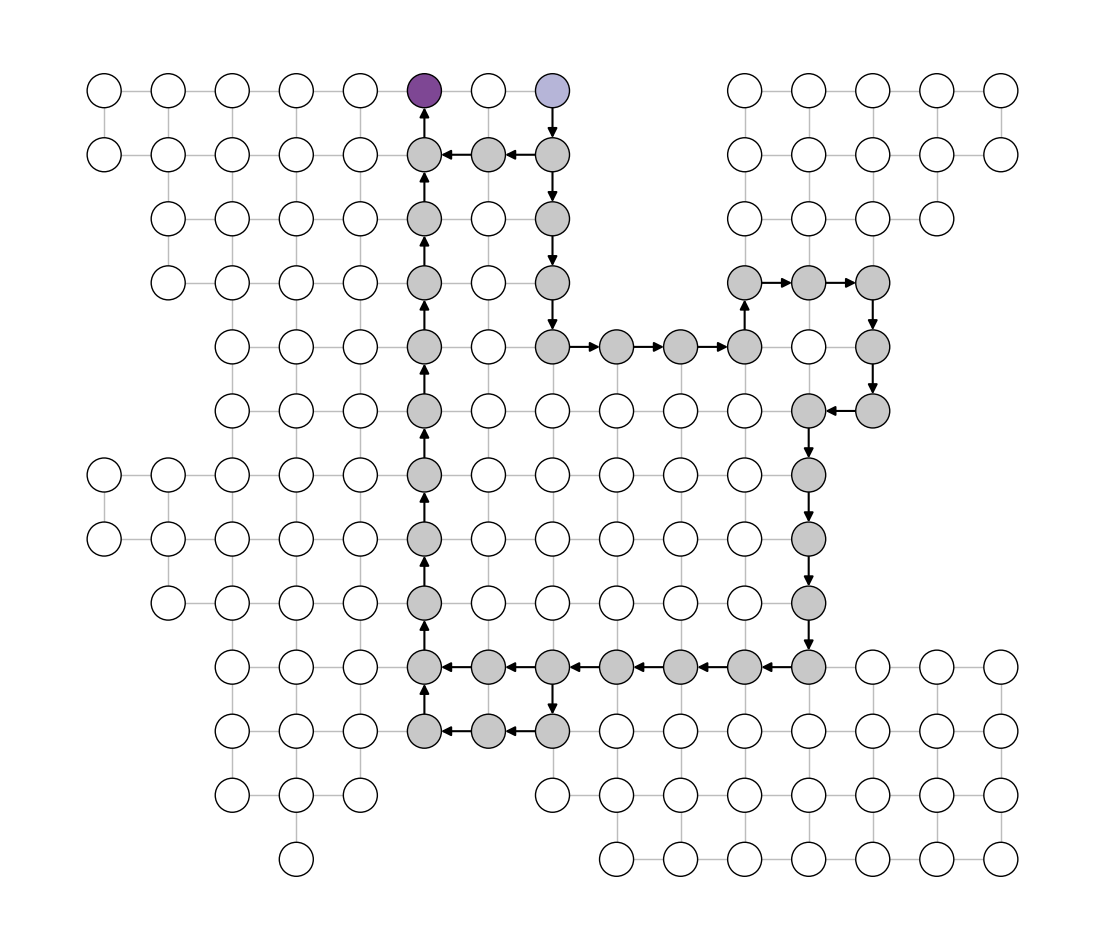}
    \caption{Lot and directions without infeasible overlaps}
    \label{fig:overlap_constraints}
\end{subfigure} 
\begin{subfigure}[b]{0.49\textwidth}
\centering
\includegraphics[width=0.49\linewidth]{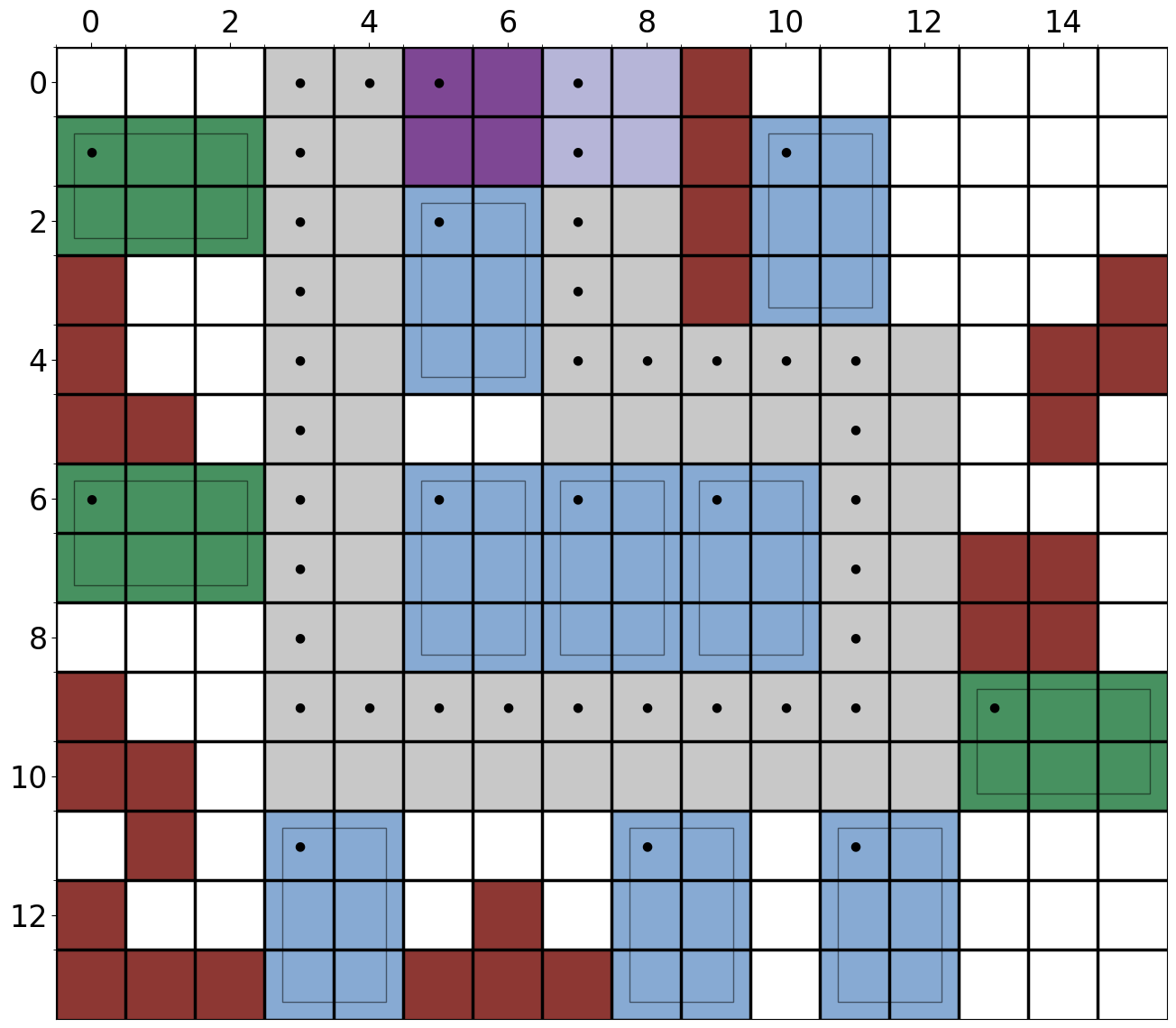}
\includegraphics[width=0.49\textwidth]{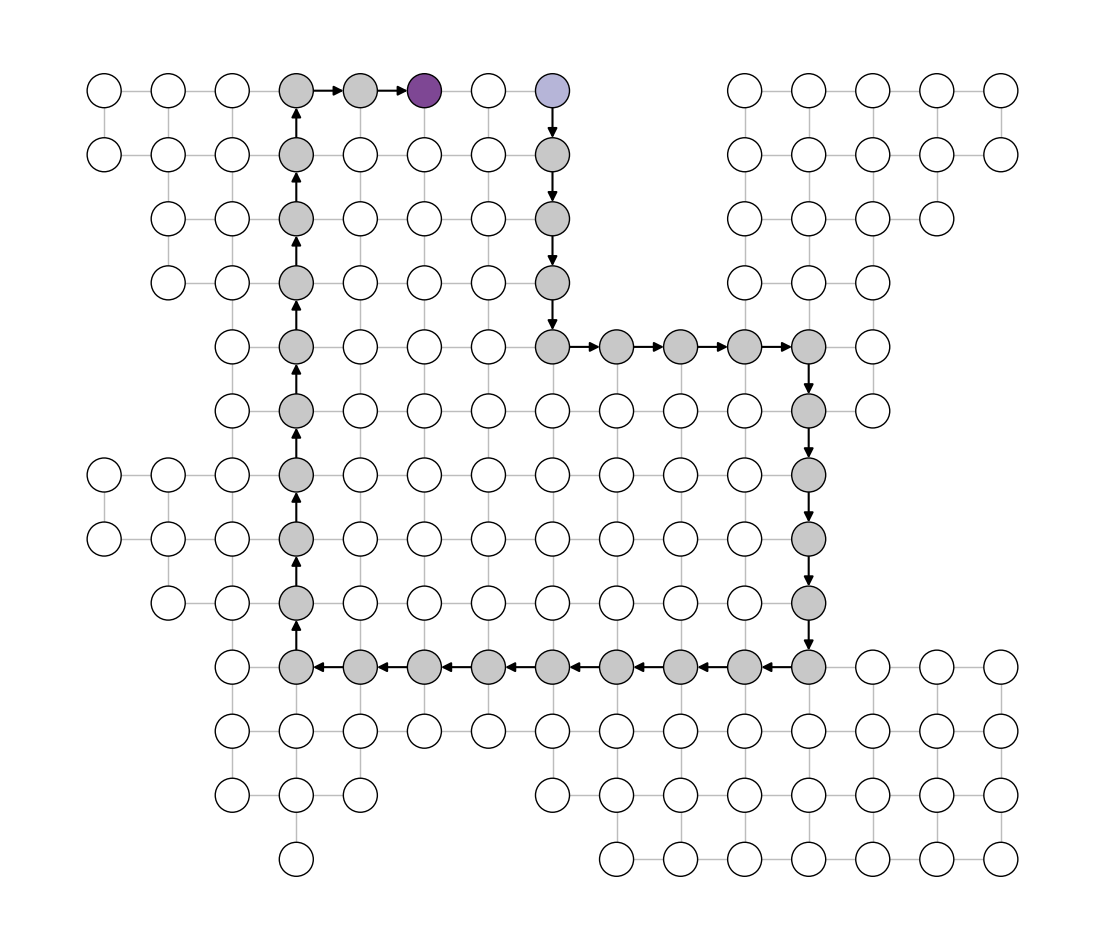}
\caption{Lot and directions with minimum length constraints}
\label{fig:one_way_turn_constraints}
\end{subfigure} 
   \caption{Results from adding turn restriction constraints}
\label{fig:turn_results_one_way}
\end{figure}

Figure \ref{fig:overlap_constraints} shows an optimal solution with these constraints. Observe that when two driving fields are next to each other, the traffic moves in the same direction. In other words, it is possible to merge these fields into larger driveways that point in the same direction. E.g., see the fields anchored at (9, 5), (9, 6), (9, 7), (10, 5), (10, 6), and (10, 7). A model that adheres to these constraints can also approximate curved driveways using a step-like sequence of driving fields, as mentioned earlier in the two-way case. However, when the size of the driving field is smaller than the vehicle length, we might want the vehicles to fit fully within two adjacent junctions. In such cases, additional restrictions are required on combinations of $z$ variables, which are illustrated in Figure \ref{fig:one_way_avoid_sharp_turns}. In these diagrams, we examine eight possible turn movements at an active driving cell. Consider the example in the upper left cell (3, 3) of the figure's left panel. Vehicles driving east turn left to head north. If $\parklength = 3$, then to avoid sharp turns, we do not allow the $z$ variables starting at (1, 3) and (2, 3) and going in the east-west direction to be active. No restrictions are placed on the $z$ variables starting at (0, 3) since the longitudinal segment from (3, 3) to (0, 3) spans three cells, which is sufficient to hold a vehicle. Generalizing this, if $z_{(i,j-1), (i,j)}$ and $z_{(i,j), (i-1,j)}$ are ones, then $z_{(i-k, j), (i-k, j-1)}$ and $z_{(i-k, j), (i-k, j +1)}$ for $k \in \{1, \ldots, \parklength-1\}$ must be zero. 

In the earlier example, this can easily be set using four constraints, one of which is $z_{(3, 2), (3, 3)} + z_{(3, 3), (2, 3)} + z_{(2, 3), (2, 4)} \leq 2$. Note that some of these inequalities may be redundant if constraints that avoid overlapping of opposite-moving traffic are included. For larger vehicles, turning might still present challenges; in such scenarios, we can modify these restrictions to either extend the distance between intersections or increase the sizes of the cells. Figure \ref{fig:one_way_turn_constraints} shows a solution satisfying these inequalities. It can be seen that driveways are more organized and facilitate smoother turns, although this results in fewer parking stalls. In both one-way and two-way cases, additional constraints may be needed at the entrance or exit, depending on the configuration of the access roads.

\begin{figure}[H]
  \centering
  \begin{subfigure}[b]{0.49\textwidth}
      \centering
    \includegraphics[width=0.49\textwidth]{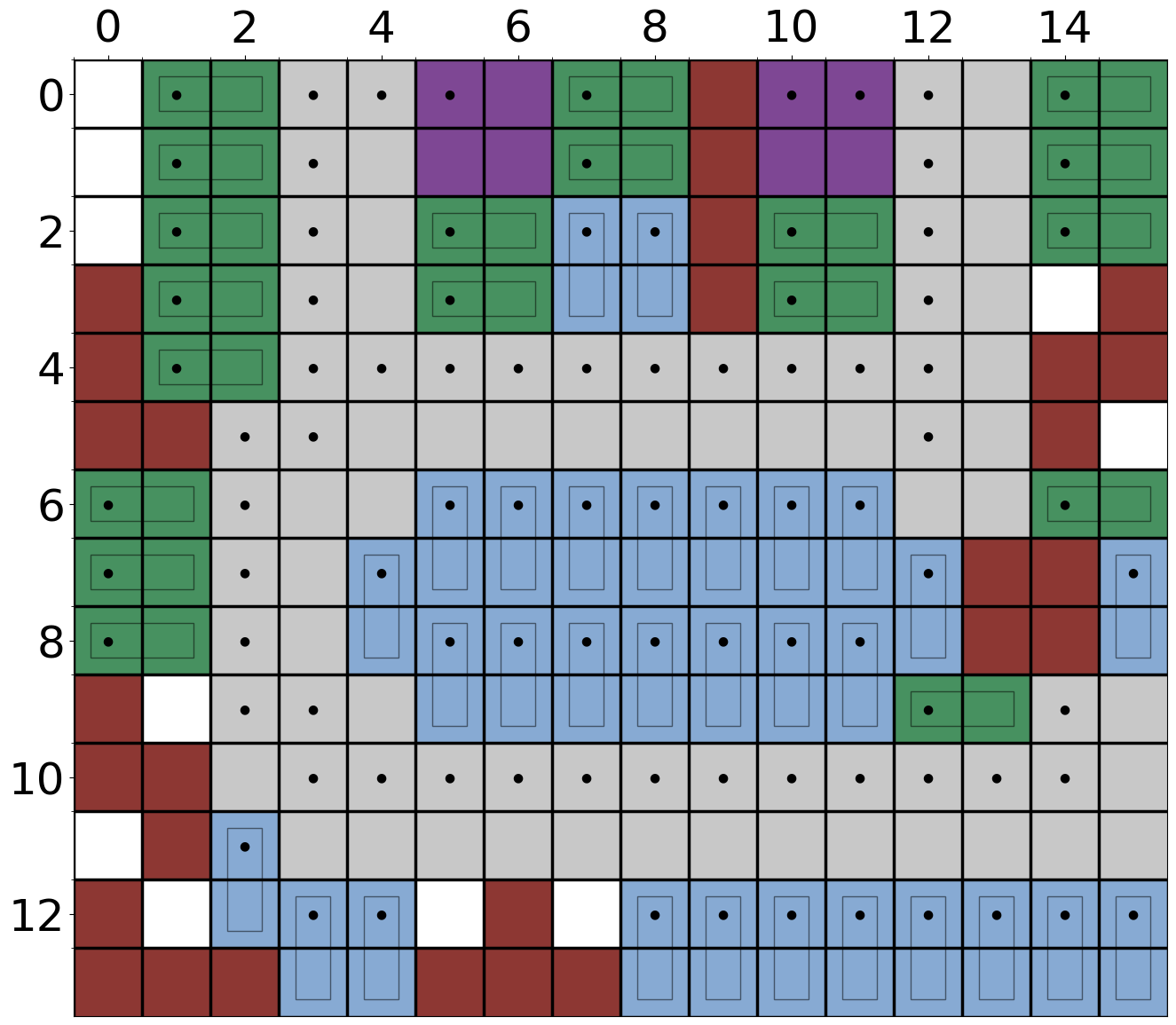}
    \includegraphics[width=0.49\textwidth]{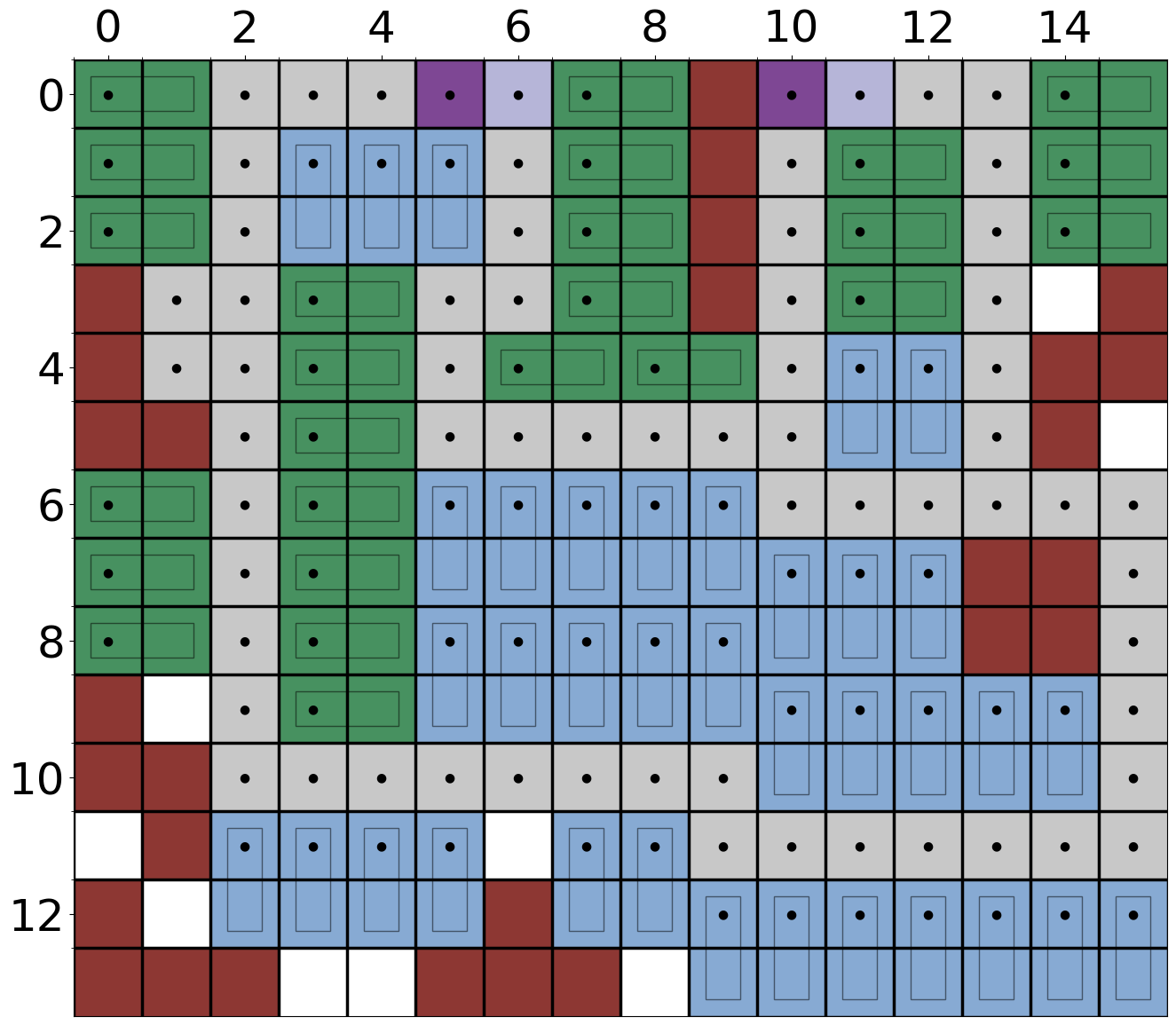}
    \caption{Connected entrances/exits}
    \label{fig:multiple_connected}
  \end{subfigure}
  \hfill
  \begin{subfigure}[b]{0.49\textwidth}
      \centering
    \includegraphics[width=0.49\textwidth]{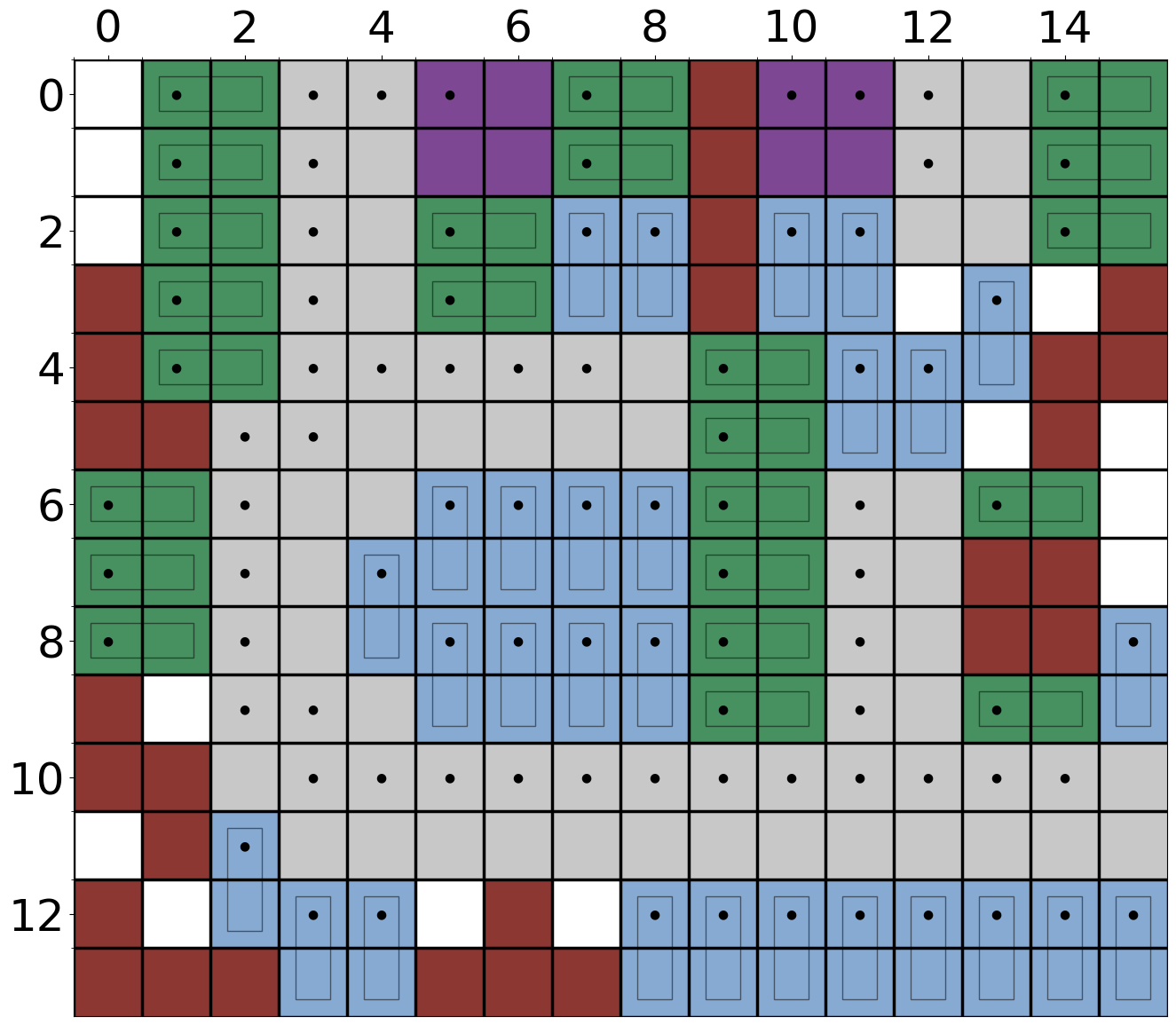}
    \includegraphics[width=0.49\textwidth]{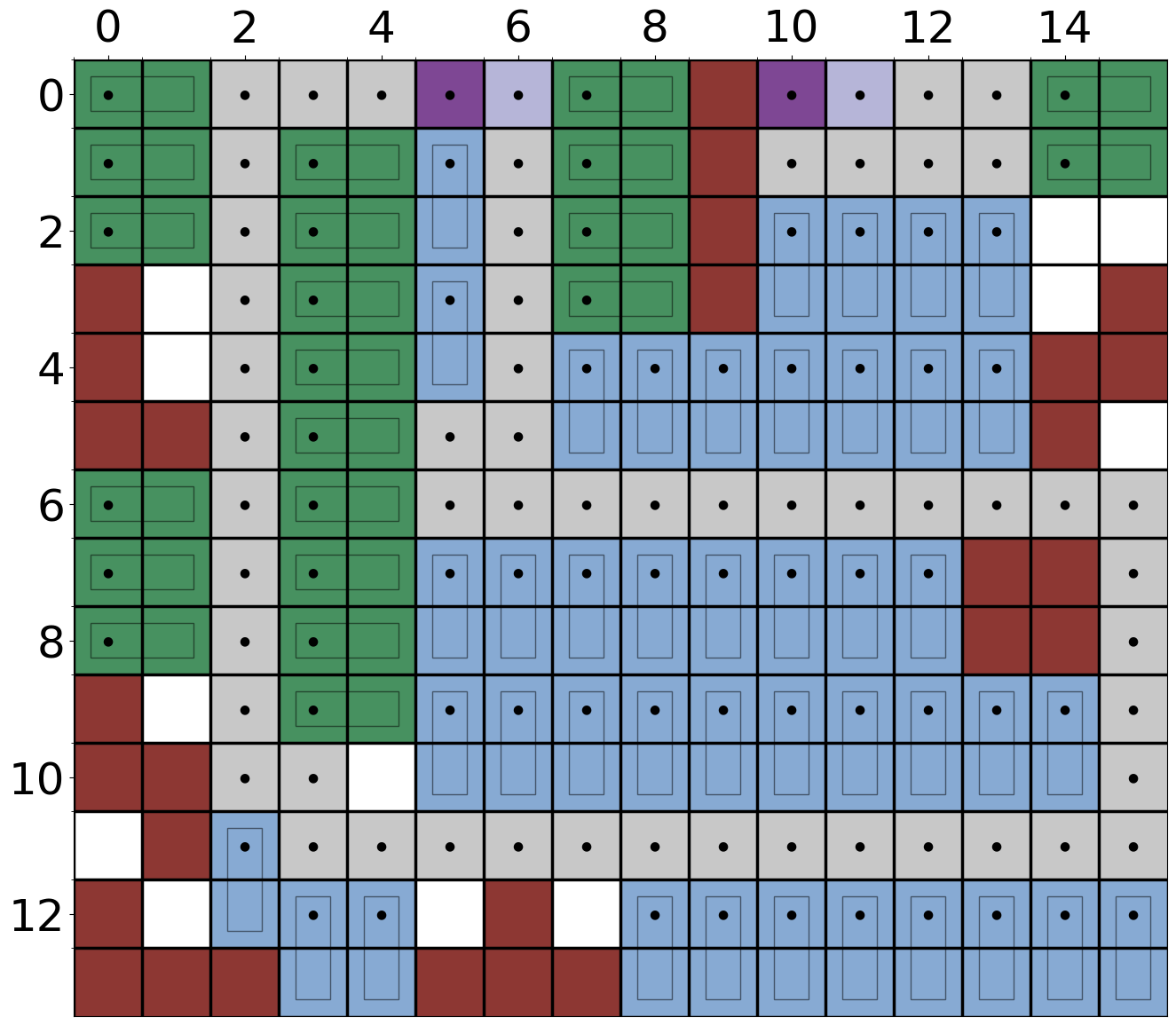}
    \caption{Potentially disconnected entrances/exits}
        \label{fig:multiple_disconnected}
  \end{subfigure}
\caption{Layouts with multiple entrances and exits}
\label{fig:multiple}
\end{figure}

\textbf{Multiple entrances and exits:}
Parking lots often have multiple entrances and exits. We can design such lots in two ways. (1) All driveways can be connected, ensuring that entrances and exits are accessible from each other. (2) We can relax this requirement and potentially create multiple disjoint lots. In the first case, we make the driving field variables active for all entrance cells and solve the problem as a single entrance/exit model with any entrance cell in the grid as the entry point. An optimal solution for the two-way and one-way case with connected entrance cells is shown in Figure \ref{fig:multiple_connected}. The number of parking stalls may be reduced compared to a single entrance/exit case due to the space allocated for entrance locations.

In the second approach, we can solve the problem by assuming that not all entrances/exits are necessarily connected, but the driveways provide access to each parking space from at least one of the entrances/exits. This type of design may be suitable for large parking lots with multiple access roads on different sides. Solving this variant again is easy and can be done by simply connecting the entrances (and the exits) in the grid graph of valid driving fields with dummy arcs. The connectivity cut constraints would then require each active driving cell to be connected to at least one entrance/exit. Figure \ref{fig:multiple_disconnected} shows an optimal solution using this technique. The driveways originating from different entrances, in this case, are disjoint. The number of parking stalls also increased compared to the connected layouts by two for both the two-way and one-way cases. Disconnected parking entrances can make the overall space behave as separate lots and increase capacity. However, the search for parking could be prolonged since if one section of the lot is full, vehicles must exit, navigate the public streets, and reenter a different area via an alternate entrance.

\textbf{Designing obstacle locations:} The models presented in this paper can optimize the number of parking stalls in cases where some cells in the parking lot's interior are blocked. Since these blocked regions often host necessary parking-related elements, locating them suitably can also be a key design objective. In this subsection, we provide two examples in which the MIP models can be extended to determine where to place obstacles while maximizing the lot capacity. The results in this section present the best lower bounds obtained by the solver within a specified time limit.

In the first example, we analyze the spatial allocation of shopping cart corrals. These are relevant in parking lots serving shopping areas where customers return empty carts to dedicated spaces after loading their purchases. Having more corrals benefits customers by reducing walking distances. However, it also reduces the number of parking stalls, which is undesirable for the lot operator or the stores. From a functional perspective, cart corrals should be evenly distributed across the lot and accessible from the driving fields. Cart corrals are typically the same size as parking fields. Thus, planners can use our MIP models to address this problem by simply fixing the locations of the corrals to those of parking fields and evaluating lot capacity for different configuration scenarios. 

\begin{figure}[H]
  \centering
  \begin{subfigure}[b]{0.24\textwidth}
      \centering
    \includegraphics[width=\textwidth]{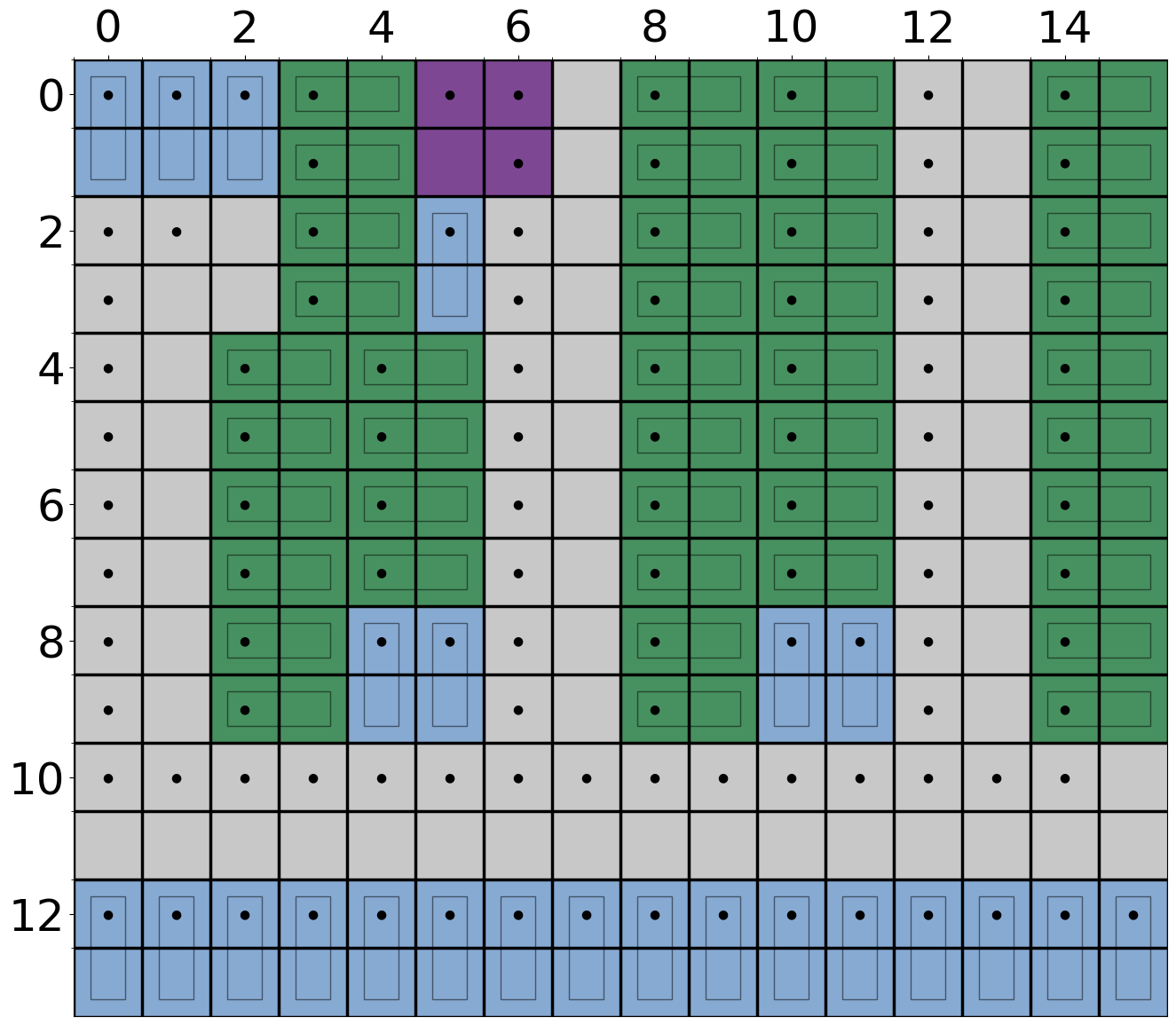}
    \caption{No obstacles (66 stalls)}
    \label{fig:no_obstacles}
  \end{subfigure}
  \hfill
  \begin{subfigure}[b]{0.24\textwidth}
      \centering
    \includegraphics[width=\textwidth]{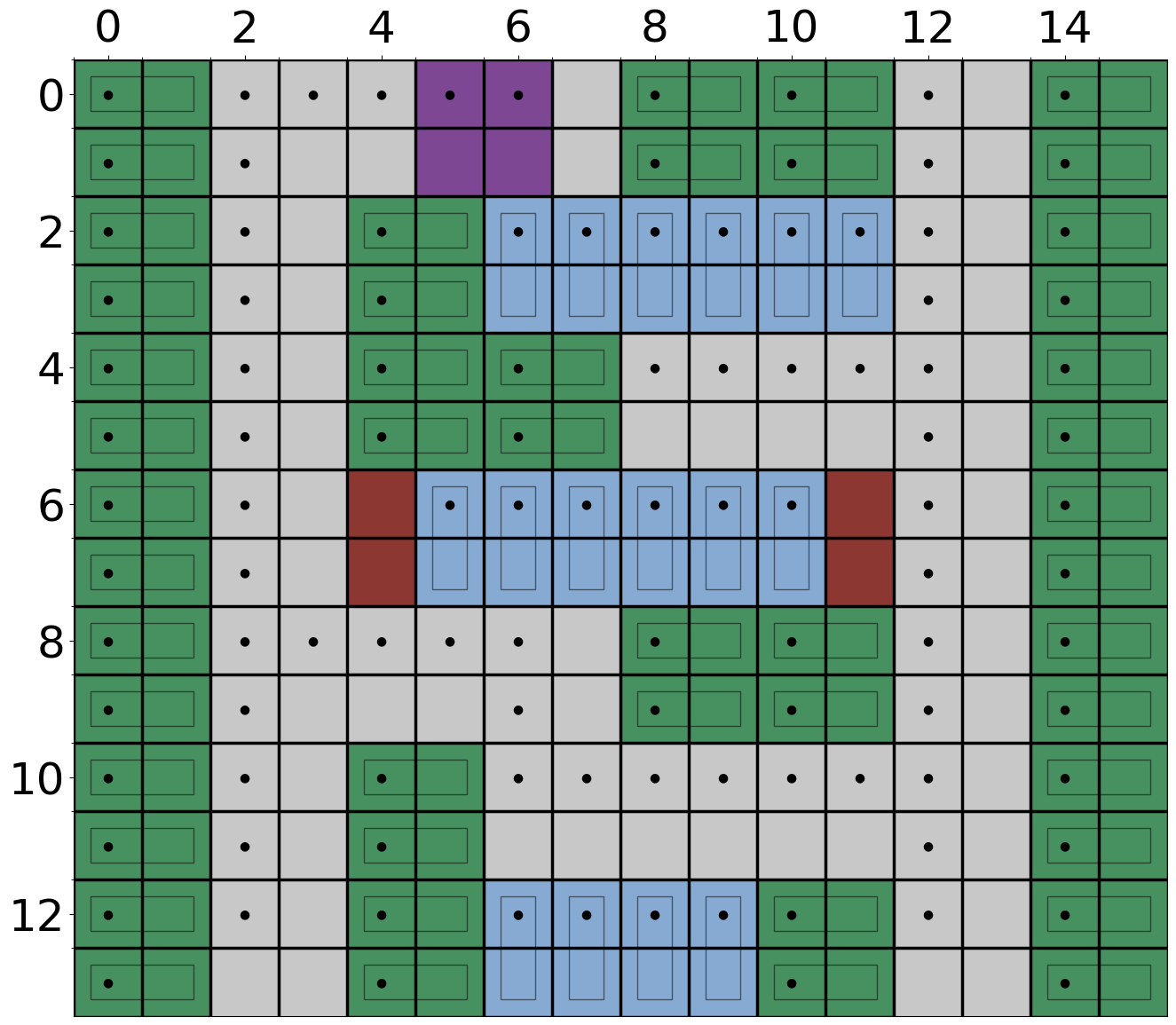}
    \caption{Two obstacles (64 stalls)}
    \label{fig:two_obstacles}
  \end{subfigure}
  \hfill
  \begin{subfigure}[b]{0.24\textwidth}
      \centering
    \includegraphics[width=\textwidth]{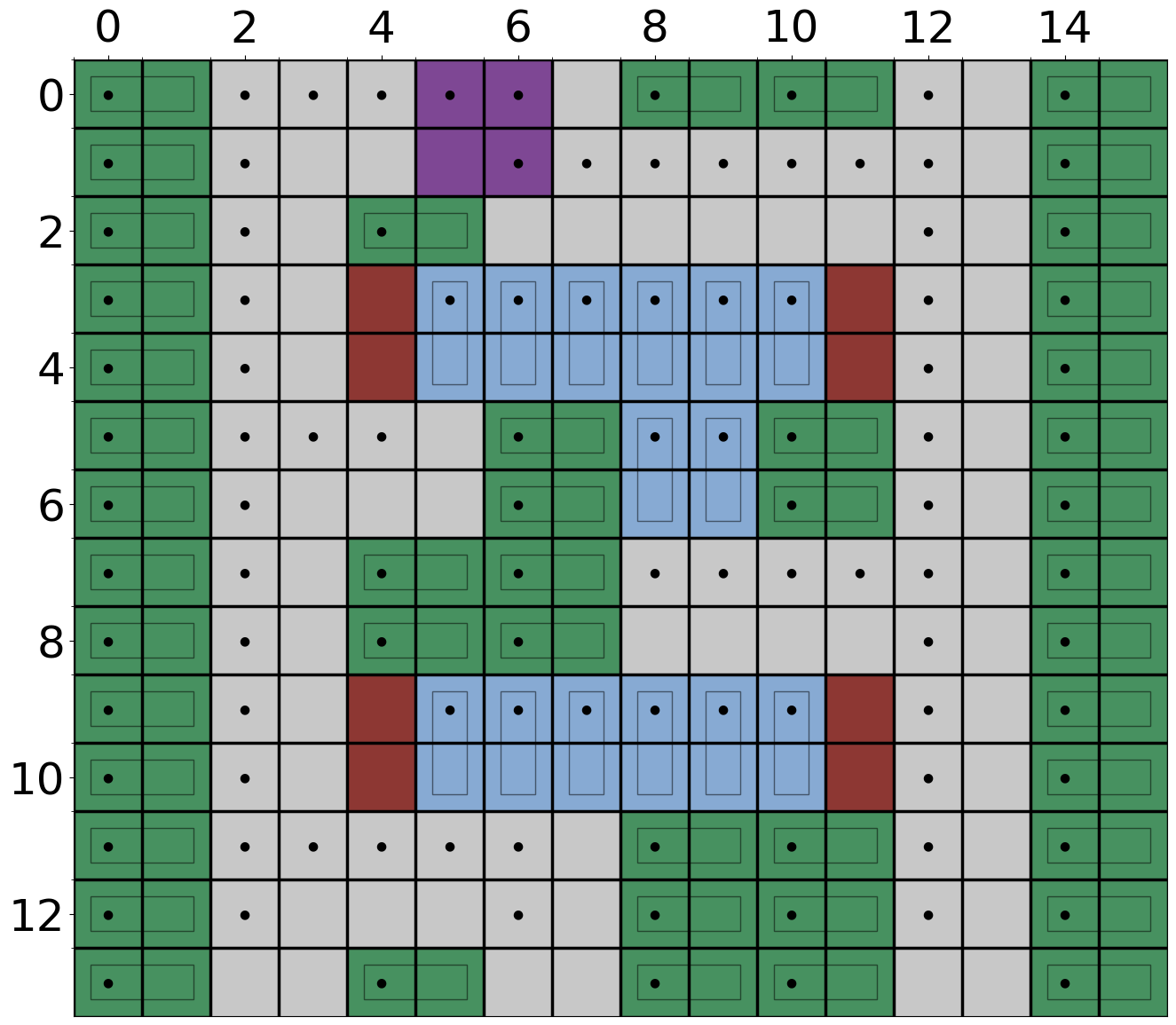}
    \caption{Four obstacles (60 stalls)}
    \label{fig:four_obstacles}
  \end{subfigure}
  \hfill
  \begin{subfigure}[b]{0.24\textwidth}
      \centering
    \includegraphics[width=\textwidth]{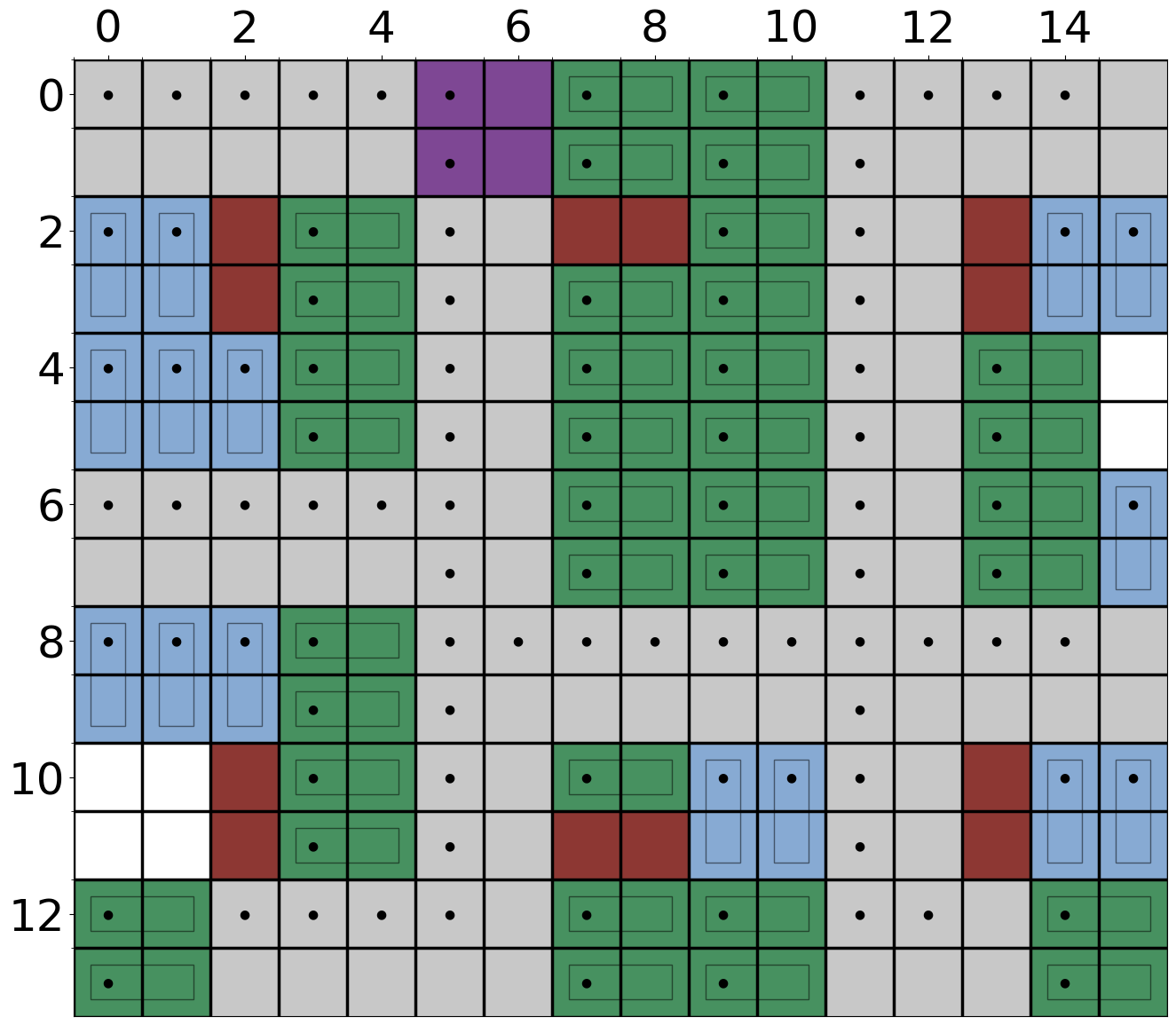}
    \caption{Six obstacles (51 stalls)}
    \label{fig:six_obstacles}
  \end{subfigure}  
  \caption{Parking lot designs with evenly distributed shopping cart corrals}
  \label{fig:corrals}
\end{figure}

For example, Figure \ref{fig:corrals} illustrates the results from a scenario analysis of a parking lot with no obstacles and a few evenly distributed ones. Figure \ref{fig:no_obstacles} shows a parking lot with 66 stalls when there are no obstacles. In Figure \ref{fig:two_obstacles}, we add constraints that force $90^{\circ}$ parking fields at (6, 4) and (6, 11), converting them into obstacles instead of using them for parking. In this case, we only lose the space assigned to the corrals, so we can accommodate 64 stalls. However, when we add four and six corrals in the configurations shown in Figures \ref{fig:four_obstacles} and \ref{fig:six_obstacles}, the number of parking spots reduces by 6 and 17, respectively. One could also select a fixed number of parking spots from the no-obstacle case in Figure \ref{fig:no_obstacles} and convert them into corrals. While this approach directly controls the reduction in lot capacity, the resulting corral locations may not be uniformly distributed.

In a second use case, we extend our model to locate planting islands inside a parking lot. Many cities around the world are facing severe urban heat island effects, and to counter them, regulations are being proposed that will influence the design of urban infrastructure. For instance, in NYC, a recent urban forest plan aims to plant trees to increase the tree canopy cover to 30\% by 2040. At the parking lot level, new rules require one tree for every eight parking spots, with at least 30\% of the trees located inside the lot for larger facilities \citep{nyc_urban_forest_2026}. These regulations, however, do not specify exactly where planting islands must be located.

We can optimize the locations of planting islands and the number of parking spots by introducing new variables and constraints. Consider the R2 resolution in which each cell is $3~\text{m}\times3~\text{m}$. Suppose $b_{ij}$ represents an obstacle location decision variable that is 1 if cell $(i, j)$ is blocked for a planting island. Let $a_{ij}$ be another decision variable that is 1 if cell $(i, j)$ is covered by the tree canopy when viewed from above. The crown sizes of mature trees vary by species. Suppose that, on average, trees planted in the lot have a crown diameter of $2\tau + 1$ cells. Thus, a cell $(i, j)$ is covered by the canopy only if a tree is planted at one of the cells in $T_{ij} = \{(i + k, j + l) \colon k, l \in \{-\tau, \ldots, \tau\} \} \cap (L \setminus B)$. Observe that $T_{ij}$ also represents the set of cells that are under the crown of a tree planted at $(i, j)$. Hence, we can add constraints \eqref{eq:canopy1} and \eqref{eq:canopy2} to connect the obstacle and canopy variables. Both sets of variables can be continuous as long as they are bounded in $[0, 1]$.
\begin{align}
    a_{ij} & \geq b_{kl} &&  \qquad \forall~(i, j) \in L \setminus B, (k, l) \in T_{ij} \label{eq:canopy1} \\
    a_{ij} & \leq \sum_{(k, l) \in T_{ij}} b_{kl} &&  \qquad \forall~(i, j) \in L \setminus B \label{eq:canopy2}
\end{align}
We also modify the single-purpose constraint (3) by adding the $b$ variables to the left-hand side. This allows a cell to be used for driving, parking, or obstacle placement, or it can simply be left empty.
\begin{align}
& b_{ij} + \drivevar{m}{n} + \sum_{(k,l) \in \invzeroparkset{i}{j}} \zeroparkvar{k}{l} + \sum_{(k,l) \in \invninetyparkset{i}{j}} \ninetyparkvar{k}{l} \leq 1  \, \, \, \, \, \, \, \, \, \forall\,  (m, n) \in  \invdriveset{i}{j}, (i,j) \in \lotset \setminus \blockedset \label{eq:singlepurpose_zero_disaggregate_obstacle} 
\end{align}
Lastly, we impose a constraint that the total canopy cover must be at least a predetermined regulatory fraction $\gamma$ of the parking lot's area. 
\begin{align}
    \sum_{(i, j) \in L \setminus B} a_{ij} \geq \gamma |L \setminus B| \label{eq:regulatory}
\end{align}
An example of this formulation is shown in Figure~\ref{fig:imposing_restrictions}. Two thresholds $\gamma = 0.3$ and $\gamma = 0.5$ are considered, and each tree is assumed to cover a square region comprising 9 cells (i.e., $\tau = 1$). The red cells indicate the locations of the planting islands, and the brown regions around them represent canopy cover. For $\gamma = 0.3$, the number of stalls decreases from 66 to 62; for $\gamma=0.5$, it reduces further to 57. The introduction of these new variables expands the problem's search space, making it more computationally challenging. These models also admit multiple configurations for the same number of parking stalls. Note that for lower values of $\gamma$, the tree cover may be absent in some portions of the lot. More uniform designs can be obtained by imposing constraints similar to \eqref{eq:regulatory} using specific subsets of $L \setminus B$ that lack coverage. These can also be dynamically added using solver callbacks. 
\begin{figure}[H]
    \centering
    \begin{subfigure}[b]{0.49\textwidth}
        \centering
        \includegraphics[width=0.49\textwidth]{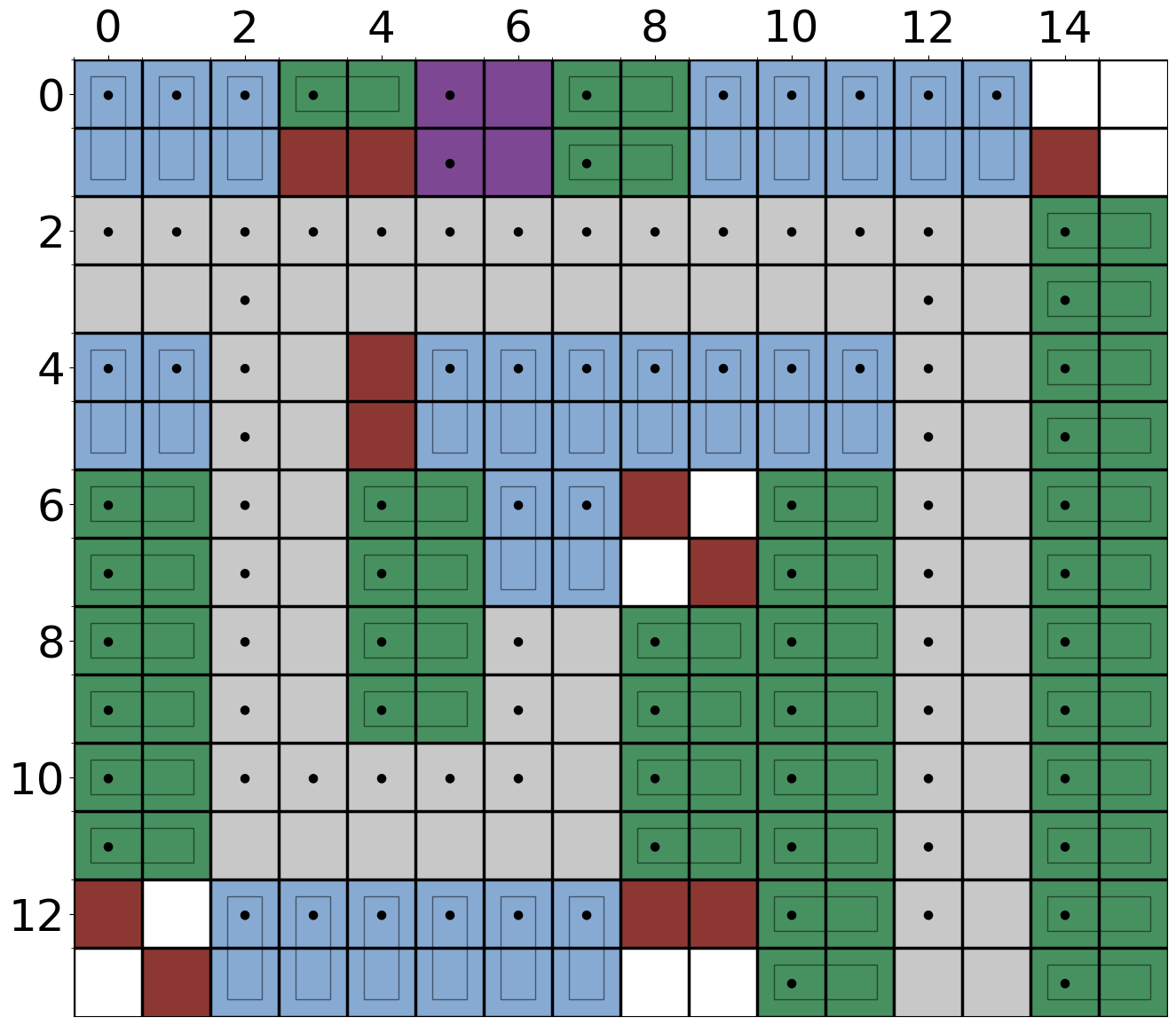}
    \includegraphics[width=0.49\textwidth]{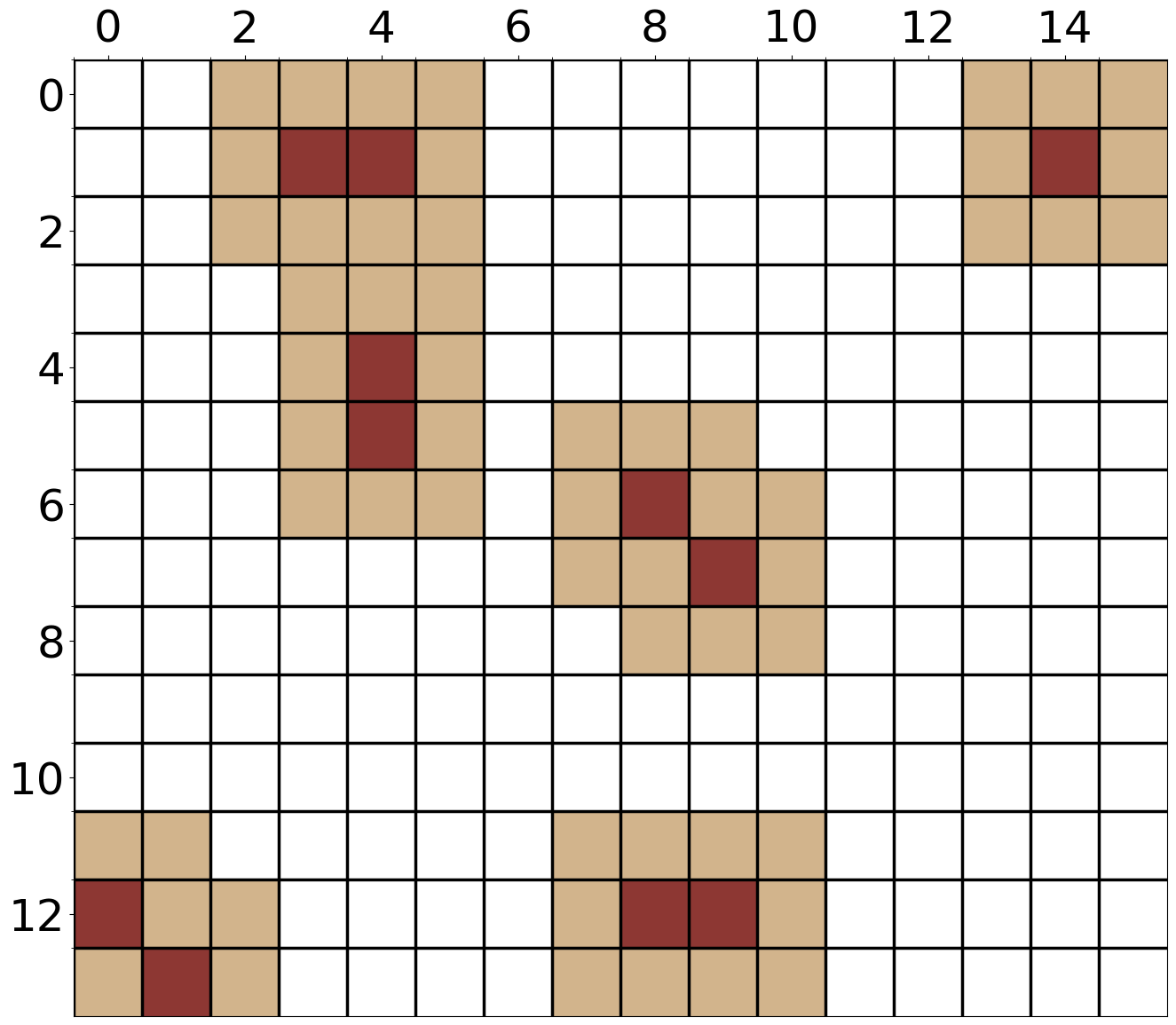}
        \caption{30\% canopy cover (62 stalls)}
        \label{fig:sharp_turn_two_way}   
    \end{subfigure}
    \begin{subfigure}[b]{0.49\textwidth}
        \centering
    \includegraphics[width=0.49\textwidth]{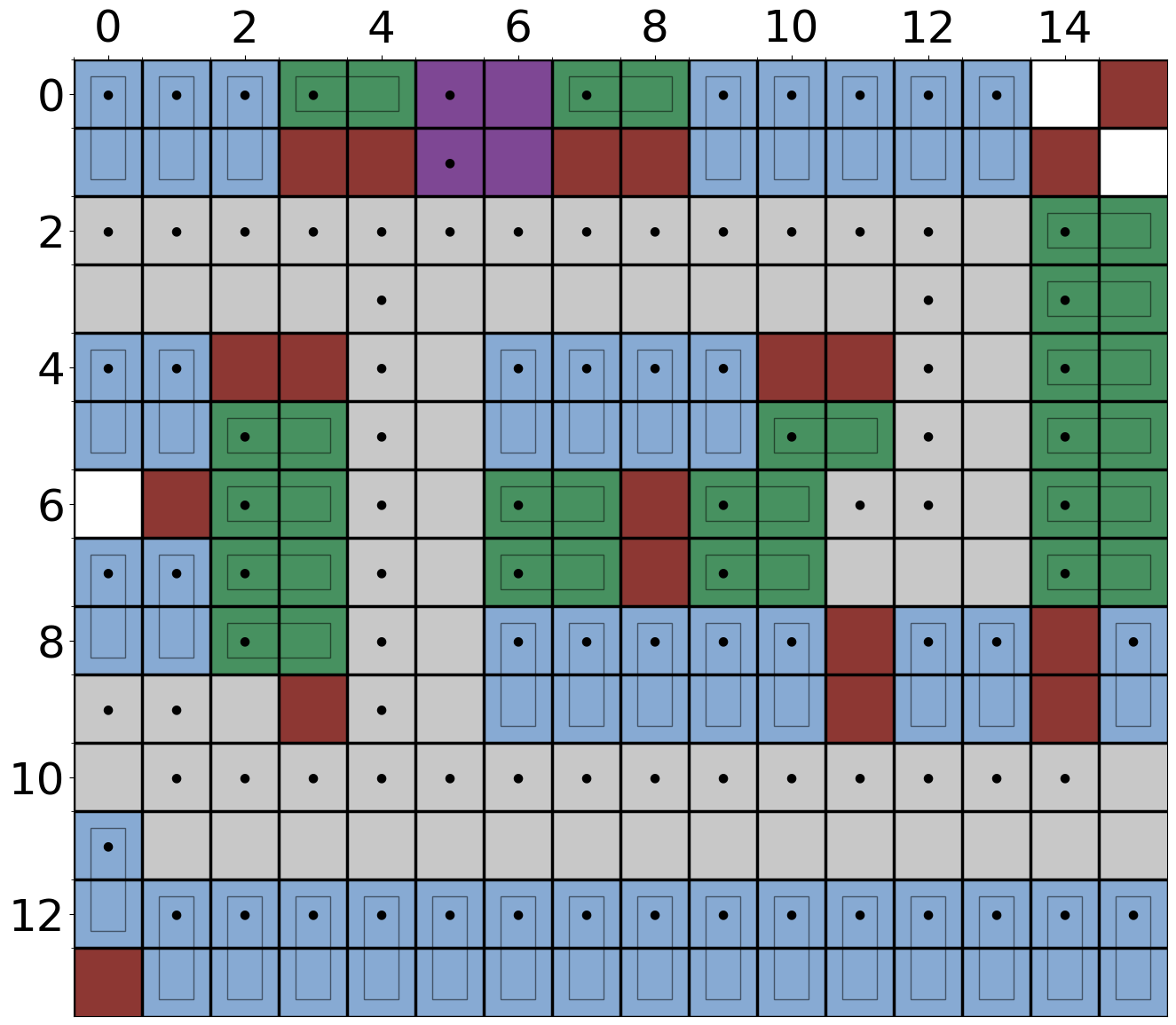}
    \includegraphics[width=0.49\textwidth]{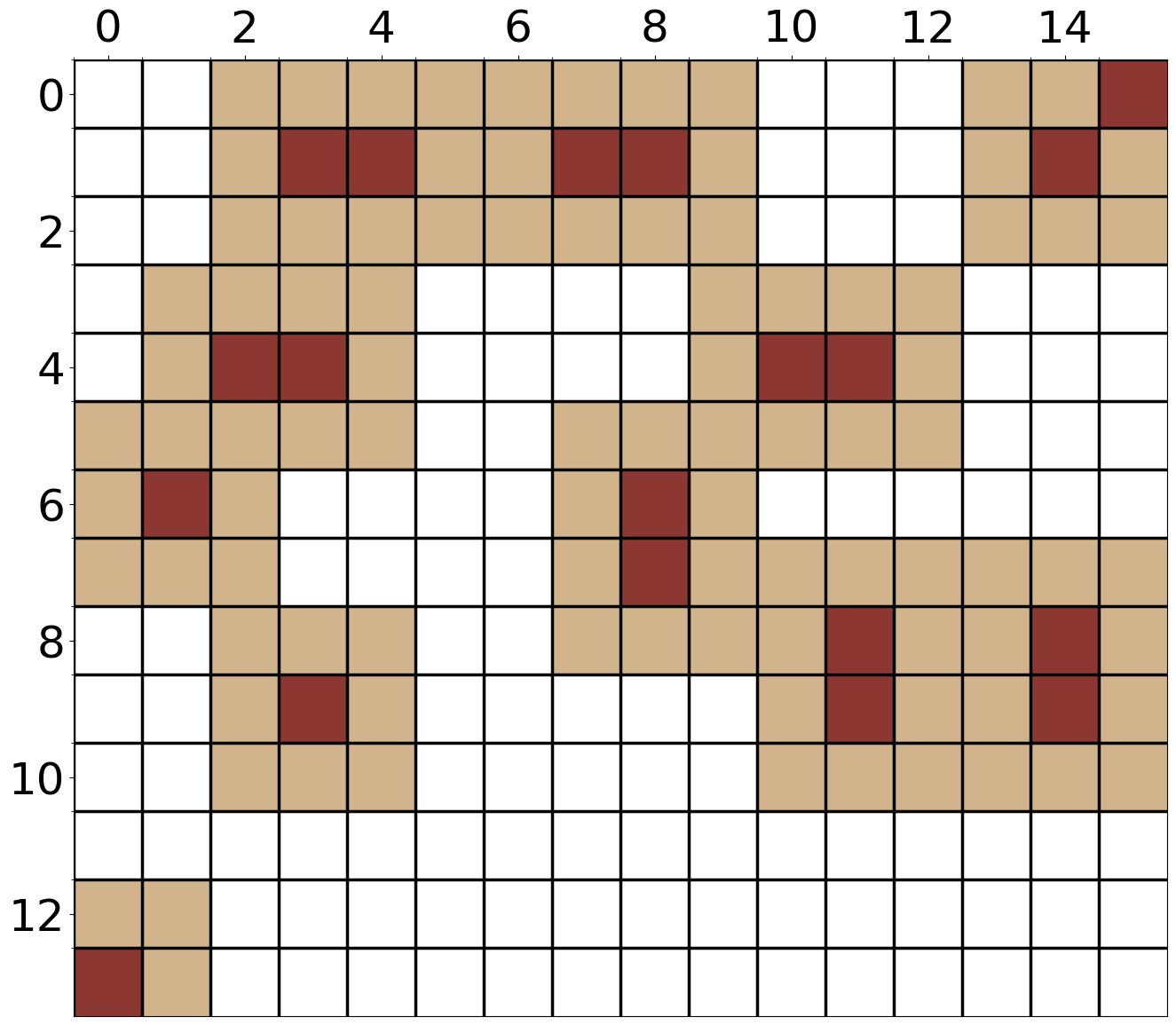}
        \caption{50\% canopy cover (57 stalls)}
        \label{fig:turn_restrictions}     
    \end{subfigure}
    \caption{Parking lot designs with planting islands. Red cells represent obstacle locations containing the tree planters, and brown (and red) cells indicate those covered by the tree canopy.}
    \label{fig:imposing_restrictions} 
\end{figure}

\section{Additional results}
\label{sec:add_stats}
This section presents additional results and instance features for the experiments conducted on the one-way and two-way versions of the problem. Recall that for each driveway configuration, three formulations were analyzed. The first baseline model solves the flow-based formulation. The second uses flow variables but includes valid inequalities, and the third optimizes lots using the CCC formulation and the branch-and-cut algorithm. 

\begin{figure}[H]
  \centering
  \begin{subfigure}[b]{0.95\textwidth} 
    \includegraphics[width=0.47\textwidth]{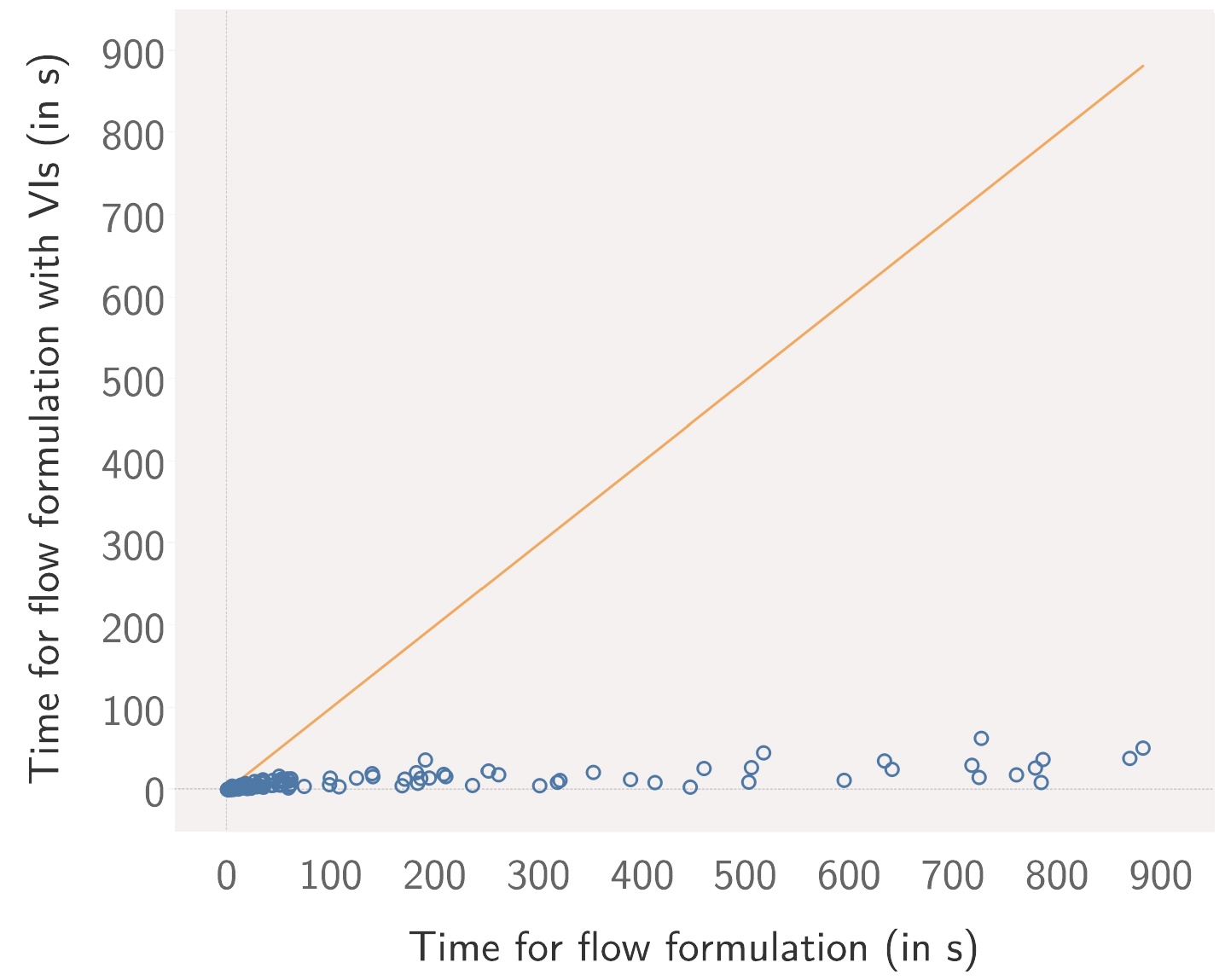}
    \hfill \hspace*{11mm}
    \includegraphics[width=0.47\textwidth]{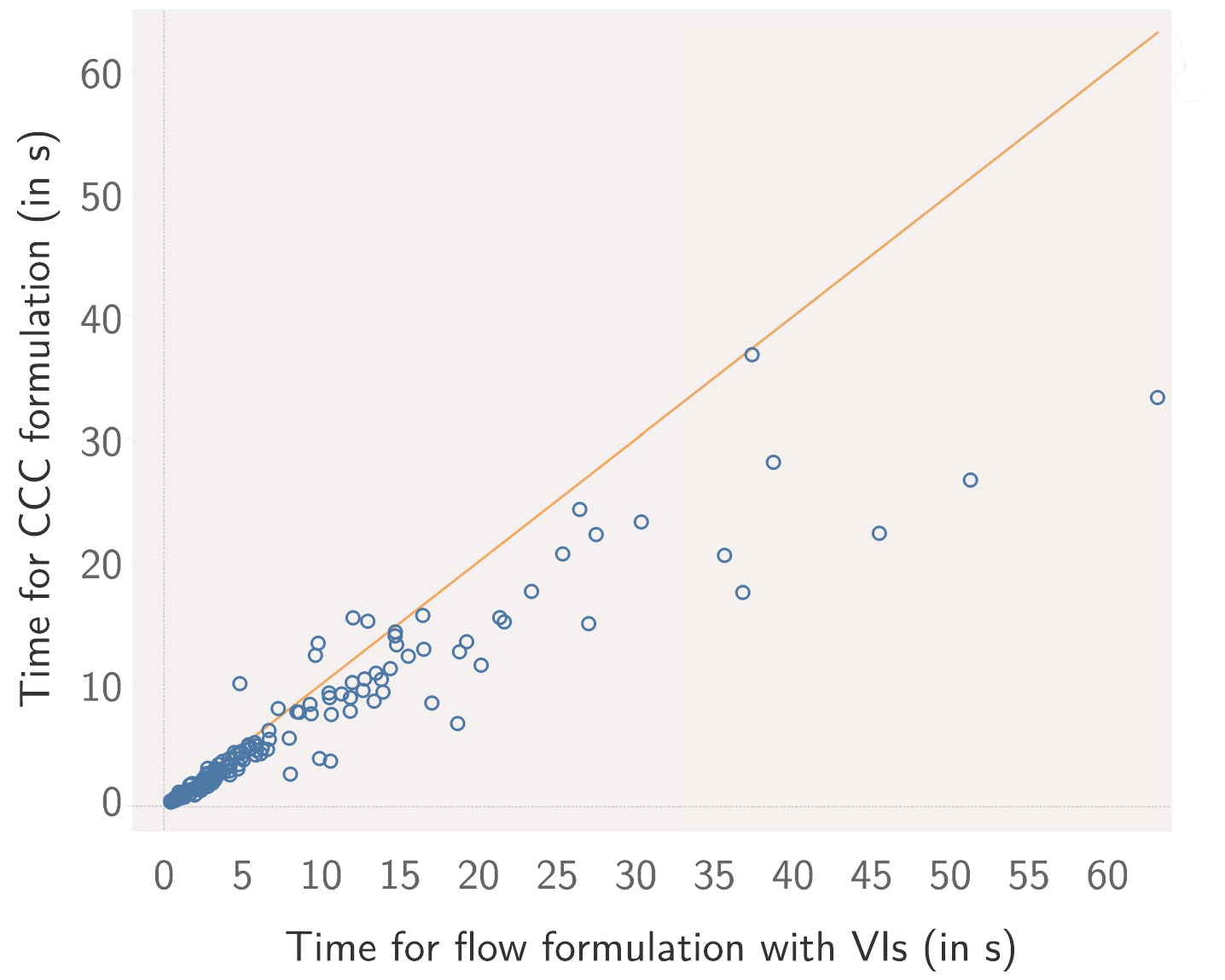}
    \caption{Runtimes for Group I instances}
    \label{fig:two_way_times}
  \end{subfigure} \\ \vspace*{5mm}
  \begin{subfigure}[b]{0.95\textwidth}
    \includegraphics[width=0.47\textwidth]{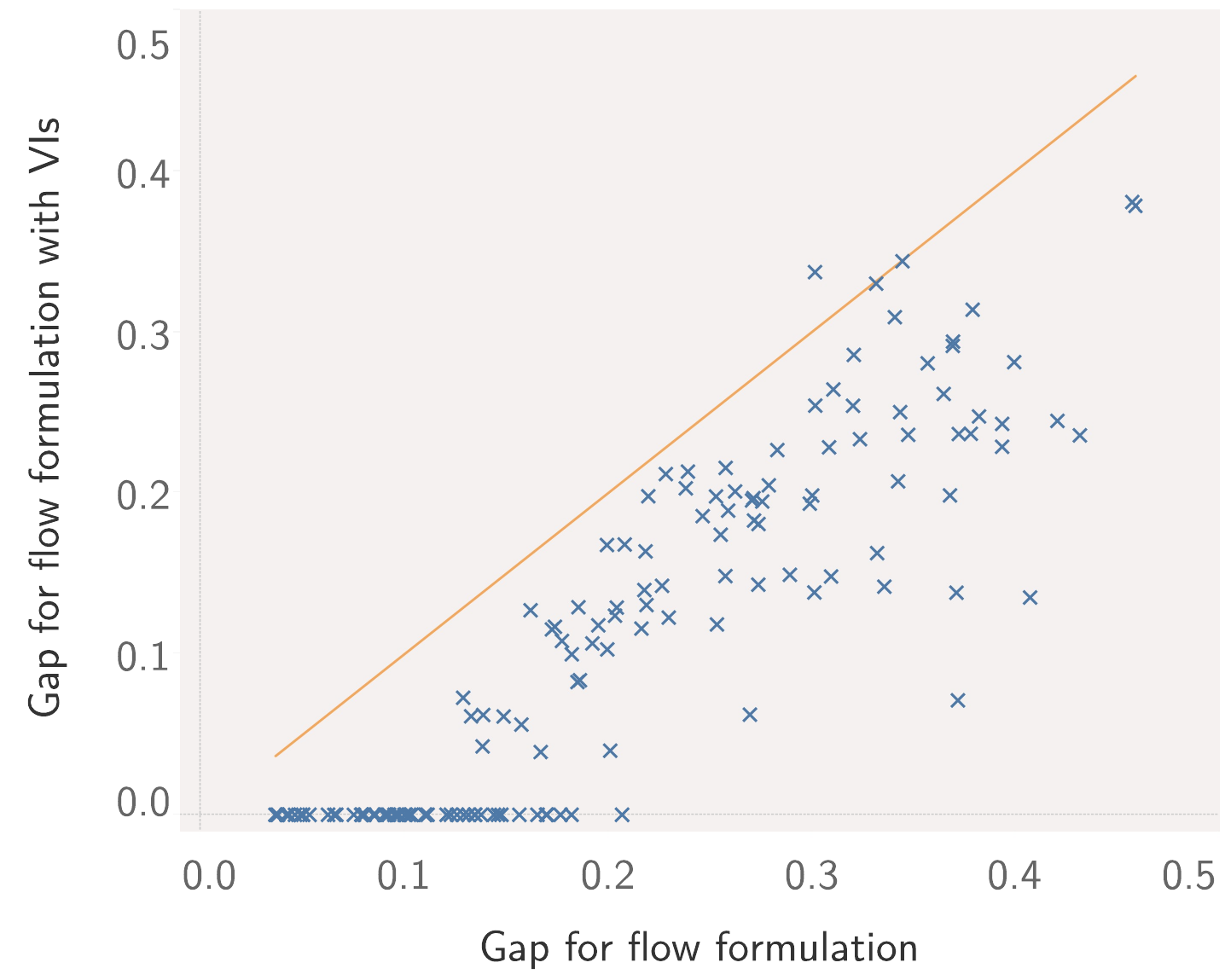}
    \hfill 
    \includegraphics[width=0.47\textwidth]{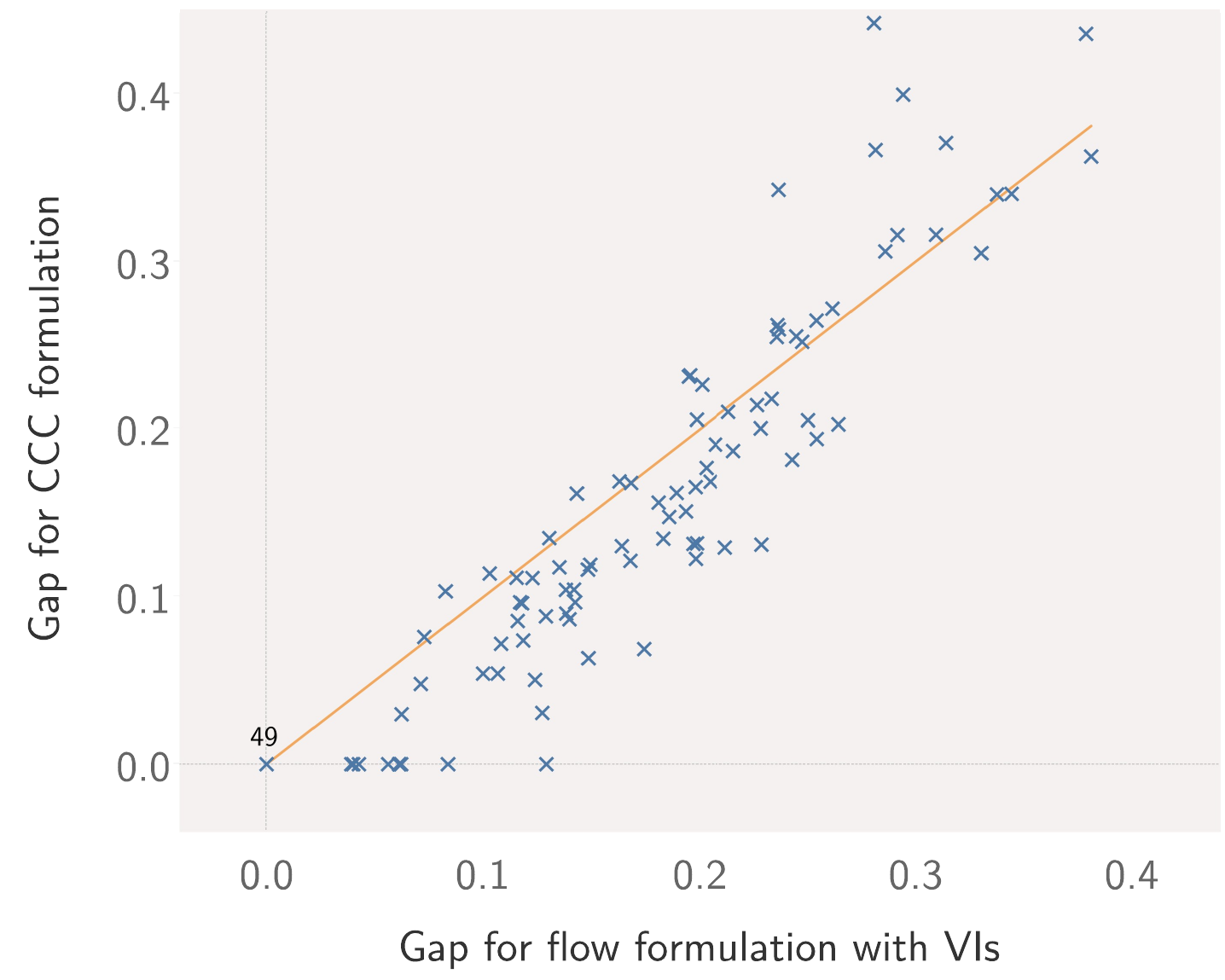}
    \caption{Optimality gaps for Group II instances}
    \label{fig:two_way_gaps}
  \end{subfigure}
  \caption{Computational results for the two-way experiments}
  \label{fig:two_way_computational}
\end{figure}

\textbf{Runtimes and gaps for experiments with a 15-minute limit:}
A scatter plot of the runtimes between the different formulations for the 325 instances is shown in Figure \ref{fig:two_way_times}. Adding hop and bidirectional hop inequalities to the flow formulation significantly improved runtimes, as shown in the left panel. The CCC approach was also faster than the flow formulation. However, as seen in the right panel, it did not offer significant advantages over the flow formulation with valid inequalities. The trends in the gaps for instances that did not find an optimal solution within the 15-minute time limit are also similar and can be seen in Figure \ref{fig:two_way_gaps}.


\begin{figure}[H]
  \centering
  \begin{subfigure}[b]{0.95\textwidth} 
    \centering
    \includegraphics[width=0.47\textwidth]{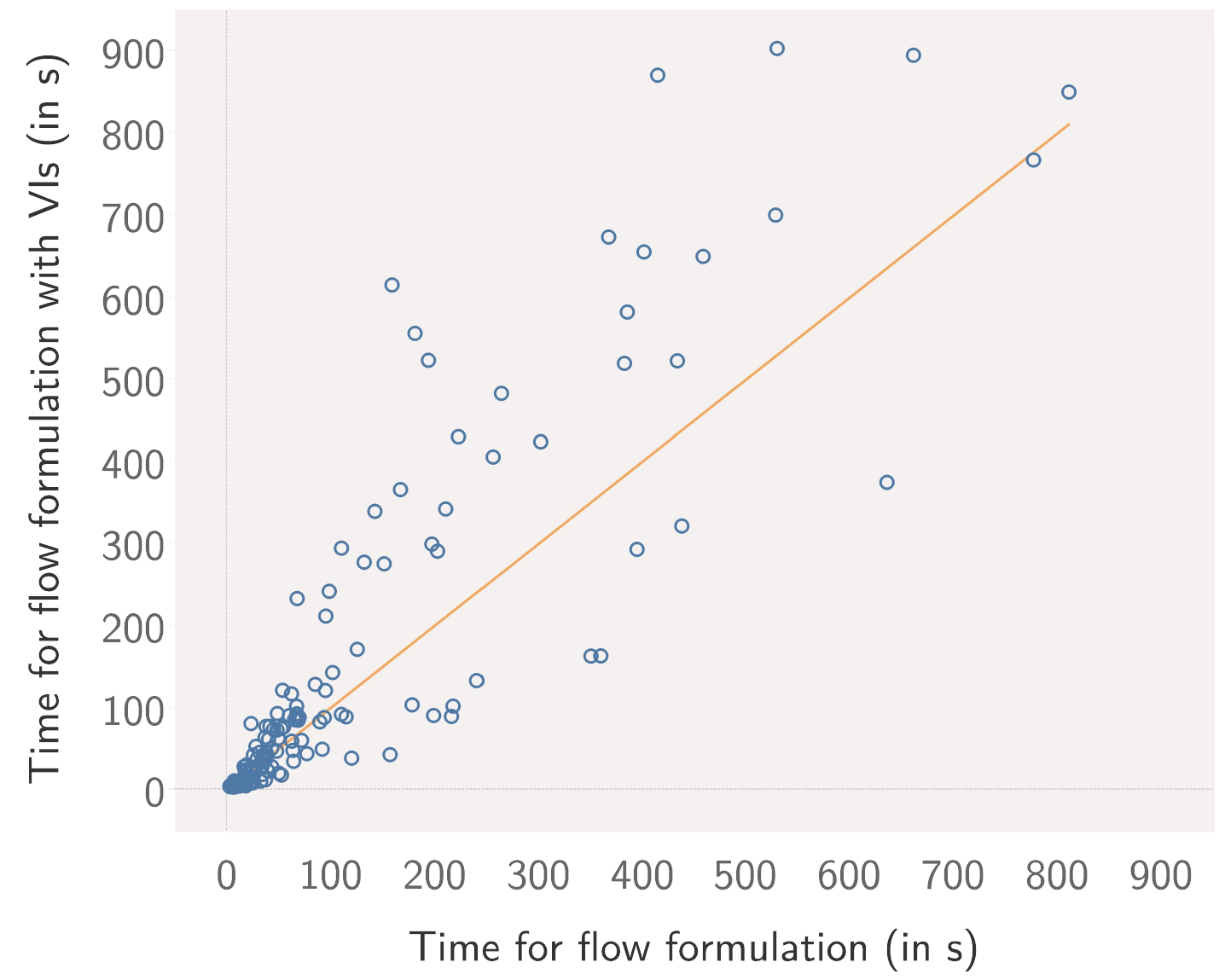}
    \hfill 
    \includegraphics[width=0.47\textwidth]{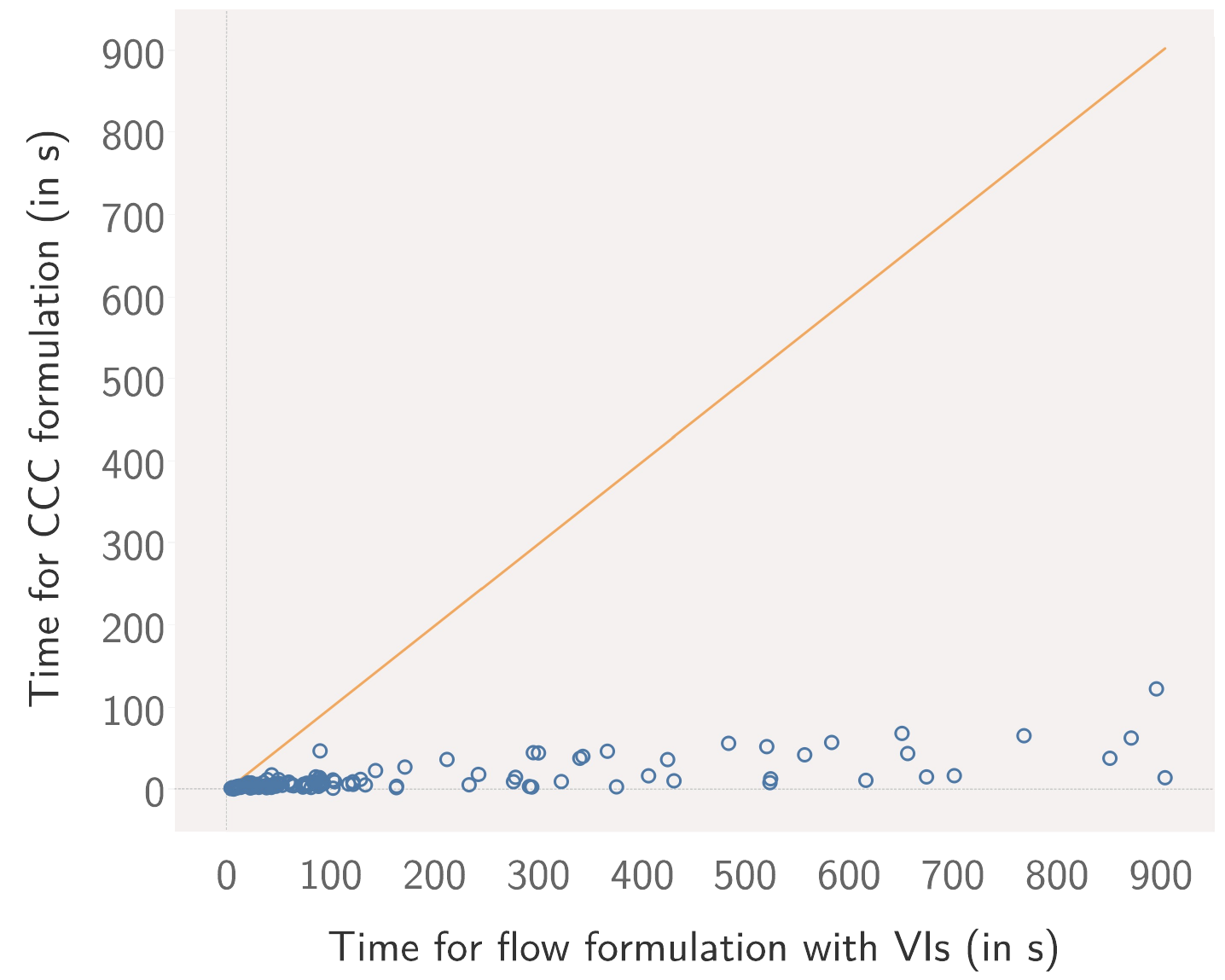}
    \caption{Runtimes for Group I instances}
    \label{fig:one_way_times}
  \end{subfigure} \\ \vspace*{5mm}
  \begin{subfigure}[b]{0.95\textwidth}
    \centering
    \includegraphics[width=0.47\textwidth]{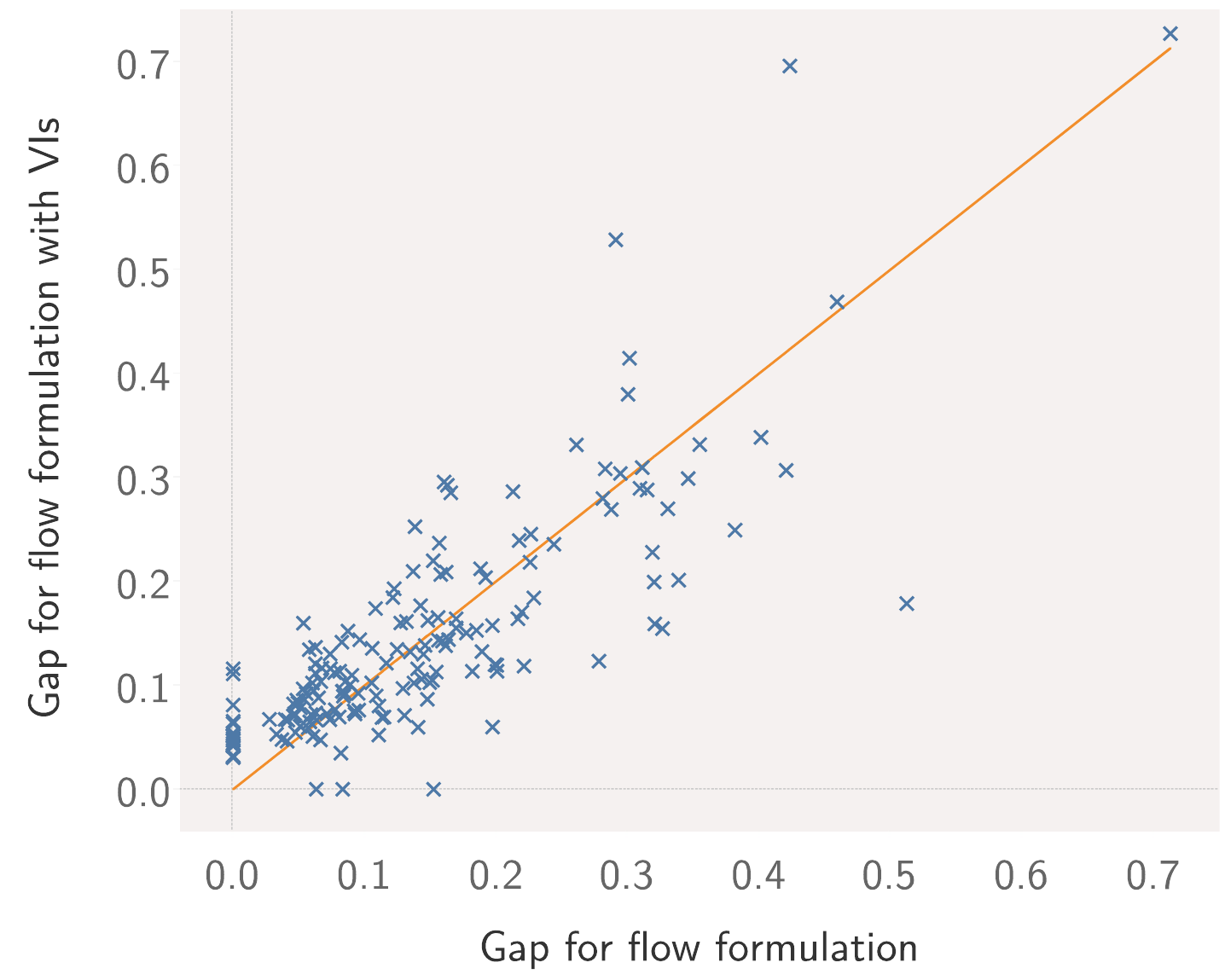}
    \hfill
    \includegraphics[width=0.47\textwidth]{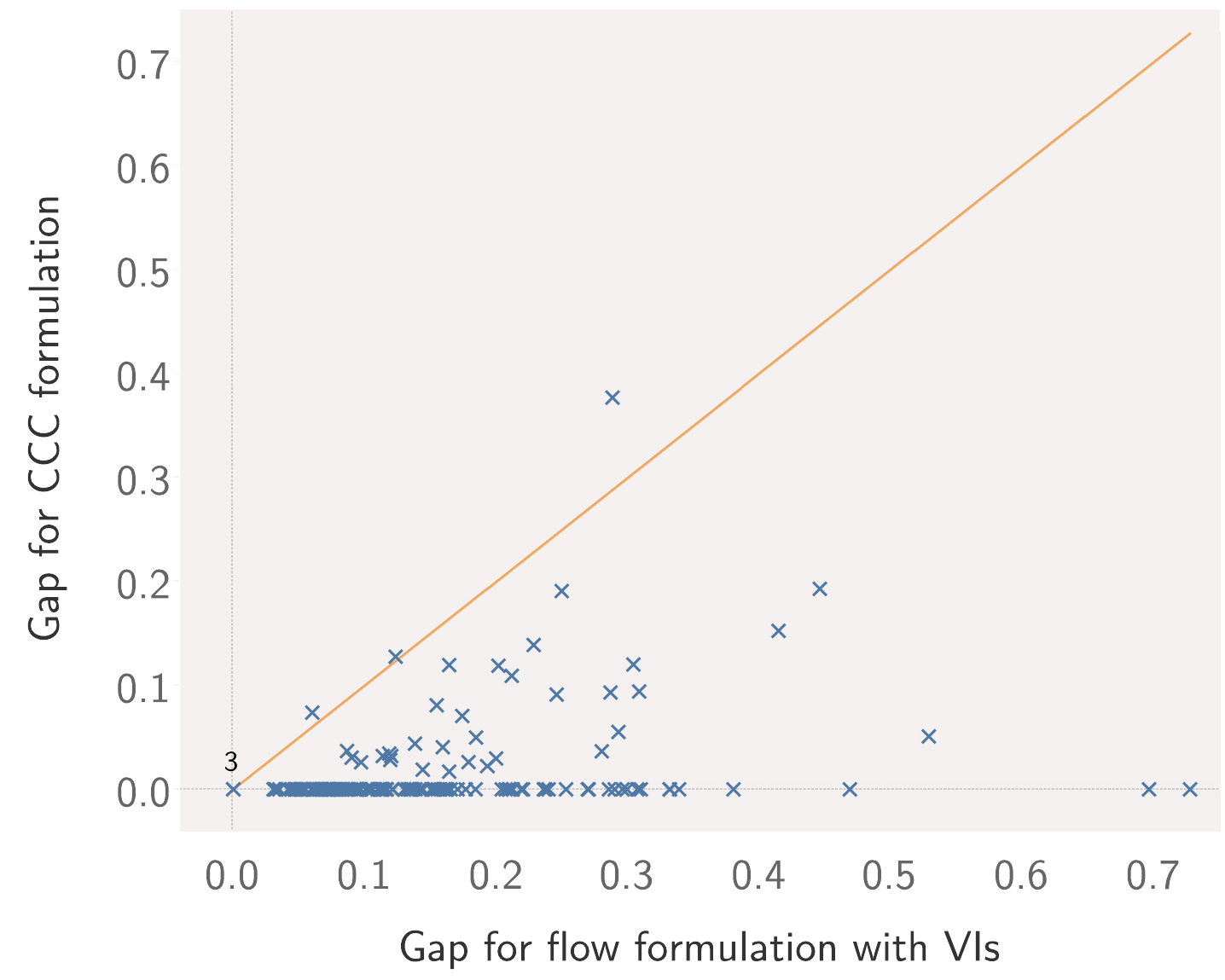}
    \caption{Optimality gaps for Group II instances}
    \label{fig:one_way_gaps}
  \end{subfigure}
  \caption{Computational results for the one-way experiments}
  \label{fig:one_way_computational}
\end{figure}

Figures \ref{fig:one_way_times} and \ref{fig:one_way_gaps} compare the times and gaps for the one-way case, respectively. Unlike the two-way scenario, adding valid inequalities to the flow formulation was counterproductive (see left panels). However, the branch-and-cut algorithm fared very well, as seen from the right panels. Five instances terminated with a very high gap when using the flow formulation and were excluded from these plots for better readability. 

\textbf{Properties of instances used for experiments with a two-hour limit:}
The properties of the instances for the extended analysis, where the model was run with a two-hour time limit, are summarized in Table \ref{tab:instance_properties}. The instance column displays abbreviated IDs, with the hyphen representing the sequence 00000. For example, the actual identifier for the instance 125-068 in the dataset is 12500000068. The area is in sq ft, and the proportion of blocked cells is expressed as a percentage of the overall number of cells in the lot. (The area measurements exclude the blocked cells.) The column $\mu \nu - |B|$ denotes the number of grid cells available for parking or driving after excluding blocked cells. The columns 2W$|D|$ and 1W$|D|$ refer to the number of valid driving cells in the two-way and one-way cases, respectively. Notice that compared to the two-way case, valid one-way driving cells are, on average, about 24\% higher and do not exhibit a multiplicative effect.

\begin{figure}[H]
    \centering
    \includegraphics[width=0.95\linewidth]{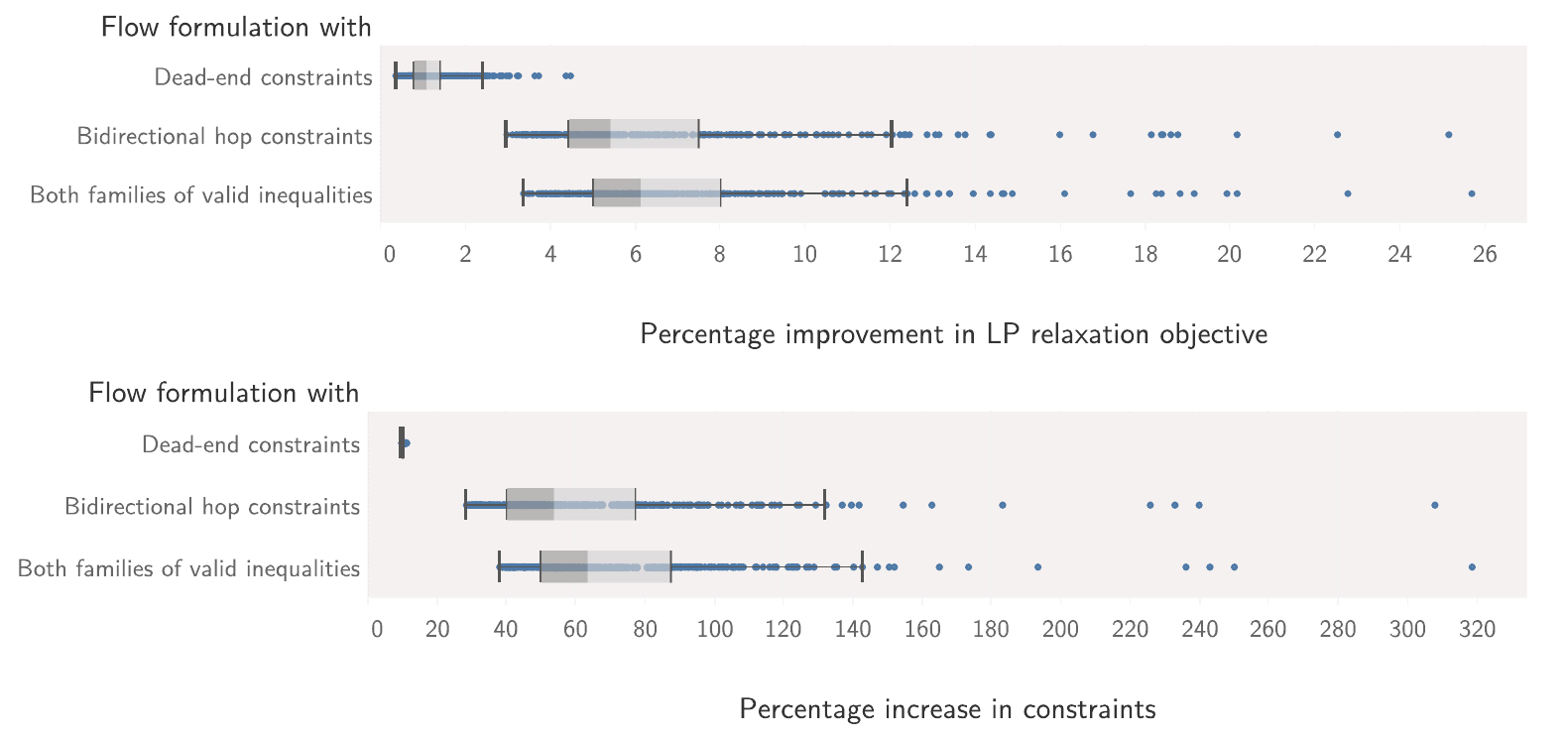}
    \caption{Comparison of the strength of valid inequalities for the one-way scenario}
    \label{fig:box_plots}
\end{figure}

\textbf{Comparison of the effectiveness of cut families:}
Recall that in the one-way scenario, two families of valid inequalities were introduced: dead-end and bidirectional hop inequalities. To test their relative effectiveness, we analyzed the improvement in the LP relaxation objectives of the flow-based formulation with and without each of these constraint families. We also examined the effect of including both sets of valid inequalities. The results for the 325 instances used in the paper are summarized in Figure \ref{fig:box_plots}. 


\begin{table}[H]
\small
  \centering
  \caption{Summary of instance properties for experiments with a two-hour limit}
    \begin{tabular}{lrrrcrrrrr}
    \hline
    \textbf{Instance ID} & \textbf{Area} & \textbf{$\mu$} & \textbf{$\nu$}  & \textbf{$\mu \nu - |B|$} & \textbf{$|B|$\%} & \textbf{$|P^0|$} & \textbf{$|P^{90}|$} & \textbf{2W$|D|$} & \textbf{1W$|D|$} \\
        \hline
125-068  & 6431  & 11    & 7     & 73    & 5.19  & 62    & 66    & 56    & 73 \\
        125-140  & 7596  & 15    & 15    & 109   & 51.56 & 93    & 93    & 78    & 109 \\
        125-299  & 9396  & 12    & 15    & 125   & 30.56 & 111   & 110   & 97    & 125 \\
        125-425  & 10941 & 18    & 18    & 152   & 53.09 & 134   & 134   & 117   & 152 \\
        125-116  & 13073 & 9     & 20    & 146   & 18.89 & 135   & 126   & 116   & 146 \\
        125-026  & 14359 & 23    & 12    & 168   & 39.13 & 145   & 156   & 134   & 168 \\
        125-269  & 16099 & 12    & 16    & 185   & 3.65  & 173   & 169   & 158   & 185 \\
        125-158  & 17726 & 11    & 33    & 239   & 34.16 & 221   & 206   & 189   & 239 \\
        125-355  & 19953 & 22    & 23    & 248   & 50.99 & 224   & 220   & 197   & 248 \\
        125-372  & 20847 & 30    & 14    & 264   & 37.14 & 234   & 249   & 220   & 264 \\
        215-145  & 22387 & 18    & 33    & 277   & 53.37 & 259   & 243   & 226   & 277 \\
        125-094  & 25118 & 23    & 15    & 291   & 15.65 & 268   & 276   & 254   & 291 \\
        125-407  & 26268 & 17    & 24    & 304   & 25.49 & 283   & 280   & 260   & 304 \\
        125-112  & 28108 & 20    & 21    & 327   & 22.14 & 306   & 306   & 286   & 327 \\
    \hline
    \end{tabular}%
  \label{tab:instance_properties}%
\end{table}%

The box plots show the percentage improvement in the LP objectives and the percentage increase in the number of constraints relative to the model without any valid inequalities. The median is indicated by the segment separating the gray portions, and the whiskers represent 1.5 times the interquartile range (IQR). The dead-end constraints were found to be less effective than the bidirectional hop inequalities, but are fewer in number. On average, adding dead-end constraints improves the LP relaxation objective by 1.24\%, compared with 6.46\% for bidirectional hop inequalities. However, relative to the size of the flow-based formulation, dead-end constraints add only 9.98\% more inequalities, compared to 62.45\% in the bidirectional hop case.

\newpage
\bibliographystyle{chicago}
\bibliography{references}

\end{document}